\documentclass[3p,final]{elsarticle}
\usepackage{ifthen}
\usepackage[ansinew]{inputenc}
\usepackage[T1]{fontenc}
\RequirePackage{geometry}

\usepackage{color}
\usepackage{ragged2e}
\usepackage{newlfont}
\usepackage{pifont}
\usepackage{yfonts}
\usepackage{dsfont}
\usepackage{stmaryrd} 
\usepackage{bbm} 
\usepackage{xspace}
\usepackage{amsmath,amssymb,amsfonts,textcomp}
\RequirePackage{amsthm}
\usepackage[mathscr]{eucal}
\usepackage{eufrak}
\usepackage{yhmath}
\RequirePackage{multirow}
\usepackage{pgf} 
\usepgfmodule[matrix]
\usepackage{tikz}

\usetikzlibrary{%
  calc,%
  arrows,%
  trees,%
  shapes.misc,
  shapes.arrows,%
  shapes.geometric,%
  chains,%
  matrix,%
  positioning,
  scopes,%
decorations.markings,%
  decorations.pathreplacing,%
  decorations.pathmorphing,
  shadows,%
backgrounds%
}

\makeatletter
\tikzset{join/.code=\tikzset{after node path={%
\ifx\tikzchainprevious\pgfutil@empty\else(\tikzchainprevious)%
edge[every join]#1(\tikzchaincurrent)\fi}}} \makeatother

\tikzset{
grph/.style ={>=latex', shorten <=-.3em, shorten >=-.3em},%
morph/.style ={>=stealth',shorten <=-.3em, shorten >=-.3em},%
every on chain/.append style={join},%
every join/.style={->},%
/seqStyle/.style ={>=stealth'},
}

\pgfkeys{
/Arrowstyle/.code=\pgfkeys{/Arrowstyle/.initial=#1},/Arrowstyle=morph,
/Dir/.code=\pgfkeys{/Dir/.initial=#1},/Dir=->,
/Dist/.code=\pgfkeys{/Dist/.initial=#1},/Dist=2,
/Bline/.code=\pgfkeys{/Bline/.initial=#1}, /Bline=0.0em,
/Lab_level/.code=\pgfkeys{/Lab_level/.initial=#1},/Lab_level=0,
/arrowstyle/.initial=morph,/dir/.initial=->,/above/.initial={},
/below/.initial={},/lab_level/.initial=0,/dist/.initial=2,/bline/.initial=-0.0em,
/isomark/.initial={\sim},/mark_level/.initial=-0.12,/side/.initial=above,
/style/.initial=, /syst/.initial={}%
}

\newcommand*{\midsloppy}{%
  \tolerance 5000%
  \hbadness 4000%
  \emergencystretch 1.5em%
  \hfuzz .1  
  \vfuzz\hfuzz}

\newtheorem{Thm}{Theorem}[section]
\newtheorem{Cor}[Thm]{Corollary}
\newtheorem{Lem}[Thm]{Lemma}
\newtheorem{Prop}[Thm]{Proposition}

\theoremstyle{definition}
\newtheorem{Defn}{Definition}[section]
\newtheorem{Rem}[Defn]{Remark}

\newtheorem{Expl}[Defn]{Example}
\newtheorem{Quest}{Question}
\newtheorem{Conj}[Quest]{Conjecture}

\newcommand*{\finprv}{\nobreak \hspace{\stretch{1}}\qed} 
\newenvironment{prv}{
\begin{proof}}{\end{proof} } 

\newcommand{\Par}[1]{ \mbox{ } \\ \smallskip \noindent\textbf{#1}.~}
\newcommand{\ParIt}[1]{ \mbox{ } \\ \smallskip \noindent\textit{#1}.~}

\numberwithin{equation}{section}

\newcommand*{\mb}[1]{\mathbb{#1}}

\newcommand*{\mbm}[1]{\mathbbm{#1}}
\newcommand*{\mrm}[1]{\mathrm{#1}}
\newcommand*{\msf}[1]{\mathsf{#1}}
\newcommand*{\mfr}[1]{\mathfrak{#1}}
\newcommand*{\scr}[1]{{\scriptstyle{#1}}}
\newcommand*{\sscr}[1]{{\scriptscriptstyle{#1}}} 

\newcommand*{\bbar}[1]{\overline{#1}} 
\newcommand*{\wtilde}[1]{\widetilde{#1}} 
\newcommand*{\what}[1]{\widehat{#1}} 

\newcommand*{\dg}{\textsf{dg}}
\newcommand*{\df}[1][\mfr{d}]{#1}
\newcommand*{\Ker}{\mrm{Ker}}

\newcommand*{\jtriv}{\msf{j}_{\msf{triv}}}

\newcommand*{\K}{\mathsf{\mrm{K}}}
\newcommand*{\Ke}{\mathsf{\mrm{K^e}}}
\newcommand*{\kk}{\mathsf{\mrm{k}}}

\newcommand*{\E}{\mathsf{\mrm{E}}}
\newcommand*{\F}{\mathsf{\mrm{F}}}
\newcommand*{\LL}{\mathsf{\mrm{L}}}

\newcommand*{\md}{\mrm{mod}}

\newcommand*{\C}{\mathcal{C}}

\newcommand*{\D}[1][]{\ind{\mathcal{D}}{\mrm{#1}}}
\newcommand*{\Db}{\mrm{D^{b}}}

\newcommand*{\T}{\cal{T}}

\newcommand*{\Grd}{\mrm{Gr}}

\newcommand*{\id}{\mathds{1}\hspace{-0.1ex}} 
\newcommand*{\restr}[1]{{{_{\mid_{#1}}}}}

\newcommand*{\Hm}{\mrm{Hom}}
\newcommand*{\End}{\mrm{End}}

\newcommand*{\Endc}[1]{\ind{\mrm{End}}{\C}(#1)}

\newcommand*{\HM}[2]{\ind{\mrm{Hom}}{#1}(#2)}

\DeclareMathOperator{\soc}{\mrm{soc}}

\newcommand*{\J}[2][]{\indu{\mrm{J}}{#2}{#1}}
\newcommand*{\Jtriv}[2][]{\indu{\mrm{J}}{#2}{#1\mrm{triv}}}

\newcommand*{\Jc}[1][]{\ind{\cal{J}}{#1}} 

\newcommand*{\Ss}{\mrm{S}}

\newcommand*{\B}{\mathrm{B}}

\newcommand{\dm}[1][n]{\msf{\mrm{d}}}

\newcommand*{\set}[1]{\left\{#1\right\}}

\newcommand*{\I}{\mathrm{I}}
\newcommand*{\N}{\mathbbm{N}}
\newcommand*{\Ns}{\ind{\N}{*}}
\newcommand*{\ggen}[1]{\langle #1\rangle} 
\newcommand*{\Psmatr}[1]{\left[\begin{smallmatrix}
#1 \end{smallmatrix} \right]} 
\newcommand*{\parr}[1]{\Bigl( #1\Bigr)}
\newcommand*{\Parr}[1]{\biggl( #1\biggr)}

\newcommand*{\ind}[2]{{#1\hspace{-0.02em}}_{#2}}
\newcommand*{\inddl}[3]{{_{#2\hspace{-0.02em}}}{\mfr{#1}\hspace{-0.02em}}_{#3}}
\newcommand*{\inds}[2]{{#1\hspace{-0.02em}}_{{\scriptscriptstyle{#2}}}}
\newcommand*{\indu}[3]{{#1\hspace{-0.02em}}_{#2}^{#3}}

\newcommand*{\indd}[3]{{_{#1\hspace{-0.02em}}}{#2}_{\hspace{-0.02em}#3}}
\newcommand*{\indds}[3]{{_{{\scriptscriptstyle{#1}}\hspace{-0.02em}}{#2\hspace{
-0.02em}}_{{\scriptscriptstyle{#3}}}}}
\newcommand*{\inddu}[4]{{_{#1\hspace{-0.02em}}{#2\hspace{-0.02em}}_{#3}^{#4}}}
\newcommand*{\inddsu}[4]{{_{{\scriptscriptstyle{#1}}\hspace{-0.02em}}{#2\hspace{
-0.02em}}_{{\scriptscriptstyle{#3}}}^{{\scriptscriptstyle{#4}}}}}

\newcommand{\myeqar}[2][rl]{\begin{array}{#1} #2 \end{array}}

\newcommand*{\som}[3]{\sum\limits_{#1}^{#2}\hspace{-0.12ex}{#3}}
\newcommand*{\Som}[2][]{\sum\limits_{#2}^{#1}\hspace{-0.12ex}}

\newcommand*{\somd}[2][]{\overset{\sscr{#1}}{\underset{\sscr{#2}}{\varoplus}
\hspace{-0.07ex}}}
\newcommand*{\prd}[2]{\raisebox{-0.15em}{$\underset{#1}{{\rule{0.14ex}{2ex}
\hspace{0.7ex}{\rule{0.14ex}{2ex}}
\hspace{-1.16ex}{\raisebox{2ex}{\rule{1.32ex}{0.18ex}}}}}$} #2}  
\newcommand*{\tens}[2][]{\inds{\overset{#1}{\otimes}}{#2}}

\newcommand*{\tenss}[1][]{\inds{{\scr{\varotimes}}}{#1}}

\newcommand*{\Union}[3][]{\raisebox{+0.2ex}{\text{${\scriptstyle{
\bigcup\limits_{#2}^{#1}}}$}}\hspace{-0.05ex}#3}
\newcommand*{\inter}[3]{\raisebox{+0.2ex}{\text{${\scriptstyle{
\bigcap\limits_{#1}^{#2}}}$}}\hspace{-0.05ex}#3}
\newcommand*{\Slash}{\raisebox{-0.25em}{\begin{tikzpicture}
\draw[semithick] (0em,0em) -- (0.25em,1em); \end{tikzpicture}}}

\newcommand*{\Frac}[2]{#1\hspace{-0.15em}{\Slash}{\hspace{-0.075em}#2}}

\newcommand*{\ninterv}[2][1]{\llbracket#1\hspace{0.1em},\hspace{0.1em}#2\rrbracket}

\newcommand*{\mul}{\hspace{-0.04em}{{\mathbf{\cdot}}}\hspace{-0.04em}}
\newcommand*{\sdash}{\textrm{-}}

\DeclareMathOperator{\im}{\mrm{Im}}

\newcommand*{\sminus}{\raisebox{0.1em}{{\tiny{$\smallsetminus$}}}}

\newcommand*{\sign}{\mrm{sign}}
\newcommand*{\idl}[1][]{\msf{1}_{#1}} 
\newcommand*{\ii}{\msf{i}} 

\newcommand*{\vphi}[1][{\hspace{0.015em}}]{\ind{\varphi}{#1}}

\newcommand*{\elts}[3]{\indds{#2}{#1}{#3}}

\newcommand*{\cderv}[1][]{\inds{\partial}{#1}}
\newcommand*{\rderv}[1][]{\partial_{#1}^{\text{\textsc{{\tiny{r}}}}}} 

\newcommand*{\lderv}[1][]{\partial_{#1}^{\text{\textsc{{\tiny{l}}}}}}

\newcommand*{\rgt}{\sscr{\text{\textsc{r}}}}

\newcommand*{\Right}{\text{\textsc{r}}}
\newcommand*{\lft}{\sscr{\text{\textsc{l}}}}

\newcommand*{\Left}{\text{\textsc{l}}}

\newcommand*{\Zc}[1][]{\ind{\mathcal{Z}}{#1}}
\newcommand*{\Zck}[1][\K]{\ind{\mathcal{Z}}{#1}}

\newcommand*{\tr}{\mfr{t}}
\newcommand*{\htr}[1][]{\ind{\what{\mfr{t}}}{#1}}

\newcommand*{\scprd}[1]{\scr{{\boldsymbol{\langle}}} {#1}
{\scr{{\boldsymbol{\rangle}}}}} 
\newcommand*{\z}[1][]{\inds{\mfr{z}}{#1}}
\newcommand*{\zc}[1][c]{\mfr{\inds{z}{\!#1}}}
\newcommand*{\ztc}[1][c]{\mfr{\inds{\widetilde{z}}{\!#1}}}

\newcommand*{\zb}[1][]{\inds{\bbar{\mfr{z}}}{#1}}

\newcommand*{\lperm}[1][]{\indu{\varepsilon}{\lft}{#1}}
\newcommand*{\rperm}[1][]{\indu{\varepsilon}{\rgt}{#1}}
\newcommand*{\cperm}[1][]{\indu{\varepsilon}{\mrm{c}}{#1}}

\newcommand*{\lgrad}[1][]{\triangle_{#1}^{\text{\textsc{{\tiny{l}}}}}}
\newcommand*{\rgrad}[1][]{\triangle_{#1}^{\text{\textsc{{\tiny{r}}}}}}
\newcommand*{\bilbox}[1][]{\Box}
\newcommand*{\bilf}[1][{{}{}}]{\inddl{\mfr{b}}#1} 
\newcommand*{\bilfb}[1][{{}{}}]{\inddl{\bbar{\mfr{b}}}#1}
\newcommand*{\bilfh}[1][]{\widehat{\ind{\mfr{b}}{#1}}}
\newcommand*{\bilfph}[1][]{\widehat{\ind{\mfr{b}}{#1}'}}
\newcommand*{\bilft}[1][{{}{}}]{\inddl{\widetilde{\mfr{b}}}#1}

\newcommand*{\bilff}[1][{{}{}{}}]{\inddU{\mfr{b}}#1}

\newcommand*{\inddU}[4]{{_{#2\hspace{-0.02em}}{#1\hspace{-0.02em}}_{#3}^{#4}}}

\newcommand*{\vlt}[1][{d}]{\msf{#1}}
\newcommand*{\vl}[2][d]{\ind{\mrm{#1}}{#2}}
\newcommand*{\vld}[2][d]{\sdual{\vl[#1]{#2}}}

\newcommand*{\val}[3][d]{\inddl{\mrm{#1}}{#2}{#3}}
\newcommand*{\vald}[3][d]{\sdual{\val[#1]{#2}{#3}}}
\newcommand*{\valp}[3][d]{\val[#1]{#2}{#3}'}
\newcommand*{\valpd}[3][d]{\sdual{\valp[#1]{#2}{#3}}}

\newcommand*{\n}[1][n]{\msf{#1}}

\newcommand*{\Bmd}[3][]{\inddu{#2}{B}{#3}{#1}}

\newcommand*{\Bmdd}[2]{\dual{\indd{#1}{B}{#2}}}

\newcommand*{\op}[1]{#1^{{\text{{\tiny{$\circ$}}}}}}

\newcommand*{\dual}[1]{#1^{\star}}
\newcommand*{\sdual}[1]{#1^{\hspace{-0.0625em}{{\sscr{\star}}}}}
\newcommand*{\ldual}[1]{{^{\raisebox{+0.16ex}{$\lft$}\hspace{-0.2em}}#1}}
\newcommand*{\rdual}[1]{#1^{\raisebox{+0.16ex}{$\rgt$}}}
\newcommand*{\ddual}[2][\diamond]{#2^{#1}}
\newcommand*{\shdual}[2][\sharp]{#2^{{\hspace{-0.1em}\tiny{#1}}}}

\newcommand*{\Q}{\msf{\mrm{Q}}}

\newcommand*{\cyc}[2][]{\inds{#2}{\msf{cyc}#1}}  
\newcommand*{\mQ}{{\mathcal{Q}}} 
\newcommand*{\mQb}{\bbar{\mathcal{Q}}} 
\newcommand*{\mQh}{\widehat{\mathcal{Q}}} 

\newcommand*{\mQk}{\ind{\mathcal{Q}}{\kk}} 

\newcommand*{\lmQ}[1][1]{\indd{\Left\hspace{-0.1em}}{\mQ}{\sscr{#1}}}
\newcommand*{\rmQ}[1][1]{\indd{\Right\hspace{-0.1em}}{\mQ}{\sscr{#1}}}
\newcommand*{\ldmQ}[1][1]{\dual{\lmQ[#1]}}
\newcommand*{\rdmQ}[1][1]{\dual{\rmQ[#1]}}
\newcommand*{\lmQp}[1][1]{\indd{\Left\hspace{-0.1em}}{\mQ}{\sscr{#1}}'}
\newcommand*{\rmQp}[1][1]{\indd{\Right\hspace{-0.1em}}{\mQ}{\sscr{#1}}'}
\newcommand*{\ldmQp}[1][1]{\dual{\lmQp[#1]}}
\newcommand*{\rdmQp}[1][1]{\dual{\rmQp[#1]}}
\newcommand*{\lmQc}[1][1]{\what{\indd{\Left\hspace{-0.1em}}{\mQ}{\sscr{#1}}}}
\newcommand*{\rmQc}[1][1]{\what{\indd{\Right\hspace{-0.1em}}{\mQ}{\sscr{#1}}}}
\newcommand*{\ldmQc}[1][1]{\what{\dual{\lmQ[#1]}}}
\newcommand*{\rdmQc}[1][1]{\what{\dual{\rmQ[#1]}}}
\newcommand*{\lmQpc}[1][1]{\what{\indd{\Left\hspace{-0.1em}}{\mQ}{\sscr{#1}}'}}
\newcommand*{\rmQpc}[1][1]{\what{\indd{\Right\hspace{-0.1em}}{\mQ}{\sscr{#1}}'}}
\newcommand*{\ldmQpc}[1][1]{\what{\dual{\lmQp[#1]}}}
\newcommand*{\rdmQpc}[1][1]{\what{\dual{\rmQp[#1]}}}

\newcommand*{\rmQb}[1][1]{\indd{\Right\hspace{-0.1em}}{\mQb}{\sscr{#1}}}

\newcommand*{\rmQbc}[1][1]{\what{\indd{\Right\hspace{-0.1em}}{\mQb}{\sscr{#1}}}}

\newcommand*{\rdmQbc}[1][1]{\what{\dual{\rmQb[#1]}}}
\newcommand*{\BmQ}[3][]{\inddu{#2\hspace{-0.1em}}{\mQ}{#3}{#1}}
\newcommand*{\BmdQ}[3][]{\inddu{#2\hspace{-0.1em}}{\mQ}{#3}{\ast #1}}
\newcommand*{\BmQc}[2]{\widehat{\indd{#1\hspace{-0.1em}}{\mQ}{#2}}}
\newcommand*{\BmdQc}[2]{\widehat{\inddu{#1\hspace{-0.1em}}{\mQ}{#2}{\ast}}}

\newcommand*{\kQ}{\kk\mQ}
\newcommand*{\kQd}{\kk\dual{\mQ}}
\newcommand*{\kQc}{\what{\kk\mQ}}
\newcommand*{\kQdc}{\what{\kk\dual{\mQ}}}
\newcommand*{\kQpc}{\what{\kk\mQ'}}

\newcommand*{\kQshc}{\widehat{\kk\shdual{\mQ}}}
\newcommand*{\kQh}{\kk\mQh}
\newcommand*{\kQhc}{\widehat{\text{$\msf{\kk}$}\widehat{\mathcal{Q}}}}

\newcommand*{\kQtc}{\widehat{\text{$\msf{\kk}$}\widetilde{\mathcal{Q}}}}

\newcommand*{\kQbc}{\what{{\kk\bbar{\mQ}}}}
\newcommand*{\kQbpc}{\what{\kk\msf{\bbar{\mQ}'}}}
\newcommand*{\kQk}{\kk\mQk}
\newcommand*{\kQkc}{\what{\kk\mQk}}

\newcommand*{\Red}{\mrm{red}}
\newcommand*{\red}[1]{\ind{#1}{\mrm{red}}}
\newcommand*{\triv}[1]{\ind{#1}{\mrm{triv}}}
\newcommand*{\skc}[1]{\mrm{skew}[#1]}
\newcommand*{\skcb}[1]{\mrm{skew}\hspace{-0.3em}\set{#1}}

\newcommand*{\sr}{\mathsf{s}} 
\newcommand*{\tg}{\mathsf{t}} 

\newcommand*{\lin}[2][]{\pgfkeys{/innersep/.initial=0.1,#1}
\begin{tikzpicture}[grph,baseline=(a.base)]
\node  (a) at (0,0) [inner sep=\pgfkeysvalueof{/innersep}em]
{$#2$}; \path[ultra thin,->,shorten <=-0.1em,shorten
>=-0.1em,densely dashed] (a.north west) edge  (a.north east);
\end{tikzpicture}}

\renewcommand*{\to}[1][]{\hspace{-0.4em}
\pgfkeys{/arrowstyle=\pgfkeysvalueof{/Arrowstyle},/dir=\pgfkeysvalueof{/Dir},
/style=,/dist=\pgfkeysvalueof{/Dist},/lab_level=\pgfkeysvalueof{/Lab_level},/above={},
/below={},/mark_level=-0.12,/bline=-0.25em,/isomark={},/side=above,#1}
\begin{tikzpicture}[\pgfkeysvalueof{/arrowstyle},baseline=\pgfkeysvalueof{/bline},
node distance=\pgfkeysvalueof{/dist}em]
\node (i) {}; \node (j)[right=of i]  {}; \node (label_pos) [inner
sep=0.12em] at ($ (i.east)!0.5! (j.west) $) {};
\path[\pgfkeysvalueof{/dir},\pgfkeysvalueof{/style}] (i) edge (j);
\node [inner sep=0em,above=\pgfkeysvalueof{/lab_level}em of
label_pos] {{$\scriptstyle{\pgfkeysvalueof{/above}}$}};
\node [inner sep=0em,below=\pgfkeysvalueof{/lab_level}em of
label_pos] {{$\scriptstyle{\pgfkeysvalueof{/below}}$}};
\node [inner sep=0em,
\pgfkeysvalueof{/side}=\pgfkeysvalueof{/mark_level}em of
label_pos] {${\scriptscriptstyle{\pgfkeysvalueof{/isomark}}}$};
\end{tikzpicture}
\pgfkeys{/above={},/below={},/mark_level=-0.12}\hspace{-0.4em}}

\newcommand*{\path}[2][]{
\pgfkeys{/Arrowstyle=grph,/Dir=->,/Dist=1,/Lab_level=-0.05,/bline=-0.0em,#1}
\makebox{\ensuremath{#2}}
\pgfkeys{/Arrowstyle=morph,/Dist=2,/Dir=->,/Lab_level=0}}
\newcommand*{\edge}[3][]{
\path{#2 \to[#1] #3}\pgfkeys{/below={}}}

\newcommand*{\graph}[3][]{\hspace{-0.5em}
\pgfkeys{/Dist=1,/dist/.initial=\pgfkeysvalueof{/Dist},clsep/.initial=1,
/bline=-0.2em,Lab_level=-0.08,
/lab_level/.initial=\pgfkeysvalueof{/Lab_level},#1}
\begin{tikzpicture}[grph,baseline=\pgfkeysvalueof{/bline}]
\matrix (m) [matrix of math nodes,ampersand replacement=\&, row
sep=\pgfkeysvalueof{/dist}em,column sep=\pgfkeysvalueof{/clsep}em,
text height=0.5em, text depth=0.0em]
{#2 \\}; #3;
\end{tikzpicture}
\pgfkeys{/Dist=2,/Lab_level=0}\hspace{-0.5em} }

\newcommand*{\Seq}[2][]{
\pgfkeys{/Arrowstyle=morph,/Dist=1.5,/Lab_level=0.0,Dir=->,/bline=-0.2em,#1}
\makebox{\ensuremath{#2}}
\pgfkeys{/Arrowstyle=morph,/Dist=2,/Dir=->,/Lab_level=0}}
\newcommand*{\morph}[3][]{ 
\Seq{#2 \to[/dist=2,#1] #3}}
\newcommand*{\diagram}[3][]{
\pgfkeys{/Dist=2,/dist=1,clsep/.initial=2,/bline=0,#1}
\begin{tikzpicture}[morph,baseline=\pgfkeysvalueof{/bline}]
\matrix (m) [matrix of math nodes,ampersand replacement=\&, row
sep=\pgfkeysvalueof{/dist}em,column
sep=\pgfkeysvalueof{/clsep}em,text height=1em, text depth=0.4em]
{#2 \\}; #3;
\end{tikzpicture} \pgfkeys{/Dist=2}
}

\newcommand*{\cyclF}[7][]{\pgfkeys{/dir=->,
/dist=1,/lab_level=-0.1,/bline=-0.25em,/scale/.initial=1,#1}
{\begin{tikzpicture}[grph,baseline=\pgfkeysvalueof{/bline},node
distance=\pgfkeysvalueof{/dist}em,scale=\pgfkeysvalueof{/scale}]
\node (a) at (0,0) {$#2$}; \node (b) at (0,1) {$#3$}; \node (c) at
(1,1) {$#4$}; \node (d) at (1,0) {$#5$};
\path[\pgfkeysvalueof{/dir}] (a)  edge node[] {} (b)
              (b)  edge node[above=\pgfkeysvalueof{/lab_level}]{$\scr{#6}$} (c)
              (c)  edge  node[] {} (d)
              (d)  edge  node[midway,below=\pgfkeysvalueof{/lab_level}]
{$\scr{#7}$} (a);
\end{tikzpicture}}
\pgfkeys{/dist=\pgfkeysvalueof{/Dist},/lab_level=\pgfkeysvalueof{/Lab_level}}}
\newcommand*{\cyclFd}[8][]{\pgfkeys{/dir=->,
/dist=1,/lab_level=-0.1,/bline=-0.25em,/scale/.initial=1,#1}
{\begin{tikzpicture}[grph,baseline=\pgfkeysvalueof{/bline},node
distance=\pgfkeysvalueof{/dist}em,scale=\pgfkeysvalueof{/scale}]
\node (a) at (0,0) {$#2$}; \node (b) at (0,1) {$#3$}; \node (c) at
(1,1) {$#4$}; \node (d) at (1,0) {$#5$};
\path[\pgfkeysvalueof{/dir}]
              (a)  edge node[] {} (b)
              (b)  edge node[above=\pgfkeysvalueof{/lab_level}]{$\sscr{#6}$} (c)
              (c)  edge
node[sloped,above=\pgfkeysvalueof{/lab_level}]{$\sscr{#7}$} (a)
              (a)  edge  node[midway,below=\pgfkeysvalueof{/lab_level}]
{$\sscr{#8}$} (d)
              (d)  edge  node[] {} (c);
\end{tikzpicture}}
\pgfkeys{/dist=\pgfkeysvalueof{/Dist},/lab_level=\pgfkeysvalueof{/Lab_level}
}}

 \newcommand*{\graphI}[6][]{\pgfkeys{/dir1/.initial=->,/dir2/.initial=->,
/dist=1,/lab_level=-0.1,/bline=-0.25em,/scale/.initial=1,#1}
{\begin{tikzpicture}[grph,baseline=\pgfkeysvalueof{/bline},node
distance=\pgfkeysvalueof{/dist}em,scale=\pgfkeysvalueof{/scale}]
#2; \path[\pgfkeysvalueof{/dir1}] (a)  edge node
[sloped,above=\pgfkeysvalueof{/lab_level}] {$\sscr{#3}$} (b)
               (b)  edge node [sloped,above=\pgfkeysvalueof{/lab_level}]
{$\sscr{#4}$} (c)
               (c)  edge  node [midway,below=\pgfkeysvalueof{/lab_level}]
{$\sscr{#5}$} (a); \path[\pgfkeysvalueof{/dir2}] (c)  edge node
[above=\pgfkeysvalueof{/lab_level}] {$\sscr{#6}$} (d);
\end{tikzpicture}}
\pgfkeys{/dist=\pgfkeysvalueof{/Dist},/lab_level=\pgfkeysvalueof{/Lab_level}}}

\newcommand*{\graphIi}[6][]{\pgfkeys{/dir1/.initial=->,/dir2/.initial=->,
/dist=1,/lab_level=-0.1,/bline=-0.25em,/scale/.initial=1,#1}
{\begin{tikzpicture}[grph,baseline=\pgfkeysvalueof{/bline},node
distance=\pgfkeysvalueof{/dist}em,scale=\pgfkeysvalueof{/scale}]
#2; \path[\pgfkeysvalueof{/dir1}]  (a)  edge node
[above=\pgfkeysvalueof{/lab_level}] {$\sscr{#3}$} (b);
\path[\pgfkeysvalueof{/dir2}]   (b)  edge node
[sloped,above=\pgfkeysvalueof{/lab_level}] {$\sscr{#4}$} (c)
               (c)  edge node [sloped,above=\pgfkeysvalueof{/lab_level}]
{$\sscr{#5}$} (d)
               (d)  edge  node [midway,below=\pgfkeysvalueof{/lab_level}]
{$\sscr{#6}$} (b);
\end{tikzpicture}}
\pgfkeys{/dist=\pgfkeysvalueof{/Dist},/lab_level=\pgfkeysvalueof{/Lab_level}}}

\newcommand*{\NodeI}[4]{
\node (a) at (-1.25,0) {$#1$}; \node (b) at (0,1) {$#2$}; \node
(c) at (1.25,0) {$#3$}; \node (d) at (2.25,0) {$#4$};}
 
\newcommand*{\NodeIi}[4]{
\node (a) at (-2.25,0) {$#1$}; \node (b) at (-1.25,0) {$#2$};
\node (c) at (0,1) {$#3$}; \node (d) at (1.25,0) {$#4$};}

\newcommand*{\trr}{{\scriptscriptstyle{\triangleright}}}

\allowdisplaybreaks[1]

\usepackage[pdftex]{hyperref}
\hypersetup{backref,implicit=true,bookmarks=true,
colorlinks=true, linkcolor=blue,filecolor=blue,
pagecolor=blue,urlcolor=blue,pdfpagemode=UseOutlines}

   \journal{\, } 

\begin{document}
\midsloppy4pt \fontfamily{ptm}\fontseries{m}\selectfont

\begin{frontmatter}
\title{Potentials and Jacobian algebras for tensor algebras of bimodules}

\author{Bertrand Nguefack\fnref{fnote}}


\fntext[fnote]{The  author acknowledges  the financial support from
{\texttt{CIDA}} for his PhD's studies during which part of this work started}




\begin{abstract}
We introduce and study potentials, mutations and Jacobian algebras
in the  framework of tensor algebras associated with
symmetrizable dualizing pairs of bimodules on a symmetric  algebra
over any commutative ground ring. The graded context is also considered by starting from graded bimodules, and the classical  non simply-laced context of modulated quivers with potentials is  a particular case. The study of potentials in this framework  is   related to symmetrically separable algebras, and we have two kinds of potentials: the symmetric  and the non symmetric ones.
When the Casimir ideal of the symmetric algebra coincides with its
center,  all potentials appear as symmetric potentials and their manipulation mimics the simply laced study of quivers with potentials. This useful information suggests that, for applications to cluster algebras theory and related fields, one may restrict a further study of modulated quivers with potentials to the setting where the ground symmetric algebra is separable over a field. Associated with this work is a generalized construction of Ginzburg dg-algebras and cluster categories associated with graded modulated  quivers with potentials.
\end{abstract}


\end{frontmatter}

\tableofcontents

\midsloppy4pt

\section{Introduction}

The main purpose in this paper is to extend  to a suitable general
framework some recent  aspects of the theory of  quivers with
potentials and corresponding Jacobian algebras started in
\cite{DWZ}. The First motivation of this work  is a result of \cite{BIRSm} relating   the mutation of cluster tilting objects in $2$-Calabi-Yau
categories  to the mutation of quivers with potentials. 
In the simply laced case, the theory of quivers with potentials
was motivated by several sources:  superpotentials in physics \cite{DM, BD,B2003},
Calabi-Yau algebras \cite{Boc,Ginzburg,IR,KR2}, cluster algebras. The
original motivation for the study of quivers with potentials comes
from the theory of cluster algebras introduced and studied in a
series of papers \cite{FZ1,FZ2,BFZ,FZ6} by S. Fomin and A.
Zelevinsky. The  underlying combinatorics of the theory of cluster
algebras is embodied in skew-symmetrizable integer matrices and
their mutations, or equivalently, in valued quivers without loops and their mutations. However, most of the time, recent categorifications of cluster theory
restrict to the simply-laced case, that is the one corresponding
to skew-symmetric matrices or equivalently to $2$-acyclic
quivers without loops. 

\subsection*{The present framework and the method}
 In this introductory discussion, we do not provide explicit definitions for some notions announced  here and kindly refer the reader to the text for full detailed definitions. The  general framework considered here is based on the existence of the so-called trace maps on simple algebras \cite[\textsection{22}]{Draxl}. We let
$\kk$ be  any   commutative ring and $(\K,\tr)$ a symmetric
$\kk$-algebra, finitely  generated projective as  $\kk$-module, here  $\tr\in\Hm_{\kk}(\K,\kk)$ is a strongly non-degenerate trace map for $\K$, that is, $\tr$ induces an isomorphism of $\K$-bimodules $\morph[/above={\sim}]{\K}{\Hm_{\kk}(\K,\kk)}$ taking each $a\in \K$ to the $\kk$-linear map $\Seq{\tr(a\mul\sdash ): \, b\mapsto \tr(ab) }$. Let  $B$ be a $\K$-bimodule, finitely generated projective as  left $\K$-module and  aright
$\K$-module. Then $B$ appears as part of a data
$\set{B,\dual{B},\bilf}$ which we call a \emph{symmetrizable
dualizing pair of $\K$-bimodules}, here     $\morph[/above={\bilf}]{B\otimes \dual{B} \oplus \dual{B}\otimes B}{\K}$     is a \emph{strongly non-degenerate bilinear form} and,   
$\tr$ is a \emph{symmetrizing map} for $\bilf$, that is, $\tr\bilf(x\otimes \xi)=\tr\bilf(\xi\otimes x)$ for all $x\in B$ and $\xi\in\dual{B}$, 
see Definition~\ref{defn.dualizing-bim}. The data $\mQ=(B,\K,\tr)$ is is called a
 \emph{$\kk$-modulated quiver} having $B$ as \emph{arrow bimodule}.  The 
\emph{path algebra} $\kQ$ of   $\mQ$ (or the path algebra of  $B$)  is  the tensor algebra of $B$ over $\K$; thus $\kQ=\mrm{T}_{\K}(B)=\somd{l\geq 0}{B^{(l)}}$ where $B^{(l)}=\kk\mQ_l$ is the $l$-fold tensor product of $B$ over $\K$ (referred to as the bimodule  generated by all length-$l$ paths in $\mQ$), with $B^{(0)}=\K$. The \emph{complete path algebra} of $\mQ$ is given by $\what{\kk\mQ}=\prd{l\geq 0}{B^{(l)}}$. Write $\kQc_{(d)}=\prd{l\geq d}{B^{(l)}}$ for all natural number $d\geq 1$ and let $\J{\kQc}=\kQc_{(1)}$. Then $\kQc$ is a topological algebra with $\J{\kQc}$-adic topology and $\J{\kQc}$ is referred as the \emph{closed arrow ideal} of $\mQ$. 
Observe that  the classical non simply-laced context  is recovered when $\K$ is a  direct product $\prd{i\in\ninterv{n}}{\kk_i}$ of division algebras over a  field $\kk$, 
here  $n\geq 1$ is a natural number, $\ninterv{n}=\set{1,\dotsc,n}$ and each $\kk_i$ is viewed  as  subfield  in $\K$ with unit $\idl[i]$. On the other hand,  the simply-laced context is obtained when  $\K$ occurs as  elementary semisimple algebra $\kk^n=\prd{i\in\ninterv{n}}{\kk_i}$ over a field $\kk$, here 
  $\kk_i=\kk$ for all $i\in\ninterv{n} $; in this case $B$ is a central $\kk^n$-bimodule and the data $\Q=(B,\kk^n)$ may be referred to as a \emph{$\kk$-quiver}, the arrows of $\Q$ correspond bijectively to the union of  $\kk$-bases of $\idl[i]\mul B\mul \idl[j]$ with $i ,j\in \ninterv{n}$. For a $\kk$-quiver $\Q=(B,\kk^n)$,  a \emph{potential} $W$ on $\Q$ was defined as a possibly infinite sum of cyclic elements  in $\what{\kk\Q}_{(2)}$;  to $W$ is associated a closed ideal $\J{W}$, called the \emph{Jacobian ideal} of $W$ and generated by the cyclic derivatives of $W$ with respect to the arrows of $\Q$,  the quotient  algebra $\Jc{(\Q,W)}:=\Frac{\what{\kk\Q}}{\J{W}}$ is called the \emph{Jacobian algebra} \cite{DWZ}.  Next,  we enrich the framework just described by  starting with $G$-graded $\K$-bimodules $B$ for an abelian group $G$ and considering potentials of homogeneous  degree with respect to $G$-grading.

In the present  framework, in order to get  an appropriate notion of potential with
respect to cyclic derivatives  we must  lift
ordinary permutations of arrows from simply laced path algebras to a
kind of skew permutations for tensor algebras $\kQc$.  This can be achieved  in two complementary ways. Let us describe  the most general and intrinsic
method of our study.  Given a symmetrizable dualizing pair of $\K$-bimodules $\set{M,M', \beta}$, we observe that the induced non-degenerate bilinear forms $\morph{M\tenss M'}{\K}$ and $\morph{M' \tenss M}{\K}$ are
\emph{dualizing morphisms} and their dual morphisms give rise the
following \emph{ Casimir morphisms} $\morph{\z[M\otimes M']:
\K}{M\otimes M'}$ and $\morph{\z[M'\otimes M] : \K}{M'\otimes M}$ 
(see subsection~\ref{subsubsec:basemutuellementduales}). These Casimir morphisms enjoy surprisingly nice properties and
are fundamental for a notion of skew permutation inside tensor
path algebras: the \emph{left permutation} and the
\emph{right permutation} of $\z[M\otimes M']$ coincide with
$\z[M'\otimes M]$ and reciprocally, the left permutation and the
right permutation of $\z[M'\otimes M]$ coincide with $\z[M\otimes
M']$, so that the complete cyclic permutation of each of above
Casimir morphisms stays invariant. Referring to  the last property
we say that each Casimir morphism $\z\in\set{\z[M\otimes M'],
\z[M'\otimes M]}$ is \emph{cyclically stable}.  Thanks to some crucial properties of Casimir morphisms, potentials for modulated quivers we easily defined as  \emph{morphisms of
$\K$-bimodules $\morph[/above={\mfr{m}}]{\K}{\kQc_{(2)}}$},
equivalently potentials  correspond to $\K$-central elements in $\kQc_{(2)}$.
For the second but complementary approach  of our study, we restrict to
\emph{symmetric potentials}: they can be obtained from  elements
of the central $\Zc(\K)$-bimodule $\kQc \otimes_{\K^e} \K$ where $\Zc(\K)$ is   the center of the algebra $\K$ and
$\K^e=\K\otimes_{\kk} \op{\K}$ is the enveloping $\kk$-algebra of $\K$.
Indeed, the ordinary cyclic permutation of cyclic tensor elements
from simply laced path algebras appears to be well-defined on
$\kQc  \otimes_{\K^e} \K$, and the manipulation of symmetric
potentials becomes less technical. In particular, if the $\kk$-algebra $\K$ is
separable over a ground field $\kk$, then by a result of Donald G.
Higman \cite{Higman1955a}, the Casimir ideal $\zc(\K)$ of $\K$ coincides with  the
center $\Zc(\K)$ of $\K$ and potentials
on $\mQ$ coincides with symmetric potentials, the latter also
holds when $\K$ is a symmetrically (or strongly) separable
algebra  over any commutative ring. The special treatment  of
 symmetric potentials in this work  is motivated  by a recent work of B. Keller on deformations
of Calabi-Yau differential graded categories and on Ginzburg
differential graded categories, in which the author considers
potentials in a path category $\cal{A}$ over a simply laced
 discrete category $\cal{R}$ as elements of
$\cal{A} \otimes_{\cal{R}^e} \cal{A}$.  We also point out that
when the Casimir ideal of $\K$ does not coincide with
$\Zc(\K)$, the class of Jacobian algebras obtained from skew
permutations and cyclic skew derivatives strictly contains the
class of Jacobian algebras obtained from symmetric potentials.

The next challenge  is to prove the following reduction process. Let 
$\mfr{m}=(\mfr{m}_l)_{l\geq 2}$ be a
potential on  $\mQ$, here  $\mfr{m}_{l} \in B^{(l)}=B\otimes B$ is the degree-$l$ component 
of $\mfr{m}$. We refer to $\mfr{m}_2$ as the \emph{trivial part of $\mfr{m}$}; the \emph{trivial part
$\triv{B}$} of $B$  is the image of $\mfr{m}_{2}$ under cyclic derivative, and it is
assumed that $\mfr{m}_{2} \in \triv{B}^{(2)}$;  the
\emph{reduced part of $B$} is $\red{B}:=\Frac{B}{\triv{B}}$.  Under some splitting conditions, we
have naturally induced symmetrizable dualizing pairs of bimodules
$\set{\triv{B},\triv{B},\beta}$ and $\set{\red{B}, \dual{\red{B}},
\bbar{\beta}}$,  yielding  a trivial modulated quiver with
potential $(\triv{\mQ},\mfr{m}_{2})$  and a modulated quiver
$\red{\mQ}$, where $\triv{B}$ is the arrow bimodule of
$\triv{\mQ}$ while $\red{B}$ is the arrow bimodule of $\red{\mQ}$.
We consider potentials $\mfr{m}$ such that $\mfr{m}_{2}$
appears as a Casimir morphism
$\mfr{m}_{2}=\triv{\mfr{m}}=\z[U\otimes V]$ with $\triv{B}=U\oplus V$, in this
case $\mfr{m}$ is called $2$-loop free and the pair $(\mQ,\mfr{m})$ is called a \emph{modulated quiver with potential}.  Now, the \emph{reduction
process} consists in constructing another modulated quiver with
potential $(\red{\mQ},\red{\mfr{m}})$  whose trivial part is zero and such that
 along some appropriate epimorphism of topological  path
algebras from $\kQc$ into $\what{\kk\red{\mQ}}$ (or along  some
appropriate automorphism of the topological path algebra $\kQc$),
the Jacobian algebra of $(\mQ,\mfr{m})$   coincides with the
Jacobian algebra of $(\red{\mQ},\red{\mfr{m}})$. Dealing with this
reduction problem, one of the technical and crucial point is to
prove that Jacobian ideals are preserved along some special
continuous isomorphisms of topological path algebras and to this
last end we must find a way to lift to the  framework of tensor
path algebras  the  "cyclic Leibniz rule" and
"chain-rule" used in \cite{DWZ} for the same purpose. For a simply-laced path algebra of a quiver $\Q$, the cyclic Leibniz rule is an
easy consequence of the existence of a canonical $\kk$-basis of
$\kk\Q$  induced  by the arrows of the quiver, relatively to
which, the ordinary cyclic permutation of cyclic elements in $\kk
\Q$ reduces to the cyclic permutation of arrows of $\Q$. However,
in the tensor path algebra $\kQ$, controlling cyclic skew
permutations of a given homogeneous potential is rather a complex
matter.

\ParIt{On the obstruction to the reduction}   We  must draw the attention of the anonymous   reader that the obstruction which arises when trying to reduce a modulated quiver with potential is of the same nature as the obstruction to the generalization of the  well-known Gabriel's theorem for presentation of finite dimensional algebras. Gabriel's theorem states that any finite dimensional algebra over an algebraically closed field admits a presentation by a quiver with relations; whereas  the non simply-laced analogue of this result states that any finite dimensional algebra $A$ over a field, with Jacobson radical $\J{A}$, admits a presentation by a modulated quiver with relations provided $A$ can be given a structure of an $(A/\J{A})$-bimodule such that  the inclusion
$\morph[/dir=right hook->]{\J[2]{A} }{\,\J{A}}$  splits as  morphism of $(A/\J{A}$)-bimodules.  The latter  splitting condition is satisfied if the ground  field is perfect. For an arbitrary  ring $\kk$,   the trivial
bimodule $\triv{B}$ needs not be a direct summand in the arrow
 bimodule $B$, the latter happens especially  when the symmetric
 enveloping  $\kk$-algebra $\K^e$ is not semisimple. Thus, if $\K$ is separable over a  field $\kk$, then, as a tensor product over a field  of  two separable $\kk$-algebras,  $\K^e$ is a separable $\kk$-algebra  and hence semisimple (see {\cite[Cor~11.6.8]{Cohn2003}}), in this case the obstruction to the reduction of modulated quivers with potentials
disappears exactly as in the case of presentation of finite dimensional algebras by modulated quivers with relations.

\subsection*{Description of main results and organization  of the paper}
The first main  result of this work is the reduction Theorem~\ref{theo.red-mQp}; it
establishes the reduction  up to \emph{weak
right-equivalences}. Focusing on symmetric potentials,  the reduction process is refined in Theorem~\ref{theo.sympot} under some natural splitting conditions, it can be obtained up to right-equivalences as in the simply-laced case. Now, whenever the
reduction is possible, for each central idempotent $e\in \K$
satisfying some natural condition, we   define the mutation
of a modulated quiver with potential at "$e$" up to weak
right-equivalences (or right-equivalences if the Casimir
ideal $\zc(\K)$ coincides with $\Zc(\K)$), and
Theorem~\ref{theo.mut-involutif} states  that the mutation at $e$ is a well-defined
involution on the set of (weak) right-equivalence classes of
modulated quivers with potentials.
Of a special interest, we  deduce (in Corollary~\ref{cor.sympot-perfectfield}) that in the setting of a separable algebra $\K$ over a field, all potentials are symmetric ones and the study of modulated quivers with potentials in such a context mimics
the simply-laced  study: cyclic (left or right) permutations are images
under a Casimir operator of corresponding ordinary permutations. This is  indeed a useful   information: for applications to cluster algebras theory, one  may  restrict a further non simply-laced study of modulated quivers with potentials
to the setting of a  perfect ground field where things behave smoothly. Next, considering graded modulated quiver with homogeneous potentials, we extend our mains results to  the graded context.

The paper is organized as follows. Section~\ref{sec.morptr} is
dedicated to  the main technical tools about symmetrizable
dualizing pairs, in Section~\ref{sec2:chp3}
we begin the discussion of  the general  approach
to potentials and cyclic derivatives.
Section~\ref{sec.reduction} deals with the
reduction problem in full generality and in
section~\ref{sec:sympot} we focus on symmetric potentials and
sharpen the main result from section~\ref{sec.reduction}. Examples  and illustrations of reduction are postponed to Section~\ref{secexample}. Then, after  discussing on mutations of modulated quivers with potentials in Section~\ref{sec:mutationsmqp}, more examples of sequences of mutations and reductions in the Dynkin type $\mb{F}_4$  appear in Section~\ref{sec:explmutationreduction}. In section~\ref{sec:gradedcontext} we consider graded modulated  quiver with potentials of homogeneous degree and provide arguments showing that the results of preceding sections hold in the graded context.  In guise of application, in the last section we introduce non simply-laced Ginzburg \dg-algebras and cluster categories associated with graded modulated quivers with potentials, generalizing the construction of cluster categories associated with graded quivers with potentials from \cite[Def~3.5]{Amiot2009} and  \cite[\S~4]{KY2011}. 

\subsection*{Some  perspectives}
In the present work, we  have not investigated rigid modulated
quivers with potentials and non-degeneracy of mutation as done in
\cite[\S6,7,8]{DWZ}; also we have not studied decorated
representations of modulated quivers with potentials. However, at least in the presence of separability over a base field or more specially in the setting  of a perfect ground field, we believe a general study of   mutations of
decorated representations of modulated quivers with potentials should be affordable. 
In order to  understand the cluster categories associated with modulated quiver with (nonzero) potentials, it is natural to prove  the following  about   non simply-laced  Ginzburg \dg-algebras.\\[0.1em]
\textbf{Conjecture}: Keller's result on the $3$-Calabi-Yau property of  simply-laced Ginzburg \dg-categories holds in the general framework,  at least when separability over a 
 field $\kk$ is assumed.

\subsubsection*{Conventions, matrix mutation and valued quiver mutation} 
 We let $\kk$  be a
 commutative ring and   $\K$   a   $\kk$-algebra  assumed to be  finitely
generated  projective as $\kk$-module. Tensor (path) algebras
occur as tensor algebras of  $\K$-bimodules $B$, with $B$ assumed
to be finitely generated and projective as  left $\K$-module and
  right  $\K$-module.  The tensor product $B\tenss[\K]B'$ of two
$\K$-bimodules $B$ and $B'$ is  denoted by $B\otimes B'$ or
simply by $B B'$. The composition of any two morphisms $\morph{f:
X}{Y}$ and $\morph{g: Y}{Z}$ in  a given category is written
either  as $g\circ f$, $g\mul f$ or as $g f$.
We shall sometimes   deal with infinite linear combinations which
 naturally occur:  thus each element
$x=(x_{\lambda})_{\lambda\in\Sigma}$ of a direct product
$\prd{\lambda\in\Sigma}{B_{\lambda}}$  of left or right
$\K$-modules   appears naturally  as an infinite sum
$x=\Som{\lambda\in\Sigma}{x_{\lambda}}$.

\ParIt{Matrix mutation} Let $n\in \N$ be a nonzero integer, we write $\ninterv{n}=\set{1,\dotsc,n}$.  Let 
$\B=(b_{i,j})_{1\leq i,j\leq n}$ be a matrix  with integer entries,  
 $\B$ is assumed  skew-symmetrizable, that is, there exists a  diagonal $n\times n$-matrix $\underline{\mbm{n}}=(\n[1],\dotsc,\n_n)$ of nonzero positive integers such that $b_{i,j}\n_j=-b_{j,i}\n_i$ for all $i,j \in \ninterv{n}$.  The
\emph{mutation of $\B$ at direction $k\in\ninterv{n}$}  is the  
skew-symmetrizable matrix $\B'=\mu_k(\B)$  described as  follows: 
define the common sign of  each pair  $a,b\in\mb{Z}$  by $\sign(a,b)=\sign(\sign(a)+\sign(b))$ where $\sign(0)=0$; then $\B'=(b_{i,j}')$ is  given by the
following mutation rule: if $k\in\set{i, j}$ then 
$b_{i,j}'=-b_{i,j}$, otherwise we have 
$b_{i,j}'=b_{i,j} + \sign(b_{i,k},b_{k,j}) b_{i,k}b_{k,j}$.

\ParIt{Valued quiver mutation} 
An arbitrary \emph{(locally finite) valued quiver} $\Q$ with valuation $\vlt$ consists
of a set of points $\I=\Q_0$, and disjoint finite sets $\Q_1(i,j)$ of valued arrows from $i$ to $j$,  where  the valuation of each $\alpha\in \Q_1(i,j)$ is a pair of natural numbers $\vlt(\alpha)=(\vl{\alpha},\vld{\alpha})$ and $\alpha$ may be pictured as $\edge[/above={\alpha}]{i}{j}$ or as
$\edge[/dist=2,/above={\vl{\alpha}\,,\,\vld{\alpha}}]{\alpha:i}{j}$; 
the valuation $\vlt$ is required to be right (or  left) symmetrizable, where the (minimal right) symmetrizing map $\morph[/above={\n}]{ \I}{ \N}$  for $\vlt$ prescribes for each  $i\in\I$ a non-zero integer $\ind{\n}{i}\in\N$ such that $\vl{\alpha}\ind{\n}{j}= \vld{\alpha}\ind{\n}{i}$ for all $\alpha\in\Q_1(i,j)$.  Arrows with valuation $(0,0)$ are referred to as \emph{$0$-valued} (or trivially valued) arrows,  they are normally not drawn in pictures, it is understood that a $0$-valued arrow is not counted among the  arrows of the valued quiver. For an integer $m\geq 2$, the valued quiver is \emph{$m$-acyclic} if it contains no $m$-cycle,  that is, a path of length $m$ of the form $\Seq{i_1\, \to[/above={\alpha_1}] i_2 \dotsm\!\dotsm i_{m-1}\, \to[/above={\alpha_{m-1}}]\,i_m \to[/above={\alpha_{m}}] \, i_1}$. The
composite of two paths $\omega\in\Q(i,j)$ and
$\omega'\in\Q(j,t)$  is their concatenation denoted by $\omega\omega'$ or $\omega\mul\omega'$.

For a valued quiver $\Q$ over a set points $\I$ with valuation $\vlt$, define its \emph{normal form}  as the valued quiver over $\I$ without parallel arrows,  with valuation still denoted by $\vlt$ and obtained from $\Q$ by replacing each  finite set \mbox{$\Q_1(i,j)=\set{\alpha_1,\dots,\alpha_m}$} by a one-element set consisting of a single valued arrow $\edge[/dist=2,/above={\val{i}{j}\,,\,\vald{i}{j}}]{\alpha: i}{j}$ with \mbox{$(\val{i}{j},\vald{i}{j})=\som{s=1}{m}{(\vl{\alpha_s},\vld{\alpha_s})}$}. Thus  a valued quiver is \emph{normalized} whenever it coincides with its normal form.
Let $\Q$  be a  normalized  valued quiver over  $\I$ with valuation $\vlt$,  then $\Q$ is completely defined by its set of points and its valuation.
 Let $k\in \I$ be  a point not lying on a $2$-cycle in $\Q$. The  \emph{mutation  of $\Q$ at  $k$} is the normalized valued quiver
$\Q'=\mu_k(\Q)$ over $\I$ with valuation $\vlt'$   as follows:
\begin{itemize}
\item[$\msf{(a)}$]  For any valued arrows
$\edge[/above={a,\,b}]{\alpha: x}{k}$ and
$\edge[/above={c,\,d}]{\beta: k}{y}$ starting or ending at $k$  in
$\Q$, there  are  corresponding  valued  arrows
$\edge[/dir=<-,/above={a,b}]{\dual{\alpha}: x}{k}$ and
$\edge[/dir=<-,/above={c,d}]{\dual{\beta}: k}{y}$  in  $\Q'$.
 \item[$\msf{(b)}$] For each pair  $ i,j\in\I\sminus\set{k}$ we have:
 $\valp{i}{j}=\!\max(\val{i}{k}\mul \val{k}{j}\!\!-\!\!\vald{j}{i},0)\!+\!
   \max(\val{i}{j}\!\!-\!\!\vald{j}{k}\mul\vald{k}{i},0)$, equivalently, 
\mbox{$\valpd{i}{j}=\!\max(\vald{i}{k}\mul \vald{k}{j}\!\!-\!\!\val{j}{i},0)
      \!+\!\max(\vald{i}{j}\!\!-\!\!\val{j}{k}\mul\val{k}{i},0)$}.
\end{itemize}

The above description of  mutation is canonical: we never add superfluous $2$-cycles. This contrasts with ordinary quiver mutation where superfluous $2$-cycles are added  and then,  some of them are  "simplified"  in a non canonical way.
By a little abuse of language, if $\Q$ and $\Q'$ are any valued quivers over a set $\I$, we still write $\Q'=\mu_k(\Q)$  if the normal form of $\Q'$ is the mutation at   $k$ of the normal form of $\Q$. Below is an illustration of valued quiver mutations, where the first two   are respectively the normal forms of the last two ones:
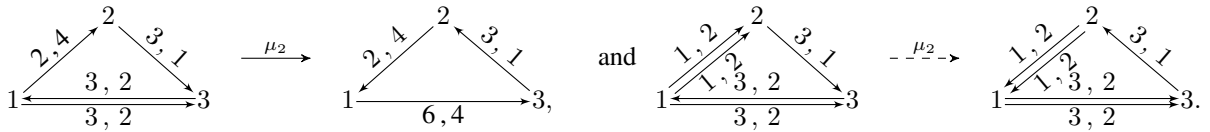
\begin{center}
$\myeqar[clcclc]{
\graph[/dist=2,/clsep=2.4]{{} \& 2 \& {} \\ 1 \& \& 3}{
\path[->]
(m-2-1) edge node[sloped,above=-0.08] {2\,,\,4} (m-1-2)
(m-1-2) edge node[sloped, above=-0.08] {$3\,,\,1$} (m-2-3)
($ (0,0.04)+(m-2-3.west) $) edge node[sloped,above=-0.08] {$3\,,\,2$} ($ (0,0.04)+(m-2-1.east) $)
($ (0,-0.04)+(m-2-1.east) $) edge node[sloped,below=-0.08] {$3\,,\,2$} ($ (0,-0.04)+(m-2-3.west) $)
} & \to[/above={\mu_2}] &
\graph[/dist=2,/clsep=2.4]{{} \& 2 \& {} \\ 1 \& \& 3,}{
\path[<-]
(m-2-1) edge node[sloped,above=-0.08] {$2\,,\,4$} (m-1-2)
(m-1-2) edge node[sloped, above=-0.08] {$3\,,\,1$} (m-2-3)
(m-2-3) edge node[below=-0.08] {6\,,\,4} (m-2-1)
} & \; \text{ and } \;
\graph[/dist=2,/clsep=2.4]{{} \& 2 \& {} \\ 1 \& \& 3}{
\path[->]
($ (0.04,-0.04)+(m-2-1.north east) $) edge node[sloped,below=-0.08] {$1\,,\,2$} ($ (0.04,-0.04)+(m-1-2.south west) $)
($ (-0.04,0.04)+(m-2-1.north east) $) edge node[sloped,above=-0.08] {$1\,,\,2$} ($ (-0.04,0.04)+(m-1-2.south west) $)
(m-1-2) edge node[sloped, above=-0.08] {$3\,,\,1$} (m-2-3)
($ (0,0.04)+(m-2-3.west) $) edge node[sloped,above=-0.08] {$3\,,\,2$} ($ (0,0.04)+(m-2-1.east) $)
($ (0,-0.04)+(m-2-1.east) $) edge node[sloped,below=-0.08] {$3\,,\,2$} ($ (0,-0.04)+(m-2-3.west) $)
} & \to[style=dashed,/above={\mu_2}] &
\graph[/dist=2,/clsep=2.4]{{} \& 2 \& {} \\ 1 \& \& 3.}{
\path[<-]
($ (0.04,-0.04)+(m-2-1.north east) $) edge node[sloped,below=-0.08] {$1\,,\,2$} ($ (0.04,-0.04)+(m-1-2.south west) $)
($ (-0.04,0.04)+(m-2-1.north east) $) edge node[sloped,above=-0.08] {$1\,,\,2$} ($ (-0.04,0.04)+(m-1-2.south west) $)
(m-1-2) edge node[sloped, above=-0.08] {$3\,,\,1$} (m-2-3)
($ (0,0.04)+(m-2-3.west) $) edge node[sloped,above=-0.08] {$3\,,\,2$} ($ (0,0.04)+(m-2-1.east) $)
($ (0,-0.04)+(m-2-1.east) $) edge[->] node[sloped,below=-0.08] {$3\,,\,2$} ($ (0,-0.04)+(m-2-3.west) $)
}
}$
\end{center}
We observe that  $2$-acyclic normalized valued quivers without loops  over a set $\I$ correspond
bijectively to skew-symmetrizable matrices with integer
coefficients indexed by $\I$, in such a way that valued quiver mutation   and matrix mutation agree:  let $\Q$ be a normalized  $2$-acyclic valued quiver without
loops over $\I$ with valuation $\vlt$, and  $(b_{i,j})_{i,j\in \I}$  the corresponding skew-symmetrizable matrix, then 
\mbox{$(b_{i,j},b_{j,i})=(\val{i}{j}-\vald{j}{i},\val{j}{i}-\vald{i}{j})$}.

\section{Trace maps and symmetrizable Dualizing pairs of
bimodules}\label{sec.morptr}

 We write $\Zc(\K)$ for the center of $\K$; the \emph{$\K$-center} $\Zck(B)$ of a
$\K$-bimodule $B$ is the $\Zc(\K)$-subbimodule   of $B$
consisting of all elements $x$ with $a x= xa$ for all $a$ in $\K$.
Recall that the \emph{left dual}
$\ldual{B}=\HM{\K}{\indd{\K}{B}{},\K}$, the \emph{$\kk$-dual
$\HM{\kk}{\indd{\K}{B}{\K},\kk}$} and the \emph{right dual}
$\rdual{B}=\HM{\K}{\indd{}{B}{\K},\K}$  of $B$ consist
respectively of left $\K$-linear maps, $\kk$-linear maps and right
$\K$-linear maps on $B$,
 with actions defined as follows: for $a,b\in \K$, $u\in
\ldual{B}$, $\xi\in\Hm_{\kk}(\indd{\K}{B}{\K},\kk)$ and $v\in
\rdual{B}$, we have $(a\mul u\mul b)(x)=u(x \mul a)\mul b$,
$(a\mul \xi\mul b)(x)=\xi(b\mul x\mul a)$ and $(a\mul v\mul
b)(x)=a\mul v(b\mul x)$ for every $x\in B$. The bimodule $B$ is 
\emph{dualizing}  if the left dual and the right dual of $B$ are
isomorphic. Recall that $\K$ is  \emph{Frobenius} if  there is an isomorphism
$\morph[/above={\sim}]{\phi:\K}{\Hm_{\kk}(\K,\kk)}$ of  left or right
$\K$-modules; if additionally  $\phi$ is  a $\K$-bimodule morphism  then $\K$ is called  a \emph{symmetric Frobenius} algebra or simply a \emph{symmetric algebra} and may be denoted  by $(\K,\tr)$ with $\tr=\phi(1)$. Symmetric algebras and traces are
related as in the following definition.
\begin{Defn}\label{defn.dualizing-bim}
\begin{enumerate}
\item[$\msf{(i)}$] A \emph{$\kk$-linear trace form} (or simply a
trace) on   $\K$ is any element $\tr$ in the $\K$-center of
$\HM{\kk}{\K,\kk}$: thus $\tr(a\mul b)=\tr(b\mul a)$ for all
$a,b\in \K$.  The \emph{radical}  of  $\tr$ is the  ideal
$R_{\tr}:=\set{a\in \K \; : \; \forall b\in \K, \tr(a b)=0}$, and
$\tr$ is   \emph{non-degenerate} if its radical is zero. When the induced $\K$-bimodule morphism  $\morph{\K}{\Hm_{\kk}(\K,\kk),\,a\mapsto (b\mapsto \tr(ab))}$ is an isomorphism,  $\tr$ is called \emph{strongly non-degenerate} and in
this case $(\K,\tr)$ is  \emph{symmetric}.
 \item[$\msf{(ii)}$]  The \emph{Casimir morphism} $\morph{\z[\K\otimes_{\kk}\K]:
\kk}{\K \otimes_{\kk} \HM{\kk}{\K,\kk} \cong \K\tens{\kk} \K}$ associated with each symmetric algebra $(\K,\tr)$
takes the unit of $\kk$ to a \emph{Casimir element}
$\Som{s\in\Lambda}{e_s \otimes \dual{e_s}} \in\Zck(\K\tens{\kk}
\K)$ characterized by the identities:
\vskip-1.5\baselineskip
\begin{equation}\label{casimir-K}
\text{for all } a\in \K,\; \quad \som{s\in\Lambda}{}{e_s
\tr(\dual{e_s}a)}= a = \som{s\in\Lambda}{}{\tr(ae_s) \dual{e_s}}.
\end{equation}
\end{enumerate}
\end{Defn}

In part $\msf{(ii)}$ above,  $\set{e_s\,: s\in\Lambda}$ is a
finite generating set for $\K$ over $\kk$ corresponding to an
epimorphism $\morph[/dir=->>,/above={p}]{\kk^{(\Lambda)}}{\K}$,
and since $\K$ is a projective $\kk$-module we
choose a right inverse
$\morph[/dir=>->,/above={q}]{\K}{\kk^{(\Lambda)}}$ for $p$
yielding a  generating  set $\set{\hat{e}_s\,: s\in\Lambda}$ for
the dual $\HM{\kk}{\K,\kk}$ to which corresponds a "dual
generating set" $\set{\dual{e}_s\,: s\in\Lambda} \subset \K$ with
$\tr(\dual{e}_s\mul \sdash) = \hat{e}_s : \K\to \kk: a\mapsto
\tr(\dual{e}_s a)=\hat{e}_s(a)$. Identities~\eqref{casimir-K} yields the following observation.

\begin{Rem} \label{rem.dualizing} Any  
$\K$-bimodule $B$ over symmetric algebra $(\K,\tr)$ is   dualizing:  the canonical maps
$\morph{\tr_{\lft}\hspace{-0.2em}=\hspace{-0.2em}\tr\hspace{-0.2em}\circ\hspace{-0.2em}
\sdash: \ldual{B}}{\HM{\kk}{B,\kk}}$ and
$\morph{\tr_{\rgt}\hspace{-0.2em}=\hspace{-0.2em}\tr\hspace{-0.2em}\circ\hspace{-0.2em}
\sdash: \rdual{B}}{\HM{\kk}{B,\kk}}$ are
 bimodule isomorphisms, and for all $v\in \HM{\kk}{B,\kk}$ we have: 
 $ \tr_{\lft}^{-1}(v):  x\mapsto
\Som{s\in\Lambda}{e_s v(\dual{e}_s x)}$ and
$ \tr_{\rgt}^{-1}(v): (x\mapsto
\Som{s\in\Lambda}{e_s v(x\dual{e}_s)})$.
\end{Rem}

We then introduce the first main tool
for the study of potentials in a general framework.
\begin{Defn} Let $B,B'$ be $\K$-bimodules together with a bimodule morphism 
$\morph{\bilf: B\otimes B' \oplus B'\otimes B}{\K}$  referred to
as the bilinear form.
\begin{itemize}
\item [$\msf{(a)}$] The data  $\set{B,B';\bilf}$ is  a
\emph{symmetrizable  pairing}  over $(\K,\tr)$ if properties
$\msf{(i)}$ and $\msf{(ii)}$  below hold.
\begin{itemize}
\item [$\msf{(i)}$] $(\K,\tr)$ is a symmetric algebra, and  $\tr$ is a symmetrizing
 trace for $\bilf$, that is, $\tr(\bilf(x\otimes x'))=\tr(\bilf(x'\otimes x))$
 for all $x\in B$ and $x'\in B'$.
\item [$\msf{(ii)}$] $\bilf$ is   non-degenerate, that is,  the \emph{adjoint morphisms}
 $\morph{ \bilf[{{}{1,\lft}}]:  B'}{\ldual{B}: x'\mapsto \bilf(\sdash\otimes x')
 }$ and $\morph{ \bilf[{{}{1,\rgt}}]: B}{\rdual{B'}: x\mapsto \bilf(x\otimes \sdash)
 }$ (or equivalently the \emph{adjoint morphisms}  $\morph{ \bilf[{{}{2,\rgt}}] : B'}{\rdual{B}: x'\mapsto \bilf(x'\otimes \sdash) }$ and $\morph{ \bilf[{{}{2,\lft}}] : B}{\ldual{B'}: x\mapsto \bilf(\sdash\otimes x)
 }$) are  injective.
\end{itemize}
\item[$\msf{(b)}$] The \emph{ordered data} $\set{B,B';\bilf}$
is  a \emph{symmetrizable  weakly dualizing pair} over $(\K,\tr)$  if 
$B$ is projective as  left and   right $\K$-module, conditions $\msf{(i)}\sdash\msf{(ii)}$ hold and  the adjoint morphism
 $\bilf[{{}{1,\lft}}]$ (or equivalently the adjoint morphism  $\bilf[{{}{2,\rgt}}]$) is  bijective. If in addition, $B$ (and thus $B'$) is finitely
generated  as  left and right $\K$-module, then the ordered
data $\set{B',B;\bilf}$
 is also a symmetrizable weakly dualizing pair over $(\K,\tr)$ and we call the (non ordered) data $\set{B,B';\bilf}$  a \emph{symmetrizable  dualizing pair} of bimodules, $\bilf$  strongly non-degenerate,  $B$ and $B'$ are called \emph{mutually dual}
 and we write: $B'=\dual{B}$ and $B=\dual{B'}$.
 \end{itemize}
\end{Defn}

Often in a weakly dualizing pair $\set{B, \dual{B};
\bilf}$ we shall omit to specify the bilinear form $\bilf$, in this case we write:
$ \scprd{x\otimes \xi}=\indd{B}{\scprd{x\otimes \xi}}{\dual{B}}=\bilf(x\otimes \xi)$ and
$ \scprd{\xi\otimes x}=\indd{\dual{B}}{\scprd{\xi\otimes x}}{B}=\bilf(\xi\otimes x)$ for all $x\in B$ and $\xi\in \dual{B}$.

Note that each symmetric algebra $(\K,\tr)$
gives rise to a natural symmetrizable dualizing pair 
$\set{\K,\K}$ with the bilinear form given by the
multiplication of $\K$.  We need the following lemma  which gives a
large class of symmetric algebras 
as well as  the existence of nonzero  traces for finite-dimensional local algebras.

\begin{Lem}\label{lem.trace-dualizing-pair}
\begin{itemize}
\item [$\msf{(a)}$] Let $\K$ be any $\kk$-algebra. If $\K$
has a  non-degenerate $\kk$-linear trace $\tr$, then the  $\Zc(\K)$-module
$\Zck(\Hm_{\kk}(\K,\kk))$  is  free  of dimension
one and  each  non-degenerate trace on $\K$ is  given by $c\mul \tr$
where $c\in \Zc(\K)$ is a central unit.
\item [$\msf{(b)}$] Suppose $\K$ is finite-dimensional over a field $\kk$, with Jacobson radical $\J{\K}$. Then there is a nonzero trace
\mbox{$\tr\in \soc(\indd{\K}{\HM{\kk}{\K,\kk}}{})\cap
\soc(\ind{\HM{\kk}{\K,\kk}}{\K})  $}.  The $\kk$-algebra $\bbar{\K}=\K/\J{\K}$  is symmetric  and each $\bbar{\K}$-bimodule $M$, finite-dimensional over $\kk$, is part of a symmetrizable dualizing pair  $\set{M,\dual{M};\bilf}$.
\end{itemize}
\end{Lem}

\begin{prv} \Par{$\msf{(a)}$} Suppose $\tr$ is a non-degenerate trace on
$\K$ and let $\tau$ be any trace on $\K$, since clearly the dual
$\HM{\kk}{\K,\kk}$ is a free left $\K$-module of dimension one;
there exists $c\in\K$ such that $\tau = c\mul \tr$. We must show
that $c\in\Zc(\K)$, thus let $a,b\in\K$: we have
$\tr(cab)=\tau(ab)=\tau(ba)=\tr(cba)=\tr(acb)$, thus
$\tr((ca-ac)b)=0$ for all $b\in\K$, so that $ca -ac \in
R_{\tr}=0$, hence $c\in\Zc(\K)$. Now suppose $\tau=c\mul\tr$ is
also non-degenerate, then we must also have $\tr=c'\tau$ for some
$c'\in\Zc(\K)$, so that $\tr=c'c\tr$ and $\tau=cc'\tau$, yielding
that $c'c=1=cc'$, this  proves part $\msf{(a)}$ of  the lemma.
\Par{$\msf{(b)}$} It is a standard result that finite dimensional
simple  algebras over a field and hence semisimple algebras are
symmetric Frobenius algebras. This can be done by invoking the
existence of the so-called reduced trace for simple algebras which
are finite-dimensional over their centers. Hence the semisimple $\kk$-algebra $\bbar{\K}$ is  symmetric for some trace $\bbar{\tr}$, and if $\morph[/dir=->>]{\pi:\,\K}{\bbar{\K}}$  is the natural projection,  then   we  get a $\kk$-linear trace $\tr=\bbar{\tr}\circ \pi$ for $\K$ with radical $R_{\tr}=\J{\K}$ and with   \mbox{$\tr \in \parr{ \soc(\indd{\K}{\HM{\kk}{\K,\kk}}{})\cap \soc(\ind{\HM{\kk}{\K,\kk}}{\K}) } $}.
The rest of the proof follows from Remark~\ref{rem.dualizing} together with the observation that the bilinear form associated with a symmetrizable dualizing pair  $\set{M,\dual{M},\bilf}$ is  induced by the corresponding non-degenerate trace form: indeed assume   $\dual{B}$ is a $\bbar{\K}$-bimodule isomorphic to one of (and thus to all) the standard duals of a  $\bbar{\K}$-bimodule $B$,  finite-dimensional over $\kk$, choose an isomorphism $\morph[/isomark={\sim}]{\phi: \dual{B}}{\HM{\kk}{B,\kk}}$; then $\phi$ yields a symmetrizable dualizing pair $\set{B,\dual{B},\bilf}$ over $(\bbar{\K},\tr)$ with  $\bilf$ given as follows:  let $x\in B, u\in\dual{B}$,  by Remark~\ref{rem.dualizing}  write  $\tr\circ u_1=\phi(u)=\tr\circ u_2$  with $u_1\in\ldual{B}$ and $u_2\in\rdual{B}$, 
then $\bilf(x\otimes u)= u_1(x)$ and $\bilf(u\otimes x)=u_2(x)$.
\end{prv}

\subsubsection*{Mutually dual projective bases and Casimir elements}
\label{subsubsec:basemutuellementduales} Assume $B$ is  part of a
symmetrizable weakly dualizing pair 
$\set{B,\dual{B};\bilf}$ over   $(\K,\tr)$. Choose a split sequence
$\Seq{\indd{\K}{B}{}\, \to[/dir=>->, /above={\pi'}]\,
\K^{(\mfr{p})}\, \to[/dir=->>,/above={\pi}] \, \indd{\K}{B}{}}$
for the left $\K$-module $B$, where $\pi$ is a split epimorphism
with right inverse $\pi'$, $\K^{(\mfr{p})}$ is  a direct sum of
copies of $\K$ indexed by a (possibly infinite) cardinal
$\mfr{p}$. We get  a \emph{left
projective basis} $(\set{x_s: \, s\in \mfr{p}}, \set{\what{x}_s,\,
s\in\mfr{p}})$ for $B$ characterized by the following property:
\mbox{$x=\Som{s\in\mfr{p}}{\what{x}_s(x)x_s}$} for all $x\in B$,
and since moreover for each $u\in \ldual{B}$, the map
$\Seq{x\mapsto \Som{s\in\mfr{p}}{(\what{x}_s u(x_s))(x)} =
\Som{s\in\mfr{p}}{(\what{x}_s(x) u(x_s))} =u(x)}$  is  a well-defined element of  $\ldual{B}$,
 $u$  naturally  occurs as  (possibly infinite) sum:
 $u=\Som{s\in\mfr{p}}{\what{x}_s u(x_s)}$. 
 We  refer to the (possible infinite) sum $\Som{s\in\mfr{p}}{\what{x}_s\otimes x_s}$  as the
 \emph{Casimir element} associated with the left projective $\K$-module $B$
 and its left dual. Thus, if $\mfr{p}$ is a finite cardinal, then under the natural isomorphism
  $\morph[/above={\sim}]{\phi:
\ldual{B}\otimes B}{\HM{\K}{\indd{\K}{B}{},\indd{\K}{B}{}}}: \,
\phi(u\otimes x)(z)=u(z)x $ (induced by the adjunction of tensor
product), the pre-image of the identity map is given by the
\emph{Casimir element}. Similarly,  the \emph{Casimir element}
$\Som{s\in\mfr{q}}{y_s\otimes \what{y}_s}$ and  the right
projective basis $( \set{y_s,:\, s\in\mfr{q}},\set{\what{y}_s:\, s
\in \mfr{q}})$ associated with the right projective $\K$-module
$B$ and its right dual  have  the following characterizing
property:{\centering{ $x=\Som{s\in\mfr{q}}{y_s\what{y}_s(x)}$}}
and {\centering{$u=\Som{s\in\mfr{q}}{u(y_s)\what{y}_s}$ for all
$x\in B$ and $u\in \rdual{B}$,}} and when $\mfr{q}$ is a finite
cardinal, the Casimir element associated with $\ind{B}{\K}$  is
the pre-image of the identity map under the natural isomorphisms
 $\morph[/above={\sim}]{\psi:B \otimes \rdual{B}}{\HM{\K}{\ind{B}{\K},\ind{B}{\K}}}: \, \psi(y\otimes
v)(z)=y v(z) $. Using the adjoint isomorphisms
$\morph{ \bilf[{{}{1,\lft}}]:  B'}{\ldual{B}
 }$  and $\morph{ \bilf[{{}{2,\rgt}}] : B'}{\rdual{B}}$, we get two pairs $(\set{x_s:\, s\in \mfr{p}},
\set{\dual{x_s} :\, s\in \mfr{p} })$ and  $(\set{y_s:\, s\in
\mfr{q}}, \set{\dual{y_s} :\, s\in \mfr{q} })$  of a left
projective basis and a right projective basis  associated with $B$
and its weak dual $\dual{B}$, having the following characterizing
identities where $x\in B,  \xi\in \dual{B}$ and the
formula expressing each $\xi$  may (naturally) appears as an
infinite sums:
\begin{equation}\label{eq.dualbases}
 \Som{s\in\mfr{p}}{\bilf(x \otimes\dual{x_s}
)x_s}=x=\Som{r\in\mfr{q}}{y_r\bilf(\dual{y_r}\otimes x )} \text{
and } \Som{s\in\mfr{p}}{\dual{x_s}\bilf(x_s \otimes
\xi)}=\xi=\Som{s\in\mfr{q}}{\bilf(\xi\otimes y_r )\dual{y_r}}.
\end{equation}
The  "elements"  $\z[\dual{B}\otimes
B] =\Som{s\in\mfr{p}}{\dual{x_s} \otimes x_s}$ and $\z[B\otimes
\dual{B}] =\Som{s\in\mfr{p}}{y_s\otimes \dual{y_s}}$ are again 
referred to as \emph{Casimir elements} associated with $\bilf$ (or
with the pair $\set{B,\dual{B};\bilf}$), note in view of 
 equations~\eqref{eq.dualbases}  that these Casimir elements are $\K$-central.

Now suppose we are given  two symmetrizable  pairing   $\set{M,\dual{M};\beta}$ and
$\set{M',\dual{M'};\beta'}$ respectively over $(\K, \tau)$ and $(\K, \tau')$. Then  for a $\kk$-linear map $\morph{f: M}{M'}$,  its \emph{left dual} $\ldual{f}$ and  its
\emph{right dual} $\rdual{f}$ (if they exist) are the
unique $\kk$-linear maps $\morph{\ldual{f},\, \rdual{f}:
\dual{M'}}{\dual{M}}$ defined by the condition:
$
\beta'(f(\sdash)\otimes \sdash)= \beta(\sdash\otimes
\ldual{f}(\sdash))$ and $\beta'(\sdash\otimes
f(\sdash))=\beta(\rdual{f}(\sdash)\otimes \sdash).$\; 
We say that $f$  is \emph{left dualizing} (respectively,
\emph{right dualizing}) when $\ldual{f}$ (respectively,
$\rdual{f}$) exists. Note that, when they exist,
$\ldual{f}$ and $\rdual{f}$ need not coincide if the symmetrizable
requirement on our pairing of bimodules is omitted. \emph{$f$ is  dualizing } if $\ldual{f}$ and $\rdual{f}$  exist and coincide, in this
case their common value  $\dual{f}$ is called the
\emph{dual} of $f$.

\begin{Lem}\label{lem.dualizing-pair}
 Let $\set{B,\dual{B}; \bilf}$ and $\set{B',\dual{B'}; \bilf'}$
  be  symmetrizable pairing  over $(\K,\tr)$, and $f\in \HM{\kk}{B,B}$.
\begin{itemize}
\item [$\msf{(1)}$] If $f$ is left dualizing then  it is  left $\K$-linear and $\ldual{f}$ is right $\K$-linear;  if $f$ is right dualizing then it is right $\K$-linear and $\rdual{f}$ is left $\K$-linear. If $f$ is a  $\K$-bimodule morphism, then it is dualizing whenever $f$ admits a left  or  right dual, in this   case  the dual of $f$ is the
unique $\K$-bimodule morphism $\morph[/dist=1.2]{\dual{f}:\dual{B'}}{\dual{B}}$   with  the following property:
\begin{equation} \label{eq.dual-morph}
 \bilf'(f(\sdash)\otimes \sdash)=\bilf(\sdash\otimes \dual{f}(\sdash))
  \;\; \text{ or equivalently  }
\;\; \bilf'(\sdash\otimes f(\sdash))=\bilf(\dual{f}(\sdash)\otimes
\sdash).
\end{equation}
\item [$\msf{(2)}$] If the data $\set{B,\dual{B}; \bilf}$
is weakly dualizing, then any
morphism $\morph[/dist=1.2]{f:B}{B'}$  of left $\K$-modules
(respectively, right $\K$-modules, $\K$-bimodules) is left
dualizing  (respectively, right dualizing, dualizing).
\end{itemize}
\end{Lem}

 \begin{prv}  For  part $\msf{(1)}$, simply apply the fact that the bilinear forms $\bilf$ and $\bilf'$ are non-degenerate and symmetrizable via the same non-degenerate trace map $\tr$.
In part $\msf{(2)}$,  the  
 the ordered data $\set{B,\dual{B}; \bilf}$  is  weakly dualizing   over $(\K,\tr)$, so that we have   adjoint bimodule isomorphisms
$\morph{ \bilf[{{}{1,\lft}}]:  B'}{\ldual{B}
 }$  and $\morph{ \bilf[{{}{2,\rgt}}] : B'}{\rdual{B}}$. Thus, when
$\morph[/dist=1.2]{f:B}{B'}$ is left $\K$-linear,
the composition map along the sequence
\begin{center}
$\Seq{ \dual{B'}\, \to[/above={\bilf[{{}{1,\lft}}]'}] \, \ldual{B'} \,
 \to[/dist=4,/above={\HM{\K}{f,\K}}]\, \ldual{B}\, \to[/above={\bilff[{{}{1,\lft}{-1}}]},
/below={\sim}] \dual{B} }$,
\end{center} is
clearly a left dual for $f$.  Similarly,  if $f$ is
right $\K$-linear then its  right dual exists. When $f$ is a bimodule
morphism, part $\msf{(1)}$ and the previous arguments prove that $f$ is  dualizing.
\end{prv}

\begin{Lem}\label{lem.dualizing-pair-Casimir}
Let $\set{B,\dual{B}; \bilf}$ and $\set{B',\dual{B'}; \bilf'}$
  be  symmetrizable dualizing pairs over $(\K,\tr)$. Then the left dual
of any left $\K$-linear  isomorphism 
$\morph[/dist=1.2]{f:B}{B'}$  yields:
$(\ldual{f^{-1}}\otimes f)(\z[\dual{B}\otimes B]) =
\z[\dual{B'}\otimes B']$. Dually, the right dual of any right $\K$-linear isomorphism $\morph[/dist=1.2]{f' : B}{B'}$ yields: $(f'\otimes \rdual{(f'^{-1})})(\z[B\otimes \dual{B}])
= \z[B'\otimes \dual{B'}]$.
\end{Lem}

\begin{prv} Let $\morph[/dist=1.2]{f:B}{B'}$ be an isomorphism of left
$\K$-modules, then in view the last part of
Lemma~\ref{lem.dualizing-pair}, $f$ and $f^{-1}$
 are left dualizing and clearly  $\ldual{f^{-1}}=(\ldual{f})^{-1}$.
 Write $\z[\dual{B}\otimes B]=\Som[n]{s=1}\dual{x_s}\otimes x_s$ for the
Casimir element in $\dual{B}\otimes B$. Now let $x'\in B, \xi'\in
\dual{B'}$, applying the characterizing properties
\eqref{eq.dualbases} for the Casimir element $\z[\dual{B}\otimes
B]$ and the definition of the left dual $\ldual{f^{-1}}$  we have:
\\
$\begin{array}{rl} x'    = & f f^{-1}(x') = f\bigl( \Som[n]{s=1}
\bilf(f^{-1}(x')\otimes \dual{x_s}) x_s \bigr) = \Som[n]{s=1}
\bilf(f^{-1}(x')\otimes \dual{x_s}) f(x_s) 
  =  \Som[n]{s=1} \bilf'(x'\otimes \ldual{f^{-1}}(\dual{x_s})) f(x_s);  \\
 \xi'   = & \ldual{f^{-1}}\ldual{f}(\xi')=
  \ldual{f^{-1}}\bigl( \Som[n]{s=1}\dual{x_s}\bilf(x_s\otimes \ldual{f}(\xi') ) \bigr)
  =  \Som[n]{s=1}\ldual{f^{-1}}(\dual{x_s})\bilf(x_s\otimes \ldual{f}(\xi') )  
  =    \Som[n]{s=1}\ldual{f^{-1}}(\dual{x_s})\bilf'(f(x_s)\otimes \xi' ),  \\
\end{array}$ \\ showing in virtue of  \eqref{eq.dualbases}
  that the element $(\ldual{f^{-1}}\otimes f)(\z[\dual{B}\otimes B])=
 \Som[n]{s=1}\ldual{f^{-1}}(\dual{x_s})\otimes f(x_s)$ is as
 claimed the Casimir element  in $\dual{B'}\otimes B'$.  The dual statement is proved in  the
 same way.
\end{prv}

\subsubsection*{Products of symmetrizable dualizing pairs}\label{prd-dual-pairs}  First,  note that if  $B$ and $B'$ are $\K$-bimodules, projective and finitely
generated as   left and  right $\K$-modules,  then their  tensor product
$B\otimes B'$ (over $\K$) is till finitely generated projective as  left and
as  right  $\K$-module. Suppose $\set{B,\dual{B}; \bilf}$ and
$\set{B',\dual{B'}; \bilf'}$  are   symmetrizable dualizing
pairs   over $(\K,\tr)$. We can form  the  product
$\set{B, \dual{B}; \bilf} \otimes \set{B', \dual{B'}; \bilf'} :=
\set{B\otimes B', \dual{B'}\otimes\dual{B}; \bilf\ast\bilf'}$ 
with the induced bilinear form  $\bilf\ast\bilf'$: for  $x\in B$, $x'\in B'$, $u\in \dual{B}$ and $u'\in \dual{B'}$ we have $(\bilf\ast\bilf')(x\otimes x'\otimes u'\otimes
u)=\bilf(x\bilf'(x'\otimes u')\otimes u)$ (and thus
$(\bilf\ast\bilf')(u'\otimes u\otimes x\otimes
x')=\bilf'(u'\bilf(u\otimes x)\otimes x')$). One defines
in the same way the product of any finite number of dualizing
pairs of bimodules. Observe that the pair
$\set{B,\dual{B}; \bilf}$ also  induces two symmetrizable
dualizing pairs  $\set{B\otimes\dual{B},B\otimes
\dual{B}}$  and  $\set{\dual{B}\otimes B,\dual{B}\otimes B}$ in
which $B\otimes\dual{B}$ and $\dual{B}\otimes B$  are self-dual
bimodules.  The next   lemma gives a simple but crucial
observation.
\begin{Lem}\label{lem.central-elt}
 \begin{itemize}
\item[$\msf{(1)}$] For a symmetrizable dualizing pair  $\set{B,\dual{B}; \bilf}$ over $(\K,\tr)$, the dual morphisms of the bilinear forms
$\morph{\bilf[{{}{1}}]:B\otimes \dual{B}}{\K}$ and  $\morph{\bilf[{{}{2}}]: \dual{B}\otimes
B}{\K}$, with $\bilf[{{}{1}}](x\otimes u)=\bilf(x\otimes u)$ and $\bilf[{{}{2}}](u\otimes x)=\bilf(u\otimes x)$ for all $x\in B, u\in \dual{B}$, coincide with the  Casimir
morphisms $\morph{\z[B\otimes\dual{B}]:\K}{B\otimes\dual{B}}$ and
$\morph{\z[\dual{B}\otimes B]:\K}{\dual{B}}$  taking the unit
element of $\K$ to the corresponding Casimir elements.
\item[$\msf{(2)}$] Suppose $\set{B,\dual{B}}\otimes
\set{B',\dual{B'}}=\set{B\otimes B',\dual{B'}\otimes \dual{B}}$ is
the product of  symmetrizable dualizing pairs over $(\K,\tr)$. Then
the  corresponding  Casimir elements    are given by
$\z[{(B\otimes B')\otimes(\dual{B'}\otimes
\dual{B})}]=\som{i=1}{q}{\som{j=1}{q'}{(y_i\otimes y_j')\otimes
(\dual{y'}_j\otimes \dual{y}_i)}}$ and  $\z[{(\dual{B'}\otimes
\dual{B})\otimes (B\otimes
B')}]=\som{s=1}{p}{\som{t=1}{p'}{(\dual{x'}_t\otimes \dual{x}_s)
\otimes (x_s\otimes x_t')}}$,  where  $\z[B\otimes
\dual{B}]=\som{i=1}{q}{y_i\otimes\dual{y_i}}$, $\z[\dual{B}\otimes
B]=\som{s=1}{p}{\dual{x_s}\otimes x_s}$, $\z[B'\otimes
\dual{B'}]=\som{t=1}{q'}{y_t'\otimes\dual{y_t'}}$  and
$\z[\dual{B'}\otimes B']=\som{j=1}{p'}{\dual{x_j'}\otimes x_j'}$.
\end{itemize}
\end{Lem}
\begin{prv} The proof is a direct application of the definition of Casimir elements and the
dual  of a morphism.
\end{prv}

\ParIt{Derivative operators} For a symmetrizable dualizing pair 
$\set{M,\dual{M};\bilf}$ over $(\K,\tr)$, let $A:=\msf{T}_{\K}(M)=\K\oplus M\oplus
(M\otimes M) \oplus (M\otimes M \otimes M)\oplus \ldots$ be the
tensor algebra of the $\K$-bimodule $M$, then  write
\begin{center}
$\morph{\lderv=\lderv[\dual{M}]:=\bilf\otimes \id: \dual{M}\otimes
M\otimes A }{A}$ and
$\morph{\rderv=\rderv[\dual{M}]:=\id\otimes\bilf: A\otimes
M\otimes \dual{M} }{A}$,
\end{center}
respectively referred to as \emph{left derivative operator} and \emph{right
derivative operator}.  We now conclude this  subsection with a property of a \emph{cyclical stability}.
\begin{Lem}\label{lem.cyc-stable} Let $\set{M,\dual{M};\bilf}$ be
 a symmetrizable dualizing pair over $(\K,\tr)$, let
 $\z=\z[M\otimes\dual{M}]$ and  $\z'=\z[\dual{M}\otimes M]$.
\begin{itemize}
\item[$\msf{(1)}$] For every bimodule morphism
$\morph[/above={\mfr{m}}]{\K}{M}$, we have:
 $\bilf(\mfr{m}(1)\otimes \sdash)=\dual{\mfr{m}}=\bilf(\sdash \otimes
\mfr{m}(1))$, and  $\mfr{m}$ is cyclically stable, that is, the bimodule morphisms $\rperm \mfr{m}:=(\id_{M}\otimes
\bilf)\circ (\id_{M}\otimes \mfr{m}\otimes \id_{\dual{M}})\circ
\z[M\otimes\dual{M}]$ and  $\lperm \mfr{m}:=(\bilf\otimes
\id_{M})\circ (\id_{\dual{M}}\otimes \mfr{m}\otimes \id_{M})\circ
\z[\dual{M}\otimes M]$ coincide with $\mfr{m}$. 
\item[$\msf{(2)}$] Consider the following morphisms referred
to as \emph{left or right  permutations} of $\z$  and  $\z'$:
 \begin{equation}
 \hspace{-2.1em}{
 \left.\begin{aligned}
 &\lperm \z\!=\!\lderv[\dual{M}](\id_{\dual{M}}\otimes \z\otimes \id_{M})\circ\z',\;
 \rperm \z\!=\!\rderv[M](\id_{\dual{M}}\otimes \z\otimes\id_{M})\circ\z', \\
 &  \lperm \z'\!=\!\lderv[M](\id_{M}\otimes \z'\otimes
\id_{\dual{M}})\circ\z, \text{ and } \rperm
\z'\!=\!\rderv[\dual{M}](\id_{M}\otimes \z'\otimes
\id_{\dual{M}})\circ\z.
      \end{aligned} \right.\!}
\end{equation}
Then $\z$ and $\z'$ are \emph{cyclically equivalent}: \mbox{$\lperm \z'=\z=\rperm \z'$} and \mbox{$\lperm \z=\z'=\rperm \z$}. 
In particular $\z$ and $\z'$ are cyclically stable:
$\lperm[2]\z:=\lperm(\lperm \z)= \z =\rperm[2]\z :=\rperm(\rperm
\z)$ and
 $\lperm[2]\z':=\lperm(\lperm \z')=\z'=\rperm[2]\z':=\rperm(\rperm
 \z')$.
\end{itemize}
\end{Lem}

\begin{prv} For part $\msf{(1)}$,
let  $\morph[/above={\mfr{m}}]{\E}{M}$ be a morphism of bimodules,
in respect to the data
$\set{M,\dual{M};\bilf}$ and $\set{\K,\K}$, and in view of
Lemma~\ref{lem.dualizing-pair}  the dual $\morph{\dual{\mfr{m}}:
\dual{M}}{\K}$ of $\mfr{m}$ exists  and satisfies the following
relation: $\bilf(\mfr{m}(1)\otimes \sdash
)=\indd{\K}{\scprd{1\otimes
\dual{\mfr{m}}(\sdash)}}{\K}=\dual{\mfr{m}}=\indd{\K}{\scprd{\dual{\mfr{m}}
(\sdash) \otimes 1}}{\K}=\bilf(\sdash \otimes \mfr{m}(1))$. In view of
Lemma~\ref{lem.central-elt}, the Casimir morphism $\z$,
sending the unit of $\K$ to the Casimir element
$\z(1)=\som{t=1}{q}{y_t\otimes\dual{y_t}}$, is
 the dual of  the bilinear form  $\morph{B\otimes\dual{B}}{\K}$.
Thus,  using the relation
 $\bilf(\mfr{m}(1)\otimes \sdash )=\bilf(\sdash \otimes
\mfr{m}(1))$ we get: $\rperm\mfr{m}(1)=(\id_{M}\otimes \bilf)\circ
(\id_{M}\otimes \mfr{m}\otimes
\id_{\dual{M}})(\som{t=1}{q}{y_t\otimes\dual{y_t}})=
\som{t=1}{q}{y_t \mul \bilf(\mfr{m}(1)\otimes
\dual{y_t})}=\som{t=1}{q}{y_t \mul \bilf(\dual{y_t} \otimes
\mfr{m}(1))}=\mfr{m}(1)$,  where the last equality follows by the
definition of projective bases and  Casimir elements, see
equations~\eqref{eq.dualbases}. Similarly,  one shows
that $\lperm \mfr{m}:=(\bilf\otimes \id_{M})\circ
(\id_{\dual{M}}\otimes \mfr{m}\otimes \id_{M})\circ
\z[\dual{M}\otimes M]=\mfr{m}$.

Part $\msf{(2)}$ is a direct  application of the definition of
Casimir morphisms and identities \eqref{eq.dualbases}. Indeed
write $\z'(1)=\som{k=1}{p}{\dual{x}_k\otimes
x_k}\in\dual{M}\otimes M$. We then have
$(\lperm\z')(1):=\parr{\lderv(\id
\otimes\z'\otimes\id)\circ\z}(1)=
\lderv\parr{\som{s=1}{q}{\som{k=1}{p}{y_s\otimes \dual{x}_k\otimes
x_k\otimes \dual{y}_s}}}
=\som{s=1}{q}{(\som{k=1}{p}{\bilf(y_s\otimes \dual{x}_k)\mul
x_k})\otimes \dual{y}_s}=\som{s=1}{q}{y_s\otimes
\dual{y}_s}=\z(1)$,
 hence $\lperm\z'=\z$. Similarly we have
  $(\rperm\z')(1)=\parr{\rderv(\id
\otimes\z'\otimes\id)\circ\z}(1)=
\rderv\parr{\som{s=1}{q}{\som{k=1}{p}{y_s\otimes \dual{x}_k\otimes
x_k\otimes
\dual{y}_s}}}=\som{s=1}{q}{y_s(\som{k=1}{p}{\dual{x}_k\mul\bilf(x_k\otimes
\dual{y}_s)})}=\som{s=1}{q}{y_s\otimes
  \dual{y}_s}=\z(1)$, hence $\rperm\z'=\z$.
In the same way one checks that $\lperm\z=\z'=\rperm\z$. 
\end{prv}

\section{Potentials  and Jacobian algebras}\label{sec2:chp3}

\subsection{Tensor path algebras of modulated quivers}
\label{sec:pathAlgmQ}

\begin{Defn}\label{def.modQ} Given a symmetrizable
dualizing pair  $\set{B,\dual{B};\bilf}$ over a symmetric algebra $(\K,\tr)$, we refer to the data  $\mQ=(B,\K,\tr)$ as  \emph{$\kk$-modulated quiver}, \index{\emph{modulated
quiver}}  $B$ is therefore referred to  as the \emph{arrow bimodule} of $\mQ$.  The \emph{dual}
of $\mQ$ is the modulated quiver $\dual{\mQ}=(\dual{B},\K,\tr)$.
\end{Defn}

Note that we may decompose $\K$ as finite direct product
$\prd{i\in\I}{\kk_i}$ of indecomposable $\kk$-algebras $\kk_i$, each $\kk_i$
is viewed as  subalgebra in $\K$ and the unity of $\K$
is $\idl=\som{i\in\I}{}{\idl[i]}$ where $\idl[i] \in \kk_i$ is the
unity of $\kk_i$, the set $\set{\idl[i]\;:\;
i\in\I}$ is then a system of central primitive orthogonal idempotents for $\K$. Thus, $(\K,\tr)$ occurs as direct product of symmetric algebras $(\kk_i,\tr_i)$ with $\tr_i=\tr\restr{\kk_i}$, $i\in \I$, and the pair $\set{B,\dual{B};\bilf}$  occurs as direct sum of
 symmetrizable dualizing pairs 
$\set{\Bmd{i}{j},\Bmdd{i}{j};\bilf[{{i}{j}}]}$ over  $(\kk_i\times\kk_j,(\tr_i,\tr_j))$ with $\Bmd{i}{j}=\idl[i]\mul B \mul \idl[i]$.  When
$\Bmd{i}{j}$ is nonzero, we have an arrow from $i$ to $j$ in $\mQ$, pictured  as $\edge[/dist=3,/above={\Bmd{i}{j}\, ,\, \Bmdd{i}{j}}]{\alpha_{i,j}: i}{j}$  or simply as  $\edge[/dist=3,/above={\Bmd{i}{j}}]{\alpha_{i,j}: i}{j}$. We say that \emph{$\mQ$ has no loop} if $\Bmd{i}{i}$ is zero
for all $i\in\I$,  we also write $\mQ_m(i,j)$ for the set of all
length-$m$ paths from $i$ to $j$, while $\mQ(i,j)$ denotes  the
set of all paths from $i$ to $j$. Observe that, if moreover
each $\kk$-algebra $\kk_i$ is a division algebra then the case of
classical modulated quivers  is recovered and the
\emph{underlying (normalized) valued quiver} of $\mQ$ over $\I$ has valued arrows
$\edge[/dist=2.2,/above={\val{i}{j}\, ,\, \vald{i}{j}}]{i}{j}$  with
$\val{i}{j}=\dim_{\kk_j}(\Bmd{i}{j})$ and
$\vald{i}{j}=\dim_{\kk_i}(\Bmd{i}{j})$  for all $i,j\in\I$.

In the sequel we assume that  $\mQ:=(B,\K,\tr)$ is a $\kk$-modulated quiver. The
\emph{tensor path algebra} $\kQ$ of $\mQ$  is by definition the tensor  algebra of
the $\K$-bimodule $B$, thus \mbox{$\kQ:=\msf{T}_{\K}(B)=\somd{m\geq 0}{\kk \ind{\mQ}{m}}$},
where $\kk\ind{\mQ}{m}=\indu{B}{}{\sscr{(m)}}$ is the $m$-fold tensor product of  $B$, referred to as the \emph{bimodule of degree-$m$ homogeneous elements} (or the
\emph{bimodule of all length-$m$ paths}), with
$\indu{B}{}{\sscr{(0)}}=\K$, \mbox{$\indu{B}{}{\sscr{(1)}}=B$} and
$\indu{B}{}{\sscr{(m+1)}}=\indu{B}{}{\sscr{(m)}}\otimes B$ for all
$m\geq 1$. We let $\kk\ind{\mQ}{(t)}=\somd{m\geq
t}{\kk\ind{\mQ}{m}}$,  the ideal $\kk\ind{\mQ}{(1)}$  is   referred to as the
\emph{arrow ideal} in $\kQ$ and we have $\Frac{\kk\mQ}{
\kk\ind{\mQ}{(1)}}=\K$. In general the arrow ideal of $\kQ$
 needs not coincide with the Jacobson radical of
$\kQ$, unless $\K$ is semisimple and  $\mQ$ is acyclic (that is,
there exists some $m\geq 1$  with $B^{(m)}=0$).
Next,  the \emph{complete tensor path algebra}  of  $\mQ$ is the direct product
 $\kQc=\prd{m\geq 0}{\indu{B}{}{\sscr{(m)}}}$,  and the
\emph{closed arrow ideal} of $\kQc$ is given by
$\J{\kQc}=\prd{m\geq 1}{\indu{B}{}{\sscr{(m)}}}$, the latter
coincides with  the Jacobson radical of $\kQc$  whenever  $\K$ is
semisimple. For all $i,j\in \I$, each $\kk_i\sdash\kk_j$-bimodule
$\idl[i]\mul\kQc \mul\idl[j]$ is referred to as the bimodule of all
elements $\xi$ with source $\sr(\xi)=i$ and target
$\tg(\xi)=j$.

\ParIt{The $\J{\kQc}$-adic topology} The $\J{\kQc}$-adic topology  on $\kQc$ admits as  system of open neighbourhoods of $0$ the family
$\set{\J[l]{\kQc}}_{{l\geq 0}}$, with $\J[l]{\kQc}=\prd{m\geq
l}{\indu{B}{}{\sscr{(m)}}}$ for each $l\geq0$.
 The closure of each subset $S\subset \kQc$ is given by
 \vskip-0.5\baselineskip
\begin{equation}\label{eq.J-adic}
  \bar{S}=\inter{l\geq 0}{}{(S+\J[l]{\kQc})}.
\end{equation}

\begin{Rem}\label{rem.formelt-(A,JA)}
$\msf{(a)}$\; $\kQc$ coincides with the projective limit
$\varprojlim_{l\geq 0}\Frac{\kQc}{\J[l]{\kQc}}$, thus the
$\J{\kQc}$-adic topology on $\kQc$ is complete and separate. Next,
let $F=\som{\lambda=(\lambda_1,\dotsc,\lambda_m)
\in\N^m}{}{a_{\lambda}t_1^{\lambda_1} \mul \dotsm\mul
t_m^{\lambda_m}}$ be any power series over $\K$ for some natural
number $m$, then for all $u=(u_1,\dotsc, u_m)$ with $u_1, \dotsc,
u_m\in \J{\kQc}$, the infinite sum
$F(u)=\som{\lambda=(\lambda_1,\dotsc,\lambda_m)
\in\N^m}{}{a_{\lambda}u_1^{\lambda_1} \mul \dotsm\mul
u_m^{\lambda_m}}$ defines a unique element in $\kQc$  given as the
limit $\lim\limits_{\lambda \to[/dist=0.5] \infty} F_{\lambda}(u)$
of the series of partial sums $F_{\lambda}(u):=\som{\theta\leq
\lambda}{}{h_{\theta}(u)}$ where
$h_{\theta}(u)=a_{\theta}u_1^{\theta_1} \mul \dotsm\mul
u_m^{\theta_m}$ for each $\theta=(\theta_1\dotsc,\theta_m)\in
\N^m$. \\[-0.5\baselineskip]
$\,\msf{(b)}$\; Let $S\subset \kQc$ be any $\kk$-submodule, then
 $\bbar{S}=\set{\Som{l\geq 0}{x_l} \; : x_l\in  S\cap
 \J[l]{\kQc}}$.
 \finprv
\end{Rem}

\subsection{Casimir morphisms and projective bases for tensor path algebras} \label{sub.syst.fleches}

Let $l\geq 0$ be a fixed natural number, in view of
Lemma~\ref{lem.central-elt} and the
 discussion preceding  it,  we have an induced
symmetrizable dualizing pair  $\set{B^{(l)},{\dual{B}}^{(l)};\bilf^l}$ over $(\K,\tr)$, where $\set{B^{(0)},{\dual{B}}^{(0)};\bilf^0}$ coincides 
the natural  dualizing pair $\set{\K,\K}$, $\bilf^0$ being the
multiplication of $\K$. We also have the following  Casimir morphisms:
\vskip-0.5\baselineskip
\begin{equation}\label{eq-projbases-casimir}
\begin{aligned}
& \z[{(l)}]:=\morph{\z[B^{(l)}\otimes{\dual{B}}^{(l)}]:\K}{B^{(l)}
\otimes{\dual{B}}^{(l)}}:\, \idl \mapsto \Som{y\in\rmQ[l]}
y\otimes\dual{y}, \text{ with }
\rmQ[0]=\set{1}=\rdmQ[0], \\[-0.4\baselineskip]
& \z[{(l)}]':= \morph{\z[{\dual{B}}^{(l)}\otimes
B^{(l)}]:\K}{{\dual{B}}^{(l)} \otimes B^{(l)} }:\, \idl \mapsto
\Som{y\in\lmQ[l]} \dual{x}\otimes x, \text{ with }
\lmQ[0]=\set{1}=\ldmQ[0].
\end{aligned}
\end{equation}
 Here  the pair $(\lmQ[l],\ldmQ[l])$, with
$\ldmQ[l]=\set{\dual{x}\,:\, x\in\lmQ[l]}$, is  a left projective
basis for the left $\K$-module $B^{(l)}$  and its dual, and
$(\rmQ[l],\rdmQ[l])$, with  $\rdmQ[l]=\set{\dual{y}\,:\,
y\in\rmQ[l]}$,  is a right projective basis for the right
$\K$-module $B$ and its dual. We get two symmetrizable weakly dualizing pairs   $\set{\kQd,\kQc,\bilfh}$ and
$\set{\kQ,\kQdc,\bilft}$ with induced bilinear forms:
\vskip-1.5\baselineskip
\begin{equation}\label{eq.dualbasesgen}
 \begin{aligned}
&  \morph{\bilfh:\; \kQc\otimes\kQd \oplus \kQd\otimes
 \kQc}{\K},\; \bilfh(x\otimes\xi)=\Som{l\geq 0}\bilf^l(x_l\otimes\xi_l ) \text{ and }  \bilfh(\xi\otimes x
)= \Som{l\geq 0} \bilf^l(\xi_l\otimes x_l ) \\[-0.5\baselineskip]
& \text{ for all } \xi=(\xi_l)_{l\geq 0} \in \kQd \text{ and }
x=(x_l)_{l\geq 0}\in \kQc.
\end{aligned}
\end{equation}
\vskip-1.5\baselineskip
\begin{equation}\label{eq.dualbasesgen+}
 \begin{aligned}
&  \morph{\bilft:\; \kQ\otimes\kQdc \oplus \kQdc\otimes
 \kQ}{\K},\; \bilft(\chi\otimes\zeta)=\Som{l\geq 0}\bilf^l(\chi_l\otimes\zeta_l) \text{ and }  \bilft(\zeta\otimes
 \chi)=  \Som{l\geq 0}\bilf^l(\zeta_l\otimes \chi_l ) \\[-0.5\baselineskip]
& \text{ for all } \zeta=(\zeta_l)_{l\geq 0}\in \kQdc \text{ and }
\chi=(\chi_l)_{l\geq 0} \in \kQ.
\end{aligned}
\end{equation}
 For all natural numbers $n,\,m$ with $m \neq 0$ and  each symbol $\text{\textsc{s}}
\in\set{\text{\textsc{l}},\, \text{\textsc{r}}}$  we put:
\vskip-0.5\baselineskip
\begin{equation}\label{eq-projbaseskQ}
\begin{aligned}
&  \BmQ{\sscr{\text{\textsc{s}}}}{m} =
\BmQ[(m)]{\sscr{\text{\textsc{s}}}}{1} :=
\underbrace{\BmQ{\sscr{\text{\textsc{s}}}}{1} \otimes
\dotsm\otimes \BmQ{\sscr{\text{\textsc{s}}}}{1}}_{m \text{ copies
}},\; \text{ and }
 \BmdQ{\sscr{\text{\textsc{s}}}}{m} =
\BmdQ[(m)]{\sscr{\text{\textsc{s}}}}{1} :=
\underbrace{\BmdQ{\sscr{\text{\textsc{s}}}}{1} \otimes
\dotsm\otimes \BmdQ{\sscr{\text{\textsc{s}}}}{1}}_{m \text{ copies
}}  \\[-0.2\baselineskip]
& \BmQc{\sscr{\text{\textsc{s}}}}{n}=\Union{l\geq
n}{\BmQ{\sscr{\text{\textsc{s}}}}{l}},\; \text{ with  dual
projective basis: }
\BmdQc{\sscr{\text{\textsc{s}}}}{n}=\Union{l\geq
n}{\BmdQ{\sscr{\text{\textsc{s}}}}{l}}.
\end{aligned}
\end{equation}
Hence, $(\ldmQc[0],\lmQc[0])$ and $(\rdmQc[0], \rmQc[0])$ are two
pairs of projective bases associated with the symmetrizable weakly
dualizing pair $\set{\kQd,\kQc;\bilfh}$ and we have the
following characterizing identities:
\vskip-\baselineskip
\begin{equation}\label{eq.dualbasesKQc}
\Som{\chi\in\lmQc[0]}{\dual{\chi}\bilfh(\chi\otimes \sdash)} =
 \id_{\kQd} = \Som{\omega\in\rmQc[0]}{\bilfh(\sdash\otimes \omega)\dual{\omega}}
 \text{ and }
 \Som{\chi\in\lmQc[0]}{\bilfh(\sdash\otimes \dual{\chi})\chi} =
 \id_{\kQc} = \Som{\omega\in\rmQc[0]}{\omega \bilfh(\dual{\omega}\otimes
 \sdash)}.
\end{equation}
 One can derive  similar conclusions for the symmetrizable weakly
dualizing pair $\set{\kQ,\kQdc;\bilft}$.

\ParIt{Continuous morphisms of path algebras}
\begin{Defn} \label{defn.morp-pathalg} Let $M$ be a $\K$-bimodule, then an algebra morphism $\morph{f: \widehat{\mrm{T}_{\K}(M)}}{\kQc}$ is called a \emph{morphism of path algebras} if $f\restr{\K}=\id_{\K}$ and $f(M) \subset \J{\kQpc}$. In this case we let $\morph{f_l: M}{B^{(l)}}$, $l\geq 1$, be the family of $\K$-bimodule morphisms such that $f\restr{M}=(f_l)_{l\geq 1}$. 
\end{Defn}

In the classical case of a semisimple algebra $\K$, one checks that an algebra morphism $f$ as above is a path algebra morphism if $f\restr{\K}=\id_{\K}$. Recall that if $A$ is a $\kk$-algebra  with a $\K$-bimodule structure such that the unity of $A$ is $\K$-central,  then any $\K$-bimodule morphism  $\morph{f_{(1)}:B}{A}$ uniquely extends to a morphism $\morph{f: \kQ}{A}$ of
$\kk$-algebras. 

\begin{Prop}\label{prop.prolong} Given any  $\kk$-modulated quiver
$\mQ'=(B',\K,\tr)$,
  the two following statements are true.
\begin{itemize}
 \item[$\msf{(a)}$] Any  family $(\phi_l)_{l\geq 1}$ of $\K$-bimodule morphisms
 $\morph{\phi_l: B}{B'^{(l)}}$  defines a  unique
continuous morphism $\morph{\phi: \kQc}{\kQpc}$  of topological path $\kk$-algebras. Furthermore, $\phi$  is an isomorphism if and only if $\morph{\phi_1: B}{B'}$ is.
 \item[$\msf{(b)}$] Any path algebra morphism  $\morph{\phi:\kQc}{\kQpc}$  
is continuous and,  if $\phi$ is also surjective  then for every subset
$I\subset\kQc$ such that $\Ker{\phi}\subset \bbar{I}$ we have
$\phi(\bbar{I})=\bbar{\phi(I)}$. Consequently  any path algebra isomorphism   $\morph{\phi:\kQc}{\kQpc}$  is an
homeomorphism of topological path $\kk$-algebras.
 \end{itemize}
\end{Prop}
\begin{prv} \Par{Statement~$\msf{(a)}$} For the first part of
$\msf{(a)}$, the existence an  extension $\phi$ follows by the
universal property of $\kQ$ and by
Remark~\ref{rem.formelt-(A,JA)},   the continuity and hence the
uniqueness  of $\phi$ follow by statement $\msf{(b)}$. For the
second part of $\msf{(a)}$,  if $\phi$ is an isomorphism of
algebras, then $\morph{\phi_1: B}{B'}$ is clearly an isomorphism
of $\K$-bimodules. Conversely, assume that $\phi_1$ is an
isomorphism of $\K$-bimodules, thus without lost of generality we
can also assume that $B'=B$ and $\phi_1=\id_{B}$. With notations
of \eqref{casimir-K} and \eqref{eq-projbaseskQ}, we take a left
projective basis $(\ldmQc[0],\lmQc[0])$ for $\kQ$ and its weak
dual $\kQc$, and a projective basis $\set{e_s,\, \dual{e_s}:\,
s\in\Lambda}$ of the symmetric algebra $\K$ over $\kk$. The system
$\Ss=\set{e_s \chi\, : s\in\Lambda, \chi\in\lmQc }$ is a
"projective basis" of $\kQc$ over $\kk$ with the corresponding
dual "projective basis" $\dual{\Ss}=\set{ \dual{\chi} \dual{e_s}\,
: s\in\Lambda, \chi\in\lmQc }$.  The elements of $\Ss$ being
ordered in an increasing order  of their degree, in view of
identities~\eqref{eq.dualbasesKQc}, each element $x\in \kQc$ is
written  as an infinite $\kk$-linear combination
$x=\Som{\chi\in\Ss}{c_{\chi} \chi}$, and the infinite matrix
representing the map $\phi$ relatively to the projective basis
$\Ss$ is a triangular matrix with the "$1$'s" on its diagonal, and
hence is invertible, consequently $\phi$  is bijective.

\Par{Statement~$\msf{(b)}$} Let $J=\J{\kQc}$,  $J'=\J{\kQpc}$ and  $\morph{\phi:\kQc}{\kQpc}$ a morphism   of path algebras. Thus $\phi\restr{\K}=\id_{\K}$ and $\phi(B) \subset J'$. 
The definition of the $J$-adic topology shows  each subset $U\subset \kQc$ containing a power of $J$  is open. By assumption, $\phi(B)\subset J'$, implying that
 $\phi(J^{l})\subset J'^{l}$ for all $l\geq 1$, so that $\phi^{-1}(J'^{l})\supset J^{l} $, showing that each pre-image $\phi^{-1}(J'^{l})$ is
an open set. Hence $\phi$ is continuous.

Now assume that $\phi$ is surjective.  The previous discussion  shows that $\phi(J^l)=J'^l$ for all $l\geq 1$,  in
particular $\phi^{-1}(J'^l)=\phi^{-1}(\phi(J^{l}))=J^{l} +
\Ker(\phi)$. Let  $I\subset \kQc$  with
$\Ker{\phi}\subset \bbar{I}$. Using  the fact that
$\phi^{-1}(V\cap V')=\phi^{-1}(V)\cap\phi^{-1}( V')$  and
$\phi^{-1}(\phi (U))=U +\Ker(\phi)$ for all subsets  $V,V'\subset
\kQpc$  and $U\subset \kQc$, we have:
 $\phi^{-1}\parr{\bbar{\phi(I)}}=\phi^{-1}\parr{\inter{l=1}{\infty}{(\phi(I) +
J'^{l})}}=\inter{l=1}{\infty}{\phi^{-1}(
 \phi(I) + J'^{l})}=\inter{l=1}{\infty}{(I + \Ker(\phi) + J^l)}\subset
 \inter{l=1}{\infty}{(\bbar{I} + J^l)} = \bbar{\bbar{I}}=\bbar{I}$; but
 also $\bbar{I}=\inter{l=1}{\infty}{(I + J^l)} \subset
\inter{l=1}{\infty}{(I+\Ker(\phi) + J^l)}$.  Consequently,
$\phi^{-1}\parr{\bbar{\phi(I)}}=\bbar{I}$, so that $\bbar{\phi(I)}
= \phi(\bbar{I})$.
\end{prv}

\begin{Defn}[{\cite[2.5]{DWZ}}]\label{defn.morph-unitr} Let $\phi$ be a path algebra
automorphism of  $\kQc$. Then $\phi$ is   a
\emph{change of arrows} if $\phi_{(2)}:=(\phi_l)_{l\geq 2}=0$. If
  $\phi_1=\id_{B}$  then $\phi$ is    an
\emph{unitriangular automorphism}. We say that $\phi$  has
\emph{depth} $d\geq 1$ if $\phi_k=0$ for all $1\leq k\leq d$, in
this  case $\phi(u)-u\in\J[l+d]{\kQpc}$ for all
$u\in\J[l]{\kQc}$.
\end{Defn}

\subsection{Non simply-laced generalization of potentials and their Jacobian ideals}
\label{subsec:gnerPotJcIdeal}
 \begin{Defn}\label{def.potential} A \emph{potential} on $\mQ$ is any  $\K$-bimodule morphism  $\morph[/above={\mfr{m}}]{\K}{\kQc_{(2)}}$; thus the
bimodule of potentials on $\mQ$ identifies  with the
$\Zc(\K)$-bimodule $\Zck(\kQc_{(2)})$ of all $\K$-central elements
in $\J[2]{\kQc}=\kQc_{(2)}$.
\end{Defn}

By Lemma~\ref{lem.cyc-stable}, $\mfr{m} \in \Zck(\kQc_{(2)})$ is
\emph{cyclically stable}, that is, homogeneous
components of $\mfr{m}$ are cyclically stable.
\begin{Expl} Suppose $\K=\E\times \F\times \LL$ as product of indecomposable symmetric $\kk$-algebras. Below, the second modulated quiver  is
obtained from the first  by a transformation  named latter on as mutation.
\\[-1.5\baselineskip]
\begin{Center}
\diagram[/dist=0.1,/clsep=3.5]{\& \F \& \\ \E \& \& \LL }{ \path[<-] (m-2-1)
edge node[sloped,above=-0.08] {$\Bmd{1}{2}, \Bmdd{1}{2}$} (m-1-2)
(m-1-2) edge node[sloped,above=-0.08] {$\Bmd{2}{3}, \Bmdd{2}{3}$}
(m-2-3);} $\to[/style=dashed,/dist=2.5]$ \diagram[/dist=0.1,/clsep=4.5]{\& \F \&
\\ \E \& \& \LL }{ \path[->] (m-2-1) edge node[sloped,above=-0.08]
{$\Bmd{1}{2}, \Bmdd{1}{2}$} (m-1-2) (m-1-2) edge
node[sloped,above=-0.08] {$\Bmd{2}{3}, \Bmdd{2}{3}$} (m-2-3)
(m-2-3) edge node[below=-0.08] {$\Bmd{1}{2}\otimes\Bmd{2}{3},
\Bmdd{2}{3}\otimes\Bmdd{1}{2}$} (m-2-1);}
\end{Center}
Then a potential on the second modulated quiver is given by the Casimir
$\z[(\Bmdd{2}{3}\otimes\Bmdd{1}{2})\otimes(\Bmd{1}{2}\otimes\Bmd{2}{3})]$.
\end{Expl}

 For the symmetrizable weakly dualizing pair  $\set{\kQd,\kQc,\bilfh}$, the associated \emph{left derivative operator}  and \emph{right
derivative operator} are respectively 
$\morph{\lderv:=\bilfh\otimes \id:\,\kQd\otimes \kQc}{\kQc}$ and
$\morph{\rderv:=\id\otimes\bilfh:\, \kQc\otimes \kQd}{\kQc}$, they are explicitly described as follows: for all $\xi \in {\dual{B}}^{(l)} $,
$v\in B^{(d)}$ with $d<l$,   $x\in B^{(l)}$ and $ u\in \kQc$, we have
\vskip-0.5\baselineskip
\begin{equation}\label{eq-derv}
\begin{aligned}
&\lderv(\xi\otimes v)=0=\rderv(v\otimes \xi),\; \lderv(\xi\otimes xu)=\lderv(\xi\otimes x)u=\bilfh(\xi\otimes x)u=\bilf^l(\xi\otimes x)u \text{ and } \\
& \rderv(ux\otimes\xi )=u\rderv(x\otimes\xi)=u\bilfh(x\otimes\xi)=u\bilf^l(x\otimes\xi).
\end{aligned}
\end{equation}
We observe that the left derivative operator is a morphism of
$\K\sdash\kQc$-bimodules, while the right derivative operator is
morphism of $\kQc\sdash\K$-bimodules. The following observations
are direct generalizations  of identities \eqref{eq.dualbases} and
\eqref{eq.dualbasesKQc}.

\begin{Rem}\label{rem.der-dualbasis} For all natural number
$l\in\N$,
\vskip-1.5\baselineskip
\begin{equation}\label{eq.dervdualbasesKQc}
 \Som{\chi\in\lmQc[l]}{\dual{\chi}\rderv(\chi\otimes \sdash)} =
 \id_{\kQd_{l}} = \Som{\omega\in\rmQc[l]}{\lderv(\sdash\otimes
 \omega)\dual{\omega}}  \text{ and }
 \Som{\chi\in\lmQc[l]}{\rderv(\sdash\otimes \dual{\chi})\chi} =
 \id_{\kQc_{(l)}} = \Som{\omega\in\rmQc[l]}{\omega\lderv(\dual{\omega}\otimes
 \sdash)}.
\end{equation}
\vskip-0.5\baselineskip
 Moreover, taking a componentwise composition, on $\kQd\otimes
\kQc \otimes \kQd$ we have $\lderv\rderv = \rderv\lderv$:
\vskip-0.25\baselineskip
\begin{center}
$\lderv\rderv(\xi\otimes x\otimes \zeta)=\lderv(\xi\otimes
\rderv(x\otimes \zeta))= \rderv(\lderv(\xi\otimes
x)\otimes\zeta)=\rderv\lderv(\xi\otimes x\otimes \zeta)$ for all
$x\in\kQc$, $\xi,\zeta\in \kQd$. \finprv
\end{center}
\end{Rem}

Now let $\morph[/above={\mfr{m}}]{\K}{\kQc_{(2)}}$ be a 
potential on $\mQ$.  The action of the  left derivative 
operator on $\mfr{m}$ yields the bimodule  morphism 
$\morph{\lderv\mfr{m}=\lderv\circ(\id\otimes
\mfr{m}):\, \kQd}{\kQc}$; and the action the right derivative operator yields the bimodule  morphism
$\morph{\rderv\mfr{m}=\rderv\circ(\mfr{m}\otimes
\id):\, \kQd}{\kQc}$. Thus, when  $\mfr{m}$ is identified with
$\mfr{m}(1)$,  for each $\xi\in\kQd$ we have:
\vskip-0.5\baselineskip
\begin{center}
$\lderv[\xi]\mfr{m}:=(\lderv\mfr{m})(\xi)=\lderv(\xi\otimes\mfr{m})$
and
$\rderv[\xi]\mfr{m}:=(\rderv\mfr{m})(\xi)=\rderv(\mfr{m}\otimes\xi)$.
\end{center}
Note that we have Casimir morphisms: $
\morph{\z[{(l)}] \! :\K}{B^{(l)} \otimes{\dual{B}}^{(l)}}$ and $ \morph{\z[{(l)}]' \! :\K}{{\dual{B}}^{(l)}
\otimes B^{(l)} }$ described by \eqref{eq-projbases-casimir}.
\begin{Defn}[\emph{Skew permutations}]\label{def.skewperm}
Let $l\in \N$. The \emph{left} and the \emph{right  permutation operators}
of order $l$ are morphisms
$\morph{\lperm[l],\rperm[l]: \Zck(\kQc_{(2)})}{\Zck(\kQc_{(2)})}$
acting on   potentials
$\morph[/above={\mfr{m}}]{\K}{B^{(d)}}$ by:
$\morph{\lperm[l]\mfr{m}:=\lderv(\id\otimes\mfr{m}\otimes\id)\circ\z[l]':\,
\K }{B^{(d)}}$ and
$\morph{\rperm[l]\mfr{m}:=\rderv(\id\otimes\mfr{m}\otimes\id)\circ\z[l]:\,
\K }{B^{(d)}}$.
\end{Defn}

Thus $\lperm[0]=\id_{\Zck(\kQc_{(2)})}=\rperm[0]$, and  $\lperm[d]$ and $\rperm[d]$ act as identity maps on  homogeneous potentials of  degree $d$.

\begin{Prop} \label{prop.der-perm-cycderv} {\;}
 \begin{itemize}
   \item[$\msf{(1)}$]
For all potential $\mfr{m}$ on $\mQ$, the  action
$\morph{\rderv(\lperm[l](\mfr{m})):\, {\dual{B}}^{(l)}\,}{\,\kQc}$
of the right derivative operator
 on the left permutation of order $l$ of $\mfr{m}$ is equal to the action
 $\morph{\lderv(\mfr{m}):\, {\dual{B}}^{(l)}\,}{\,\kQc}$ of the left derivative operator on $\mfr{m}$.
Likewise, the action
$\morph{\lderv(\rperm[l](\mfr{m})):\, {\dual{B}}^{(l)}\,}{\,
\kQc}$ of the left derivative operator
 on the right permutation of order $l$ of $\mfr{m}$ is equal to the action
 $\morph{\rderv(\mfr{m}):\, {\dual{B}}^{(l)}\, }{\, \kQc}$ of the
  right derivative operator on $\mfr{m}$.
\item[$\msf{(2)}$] For every $l\in\N$, we have
 $\lperm[l]\circ \rperm[l] =\id_{\Zck(\kQc_{(2)})}=\rperm[l]\circ
\lperm[l]$.
 \item [$\msf{(3)}$] We have a \emph{cyclic permutation} operator
$\morph{\cperm:\, \Zck(\kQc_{(2)})}{\Zck(\kQc_{(2)})}$ 
  defined on    homogeneous potentials
$\mfr{m}$  of degree $d+1$ by: $\cperm
\mfr{m}=\som{t=0}{d}{\lperm[t] \mfr{m}}=\som{t=0}{d}{\rperm[t]
\mfr{m}}$. Consequently, there is a \emph{cyclic derivative}
operator $\morph{\cderv:\,\kQd\tens{\Zc(\K)}
\Zck(\kQc_{(2)}) \oplus \Zck(\kQc_{(2)})\tens{\Zc(\K)}\kQd
}{\kQc}$  such that
$\cderv(\xi\otimes \mfr{m}) =\lderv(\xi\otimes \cperm\mfr{m})
=\rderv(\cperm\mfr{m}\otimes \xi)$ for all $\xi\in\kQd$ and
$\mfr{m}\in \Zck(\kQc_{(2)})$. Hence the action of the cyclic
derivative on each potential $\mfr{m}$ is the bimodule morphism  
$\cderv \mfr{m}: \morph{\kQd}{\kQc}$ with $\cderv
\mfr{m}=\lderv \cperm \mfr{m}=\rderv\cperm \mfr{m}$.
\item[$\msf{(4)}$] Let  $\xi\in{\dual{B}}^{(s)}$ and
$\zeta\in{\dual{B}}^{(t)}$ with $1\leq s,t\in\N$. On the $\Zc(\K)$-module
 of potentials we have:\\
$\begin{array}{rl}\cderv((\xi\otimes\zeta)\otimes\sdash
)=\lderv((\xi\otimes\zeta)\otimes\cperm(\sdash)
)=\lderv(\xi\otimes \cderv(\zeta\otimes \sdash) ) = &
\rderv(\cderv(\xi\otimes \sdash)\otimes \zeta) \\
 = & \rderv(\cperm(\sdash)\otimes (\xi\otimes\zeta))=
\cderv(\sdash\otimes (\xi\otimes\zeta)).
\end{array}$
\end{itemize}
\end{Prop}

\begin{prv} \Par{Statement
$\msf{(1)}$} Let $\mfr{m}$ be a  potential on $\mQ$ and $l\in
\N$. For Casimir morphisms $
 \morph{\z[{(l)}]:\K}{B^{(l)}
\otimes{\dual{B}}^{(l)}}$ and $ 
\morph{\z[{(l)}]' :\K}{{\dual{B}}^{(l)}
\otimes B^{(l)} }$,  as in \eqref{eq-projbases-casimir} we have $\z[l](1)=\Som{y\in\rmQ[l]}
y\otimes\dual{y} $ and $\z[l]'(1)=\Som{x\in\lmQ[l]}
\dual{x}\otimes x $, where $(\lmQ[l],\ldmQ[l])$ and
$(\lmQ[l],\ldmQ[l])$ are left and right projective bases for the
bimodule $B^{(l)}$ and its dual. By definition
$\rderv(\lperm[l](\mfr{m}))=\rderv\circ(\lperm[l](\mfr{m})\otimes
\id)$,   Definition~\ref{def.skewperm}  shows  that
$\lperm[l](\mfr{m})=\lderv \circ (\id\otimes\mfr{m} \otimes
\id)\circ \z[l]'$. Using \eqref{eq-derv} and
identities~\eqref{eq.dervdualbasesKQc} from
Remark~\ref{rem.der-dualbasis}, for all $\xi \in {\dual{B}}^{(l)} $ we have:\\
$
\begin{array}{rl} \bigl(\rderv(\lperm[l](\mfr{m}))\bigr)(\xi)=
& \rderv\bigl(\lperm[l](\mfr{m})(1)\otimes
\xi\bigr)=\rderv\Parr{\lderv\parr{\Som{x\in
\lmQ[l]}{\dual{x}\otimes \mfr{m}(1) \otimes x}}\otimes
\xi} =  \lderv\parr{\Som{x\in \lmQ[l]}{\dual{x}\otimes \mfr{m}(1)
\mul\rderv(x\otimes \xi)}} \\
= &  \lderv\parr{\Som{x\in
\lmQ[l]}{\dual{x}\mul\rderv(x\otimes \xi)
\otimes \mfr{m}(1) }}\hspace{-0.2em}=\hspace{-0.2em} 
\lderv\Parr{\parr{\Som{x\in \lmQ[l]}{\dual{x}\mul\rderv(x\otimes
\xi)}} \otimes \mfr{m}(1)}\hspace{-0.2em} = \hspace{-0.2em}\lderv(\xi \otimes
\mfr{m}(1))\hspace{-0.2em}=\hspace{-0.2em}(\lderv \mfr{m}) (\xi).
\end{array}
$ \\ Thus  $\rderv(\lperm[l](\mfr{m}))=\lderv{\mfr{m}}$ on ${\dual{B}}^{(l)}$. Similarly, one proves using   $\z[l]$  that
$\lderv(\rperm[l](\mfr{m}))=\rderv{\mfr{m}}$ on
${\dual{B}}^{(l)}$.
 \Par{Statement $\msf{(2)}$} Let $\mfr{m}$ be any potential on $\mQ$. By statement
$(1)$ and identities~\eqref{eq.dervdualbasesKQc} we
we have:\\
$\begin{array}{rl} \rperm[l]\circ\lperm[l] (\mfr{m})= &
\rderv\circ(\id \otimes \lperm[l](\mfr{m}) \otimes \id)\circ
\z[l]=(\id\otimes \rderv\circ(\lperm[l](\mfr{m})\otimes
\id))\circ\z[l]=(\id\otimes
\rderv\lperm[l](\mfr{m}))\circ\z[l] \\
= & (\id\otimes \lderv \mfr{m} )\circ \z[l]=
\Som{y\in\rmQ[l]}{y\lderv(\dual{y}\otimes \mfr{m})}=\mfr{m}.
\\[.1em]
\lperm[l]\circ\rperm[l] (\mfr{m})= & \lderv\circ(\id  \otimes
\rperm[l](\mfr{m}) \otimes \id)\circ \z[l]'=( \lderv\circ(\id
\otimes \rperm[l](\mfr{m}))\otimes \id
)\circ\z[l]'=((\lderv\rperm[l] \mfr{m})\otimes \id)\circ\z[l]' \\
= & (\rderv \mfr{m} \otimes \id)\circ \z[l]'=
\Som{x\in\lmQ[l]}{\lderv(\mfr{m}\otimes\dual{x})x}=\mfr{m}.
\end{array}$
\Par{Statement $\msf{(3)}$}  To show that the cyclic permutation
operator $\cperm$ is properly defined, it suffices  to consider
the case of an homogeneous potential $\mfr{m}$ of degree $d+1$
with $d\geq 1$. Statement $\msf{(2)}$ and  the cyclical
stability of $\mfr{m}$  show that
$\lperm[d+1]\mfr{m}=\mfr{m}=\rperm[d+1] \mfr{m}$ and
$\som{t=0}{d}{\lperm[t] \mfr{m}}=\som{t=0}{d}{\lperm[t]
\rperm[d+1](\mfr{m})}=\som{t=0}{d}{\rperm[d+1-t]
\mfr{m}}=\som{t=0}{d}{\rperm[t] \mfr{m}}$,  thus $\cperm$
is properly defined. For the existence of the cyclic derivative
operator, consider an arbitrary potential $\mfr{m}$. Observe that
 the cyclically stability also shows that $\rperm[l]\cperm =
\cperm =\lperm[l]\cperm $ for all natural number $l$. Let
$\xi=\Som[n]{l=1}{\xi_l} \in\kQd=\K\oplus \dual{B} \oplus
{\dual{B}}^{(2)} \oplus\dotsm $ with $\xi_l\in{\dual{B}}^{(l)}$.
Applying  statement $\msf{(1)}$, we have: $\lderv(\xi\otimes
\cperm\mfr{m})=\Som[n]{l=0}{\lderv(\xi_l \otimes
(\rperm[l]\cperm\mfr{m}))}=\Som[n]{l=0}{\rderv(\cperm\mfr{m}
\otimes \xi_l)}=\rderv(\cperm\mfr{m}\otimes \xi)$, establishing
the existence and the definition of the cyclic derivative
operator. \Par{Statement $\msf{(4)}$} We apply the  definition of
the cyclic derivative and the fact the left derivative and the
right derivative pointwise commute. Let
$\xi\in{\dual{B}}^{(s)},\, \zeta\in{\dual{B}}^{(t)}$, on the
$\Zc(\K)$-module  of potentials we have:\\
 $
\begin{array}{rl}
\lderv((\xi\otimes\zeta)\otimes\cperm(\sdash) ) & =
\lderv(\xi\otimes\lderv(\zeta \otimes\cperm(\sdash)) ) =
\lderv(\xi\otimes\cderv(\zeta \otimes \sdash)
)=\lderv(\xi\otimes\rderv( \cperm(\sdash)\otimes \zeta)
) \\
& =\lderv\rderv(\xi\otimes \cperm(\sdash)\otimes \zeta )
=\rderv\lderv(\xi\otimes \cperm(\sdash)\otimes \zeta )=
\rderv(\lderv(\xi\otimes \cperm(\sdash))\otimes \zeta )\\
&  = \rderv(\cderv(\xi\otimes \sdash)\otimes \zeta ) =
\rderv(\cderv(\sdash\otimes \xi)\otimes \zeta )
=\rderv(\rderv(\cperm(\sdash)\otimes \xi)\otimes \zeta) \\
& =\rderv(\cperm(\sdash)(\xi\otimes \zeta)).
\end{array}
$

This establishes the identities of statement $\msf{(4)}$ and
completes the proof of the proposition.
\end{prv}

In the sequel, for  $\xi\in\kQd$ and $x\in\kQc$  we put:
$\lderv[\xi]x=\lderv(\xi\otimes x)$, $\rderv[\xi]x=\rderv(x\otimes
\xi)$ and
$\cderv[\xi]\!=\!\cderv(\xi\otimes\sdash)\!=\!\cderv(\sdash\otimes
\xi)$.

\begin{Defn}\label{defn.JacAlg}
For a  potential $\mfr{m}\in\Zck(\kQc_{(2)})$, the
closure $\J{\mfr{m}}:=\bbar{\ggen{(\cderv \mfr{m})(\dual{B})}}$ of
the ideal $\ggen{(\cderv \mfr{m})(\dual{B})}$ in $\kQc$ is  
the \emph{Jacobian ideal} of   $\mfr{m}$, the corresponding \emph{Jacobian algebra} is  
$\Jc[\mfr{m}]=\Jc(\mQ,\mfr{m})=\Frac{\kQc}{\J{\mfr{m}}}$.
\end{Defn}

\begin{Defn} The $\Zc(\K)$-module $\skc{\Zck(\kQc_{(2)}),\Zck(\kQc_{(2)})}$ of \emph{skew commutators}  consists of finite sums of  potentials of the form
$\mfr{m}-\lperm\mfr{m}$ (or equivalently, of
the form $\mfr{m}-\rperm\mfr{m}$). The \emph{closed
$\Zc(\K)$-module  of skew commutators} is  the closure
$\skcb{\Zck(\kQc_{(2)}),\Zck(\kQc_{(2)})}$ of $\skc{\Zck(\kQc_{(2)}),\Zck(\kQc_{(2)})}$. Two potentials $\mfr{m}$ and $\mfr{m}'$ are  \emph{cyclically equivalent} whenever
\mbox{$\mfr{m}-\mfr{m}'$} lies in $\skcb{\Zck(\kQc_{(2)}),\Zck(\kQc_{(2)})}$, in  this case  we have $\cderv[\xi]\mfr{m}=\cderv[\xi] \mfr{m}'$ for all
$\xi\in\dual{B}$,  $\J{\mfr{m}}=\J{\mfr{m}'}$ and
$\Jc(\mQ,\mfr{m})=\Jc(\mQ,\mfr{m}')$.
\end{Defn}

\subsection{Intrinsic description of potentials, modulated quivers with potentials}

\begin{Lem} \label{lem.wdualizing-pair-split}
Let $\set{M,\dual{M};\bilf}$ be a symmetrizable weakly  dualizing
pair, $\set{U,\dual{U};\rho}$ and
$\set{\bbar{U},\dual{\bbar{U}};\bbar{\bilf}}$  two symmetrizable
pairing  over the  symmetric algebra $(\K,\tr)$.
Let $\Seq{(\varepsilon):\;  0 \to U \to[/above={f}] \, M \,
\to[/above={g}]\, \bbar{U} \to 0}$ and
$\Seq{(\dual{\varepsilon}):\;  0 \to \dual{\bbar{U}}
\to[/above={\dual{g}}] \, \dual{M} \, \to[/above={\dual{f}}]\,
\dual{U} \to 0}$  be mutually dual  exact sequences of $\K$-bimodule morphisms
 such that $(\varepsilon)$ splits as  sequence of
    left $\K$-linear maps. Then, the dual sequence $(\dual{\varepsilon})$
      splits as  sequence of  right $\K$-linear maps and   the following holds. Let
    $\Seq{ U \to[/dir=<-,/above={f'}] \, M \,
\to[/dir=<-,/above={g'}]\, \bbar{U}}$ be left
$\K$-linear morphisms such that the map $\morph[/above={\sim}]{h:=[
f,g']:\, U \oplus \bbar{U}}{M}$ is a left $\K$-linear isomorphism with inverse
 $\morph[/above={\sim}]{ h^{-1}= \Psmatr{ f' \\
   g}:\, M}{U \oplus \bbar{U} }$, then  $f'$ and $g'$ are left dualizing, 
$\morph[/above={\sim}]{ \ldual{h}= \Psmatr{ \dual{f} \\
   \ldual{g'}}:\, \dual{M}}{\dual{U} \oplus \dual{\bbar{U}}}$ and
\mbox{$\ldual{h^{-1}} = [ \ldual{f'},\dual{g}]$}.
Moreover,   if   $(\varepsilon)$ also splits as  
sequence of right $\K$-linear maps  then $\set{U,\dual{U};\rho}$ and
$\set{\bbar{U},\dual{\bbar{U}};\bbar{\bilf}}$ are weakly dualizing pairs.
\end{Lem}

\begin{prv} By assumption on $(\varepsilon)$, let  $\Seq{ U \to[/dir=<-,/above={f'}] \, M \, \to[/dir=<-,/above={g'}]\, \bbar{U}}$ be left $\K$-linear morphisms with $f'f=\id_U$, $gg'=\id_{\bbar{U}}$ and $f'g'=0$.
Part $\msf{(2)}$ of Lemma~\ref{lem.dualizing-pair} states that 
$f'$ is left dualizing, and clearly  $\dual{f}\circ
\ldual{f'}=\id_{\dual{U}}$, implying (using basic module theory) that  the dual
sequence $(\dual{\varepsilon})$ splits as  sequence of right $\K$-linear morphisms and  there is a unique right $\K$-linear morphism
$\morph{g'':\dual{M}}{\dual{\bbar{U}}}$  such $\id_{\dual{M}}= \ldual{f'}\dual{f} + \dual{g} g''$,  $g''\dual{g}=\id_{\dual{\bbar{U}}}$ and $g''
\ldual{f'}=0$.
Next, we must check  that  $\morph{g'':\dual{M}}{\dual{\bbar{U}}}$ also serves as left dual for $\morph{g':\bbar{U} }{\dual{M}}$. Let $\xi\in\dual{M}$ and
$\bbar{x}=g(x)\in \bbar{U}$ with $x\in M$. Then, we have
$\xi=\id_{\dual{M}}(\xi)=\dual{g}g''(\xi)+\ldual{f'}\dual{f}(\xi)$,
thus:\\
$\begin{array}{rl} \bilf(g'(\bbar{x}) \otimes \xi)= &
\bilf(g'(\bbar{x})\otimes
(\dual{g}g''(\xi)+\ldual{f'}\dual{f}(\xi))) =
\bilf(g'(\bbar{x})\otimes \dual{g}g''(\xi)) +
\bilf(g'(\bbar{x}) \otimes \ldual{f'}\dual{f}(\xi)) \\
=&  \bbar{\bilf}(gg'(\bbar{x})\otimes g''(\xi)) + \rho(
f'g'(\bbar{x}) \otimes \dual{f}(\xi)) = \bbar{\bilf}(\bbar{x}
\otimes g''(\xi)) + \rho(
0 \otimes \dual{f}(\xi)) \\
= & \bbar{\bilf}(\bbar{x}\otimes g''(\xi)),
\end{array}$ \\
showing that $g''$ is indeed the left dual of $g'$.

Clearly, $\morph[/above={\sim}]{
\ldual{h^{-1}} = [ \ldual{f'},\dual{g}]:\, \dual{U} \oplus
\dual{\bbar{U}}}{\dual{M}}$. The adjoint map $\morph{\bilf[{{}{1,\lft}}]:\,
\dual{M}}{\ldual{M}, \xi \mapsto \bilf(\sdash\otimes \xi)}$ is a
bimodule isomorphism,  let be  $\morph[/above={\sim}]{\what{h}:=\Psmatr{\what{f} \\
\what{g}'}:\,\ldual{M}}{ \ldual{U} \oplus \ldual{\bbar{U}}}$ be the right $\K$-linear isomorphism with
 $ \what{f}(\alpha)= \alpha\circ f$ and  $\what{g}'(\alpha)=\alpha\circ g'$ for all $\alpha \in \ldual{M}$; note that $\what{f}$ is a  bimodule morphism while $\what{g}'$ is right $\K$-linear.  We want to compute the composition:
 $\Seq{\theta:=\what{h} \circ \bilf[{{}{1,\lft}}] \circ \ldual{h^{-1}}\, : \,
 \dual{U} \oplus \dual{\bbar{U}}\, \to[/above={\sim}] \, \ldual{U} \oplus
\ldual{\bbar{U}} }$. 
Writing each element in $\dual{U} \oplus \dual{\bbar{U}}$ as 
formal sum $\upsilon + \zeta$ with $\upsilon \in \dual{U}$ and $\zeta\in \dual{\bbar{U}}$, we get:\\
$\begin{array}{rl} \theta(\upsilon + \zeta) = & \Psmatr{\what{f} \\
\what{g}'}\Bigl(\bilf(\sdash \otimes (\ldual{f'}(\upsilon) +
\dual{g}(\zeta) ) )\Bigr) = \what{f}(\bilf(\sdash \otimes
(\ldual{f'}(\upsilon) + \dual{g}(\zeta) ))  +
\what{g}'(\bilf(\sdash
\otimes (\ldual{f'}(\upsilon) + \dual{g}(\zeta) )) \\
= & \bilf(f(\sdash) \otimes (\ldual{f'}(\upsilon) +
\dual{g}(\zeta) )) + \bilf(g'(\sdash) \otimes
(\ldual{f'}(\upsilon) + \dual{g}(\zeta)
)) \\
= & \bilf(f(\sdash) \otimes \ldual{f'}(\upsilon)) +
\bilf(f(\sdash) \otimes  \dual{g}(\zeta)) + \bilf(g'(\sdash)
\otimes \ldual{f'}(\upsilon)) + \bilf(g'(\sdash) \otimes
\dual{g}(\zeta)) \\
= & \rho(f'f(\sdash) \otimes \upsilon) + \bbar{\bilf}(g f(\sdash)
\otimes \zeta) + \rho(f'g'(\sdash) \otimes \upsilon) +
\bbar{\bilf}(gg'(\sdash) \otimes \zeta) \\
= & \rho(\sdash \otimes \upsilon) + 0 + 0+ \bbar{\bilf}(\sdash
\otimes \zeta). \\
\end{array} $ \\
Thus   $\theta$ is the direct sum of $\morph{\inds{\rho}{1,\lft}:\, \dual{U}}{\ldual{U}, \upsilon
\mapsto \rho(\sdash\otimes \upsilon)}$ and
$\morph{\bilfb[{{}{1,\lft}}] :\,
\dual{\bbar{U}}}{\ldual{\bbar{U}}, \zeta \mapsto \bilfb(\sdash\otimes
\zeta)}$, implying  that these adjoint maps are isomorphisms. Since 
$M$ is projective as left and right $\K$-module, so are $U$ and $\bbar{U}$. 
Thus, for the last statement of the lemma,
if $(\varepsilon)$ also splits as  sequence of
    right $\K$-linear maps, then  $U$ and $\bbar{U}$ are also
    projective as right $\K$-modules, and  we deduce  that  $\set{U,\dual{U};\rho}$ and
$\set{\bbar{U},\dual{\bbar{U}};\bbar{\bilf}}$ are weakly dualizing.
\end{prv}

The following result shows  that  potentials appear as
Casimir elements provided some splitting conditions hold.
\begin{Prop} \label{prop.dec.pot-mQ}
\begin{itemize}
\item[$\msf{(1)}$]
 Let $\set{M,\dual{M};\beta}$ and
$\set{M',\dual{M'};\beta'}$ be  symmetrizable dualizing pairs
 over $(\K,\tr)$, and $\morph{\mfr{m}: \K }{M\otimes
M'}$ an homogeneous potential. Then the right derivative
$\morph{\rderv \mfr{m} : \dual{M'} }{M}$  and the left derivative
$\morph{\lderv \mfr{m} : \dual{M}}{M'}$ are mutually dual
morphisms; they induce three symmetrizable pairings $\set{U,V;\gamma}$,
$\set{\bbar{U},\dual{\bbar{U}};\bbar{\beta}}$ and
$\set{\bbar{V},\dual{\bbar{V}}; \bbar{\beta}'}$ with \mbox{$U=\im
(\rderv \mfr{m})$},  \mbox{$V=\im(\lderv \mfr{m})$},
\mbox{$\dual{\bbar{U}}=\Ker(\lderv \mfr{m})$},
\mbox{$\bbar{U}=\Frac{M}{U}$}, \mbox{$\dual{\bbar{V}}=\Ker(\rderv
\mfr{m})$} and \mbox{$\bbar{V}=\Frac{M'}{V}$},  together with the
following pairs of mutually dual  exact sequences of canonical
injections and   projections:
\\[-1.5\baselineskip]
 \begin{center}
 $\begin{array}{l}
   (\varepsilon):\; \Seq[/Dist=1.2]{\dual{\bbar{U}}  \to[/dir=right hook->]
\dual{M} \to[/dir=->>,/above={\msf{p}}] V,} \;
(\dual{\varepsilon}):\;  \Seq[/Dist=1.2]{ U \to[/dir=right hook->]
M  \to[/dir=->>,/above={\bbar{\msf{p}}}]\bbar{U} }
  \text{ and }
 (\varepsilon'):\;
\Seq[/Dist=1.2]{\dual{\bbar{V}} \to[/dir=right hook->] \dual{M'}
 \to[/dir=->>,/above={\msf{p}'}] U,}  \; (\dual{\varepsilon'}): \;
  \Seq[/Dist=1.2]{V \to[/dir=right hook->]
M'  \to[/dir=->>,/above={\bbar{\msf{p}}'}] \bbar{V}}.
\end{array} $
\end{center}
Moreover, if $(\dual{\varepsilon})$ and $(\dual{\varepsilon'})$
split as sequences of left and right $\K$-linear maps,
then so does $(\varepsilon)$ and $(\varepsilon')$ and, $\set{U,V;\gamma}$,
$\set{\bbar{U},\dual{\bbar{U}};\bbar{\beta}}$ and
$\set{\bbar{V},\dual{\bbar{V}}; \bbar{\beta}'}$  are
dualizing. In  this case,
 let $\morph[/above={\sim}]{h:=[\id, \bbar{\msf{q}}]:\, U \oplus
\bbar{U}}{M}$ be  a right $\K$-linear isomorphism  with inverse
 \mbox{$h^{-1}= \Psmatr{ \rho \\ \bbar{\msf{p}}}$}  and  $\morph[/above={\sim}]{h':=[\id,\bbar{\msf{q}}']:\, V \oplus \bbar{V}}{M'}$ a  left $\K$-linear isomorphism with inverse  \mbox{$ h'^{-1}= \Psmatr{ \rho' \\
   \bbar{\msf{p}}'}$}, then $ \rdual{(h^{-1})} = [ \rdual{\rho},\id]$,
$\ldual{h'^{-1}} = [ \ldual{\rho'},\id]$ and we have:
\\[-1.5\baselineskip]
\begin{center}
 $  \z[M\otimes \dual{M}] =
(h\otimes \rdual{(h^{-1})})(\z[(U \oplus \bbar{U})\otimes (V
\oplus \dual{\bbar{U}})]) = (\id\otimes\rdual{\rho})(\z[U\otimes
V]) + (\bbar{\msf{q}}\otimes \id)(\z[\bbar{U}\otimes
\dual{\bbar{U}}])$, \\
$ \z[\dual{M'} \otimes M'] =
(\ldual{h'^{-1}}\otimes h')(\z[(U \oplus \dual{\bbar{V}})\otimes (V
\oplus \bbar{V})]) = (\ldual{\rho'}\otimes \id)(\z[U\otimes V]) +
(\id \otimes \bbar{\msf{q}}')(\z[\dual{\bbar{V}}\otimes \bbar{V}])$ and  
$\mfr{m}=\z[U\otimes V].$
\end{center}
\item[$\msf{(2)}$] Let $W\in B^{(2)}$ be a  potential.
Then the  cyclic derivative $\morph[/above={\cderv
W}]{\dual{B}}{B}$   is a self-dual morphism  inducing  two
symmetrizable pairings  $\set{B_0,B_0}$ and
$\set{\bbar{B}, \dual{\bbar{B}}}$ with $B_0=\im(\cderv W)$,
$\bbar{B}=\Frac{B}{B_0}$ and $\dual{\bbar{B}}=\Ker(\cderv W)$,
together with mutually dual  exact sequences of canonical
injections and projections:
$\Seq[/Dist=1.2]{(\vartheta):\, \dual{\bbar{B}}
\to[/dir=right hook->] \dual{B}
 \to[/dir=->>,/above={\msf{p}}] B_0,}$ 
$\Seq[/Dist=1.2]{ (\dual{\vartheta}):\, B_0 \to[/dir=right
hook->,/above={\dual{\msf{p}}}] B
 \to[/dir=->>,/above={\bbar{\msf{p}}}] \bbar{B}}$.
Moreover, if $B_0$ is a direct summand in $B$ then
$\mQ=\triv{\mQ} \oplus \bbar{\mQ}$ as direct sum of two modulated
quivers, with $\bbar{\mQ}=(\bbar{B},\K,\tr)$,
$\triv{\mQ}=(B_0,\K,\tr)$  and  we have $\cperm W=\z[B_0\otimes B_0] \in B_0^{(2)}$.
\end{itemize}
\end{Prop}

\begin{prv}  Since part $\msf{(2)}$   is a direct
application of part $\msf{(1)}$, we only need to prove part
$\msf{(1)}$. But then, in view of
Lemma~\ref{lem.dualizing-pair-Casimir} and
Lemma~\ref{lem.wdualizing-pair-split}, it suffices to establish
the first part of $\msf{(1)}$ and the identity
$\mfr{m}=\z[U\otimes V]$ in the second part of $\msf{(1)}$. To
start, we must show that  the left
derivative $\morph{f:=\lderv \mfr{m} : \dual{M} }{M'}$ and the
right derivative $\morph{f':=\rderv \mfr{m} : \dual{M'}}{M}$ are
mutually dual. Once again by Lemma~\ref{lem.dualizing-pair}, $f$
and $f'$ are already dualizing morphisms. Let $\xi \in \dual{M}$,
$\xi'\in\dual{M'}$, we have:
$\beta(\xi\otimes\dual{f}(\xi'))=\beta'(f(\xi)\otimes
\xi')=\rderv[\xi']f(\xi) =\rderv[\xi']\lderv[\xi]\mfr{m} =
\lderv[\xi]\rderv[\xi']\mfr{m}= \beta(\xi,\rderv[\xi']\mfr{m})$,
showing in view of Lemma~\ref{lem.dualizing-pair} that
$\dual{f}=\rderv[\xi']\mfr{m}=f'$. Next, the pair
$( f,  \dual{f})$ clearly induces a well-defined bilinear form
$\morph{\gamma: U\otimes V \oplus V\otimes U}{\K}$ such that:
$\gamma(f(\xi)\otimes \dual{f}(\xi')):=\beta(\xi\otimes
\dual{f}(\xi'))=\beta'(f(\xi)\otimes \xi')$ and
$\gamma(\dual{f}(\xi')\otimes f(\xi)):=\beta(\dual{f}(\xi')
\otimes\xi)=\beta'(\xi'\otimes f(\xi))$ for all $\xi \in
\dual{M},\, \xi'\in\dual{M'}$, and $\gamma$ is also non-degenerate
since $\beta$ and $\beta'$ are (strongly) non-degenerate. We get
that the data $\set{U,V,\gamma}$ is a symmetrizable pairing  over $(\K,\tr)$. 
In the same way,  that there
are canonically induced symmetrizable pairing 
$\set{\bbar{U},\dual{\bbar{U}};\bbar{\beta}}$ and
$\set{\bbar{V},\dual{\bbar{V}}; \bbar{\beta}'}$ over $(\K,\tr)$. Now, we want
to show that   the sequence $\Seq[/Dist=1.2]{(\varepsilon):\;
\dual{\bbar{U}} \to[/dir=right hook->] \dual{M}
\to[/dir=->>,/above={\msf{p}}] V,}$, where  $\msf{p}$ is the
projection defined by $f$,  is dualizing and $(\dual{\varepsilon})$ is the short exact sequence
 $ \Seq[/Dist=1.2]{U \to[/dir=right hook->] M
\to[/dir=->>,/above={\bbar{\msf{p}}}]\bbar{U} }$ defined by the subbimodule $U\subset M$.  
For all  $u \in U \subset M$ and $\xi\in \dual{M}$, we have: $\gamma(u\otimes
\msf{p}(\xi))=\gamma(u\otimes f(\xi))=\beta(u\otimes \xi)$,
showing that  the  inclusion  $\morph[/dir=right
hook->]{U}{M}$ is the right dual of $\msf{p}$ and hence 
the dual of $\msf{p}$ in view of Lemma~\ref{lem.dualizing-pair}. For all
$\xi_0\in \dual{\bbar{U}}$ and $x\in M$, writing $\bbar{x}$ for
the coset of $x$ in $\bbar{U}=\Frac{M}{U}$, by definition  we
have: $\bbar{\beta}(\xi_0\otimes\bbar{x})=\beta( \xi_0\otimes x)$,
showing as before that the inclusion $\morph[/dir=right
hook->]{\dual{\bbar{U}}}{\dual{M}}$ is dualizing and its dual is the canonical projection
$\morph[/dir=->>]{\bbar{\msf{p}}: M}{\Frac{M}{U}}$. In the same
way, on can check that there are mutually dual exact sequences  $
\Seq[/Dist=1.2]{(\varepsilon'):\;  \dual{\bbar{V}} \to[/dir=right
hook->] \dual{M'} \to[/dir=->>,/above={\msf{p}'}] U,} \;
  \Seq[/Dist=1.2]{(\dual{\varepsilon'}): \; V \to[/dir=right hook->]
M'  \to[/dir=->>,/above={\bbar{\msf{p}}'}] \bbar{V}}$.

As said   before, the rest of the proof
 of  $\msf{(1)}$, except for the relation
$\mfr{m}=\z[U\otimes V]$,  is given by
Lemmas~\ref{lem.dualizing-pair-Casimir}~and~\ref{lem.wdualizing-pair-split}. But, writing
$\z[M\otimes\dual{M}]=\Som[n]{r=1}y_r\otimes\dual{y_r}$ and
$\z[\dual{M'}\otimes M']=\Som[n']{s=1}\dual{x_s'}\otimes x_s'$ for
the Casimir elements in $M\otimes\dual{M}$ and 
$\dual{M'}\otimes M'$, by \eqref{eq.dervdualbasesKQc} we
have $\mfr{m}=\Som[n]{r=1}y_r\otimes
\lderv[\dual{y_r}]\mfr{m}=(\id\otimes
\msf{p})(\z[M\otimes\dual{M}])$ and
$\mfr{m}=\Som[n']{s=1}(\rderv[\dual{x_s'}]\mfr{m})\otimes x_s'
=(\msf{p}'\otimes \id)(\z[\dual{M'}\otimes M'])$.  Hence, the
equality $\mfr{m}=\z[U\otimes V]$ follows from the  relations
$\z[M\otimes \dual{M}] = (\id\otimes\rdual{\rho})(\z[U\otimes V])
+ (\bbar{\msf{q}}\otimes \id)(\z[\bbar{U}\otimes
\dual{\bbar{U}}])$ and $
 \z[\dual{M'} \otimes M'] =(\ldual{\rho'}\otimes \id)(\z[U\otimes V]) + (\id \otimes
\bbar{\msf{q}}')(\z[\dual{\bbar{V}}\otimes \bbar{V}])$. 
\end{prv}

For a potential $\mfr{m}$ on $\mQ$ with degree-$2$ component
$\mfr{m}_{2}\in B\otimes B$, let
$\triv{B}=\cderv(\dual{B}\otimes \mfr{m}_{2})$.  Then we have an induced symmetrizable pairing  $\set{U,V}$ with $U=(\rderv \mfr{m}_{2})(\dual{B})$, $V=(\lderv \mfr{m}_{2})(\dual{B})$  and $\triv{B}\subseteq U+V$. Thus, if  $U$ and $V$  are also projective as left  and  right $\K$-modules then the pair $\set{U,V}$ is  dualizing. 
\begin{Defn} \label{defn.mQp} The potential
 $\mfr{m}$ is   \emph{$2$-loop free}  if as left and right $\K$-module $\triv{B}$ is a direct summand in $B$  and $U\cap V=0$.  In this case,  the paring $\set{U,V}$  is dualizing, $\mfr{m}_{2}=\z[U\otimes V] \in U\otimes V$ is a Casimir element and $\triv{B}=U\oplus V$;   the data $(\mQ,\mfr{m})$   is called a \emph{modulated quiver with potential}.
\end{Defn}

\section{Reduction of modulated quivers with potentials}~\label{sec.reduction}
The main  results in section require some  preparation about Jacobian ideals.  

\subsection{The cyclic Leibniz rule and the chain-rule. }\label{sec:Leibnizrule}

In the  study  of quivers with potentials, the \emph{cyclic Leibniz rule} is an easy
consequence of the  fact that any simply-laced path algebra has a "symmetric" path $\kk$-basis and
the computation of cyclic derivatives only requires the ordinary
cyclic permutation  of arrows in the quiver. 
However, such a  symmetry is generally absent in the present  framework.
Thanks to  properties of symmetrizable dualizing pairs, the following result
controls skew permutations of potentials along morphisms of $\K$-bimodules.

\begin{Prop}\label{prop.morph-et-perm}  Let $\morph{f: \widehat{\mrm{T}_{\K}(U)}}{\kQc_{(1)}}$ and
 $\morph{h: \widehat{\mrm{T}_{\K}(V)}}{\kQc_{(1)}}$
 be path algebra  morphisms, $\set{U,\dual{U};\beta}$ and  $\set{V,\dual{V};\mu}$ 
 symmetrizable dualizing pairs over $(\K,\tr)$. Let  $W=\Som[q]{k=1}{y_k\otimes
v_k}$ and $S=\Som[p]{k=1}{u_k\otimes x_k}$ be potentials  with $y_k\in U,\,x_k\in V,\, u_k\in
\widehat{\mrm{T}_{\K}(U)},\, v_k\in \widehat{\mrm{T}_{\K}(V)}$.
Then, for all $l\geq 1$ and potentials \mbox{$(f_l\otimes h)(W)
=\Som[q]{k=1}{f_l(y_k)\otimes h(v_k)}$}, \mbox{$(f\otimes
h_l)(S) =\Som[p]{k=1}{f(u_k)\otimes h_l(x_k)}$}, $(f\otimes h)(W)$ and $(f\otimes h)(S)$ we have:
$(h\otimes f_l)(\lperm W)=\lperm[l] (f_l\otimes h)(W)$,
$(h_l\otimes f)(\rperm S) = \rperm[l] (f\otimes h_l)(S)$ and 
 \\[-1.5\baselineskip]
\begin{align}\label{eq.lperm-morph-control}
& (h\otimes f)(\lperm
W)=\Som{x\in\lmQc[0]\,}{\Som[q]{k=1}\bilfh(\dual{x}\otimes
f(y_k))h(v_k)\otimes
x}=\Som{l\in\Ns}{\lperm[l]\Som[q]{k=1}{f_l(y_k)\otimes h(v_k)}} =
\Som{l\in\Ns}{\lperm[l] (f_l\otimes h)(W) }
\\[-1\baselineskip] \label{eq.rperm-morph-control}
& (h\otimes f)(\rperm S)=\Som{y\in\rmQc[0]\,}{ y\otimes
\Som[p]{k=1}{f(u_k) \bilfh(h(x_k)\otimes
\dual{y})}}=\Som{l\in\Ns}{ \rperm[l] \Som[p]{k=1}{f(u_k)\otimes
h_l(x_k)}} =\Som{l\in\Ns}{ \rperm[l] (f\otimes h_l)(S)}.
\end{align}
Moreover, any  morphism $\morph{\phi: \kQc}{\kQpc}$ of path
algebras  over the same symmetric algebra $(\K,\tr)$
sends cyclically equivalent potentials to cyclically equivalent
ones.
\end{Prop}

 \begin{prv}  Fix a natural number $l\geq 1$, in view of
 \eqref{eq-projbases-casimir} we have Casimir
 morphisms
$ \morph{\z[{(l)}]: \K}{B^{(l)} \otimes{\dual{B}}^{(l)}}$,
 $\morph{\z[{(l)}]': \K}{{\dual{B}}^{(l)} \otimes B^{(l)} }$  with
 $\z[{(l)}]\equiv \Som{y\in\rmQ[l]} y\otimes\dual{y}$  and
$\z[{(l)}]' \equiv \Som{y\in\lmQ[l]} \dual{x}\otimes x.$
We also consider the Casimir elements  $\z[\dual{U}\otimes
U]\equiv\Som{u\in\Ss}{u\otimes\dual{u}}$ and
$\z[V\otimes\dual{V}]\equiv \Som{v\in\Ss'}{\dual{v}\otimes v}$ associated with the  dualizing
pairs $\set{U,\dual{U};\beta}$ and $\set{V,\dual{V};\mu}$.
Let us  prove that \eqref{eq.lperm-morph-control} holds. By
Lemma~\ref{lem.dualizing-pair},  each bimodule morphism
$\morph{f_l: U}{B^{(l)}}$ is dualizing and its dual 
$\morph{\dual{f_l}: {\dual{B}}^{(l)}}{\dual{U}}$ is characterized by the relations: 
$\beta(u\otimes \dual{f}(\xi))= \bilf^l(f(u)\otimes \xi)
$ and $\beta(\dual{f}(\xi)\otimes u) =\bilf^l(\xi\otimes f(u))$
for all $u\in U$ and $\xi\in {\dual{B}}^{(l)}$. Let
$\mfr{m}_l:=\Som[q]{k=1}{f_l(y_k)\otimes h(v_k)}$, note that
$\Som[q]{k=1}{y_k\otimes h(v_k)}$ is   $\K$-central   as  image of
$\mfr{m}_l$ by the bimodule  morphism $\id_{U}\otimes h$. By
definition,  $\lperm W= \lderv(\id\otimes W\otimes \id)\circ
\z[\dual{U}\otimes U]=
\Som{u\in\Ss}{\Som[q]{k=1}{\beta(\dual{u}\otimes y_k)\mul
v_k\otimes u}}$.  In the following computation of
 $\lperm[l]\mfr{m}_{l}$,  use \eqref{eq.dualbases} for $\dual{f_l}(\dual{x}) $ in the second row and $f_l(u)$ in the height row:\\
$\begin{array}{rl} \lperm[l]\mfr{m}_{l} & :=\hspace{-0.2em}\lderv(\id\otimes
\mfr{m_l}\otimes\id)\circ\z[l]' =\hspace{-0.4em}
\Som{x\in\lmQ[l]}\Som[q]{k=1}{\bilf^l(\dual{x}\otimes
f_l(y_k))\otimes h(v_k) \otimes x}  
= \hspace{-0.4em}\Som{x\in\lmQ[l]}\Som[q]{k=1}{\beta(\dual{f_l}(\dual{x})
\otimes y_k)\otimes h(v_k) \otimes x} \\
& = \hspace{-0.4em}\Som{x\in\lmQ[l]}\Som[q]{k=1}{\beta(\Som{u\in\Ss}{\dual{u}
\beta(u\otimes\dual{f_l}(\dual{x}))}\otimes y_k)\otimes h(v_k)
\otimes x} \\ 
& =  \Som{u\in\Ss}\Som{x\in\lmQ[l]}\Som[q]{k=1}{\beta(\dual{u}
\otimes {\beta(u\otimes\dual{f_l}(\dual{x}))} y_k)\otimes h(v_k)
\otimes x} \\
& =   \Som{u\in\Ss}\Som{x\in\lmQ[l]} \lderv\parr{
\dual{u}\otimes \beta(u\otimes\dual{f_l}(\dual{x}))\mul
\Som[q]{k=1}{y_k \otimes h(v_k)}}\otimes x  \\
& = \Som{u\in\Ss}\Som{x\in\lmQ[l]} \lderv\parr{ \dual{u}\otimes
\Som[q]{k=1}{y_k \otimes h(v_k)}
\mul \beta(u\otimes\dual{f_l}(\dual{x})) }\otimes x \\
& =  \Som{u\in\Ss} \lderv\parr{ \dual{u}\otimes  \Som[q]{k=1}{y_k
\otimes h(v_k)}} \otimes \parr{\Som{x\in\lmQ[l]}
\beta(u\otimes\dual{f_l}(\dual{x}))x} \\
&  = \Som{u\in\Ss}\Som[q]{k=1}{ \beta(\dual{u}\otimes y_k)
  h(v_k)} \otimes \parr{\Som{x\in\lmQ[l]}
\beta(u\otimes\dual{f_l}(\dual{x}))x}  \\
& =  \Som{u\in\Ss} \Som[q]{k=1}{\beta(\dual{u}\otimes y_k) h(v_k)}
\otimes \parr{\Som{x\in\lmQ[l]} \beta(f_l(u)\otimes \dual{x})x}  
= \Som{u\in\Ss} \Som[q]{k=1}{\beta(\dual{u}\otimes y_k) \otimes
h(v_k)} \otimes f_l(u)  \\
& =  (h\otimes f_l)\parr{
\Som{u\in\Ss}\Som[q]{k=1}{\beta(\dual{u}\otimes y_k) \otimes v_k }
\otimes u} = (h\otimes f_l)(\lperm W).
\end{array} $ \\
 Hence, $\lperm[l]\mfr{m}_{l} =(f\otimes f_l)(\lperm W)$. In view
of \eqref{eq.dualbasesgen}, for all $\xi\in {\dual{B}}^{(l)}$,
$z=(z_l)_l\in\kQc$, with $z_l\in B^l$, we have:
$\bilfh(z\otimes\xi)=\bilf^l(z_l\otimes\xi )$ and 
$\bilfh(\xi\otimes z )=  \bilf^l(\xi\otimes z_l )$.  Note
that $\lmQ[0]=\set{1}=\ldmQ[0]$ and $\bilf(1\otimes
\alpha)=0=\bilf(\alpha\otimes 1)$ for all $\alpha\in B^{(t)}$ with
$t>0$. Then  letting $ T:=
\Som{x\in\lmQc[0]\,}{\Som[q]{k=1}\bilfh(\dual{x}\otimes
f(y_k))h(v_k)\otimes x} $  we obtain: \\
 $\begin{array}{rl}
 T := & \Som{x\in\lmQc[0]\,}{\Som[q]{k=1}\bilfh(\dual{x}\otimes
f(y_k))h(v_k)\otimes x}  =
\Som{x\in\lmQc[1]\,}{\Som[q]{k=1}\bilfh(\dual{x}\otimes
(\Som{l\in\Ns}f_l(y_k)))h(v_k)\otimes x}  \\ = &
\Som{l\in\Ns}\Som{x\in\lmQc[l]\,}{\Som[q]{k=1}\bilf^l(\dual{x}
\otimes f_l(y_k))h(v_k)\otimes x} =
\Som{l\in\Ns}{\lperm[l]\mfr{m}_{l}}=\Som{l\in\Ns}{(h\otimes
f_l)(\lperm W)}
\\  =& (h\otimes \Som{l\in\Ns}{f_l})(\lperm W) = (h\otimes f)(\lperm W).
\end{array} $ \\
Hence, \eqref{eq.lperm-morph-control} is proved. Dually, \eqref{eq.rperm-morph-control}  also holds.

For the last part of the proposition, let $\morph{\phi:
\kQc}{\kQpc}$ be a morphism of path algebras.  As in
Proposition~\ref{prop.prolong}, $\phi$ is continuous and induced by a family  of $\K$-bimodule morphisms  $\morph{\phi_l: B}{B'^{(l)}},\, l\geq 1$. For all
potential $\mfr{m}=\Som[q]{k=1}{y_k\otimes v_k}\in B^{(d)}$ with
$y_k\in\B$, letting
$\mfr{m}_l:=\Som[q]{k=1}{\phi_l(y_k)\otimes\phi(v_k)}$ yields that:
$\phi(\mfr{m}-\lperm
\mfr{m})=\phi(\mfr{m})-\phi(\lperm\mfr{m})=\Som{l
\in\Ns}{\mfr{m}_l}-\Som{l\in\Ns}{\lperm[l]\mfr{m}_l}=
\Som{l\in\Ns}(\mfr{m}_l-\lperm[l]\mfr{m}_l)$.
Therefore, $\phi$ sends any skew commutator to an element of the closed  $\Zc(\K)$-module of skew commutators in $\kQpc$, and since $\phi$ is continuous we conclude that $\phi$ sends the
closed $\Zc(\K)$-module  of skew commutators in $\kQc$ to the closed
$\Zc(\K)$-module of skew commutators in $\kQpc$. Hence $\phi$ sends cyclically equivalent potentials to cyclically equivalent ones.
\end{prv}

In the next step we  develop a differential calculus on
potentials. Consider the topological $\K$-bimodule
\\[-1\baselineskip]
 \begin{Center}
 $\kQc\bbar{\otimes}\kQc=\prd{d,e\geq 0}{(B^{(d)}\otimes_{\kk} B^{(e)})}$,
 \end{Center}
having as system of open neighborhoods of $0$   the
 subbimodules $\prd{d+e\geq m}{(B^{(d)}\otimes_{\kk} B^{(e)})}, m\geq 0$. 
Thus $\kQ\otimes_{\kk}\kQ$ is dense in
$\kQc\bbar{\otimes}\kQc$. When we fix projective bases
$(\lmQ[l],\ldmQ[l])$ and $(\rmQ[l],\rdmQ[l])$, we equally lift the
corresponding Casimir morphisms to the following:
 $\morph{\zb[(l)]: \kk}{ B^{(l)} \tenss[\kk]
{\dual{B}}^{(l)}}$  and $\morph{\zb[(l)]': \kk}{ {\dual{B}}^{(l)}
\tenss[\kk] B^{(l)}}$, with $\zb[(l)]\equiv\Som{y\in\rmQ[l]}{
y\tenss[\kk] \dual{y}}$ and
$\zb[(l)]'\equiv\Som{x\in\lmQ[l]}{\dual{x}\tenss[\kk] x}$.
 Left and  right derivative operator on $\kQc$ are naturally extended to derivative operators
$\morph{\lderv, \rderv: \kQd \otimes
(\kQc\bbar{\otimes}\kQc)}{\kQc\bbar{\otimes}\kQc}$ as follows: for
all $\xi\in \kQd$ and $v_1\tenss[\kk] v_2 \in
\kQc\bbar{\otimes}\kQc$  we have 
\begin{center}  $\lderv[\xi](v_1\tenss[\kk]
v_2)=\lderv(\xi\otimes(v_1\tenss[\kk] v_2)):=(\lderv[\xi]
v_1)\tenss[\kk] v_2$ and
  $\rderv[\xi](v_1\tenss[\kk] v_2)=\rderv((v_1\tenss[\kk] v_2)\otimes \xi)
  :=v_1\tenss[\kk](\rderv[\xi] v_2)$.
  \end{center}
Let $\xi\in\kQd$;  we define 
$\morph{\lgrad[\xi],\, \rgrad[\xi]:\,\kQc}{\kQc\bbar{\otimes}
\kQc}$  and
$\morph{\bilbox:\,\kQc\bbar{\otimes}\kQc
\times \, \kQc}{\kQc}$ as follows: for all $u,w,
v=(v_l)_l \in\kQc$ with $ v_l \in B^{(l)}, l\in\N$, we
have:
\\[-1\baselineskip]
\begin{equation}\label{eq-grad-box}
\begin{aligned}
& \lgrad[\xi](v)=\Som{l\in\N}{\lgrad[\xi](v_l)}
=\lderv[\xi]\Som{l\in\N}{\Som[l-1]{t=0}\Som{x\in\lmQ[t]}(\lderv[\dual{x}]v_l)\tenss[\kk]
x}  ,\; \, \rgrad[\xi]v=\Som{l\in\N}{\rgrad[\xi](v_l)}
=\rderv[\xi]\Som{l\in\N}{\Som[l-1]{t=0}\Som{y\in\rmQ[t]}y\tenss[\kk](\rderv[\dual{y}]v_l)},\\
& \displaystyle{(u\tenss[\kk] v)\bilbox w= u\mul w\mul v}.
 \end{aligned}
\end{equation}
 
\ParIt{Notations} Let $d,k\geq 1$, for a product ${u_{k,0}\mul u_{k,1}\dotsm\mul u_{k,d}}$ of elements of $\kQc$, we put: $u_{k,<0}:=1=u_{k,>d}$, 
$u_{k,<r}:=u_{k,0}\dotsm u_{k,r-1}$ and $u_{k,>r}:=u_{k,r+1}\dotsm
u_{k,d}$ for $0<r<d$. Recall that we have a weakly dualizing pair $\set{\kQd,\kQc;\bilfh}$, 
we also put $\bilfh[r]:=\bilfh$ for $r>0$; $\morph{\bilfh[0]:=\bilf^0:\K\otimes\K}{\K}$ is the multiplication of $\K$.

\begin{Lem}[{cyclic Leibniz rule}]\label{lem.cyclicLeibnizrule}
Let $\mfr{m}=\Som[n]{k=1}{u_{k,0}\mul u_{k,1}\dotsm\mul u_{k,d}}$
be a potential on $\mQ$ with $d\geq 1$ and \mbox{$u_{k,r}\in\kQc$}. 
Then for all $\xi\in\dual{B}$   the following cyclic Leibniz rule holds:
\\[-1\baselineskip]
 \begin{equation}\label{eq.cyclicLeibnizrule}
 \begin{aligned}
  \cderv[\xi]\mfr{m} & =
\Som[d]{r=0}\Som{x\in\lmQc[0]\,}\Som[n]{k=1}{
\lgrad[\xi](\bilfh[r](\dual{x}\otimes u_{k,<r})u_{k,r})\bilbox
(u_{k,>r}  x)}  \\[-0.5\baselineskip]
 & = \Som[d]{r=0}\Som{y\in\rmQc[0]\,}\Som[n]{k=1}{
\rgrad[\xi](u_{k,d-r}\bilfh[r](u_{k,>d-r}\otimes \dual{y}))\bilbox
(y u_{k,<d-r})}.
\end{aligned}
\end{equation}
\end{Lem}

\begin{prv} For \mbox{$k\in\set{1,\ldots, n}$}, let
\mbox{$u_{k,r}=\hspace{-0.4em}\Som{l\in\Ns}{u_{k,r,l}}$} with 
\mbox{$u_{k,r,l}\in B^{(l)}, l\geq 1$ }. Recall that
\mbox{$\lmQc[0]=\hspace{-0.4em}\Union{t\in\N}{\lmQ[t]}$} with \mbox{$\lmQ[0]=\set{1}=\ldmQ[0]$}
and in view \eqref{eq.dualbasesgen},  for all
\mbox{$x\in\Union{t\in\Ns}{\lmQ[t]}$} and \mbox{$u\in\kQc$} we have
\mbox{$\bilf^0(\dual{x}\otimes u)=0$}.  We will use induction on \mbox{$d\geq 1$} to
establish 
\mbox{$\; (L): \: \cderv[\xi]\mfr{m} =\hspace{-0.4em}
\Som[d]{r=0}\Som{x\in\lmQc[0]\,}\Som[n]{k=1}{
\lgrad[\xi](\bilfh[r](\dual{x}\otimes u_{k,<r})u_{k,r})\bilbox
(u_{k,>r} x)}. $}

Assume that $d=1$. Then
$\mfr{m}=\Som{l,l'\in\Ns}{\Som[n]{k=1}{u_{k,0,l}\, \mul
u_{k,1,l'}}}= \Som{p\geq 2}\mfr{m}_p$ where 
$\mfr{m}_p:=\Som{\begin{subarray}{c} l+l'=p \, \\ l,l'\in\Ns
\end{subarray}}{\Som[n]{k=1}{u_{k,0,l}\, \mul u_{k,1,l'}}}$ is a degree-$p$ homogeneous potential.
The cyclic permutation of $\mfr{m}$ is given by
\[
\begin{array}{rl}
\cperm\mfr{m} & = \Som{p\geq 2\,}\cperm\mfr{m}_p =\Som{p\geq
2\,}\Som[p-1]{t=0}\lperm[t]\mfr{m}_p = \Som{p\geq
2\,}\Som{l+l'=p\,}\Som[p-1]{t=0}\Som{x\in\lmQ[t]\,}\Som[n]{k=1}\lderv(\dual{x}\otimes
u_{k,0,l}\, \mul u_{k,1,l'} )x \\
& = \Som{p\geq
2\,}\Som{l+l'=p\,}\Som[l-1]{t=0}\Som{x\in\lmQ[t]\,}\Som[n]{k=1}\lderv(\dual{x}\otimes
u_{k,0,l}\, \mul u_{k,1,l'} )x + \Som{p\geq
2\,}\Som{l+l'=p\,}\Som[p-1]{t=l}\Som{x\in\lmQ[t]\,}\Som[n]{k=1}\lderv(\dual{x}\otimes
u_{k,0,l}\, \mul u_{k,1,l'} )x. 
\end{array}
\]
Let $\cperm\mfr{m} =S_1+S_2$ where $S_1$ is  the first term in the last line above and $S_2$ the second.   In view \eqref{eq-grad-box} defining  operators $\lgrad[\xi]$ and $\bilbox$,  we compute
 the left derivative $\lderv[\xi]S_1$ as follows.\\
$\begin{array}{rl} \lderv[\xi] S_1:= & \lderv[\xi]\Som{p\geq
2\,}\Som{l+l'=p\,}\Som[l-1]{t=0}\Som{x\in\lmQ[t]\,}\Som[n]{k=1}\lderv(\dual{x}\otimes
u_{k,0,l}\, \mul u_{k,1,l'} )x  \\ = & \lderv[\xi]
\Som{l\in\Ns\,}\Som[l-1]{t=0\,}\Som{x\in\lmQ[t]\,}\Som[n]{k=1}\lderv(\dual{x}\otimes
u_{k,0,l})\, \mul (\Som{l'\in\Ns\,} u_{k,1,l'}) x  \\
= & \lderv[\xi]
\Som{l\in\Ns}\Som[l-1]{t=0\,}\Som{x\in\lmQ[t]\,}\Som[n]{k=1}\lderv(\dual{x}\otimes
u_{k,0,l})\, \mul  u_{k,1} x =\Bigl(\Som[n]{k=1}\lderv[\xi]
\Som{l\in\Ns}\Som[l-1]{t=0}\Som{x\in\lmQ[t]\,}(\lderv(\dual{x}\otimes
u_{k,0,l})\tenss[\kk]x )\Bigr)\, \bilbox  u_{k,1}  \\
= & \Som[n]{k=1}(\lgrad[\xi] u_{k,0}) \bilbox u_{k,1} =
\Som{x\in\lmQc[0]}\Som[n]{k=1}\lgrad[\xi](\bilfh[0](\dual{x}\otimes
u_{k,<0}) u_{k,0}) \bilbox (u_{k,>0}x). \hfill{ } \hfill{(\ast)}
\end{array}$

And,  using \eqref{eq-projbaseskQ} describing   projective
bases   associated with $\mQ$,  for all $t=l+s\geq l$ with $0\leq s <l'$, we have
$\lmQ[t]=\lmQ[l]\otimes \lmQ[s]:=\set{x\otimes z\; :
x\in\lmQ[l],\, z\in\lmQ[s]}$;  thus the left derivative $
\lderv[\xi] S_2$ is
computed as follows. \\
$ \begin{array}{rl}
 \lderv[\xi] S_2:= &  \lderv[\xi] \Som{p\geq
2\,}\Som{l+l'=p\,}\Som[p-1]{t=l}\Som{x\in\lmQ[t]\,}\Som[n]{k=1}\lderv(\dual{x}\otimes
u_{k,0,l}\, \mul u_{k,1,l'} )x  \\
= & \lderv[\xi]
\Som{l'\in\Ns\,}\Som{l\in\Ns\,}\Som[l'-1]{s=0}\Som{x\in\lmQ[l]\,}\Som{z\in\lmQ[s]\,}
\Som[n]{k=1}\lderv(\dual{z}\otimes\dual{x}\otimes u_{k,0,l}\, \mul
u_{k,1,l'} )x z \\
= & \lderv[\xi] \Som{l'\in\Ns\,} \Som{l\in\Ns\,}
\Som[l'-1]{s=0}\Som{x\in\lmQ[l]\,}\Som{z\in\lmQ[s]\,}
\Som[n]{k=1}\lderv(\dual{z}\otimes\bilfh(\dual{x}\otimes
u_{k,0,l}) \mul u_{k,1,l'} )x z \\
= & \lderv[\xi] \Som{l'\in\Ns\,}
\Som[l'-1]{s=0}\Som{z\in\lmQ[s]\,}
\Som[n]{k=1}\lderv(\dual{z}\otimes
(\Som{l\in\Ns\,}\Som{x\in\lmQ[l]\,}\bilfh(\dual{x}\otimes
u_{k,0,l})) \mul u_{k,1,l'} x) z \\
= & \lderv[\xi] \Som{l'\in\Ns\,}
\Som[l'-1]{s=0}\Som{z\in\lmQ[s]\,}
\Som[n]{k=1}\lderv(\dual{z}\otimes
(\Som{x\in\lmQc[0]\,}\bilfh(\dual{x}\otimes u_{k,0})) \mul
u_{k,1,l'} x) z \\
= & \Som{x\in\lmQc[0]\,} \lderv[\xi] \Som{l'\in\Ns\,}
\Som[l'-1]{s=0}\Som{z\in\lmQ[s]\,}
\Som[n]{k=1}\lderv(\dual{z}\otimes (\bilfh(\dual{x}\otimes
u_{k,0})) \mul u_{k,1,l'} x) z  \\
= & \Som{x\in\lmQc[0]\,} \lderv[\xi] \Som{l'\in\Ns\,}
\Som[l'-1]{s=0}\Som{z\in\lmQ[s]\,}
\Som[n]{k=1}\lderv(\dual{z}\otimes (\bilfh(\dual{x}\otimes
u_{k,0})) \mul u_{k,1,l'}) x z \\
= & \Som{x\in\lmQc[0]\,} \lderv[\xi] \Som{l'\in\Ns\,}
\Som[l'-1]{s=0}\Som{z\in\lmQ[s]\,}
\Som[n]{k=1}(\lderv(\dual{z}\otimes (\bilfh(\dual{x}\otimes
u_{k,0})) \mul u_{k,1,l'})\tenss[\kk] z)\bilbox x  \\ = &
\Som{x\in\lmQc[0]\,}\Som[n]{k=1}\lgrad[\xi](
(\bilfh(\dual{x}\otimes u_{k,0})) \mul u_{k,1})\bilbox x =
\Som{x\in\lmQc[0]\,}\Som[n]{k=1}\lgrad[\xi](
(\bilfh(\dual{x}\otimes u_{k,<1})) \mul u_{k,1})\bilbox
(u_{k,>1}x). \hfill{ \quad } \hfill{(\ast\ast)}
\end{array}
$ \\
Hence, combining $(\ast\ast)$ and $(\ast)$ above, $(L)$ is proved for $d=1$. For the induction step, assume $d>1$
and the result true for $d-1$. We write: $\mfr{m} =
\Som[n]{k=1}{(u_{k,0}\mul u_{k,1})u_{k,2}\dotsm\mul u_{k,d}}$, by
a direct application of the
induction assumption and of the proof of the case "$d=1$" above we get:\\
 $
\begin{array}{rl}
\cderv[\xi]\mfr{m} & = \cderv[\xi] \Som[n]{k=1}{(u_{k,0}\mul
u_{k,1})u_{k,2}\dotsm\mul  u_{k,d}} \\
& = \Som{x\in\lmQc[0]\,}\Som[n]{k=1}{
\lgrad[\xi](\bilfh[0](\dual{x}\otimes
1)\mul(u_{k,0}u_{k,1}))\bilbox (u_{k,>1} x)} +
\Som[d]{r=2}\Som{x\in\lmQc[0]\,}\Som[n]{k=1}{
\lgrad[\xi](\bilfh[r](\dual{x}\otimes u_{k,<r})u_{k,r})\bilbox
(u_{k,>r} x)} \\
& = \Som[n]{k=1}{ (\lgrad[\xi](u_{k,0}u_{k,1}))\bilbox u_{k,>1}} +
\Som[d]{r=2}\Som{x\in\lmQc[0]\,}\Som[n]{k=1}{
\lgrad[\xi](\bilfh[r](\dual{x}\otimes u_{k,<r})u_{k,r})\bilbox
(u_{k,>r} x)} \\
 & = \Som[1]{r=0}\Som{x\in\lmQc[0]\,}\Som[n]{k=1}{
\lgrad[\xi](\bilfh[r](\dual{x}\otimes u_{k,<r})u_{k,r})\bilbox
(u_{k,>r} x)} + \Som[d]{r=2}\Som{x\in\lmQc[0]\,}\Som[n]{k=1}{
\lgrad[\xi](\bilfh[r](\dual{x}\otimes u_{k,<r})u_{k,r})\bilbox
(u_{k,>r} x)} \\
& = \Som[d]{r=0}\Som{x\in\lmQc[0]\,}\Som[n]{k=1}{
\lgrad[\xi](\bilfh[r](\dual{x}\otimes u_{k,<r})u_{k,r})\bilbox
(u_{k,>r}  x)}.
\end{array} $ \\
Hence,  $(L)$ is proved. Dually, the Leibniz rule involving only the operator
$\rgrad[\xi]$  holds.
\end{prv}

\begin{Lem}[cyclic chain-rule]\label{lem.chainrule}
Let $\morph{\phi:\kQc}{\kQpc}$ be a morphism of path algebras for
a given modulated quiver $\mQ'=(B',\K,\tr)$.  Then
for all potential $\mfr{m}$ on $\mQ$ and all $\xi\in\dual{B'}$ we
have:
\\[-1\baselineskip]
\begin{equation} \label{eq-chainrule}
\cderv[\xi]\phi(\mfr{m})=\Som{y\in\rmQ}{(\lgrad[\xi]
\phi(y))\bilbox
\phi(\cderv[\dual{y}]\mfr{m})}=\Som{x\in\lmQ}{(\rgrad[\xi]
\phi(x))\bilbox \phi(\cderv[\dual{x}]\mfr{m})}.
\end{equation}
\end{Lem}

\begin{prv} Since each potential $W$ decomposes as  sum  
of homogeneous potentials,  it suffices to prove the chain-rule for
homogeneous potentials.  Thus we may
assume that $\mfr{m}$ is homogeneous and write $\mfr{m}=
\Som[n]{k=1}{u_{k,0}\mul u_{k,1}\dotsm\mul u_{k,d}}\in\B^{(d+1)}$
for some $d \geq 1$ and   $u_{k,r}\in B $, $r=0,1\dots,d$. We have \mbox{$\phi(\mfr{m})=
\Som[n]{k=1}{\phi(u_{k,0})\mul \phi(u_{k,1})\dotsm\mul
\phi(u_{k,d})}$}. The bimodule $B'$ is part of a symmetrizable
dualizing pair  $(B',\dual{B'};\bilf')$ over $(\K,\tr)$, and in
view of \eqref{eq-projbaseskQ} we have projective bases
$(\lmQp[t],\ldmQp[t])$, $(\rmQp[t],\rdmQp[t])$,
$(\lmQpc[t],\ldmQpc[t])$ and $(\rmQpc[t],\rdmQpc[t])$ with
$t\in\N$, $\lmQp[0]=\set{1}=\ldmQp[0]$ and
$\rmQp[0]=\set{1}=\rdmQp[1]$. Now, let  $\xi\in \dual{B'}$, we
 will establish the chain-rule~\eqref{eq-chainrule} in terms of
 operator $\lgrad[\xi]$. The cyclic Leibniz rule
\eqref{lem.cyclicLeibnizrule} yields that:
\\[-0.5\baselineskip]
 \begin{equation*}
 \begin{aligned}
  \cderv[\xi]\phi(\mfr{m}) & =
\Som[d]{r=0}\Som{x'\in\lmQpc[0]\,}\Som[n]{k=1}{
\lgrad[\xi](\bilfph[r](\dual{x'}\otimes
\phi(u_{k,<r}))\phi(u_{k,r}))\bilbox (\phi(u_{k,>r})  x')} \\[-0.5\baselineskip]
& = \Som[d]{r=0}\Som{x'\in\lmQpc[0]\,}\Som[n]{k=1}{
\lgrad[\xi](\phi(\bilfph[r](\dual{x'}\otimes \phi(u_{k,<r})) \mul
u_{k,r}))\bilbox (\phi(u_{k,>r})  x')}, \qquad (L)
\end{aligned}
\end{equation*}
where  $u_{k,<r}=u_{k,0}\dotsm u_{k,r-1}$ and
$u_{k,>r}=u_{k,r+1}\dotsm u_{k,d}$ for $0<r<d$; $u_{k,<0}=1=u_{k,>d}$; 
$\bilfph[r]=\bilfph$ for $r>0$ and $\bilfph[0]=\bilf'^0:\K\otimes \K \to \K$ is the multiplication of
$\K$. By \eqref{eq.dervdualbasesKQc},
each $\bilfph[r](\dual{x'}\otimes \phi(u_{k,<r})) \mul
u_{k,r} \in B$ expresses as: $\bilfph[r](\dual{x'}\otimes
\phi(u_{k,<r})) \mul u_{k,r} = \Som{y\in\rmQ}y\mul
\lderv[\dual{y}]\bigl(\bilfph[r](\dual{x'}\otimes \phi(u_{k,<r}))
\mul u_{k,r}\bigr)$. But $(u a\tenss[\kk] v)\bilbox
w =(u \tenss[\kk] v)\bilbox aw=u a w v$ for all $u,v\in\kQc$ and
$a\in\K$. Hence, the cyclic Leibniz rule~$(L)$ yields: \\
 $\begin{array}{rl}
  \cderv[\xi]\phi(\mfr{m}) & =
\Som[d]{r=0}\Som{x'\in\lmQpc[0]\,}\Som[n]{k=1}{
\lgrad[\xi]\Bigl(\phi\bigl( \Som{y\in\rmQ}y\mul
\lderv[\dual{y}](\bilfph[r](\dual{x'}\otimes \phi(u_{k,<r})) \mul
u_{k,r})  \bigr)\Bigr) \bilbox (\phi(u_{k,>r}) x')} \\
& = \Som{y\in\rmQ} \Som[d]{r=0}\Som{x'\in\lmQpc[0]\,}\Som[n]{k=1}{
\lgrad[\xi](\phi(y)) \bilbox (
\lderv[\dual{y}](\bilfph[r](\dual{x'}\otimes \phi(u_{k,<r})) \mul
u_{k,r}) \phi(u_{k,>r}) x')} \\
& = \Som{y\in\rmQ} \lgrad[\xi](\phi(y))
 \bilbox \bigl( \Som[d]{r=0}\Som{x'\in\lmQpc[0]\,}\Som[n]{k=1}{
\lderv[\dual{y}](\bilfph[r](\dual{x'}\otimes \phi(u_{k,<r})) \mul
u_{k,r}) \phi(u_{k,>r}) x'} \bigr)\\
& = \Som{y\in\rmQ} \lgrad[\xi](\phi(y))
 \bilbox \bigl( \lderv[\dual{y}] \Som[d]{r=0}\Som{x'\in\lmQpc[0]\,}\Som[n]{k=1}{
\bilfph[r](\dual{x'}\otimes \phi(u_{k,<r})) \mul u_{k,r}
\phi(u_{k,>r}) x'} \bigr) \qquad (\ast)
  \end{array} $\\
For each $r>0$, we have dualizing pairs $\set{B^{(r)},{\dual{B}}^{(r)}}$ and $\set{B\otimes
B^{(d-r)},{\dual{B}}^{(d-r)}\otimes \dual{B}}$, and for $\K$-bimodule morphisms $\morph{f_r:=\phi\restr{B^{r}}: B^{r}}{\kQpc}$ and
$\morph{h_r:=\id_{B}\otimes (\phi\restr{B^{(d-r)}}) : B\otimes
B^{(d-r)} }{\kQpc}$ we observe that 
$(f_r\otimes h_r)(\mfr{m})=\Som[n]{k=1}{
 \phi(u_{k,<r}) \mul u_{k,r}
\phi(u_{k,>r})}$.
whence, invoking relation~\eqref{eq.lperm-morph-control} of
Proposition~\ref{prop.morph-et-perm} to control the left
permutation of $\mfr{m}$ with respect the pair  $\set{B^{(r)},{\dual{B}}^{(r)}}$, we get:
\\[-1\baselineskip]
\begin{center}
$(h_r\otimes f_r)(\lperm[r]\mfr{m}) =
\Som{x'\in\lmQpc[0]\,}\Som[n]{k=1}{ \bilfph[r](\dual{x'}\otimes
\phi(u_{k,<r})) \mul u_{k,r} \phi(u_{k,>r}) x'}$.
\end{center}  For each
$y\in\rmQ$, let us compute the term $S_{y,r}:= \lderv[\dual{y}]
\Som{x'\in\lmQpc[0]\,} \Som[n]{k=1}{ \bilfph[r](\dual{x'}\otimes
\phi(u_{k,<r})) \mul u_{k,r}
\phi(u_{k,>r}) x'}$: \\
$ \begin{array}{rl} S_{y,r} & =  \lderv[\dual{y}]
\Som{x'\in\lmQpc[0]\,} \Som[n]{k=1}{ \bilfph[r](\dual{x'}\otimes
\phi(u_{k,<r})) \mul u_{k,r} \phi(u_{k,>r}) x'} =
\lderv[\dual{y}]((h_r\otimes f_r)(\lperm[r]\mfr{m}))
\\
&= \lderv[\dual{y}] (h_r\otimes f_r)\bigl( \Som[n]{k=1}{
 \Som{z\in\lmQ[r]}\bilfh(\dual{z}\otimes u_{k,<r}) \mul u_{k,r}
u_{k,>r} z}\bigr)  =   \Som[n]{k=1}{
 \Som{z\in\lmQ[r]} \lderv[\dual{y}](\bilfh(\dual{z}\otimes u_{k,<r}) \mul
 u_{k,r}) \phi(u_{k,>r} z)} \\
 & = \Som[n]{k=1}{
 \Som{z\in\lmQ[r]} \phi(  \lderv[\dual{y}](\bilfh(\dual{z}\otimes u_{k,<r}) \mul
 u_{k,r}) \mul u_{k,>r} z)} =  \phi( \Som[n]{k=1}{
 \Som{z\in\lmQ[r]}  \lderv[\dual{y}](\bilfh(\dual{z}\otimes u_{k,<r}) \mul
 u_{k,r}) \mul u_{k,>r} z}) \\
  & =  \phi( \lderv[\dual{y}] \Som[n]{k=1}{  \Som{z\in\lmQ[r]} \bilfh(\dual{z}\otimes u_{k,<r}) \mul
 u_{k,r} \mul u_{k,>r} z}) \\ & = \phi(\lderv[\dual{y}]
 \lperm[r]\mfr{m}). \qquad (\ast\ast)
\end{array} $ \\
From $(\ast\ast)$ and $(\ast)$  we get:
\mbox{$\cderv[\xi]\phi(\mfr{m}) = \hspace{-0.4em}\Som{y\in\rmQ}
\lgrad[\xi](\phi(y)) \bilbox \bigl( \phi(\lderv[\dual{y}]
\Som[d]{r=0} \lperm[r]\mfr{m} ) \bigr) = \hspace{-0.4em} \Som{y\in\rmQ}
\lgrad[\xi](\phi(y))  \bilbox \bigl( \phi(\lderv[\dual{y}]
\cperm\mfr{m} ) \bigr) $} $ = \Som{y\in\rmQ}
\lgrad[\xi](\phi(y))  \bilbox \bigl( \phi(\cderv[\dual{y}]\mfr{m}
) \bigr).$
Dually,  the chain-rule  in terms of
 operator $\lgrad[\xi]$  also holds.
\end{prv}


\subsection{The reduction process}
Throughout this subsection, let $(\mQ,\mfr{m})$ be a modulated quiver with potential where $\mQ=(B,\K,\tr)$. We have  $\mfr{m}=(\mfr{m}_l)_{l\geq 2}$ with $\mfr{m}_l\in B^{(l)}$; put  $\mfr{m}=\mfr{m}_2 + \mfr{m}_{(3)}$ with  $\mfr{m}_{(3)}=(\mfr{m}_l)_{l\geq 3}$.
Let $U=(\rderv \mfr{m}_{2})(\dual{B})$, $V=(\lderv \mfr{m}_{2})(\dual{B})$, $\triv{B}=(\cderv \mfr{m}_{2})(\dual{B})$,  $\bbar{B}=B/\triv{B}$ and  
\mbox{$\dual{\bbar{B}}=\ker(\cderv\mfr{m}_{2}) \subset \dual{B}$}. Applying 
Proposition~\ref{prop.dec.pot-mQ}  and the fact that $\mfr{m}$ is  $2$-loop free 
(Definition~\ref{defn.mQp}), we get  the following. 
\begin{Rem}\label{rem.bilf-trivpart-redpart} $U\cap V=0$, $\triv{B}=U\oplus V$ and over $(\K,\tr)$ we have induced symmetrizable dualizing pairs $\set{U,V}$, $\set{\triv{B},\triv{B}}$ and $\set{\bbar{B}, \dual{\bbar{B}};\bilfb}$, together 
with mutually dual  canonical exact sequences
$\Seq[/Dist=1]{(\vartheta):\, \dual{\bbar{B}}
\to[/dir=right hook->] \dual{B}
 \to[/dir=->>,/above={\msf{p}}] \triv{B}}$ and $ 
 \Seq[/Dist=1]{ (\dual{\vartheta}):\, \triv{B} \to[/dir=right
hook->,/above={\dual{\msf{p}}}] B
 \to[/dir=->>,/above={\rho}] \bbar{B}} $ which split as sequences of left and right  $\K$-modules.
 For all $\xi\in\dual{\bbar{B}}, x\in\triv{B}$ we have: $\bilf(\xi\otimes x)=\bilfb(\xi\otimes \rho(x))=0=\bilfb(\rho(x)\otimes \xi)=\bilf(x\otimes \xi)$.
\end{Rem}

\begin{Defn}\label{defn.trivpart-reducedpart}  
 The bimodule $\triv{B}$ is the  \emph{trivial part} of $B$, $\triv{\mQ}:=\set{\triv{B},\K,\tr}$  the \emph{trivial part} of $\mQ$ and   $\triv{(\mQ,\mfr{m})} :=(\triv{\mQ},\mfr{m}_{2})$  the \emph{trivial part} of $(\mQ,\mfr{m})$.  The bimodule $\bbar{B}:=B/\triv{B}$ is the \emph{reduced part} of $B$ and  $\red{\mQ}=\bbar{\mQ}:=(\bbar{B},\K,\tr)$ the \emph{reduced  part} of $\mQ$. The modulated quiver with potential $(\mQ,\mfr{m})$
 is   \emph{reduced}  if $\mfr{m}$ belongs to $\J[3]{\kQc}$, that is, if
$\mfr{m}_{2}=0$.  $(\mQ,\mfr{m})$ is   \emph{trivial} if the reduced
part of $B$ is zero.
The \emph{trivial part of $(\mQ,\mfr{m})$  splits} if    $\triv{B}$ is a direct summand in $B$.
\end{Defn}


\ParIt{A note on presentations of Jacobian algebras} The first obstruction to  reduce a modulated quiver with potential is  of the same nature as the obstruction to the presentation of finite-dimensional algebras over non algebraically closed   fields by modulated quivers with relations (see \cite{Ben}). Indeed,  let  $A:=\Jc(\mQ,\mfr{m})$ and $\J{A}:=\Frac{\J{\kQc}}{\J{\mfr{m}}}$;  if $A$ admits a presentation by a modulated quiver with relations,  then   $A$  is an  $(A/\J{A})$-bimodule and  $\J[2]{A}$ is a direct summand in $\J{A}$. Now let $\what{J}=\J{\kQc}$ and observe that  $\triv{B}=B\cap(\what{J}^2+\J{\mfr{m}})$ and we have: $\Frac{\J{A}}{\J[2]{A}} =\Frac{(\what{J}/\J{\mfr{m}})}{((\what{J}^2+\J{\mfr{m}})/\J{\mfr{m}})} \cong
 \Frac{\what{J}}{(\what{J}^2+\J{\mfr{m}})} = \Frac{(B+\what{J}^2)}{(\what{J}^2+\J{\mfr{m}})} \cong 
\Frac{B}{(B\cap(\what{J}^2+\J{\mfr{m}})) } = \Frac{B}{\triv{B}} =\bbar{B}$, illustrating  $\msf{(1.ii)}$ below.

\begin{Defn}\label{defn.trivpart-reduction}
 A  \emph{reduction} on $(\mQ,\mfr{m})$    is a path algebra epimorphism $\morph{\phi: \kQc }{\kQpc}$ (with $\mQ'=(B',\K,\tr)$) such that $\phi(\mfr{m})$ is reduced and  the following conditions hold.  
 \begin{itemize}
 \item [$\msf{(1.i)}$] $\Ker(\phi)$  is the closed ideal generated by the image of a $\K$-bimodule morphism  $\morph{f=\Psmatr{\id \\ f'}: \triv{B}}{\triv{B}\oplus \J[2]{\kQc}}$ with $\im(f)\subset(\cderv \mfr{m})(\dual{B})$.
 \item [$\msf{(1.ii)}$] Let $\morph{\pi: \kQc}{\kQc/\Ker(\phi)}$ be the natural projection and $\morph{\bbar{\rho}: \Frac{(B+\Ker(\phi))}{\Ker(\phi)}}{\bbar{B}}$ the $\K$-bimodule epimorphism  with $\rho=\bbar{\rho}\circ \pi\restr{B}$. Then $\bbar{\rho}$  has a right inverse   $\morph{ \bbar{\rho}':  \bbar{B}}{(B+\Ker(\phi))/\Ker(\phi)}$ which lifts to a left (respectively, right) $\K$-linear  map $\morph{ \rho':  \bbar{B}}{B}$  such that $ \phi\circ \rho' =\phi_1\circ \rho'$. 
\end{itemize}
A reduction $\phi$ splits whenever $\morph{ \rho':  \bbar{B}}{B}$  is a bimodule morphism; in this case $\triv{(\mQ,\mfr{m})}$ splits.
\end{Defn}

\begin{Lem} \label{lem.trivIdeal-closedIdeal}
\begin{itemize}
\item[$\msf{(a)}$]
Let  $J_0$ be any closed ideal in $\kQc$  satisfying $\msf{(1.i)}$ above. Then
$(\cderv \mfr{m})(\dual{B}) \subset f(\triv{B})+(\cderv \mfr{m})(\dual{\bbar{B}}) $ and consequently,
$\J{\mfr{m}}=\bbar{J_0+\ggen{ (\cderv \mfr{m})(\dual{\bbar{B}}) }}$.
\item[$\msf{(b)}$]  Let $I$ and $I'$ be two $\kk$-modules in  a path algebra $A$ with $J$-adic topology, where $J$ is the complete arrow ideal in $A$.  Then the closed module $\bbar{\bbar{I}+\bbar{I'}}$ coincides with $\bbar{I+I'}$. If $I$ and $I'$ are   ideals in $A$  then $\bbar{I}\,\bbar{I'}\subset \bbar{I I'}$,  consequently, if $I'\subset\bbar{I+JI'+I'J}$ then  $I'\subset \bbar{I}$.
\end{itemize}
\end{Lem}
\begin{prv} Since part $\msf{(a)}$ is a direct consequence of the assumptions and the definition of   $\dual{\bbar{B}}$ and  $\J{\mfr{m}}$,
 we turn to part $\msf{(b)}$. Under  $J$-adic topology, the closure of a subset $S\subset A$ is  $\bbar{S}=\inter{l\geq 0}{}{(S+J^l)}=\inter{l\geq l_0}{}{(S+J^l)}$, with $l\in \N$.  Thus each subset $S +J^l$ is closed,  and for  two $\kk$-modules $I,I'\subset A$ we have:
$
\bbar{\bbar{I}+\bbar{I'}}= \inter{l\geq 0}{}{\bigl( \inter{s\geq l}{}{(I+J^s)}+
\inter{t\geq l}{}{(I'+J^t)} +J^l \bigr)} \subset \inter{l\geq 0}{}{(I+J^l+I'+J^l+J^l)}=\inter{l\geq 0}{}{(I+I'+J^l)}=\bbar{I+I'}$.
Thus $\bbar{\bbar{I}+\bbar{I'}} \subset \bbar{I+I'} \subset \bbar{\bbar{I}+\bbar{I'}}$, showing that $\bbar{\bbar{I}+\bbar{I'}}=\bbar{I+I'}$.
Assume  that $I$ and $I'$ are ideals in $A$. For  all $l\in \N$, note that  $\bbar{I}\,\bbar{I'}\subset (I+J^l)(I'+J^l) \subset I I' +J^l$, implying that  $\bbar{I}\,\bbar{I'} \subset \bbar{I I'}$.  Now, suppose that $I'\subset \bbar{I+JI'+I'J}$, we must  show that $I'\subset \bbar{I}$. Applying the relations just proved, we have\\[0.4em]
$\myeqar{
I'\subset & \bbar{I+JI'+I'J} \subset \bbar{I+J(\bbar{I+JI'+I'J})+(\bbar{I+JI'+I'J})J}
\subset  \bbar{I+\bbar{J(I+JI'+I'J)}+\bbar{(I+JI'+I'J)J}} \\[0.2em]
\subset & \bbar{I+J(I+JI'+I'J)+(I+JI'+I'J)J} \subset  \bbar{I+J^2I'+JI'J+I'J^2},\\
}$
thus repeating the previous procedure, for each $ l\geq 2$  we get
 $I'\subset \bbar{I+\Som[l]{s=0}(J^{l-s}I'J^{s})} \subset \bbar{I+J^l}=I+J^l$, implying that $I'\subset \bbar{I}$ as claimed. This completes the proof of the lemma.
\end{prv}


Let  $L:=\bbar{\kQc\mul\triv{B}\mul\kQc}$  be the closed ideal  generated by $\triv{B}$. We  
get the following  facts on Jacobian ideals. 
\begin{Thm}\label{theo.Jm-et-morph-Jtriv} Let $\morph{\phi:\kQc}{\kQpc}$ be  any path 
$\kk$-algebra  epimorphism  with \mbox{$\mQ'\!=\!(B',\K,\tr)$}. Then the  following statements hold.
\begin{itemize}
    \item [$\msf{(1)}$] $\J{\phi(\mfr{m})}\subset
    \phi(\J{\mfr{m}})$. If  $\phi$ is an isomorphism,
     then $\phi(\J{\mfr{m}})=\J{\phi(\mfr{m})}$,  $\dual{\phi_{1}}(\Ker(\cderv \phi(\mfr{m})))=
\Ker(\cderv \mfr{m})$  and there is an induced isomorphism $\morph[/isomark={\sim}]{\Jc(\mQ,\mfr{m})}{\Jc(\mQ',\phi(\mfr{m}))}$.
 \item [$\msf{(2)}$] If $\phi$ is a reduction on $(\mQ,\mfr{m})$ then
 $\phi(\J{\mfr{m}})=\J{\phi(\mfr{m})}$ and  $ \Jc(\mQ,\mfr{m}) \cong  \Jc(\mQ',\phi(\mfr{m}))$. Moreover there is a left (respectively, right) $\K$-linear isomorphism $\morph{[\id, \rho'] : L\oplus\kQpc}{\kQc}$ with $\rho'(B')\subset B$ and $\rho'(uv)=\rho'(u)\rho'(v)$ for all $u,v\in \J{\kQpc}$ such that $\phi\rho'=\id_{\kQpc}$.
\end{itemize}
\end{Thm}

\begin{Defn} Let $\morph{\phi:\kQc}{\kQpc}$  be a  reduction on $(\mQ,\mfr{m})$. The data  $(\mQ',\phi(\mfr{m}))$ is a  \emph{reduced modulated quiver with
potential  associated with} $(\mQ,\mfr{m})$ and $\phi$ is referred to as reduction from $(\mQ,\mfr{m})$ to $(\mQ',\phi(\mfr{m}))$. We also refer to   $\Ker(\phi)$ as \emph{trivial part} in $\J{\mfr{m}}$.
\end{Defn}

\begin{prv}[Proof of Theorem~\ref{theo.Jm-et-morph-Jtriv}] We let $\what{\mrm{J}}=\J{\kQc}$ and $\what{\bbar{\mrm{J}}}=\J{\kQbc}$. We have  dualizing pairs   $\set{B,\dual{B};\bilf}$, $\set{\bbar{B},\dual{\bbar{B}};\bilfb}$ and    $\set{B',\dual{B'};\bilf'}$.
Let us prove part $\msf{(1)}$. The chain-rule~\eqref{eq-chainrule} shows that, for all  $\xi'\in \dual{B'}$ we have   $\cderv[\xi'](\phi(\mfr{m})) =
\Som{y\in\rmQ}{(\lgrad[\xi'] \phi(y))\bilbox
\phi(\cderv[\dual{y}]\mfr{m})}$, implying by the surjectivity of $\phi$ that
$\J{\phi(\mfr{m})}\subset \phi(\J{\mfr{m}})$. Next, assume  that $\phi$ is an isomorphism. Applying
the previous observations to $\phi^{-1}$ and   $\phi(\mfr{m})$
shows that $\J{\mfr{m}}\subset \phi^{-1}(\J{\phi(\mfr{m})})$ and
then $\phi(\J{\mfr{m}})\subset
\phi(\phi^{-1}(\J{\phi(\mfr{m})}))=\J{\phi(\mfr{m})}$, so
that  $\phi(\J{\mfr{m}})=\J{\phi(\mfr{m})}$. The the degree-$1$ component  $\morph{\phi_{1} : B }{B'}$ of $\phi$ is  a bimodule isomorphism with dual $\morph{ \dual{\phi_{1}}: \dual{B'}
}{\dual{B}}$. Write  $\mfr{m}$ as sum $\Som{l\geq
2}\mfr{m}_l$ of homogeneous potentials  $\mfr{m}_l \in B^{(l)}$. We have $\phi(\what{\mrm{J}}^l)
= \what{\mrm{\bbar{J}}}^l$, $\phi(\mfr{m})=\Som{l\geq 2}\phi(\mfr{m}_l)$
and each $\phi(\mfr{m}_l)$ belongs  to $ \set{0}\cup
( \what{\mrm{\bbar{J}}}^l\setminus  \what{\mrm{\bbar{J}}}^{l+1})$. Let $\xi'\in\dual{B'}$. 
Note that:   $\cderv[\xi'](\phi(\mfr{m}))=0$  if and only if 
  $\cderv[\xi'](\phi(\mfr{m}_l))=0$  for all $l\geq 2$.

For each $l \geq 2$,  by the chain-rule we have
$\cderv[\xi'](\phi(\mfr{m}_l)) = \Som{y\in\rmQ}{(\lgrad[\xi']
\phi(y))\bilbox \phi(\cderv[\dual{y}]\mfr{m}_l)}$ with each
$\phi(\cderv[\dual{y}]\mfr{m}_l)$ lying in $\set{0}\cup
(\what{\mrm{\bbar{J}}}^{l-1}\setminus  \what{\mrm{\bbar{J}}}^{l}))$. But for every $x\in B$ we have 
$\phi(x)=\Som{d\geq 1}\phi_d(x)$ with  $\phi_d(x)\in
B^{(d)}$,  thus  $\lgrad[\xi'] \phi_{1}(x)=(\lderv[\xi']
\phi_{1}(x) \tenss[\kk] 1)=(\bilf'(\xi'\otimes
\phi_{1}(x))\tenss[\kk] 1)= (\bilf(\dual{\phi_{1}}(\xi')\otimes
x)\tenss[\kk] 1)$. We get that,  if $\cderv[\xi'](\phi(\mfr{m}_l))$ is zero then
the term $\zeta_{l,0}':=\Som{y\in\rmQ}{(\lgrad[\xi']
\phi_{1}(y))\bilbox \phi(\cderv[\dual{y}]\mfr{m}_l)}$  is also
zero.  But  
 $ 
\zeta_{l,0}'=   \Som{y\in\rmQ}{\bilf(\dual{\phi_{1}}(\xi')\otimes
y) \phi(\cderv[\dual{y}]\mfr{m}_l)} = \phi\bigl(
\Som{y\in\rmQ}{\bilf(\dual{\phi_{1}}(\xi')\otimes y)
\cderv[\dual{y}]\mfr{m}_l} \bigr)  
=   \phi\bigl(  \cderv(
\Som{y\in\rmQ}{\bilf(\dual{\phi_{1}}(\xi')\otimes y) \dual{y}
\otimes \mfr{m}_l}) \bigr)  =  \phi\bigl(  \cderv( \dual{\phi_{1}}(\xi') \otimes
\mfr{m}_l ) \bigr),$  
showing  that $\zeta_{l,0}'$ is zero if and only if
$\dual{\phi_{1}}(\xi')$ belongs to $ \Ker(\cderv \mfr{m}_l)$.  We conclude
that, if $\xi'\in \Ker(\cderv \phi(\mfr{m}))$ then
$\dual{\phi_{1}}(\xi') \in \bigcap\limits_{l \geq 2}(\Ker(\cderv
\mfr{m}_l)) = \Ker(\cderv\mfr{m})$. Thus $\dual{\phi_1}(\Ker(\cderv
\phi(\mfr{m}))) \subseteq \Ker(\cderv(\mfr{m}))$. Applying the
previous argument to $\phi^{-1}$ and   $\phi(\mfr{m})$ we get
 ${\dual{\phi_{1}}}^{-1}(\Ker(\cderv
    \mfr{m}))  \subseteq \Ker(\cderv \phi(\mfr{m}))$, implying that
     $\dual{\phi_{1}}(\Ker(\cderv \phi(\mfr{m})))= \Ker(\cderv
\mfr{m})$.

We now prove  part $\msf{(2)}$. Here, $\phi$ is by assumption a  reduction on $(\mQ,\mfr{m})$. Condition $\msf{(1.i)}$  of Definition~\ref{defn.trivpart-reduction} implies that $B\cap (\Ker(\phi)+\J[2]{\kQc})=\triv{B}$ and  $\morph{\phi_1: B}{B'}$   yields a bimodule isomorphism $\morph[/above={\sim}]{\bbar{\phi}_1: \bbar{B}}{B'}$ with $\phi_1=\bbar{\phi}_1\rho$, the projection $\morph{\rho: B}{\bbar{B}}$  induces a morphism 
$\morph{\bbar{\rho}: (B+\Ker(\phi))/\Ker(\phi)}{\bbar{B}}$ with $\bbar{\rho}(x+\Ker(\phi)) =\rho(x)$ for all $x\in B$. Still by assumption, 
$\bbar{\rho}$ has a right inverse 
$\morph{ \bbar{\rho}':  \bbar{B}}{(B+\Ker(\phi))/\Ker(\phi)}$ which lifts to a left (respectively, right) $\K$-linear  map $\morph{ \rho':  \bbar{B}}{B}$  such that $ \phi\circ \rho' =\phi_1\circ \rho'$.  Note that $\rho\rho'=\id_{\bbar{B}}$ and $\phi_1  \rho'=(\bbar{\phi}_1\rho)  \rho'=\bbar{\phi}_1$. Hence, without loss of generality, we may assume that $B'=\bbar{B}$,  $\bbar{\phi}_1=\id_{\bbar{B}}$ and  $\rho'$ is right $\K$-linear. Thus $\phi_1=\rho$, $\phi\circ \rho' =\phi_1\circ \rho' = \id_{\bbar{B}}$, we  have a right $\K$-linear  isomorphism 
$\morph{h'=[\id,\rho'] : \triv{B}\oplus: \bbar{B}}{B}$ and by duality a left $\K$-linear
 isomorphism $\Seq{\rdual{(h'^{-1})}=[\msf{j}, \id] :\triv{B} \oplus
\dual{\bbar{B}} \, \to[/above={\sim}]\, \dual{B}}$ where
$\msf{j}$ is  a left $\K$-linear right inverse for
$\morph{\msf{p}=\cderv\mfr{m}_{2}:\dual{B}}{\triv{B}}$. For some $p,q\in \N$, we have Casimir elements $\z[\triv{B}\otimes \triv{B}]=\Som[p]{k=1}y_{k}x_{k}$ and  $\z[\bbar{B}\otimes  \dual{\bbar{B}}]=\Som[q]{s=1}\bbar{y}_{p+s}\otimes\dual{y_{p+s}}$. For $k\in \ninterv{p}$ and $s\in\ninterv[p+1]{p+q}$, put  $\dual{y_{k}}=\msf{j}(x_{k}) \in \dual{B}$  and    $y_s=\rho'(\bbar{y}_s)$. Part $\msf{(1)}$ of
Proposition~\ref{prop.dec.pot-mQ} states that 
\\[-1.0\baselineskip]
\begin{equation} \label{eq.theo.Jm-et-morph-Jtriv} \tag{$\star$}
  \z[B\otimes \dual{B}] = (\id\otimes \msf{j}) (\z[\triv{B}\otimes \triv{B}])
  +(\rho'\otimes \id) (\z[\bbar{B}\otimes
 \dual{\bbar{B}}]) =\Som[p+q]{k=1}y_k\dual{y_k}.
 \end{equation}
By Remark~\ref{rem.bilf-trivpart-redpart} we have: 
$\bilf(\xi\otimes x)=\bilfb(\xi\otimes \rho(x))=0=\bilfb(\rho(x)\otimes \xi)=\bilf(x\otimes \xi)$  for  all $\xi\in\dual{\bbar{B}},x\in\triv{B}$.
Let  $\xi \in\dual{\bbar{B}}$. Using the chain-rule,
 identities~\eqref{eq.dervdualbasesKQc} as well as  previous observations, we have:\\
 $\begin{array}{rl}
\phi(\cderv[\xi]\mfr{m}) = & \phi\bigl( \lderv[\xi] (\cperm
\mfr{m}) \bigr) = \phi\bigl( \lderv[\xi] \bigl(\Som[p+q]{k=1}y_k\otimes \lderv[\dual{y_k}] (\cperm \mfr{m}) \bigr) \bigr) =    
\Som[p]{k=1}{\bilf(\xi\otimes y_k)
\phi(\cderv[\dual{y_k}]\mfr{m})} +
\Som[p+q]{k=p+1}{\bilf(\xi\otimes y_k)
\phi(\cderv[\dual{y_k}]\mfr{m})} \\
= & 0 + \Som[p+q]{k=p+1}{\bilfb(\xi\otimes \rho(y_k))
\phi(\cderv[\dual{y_k}]\mfr{m})} =
\Som[p+q]{k=p+1}{\bilf(\xi\otimes \phi_1(y_k))
\phi(\cderv[\dual{y_k}]\mfr{m})} =  \Som[p+q]{k=p+1}{\lderv[\xi] (\phi_1( y_k))
\phi(\cderv[\dual{y_k}]\mfr{m})} \\
= & \Som[p+q]{k=1}{(\lgrad[\xi] \phi_1(y_k))\bilbox
\phi(\cderv[\dual{y_k}]\mfr{m})} = \Som[p+q]{k=1}{(\lgrad[\xi] \phi(y_k))\bilbox
\phi(\cderv[\dual{y_k}]\mfr{m})}   \\
= &    (\Som[p+q]{k=1}{(\lgrad[\xi]
\phi(y_k))\bilbox \phi(\cderv[\dual{y_k}]\mfr{m})} - \Som[p]{k=1}{(\lgrad[\xi]
\phi(y_k))\bilbox \phi(\cderv[\dual{y_k}]\mfr{m})}  \\
 = & \cderv[\xi]\phi(\mfr{m}) - \Som[p]{k=1}{(\lgrad[\xi]
\phi(y_k))\bilbox \phi(\cderv[\dual{y_k}]\mfr{m})}  \in \J{\phi(\mfr{m})} + (\what{\bbar{\mrm{J}}}\phi(\J{\mfr{m}}) +\phi(\J{\mfr{m}}) \what{\bbar{\mrm{J}}}).
 \end{array}$ \\
Hence,
$ \phi((\cderv \mfr{m})(\dual{\bbar{B}}))\subset   \J{\phi(\mfr{m})} + (\what{\bbar{\mrm{J}}}\phi(\J{\mfr{m}}) +\phi(\J{\mfr{m}}) \what{\bbar{\mrm{J}}})$.  But   \mbox{Lemma~\ref{lem.trivIdeal-closedIdeal}-$\msf{(b)}$} shows  that  $\J{\mfr{m}}$ coincides with the  closure of  $\Ker(\phi)+\ggen{ (\cderv \mfr{m})(\dual{\bbar{B}}) }$, implying that  $ \phi(\J{\mfr{m}})$ is contained in the closure of $\J{\phi(\mfr{m})} + (\what{\bbar{\mrm{J}}}\phi(\J{\mfr{m}}) +\phi(\J{\mfr{m}}) \what{\bbar{\mrm{J}}})$ and, applying part $\msf{(b)}$ of Lemma~\ref{lem.trivIdeal-closedIdeal}  we get $\phi(\J{\mfr{m}})\subset \J{\phi(\mfr{m})}$. By  the chain-rule  we also have  $\J{\phi(\mfr{m})}\subset \phi(\J{\mfr{m}})$. Thus  $\phi(\J{\mfr{m}})=\J{\phi(\mfr{m})}$ and the latter also  shows that   $\phi$ induces an isomorphism of Jacobian algebras from $\Jc(\mQ,\mfr{m})$ to $\Jc(\mQ',\phi(\mfr{m}))$.  To complete the proof of $\msf{(2)}$, we will  extend  the right $\K$-linear map $\morph{ \rho':  \bbar{B}}{B}$ to a
continuous right $\K$-linear morphism   again denoted by
  $\rho': \kQbc \to \kQc$, with $\rho'\restr{\K}=\id_{\K}$. 
Recall that  $L$ is  the closed ideal in $\kQc$ generated by $\triv{B}$.  With notations of \eqref{eq.theo.Jm-et-morph-Jtriv} above, let $\rmQ[1]=\set{y_k \ : \ k\in\ninterv{p+p}}$ and $\rdmQ[1]=\set{\dual{y_k}  \ : \ k\in\ninterv{p+p}}$. Then $(\rmQ,\rdmQ)$ is a  right projective basis for the pair $\set{B,\dual{B}; \bilf}$, while
 $(\set{y_{k}\, : k\in\ninterv{p} },\set{x_{k}\, : k\in\ninterv{p} })$ and 
 $(\set{\bbar{y}_{p+k}\,:   k\in \ninterv{q}},\set{\dual{y_{p+k}}\, : 
   k\in \ninterv{q} })$ are  right projective bases for the pairs $\set{\triv{B},\triv{B}}$ and  $\set{\bbar{B},\dual{\bbar{B}};\bilfb}$ respectively.  And by definition we have: 
 $ \rho'(\bbar{y}_k) =y_k$  for all $k\in\ninterv[p+1]{p+q}$. In view of subsection~\ref{sub.syst.fleches}, we form corresponding     right projective bases $(\rdmQc[0],\rmQc[0])$ and $(\rdmQbc[0],\rmQbc[0])$ for the weakly
dualizing pairs $\set{\kQd,\kQc;\bilfh}$  and $\set{\kk\dual{\bbar{\mQ}},\what{\kk\dual{\bbar{\mQ}}};\what{\bilfb}}$ respectively. Here $\rmQc[0]=\set{1}\cup \rmQc[1]$ with $\dual{1}=1\in\K$; each   $y\in \rmQc[1]$  expresses as $y=y_{i_1}\otimes \dotsm \otimes y_{i_l}$ with $l\geq 1$ and $i_1,\dotsc,i_l \in \ninterv{q+p}$, the corresponding dual is  $\dual{y}=\dual{y_{i_l}}\otimes \dotsm \otimes \dual{y_{i_1}}$. A similar description is given for $\rmQbc[0]$. Next, put $\mrm{Y}=\rmQc[1]\cap L$, it consists of basis elements $y_{i_1}\otimes \dotsm \otimes y_{i_l}$ such that at least one of the integers $i_1,\dotsc,i_l$ belongs to  $\ninterv{p}$. Also put $\mrm{Y}'=\rmQc[1]\sminus\mrm{Y}$. Therefore, $\rho'$ is defined on each   basis element $\bbar{y}=\bbar{y}_{i_1} \otimes \dotsm \otimes \bbar{y}_{i_l}\in \rmQbc$ by: $\rho'(\bbar{y}):=y_{i_1}\otimes \dotsm \otimes y_{i_l}$. Thus, in virtue of identities~\eqref{eq.dualbasesKQc} from subsection~\ref{sub.syst.fleches}, for each $x\in \kQc$ and $\bbar{x}=\Som{\bbar{y}\in\rmQbc[0]}\bbar{y} \bilfh(\dual{ y}\otimes x)$ we have 
 $\rho'(\bbar{x})=
 \Som{\bbar{y}\in\rmQbc[0]}\rho'(\bbar{y})\what{\bilf}(\dual{y}\otimes x)$.
By construction,   $\rho': \kQbc \to \kQc$ has the desired  properties.

\end{prv}

Different presentations of Jacobian algebras by reduced modulated quivers with potentials can be compared using the following concept.

\begin{Defn}\label{defn.equivdroite} Let $(\mQ',\mfr{m}')$
be another  modulated quiver with potential  with
\mbox{$\mQ'=(B',\K,\tr)$}.
 A \emph{weak right-equivalence} between $(\mQ,\mfr{m})$ and
$(\mQ',\mfr{m}')$ is  a path algebra  isomorphism
$\morph[/isomark={\approx}]{\phi: \kQc}{\kQpc}$  
 such that $\J{\phi(\mfr{m})}=\J{\mfr{m}'}$.  If moreover $\phi(\mfr{m})$
is cyclically equivalent to $\mfr{m}'$ then $\phi$ is  a
\emph{right-equivalence}.
 \end{Defn}

Under the assumption that the trivial part of $(\mQ,\mfr{m})$, the first  main result of this work gives the existence and uniqueness of split reductions up to weak 
right-equivalences. As before, $\morph{\rho: \kQc}{\kQbc}$ is the natural projection.
\begin{Thm}[{reduction theorem}]\label{theo.red-mQp} Assume the trivial part of
$(\mQ,\mfr{m})$  splits and write $\mQ=\triv{\mQ}\oplus \red{\mQ}$.  Then there is a right-equivalence
$\phi$  from $(\mQ,\mfr{m})$  to a  direct sum
$\triv{(\mQ,\mfr{m})}\oplus (\red{\mQ},\red{\mfr{m}})$, yielding a split reduction
$\pi_{\mfr{m}}=\rho\phi$ from $(\mQ,\mfr{m})$   into a reduced modulated quiver with potential
$\Red(\mQ,\mfr{m})=(\red{\mQ},\red{\mfr{m}})$, with kernel $\Jtriv{\mfr{m}}$ such that  $\mfr{m}-\red{\mfr{m}}$ is cyclically
equivalent to an element in $ (\Jtriv{\mfr{m}})^{2}$.
 Furthermore, the split reduction process $\morph[/dir=|->]{\Red:
(\mQ,\mfr{m})}{\Red(\mQ,\mfr{m})}$ is a well-defined operation on
weak-right equivalence classes of modulated quivers with
potentials.
\end{Thm}


The proof of the first part of  Theorem~\ref{theo.red-mQp} is the object of  the discussion  from Lemma~\ref{lem.cyclequiv} to Lemma~\ref{lem.splitting-mQ}.

\begin{Lem} \label{lem.cyclequiv} Let $N$ be a direct summand in 
$B$ and $S$ a potential lying in the closed ideal $\bbar{\ggen{N}}$ generated by $N$. Then $S$ is cyclically equivalent to a potential lying in $N \mul \J{\kQc}$ and to a potential lying in $\J{\kQc}\mul N$.
\end{Lem}

\begin{prv} Elements of $\bbar{\ggen{N}}$ can be written as   possibly infinite sums of elements $u_l,\, l\geq 1$,  with $u_l \in \ggen{N}\cap B^{(l)}=\Som[l-1]{s=0}{B^{(s)} N B^{(l-s-1)}}$. Thus   $S=\Som{l\geq 1}S_l$ where $S_l$ is a potential lying in $\ggen{N} \cap B^{(l+1)}$. By assumption,  $B= N\oplus N'$  for some subbimodule $N'\subset B$. There is a corresponding  decomposition $\dual{B}= \dual{N} \oplus \dual{N'}$  such that for all $(\xi,\xi') \in \dual{N}\times
 \dual{N'}$ and $(x,x')\in N\times N'$ we have:
 $\bilf(\xi\otimes x')=0=\bilf(x'\otimes \xi)$ and $\bilf(\xi'\otimes x)=0=\bilf(x\otimes
 \xi')$. For each  $l\geq 1$ we have
 $B^{(l+1)} = (N\oplus N')B^{(l)} =
\bigl(\somd[l]{s=0}{N'^{(s)}\mul N \mul B^{(l-s)}}\bigr) \, \oplus
N'^{(l+1)}$,  and  each   $S_l\in \ggen{N} \cap
B^{(l+1)}$  expresses as sum $\Som[l]{s=0}S_{l,s}$   of
potentials with $S_{l,s}\in N'^{(s)} \mul N \mul B^{(l-s)}$. Hence the left permutation $\lperm[s]S_{l,s}$ of order $s$ of each $S_{l,s}$
belongs to $N B^{(l)}$ while the right permutation $\rperm[l-s]S_{l,s}$ of order $l-s$ of  each $S_{l,s}$ belongs to $B^{(l)}N$. Thus $S$ is cyclically equivalent to a
potential lying in $N \mul \J{\kQc}$ and to a potential lying in $\J{\kQc}\mul N$.
\end{prv}

Denoting as before the reduced part of
$B$  by  $\bbar{B}$  and using the assumption in    Theorem~\ref{theo.red-mQp}, we simply write $B=\triv{B}\oplus \bbar{B}$.  Part $\msf{(2)}$ of
Proposition~\ref{prop.dec.pot-mQ} shows that  the pair
$\set{B,\dual{B};\bilf}$ occurs as  direct sum of naturally induced
dualizing pairs $\set{\triv{B},\dual{\triv{B}}} \oplus \set{\bbar{B},\dual{\bbar{B}}}$ with
$\dual{B}=\dual{\triv{B}}\oplus \dual{\bbar{B}}$ and  we have a
right  inverse $\morph[/dist=2.5,/above={{\jtriv}}]{\triv{B}}{
\dual{B}}$   for the cyclic derivative
$\morph[/dir=->>]{\cderv\mfr{m}_{2} : \dual{B}}{\triv{B}}$  such
that $\im(\jtriv)=\dual{\triv{B}}$. In view of   Remark~\ref{rem.bilf-trivpart-redpart},
letting $U=\rderv \mfr{m}_2(\dual{B})$ and $V=\lderv \mfr{m}_2(\dual{B})$, we have
$U\cap V=0$, $\triv{B}=U\oplus V$, 
$\set{\triv{B},\dual{\triv{B}}} = \set{U,\dual{U}} \oplus
\set{V,\dual{V}}$ as  direct sum of naturally induced
 dualizing pairs   with
$\dual{U}=\jtriv(V)$ and $\dual{V}=\jtriv(U)$  and   $\morph[/dir=->>]{\cderv\mfr{m}_{2} : \dual{\triv{B}}}{\triv{B}}$ occurs  as direct sum of the partial
derivatives $\morph[/above={\sim}]{\rderv\mfr{m}_{2} : \dual{V}}{U}$
and $\morph[/above={\sim}]{\lderv\mfr{m}_{2} : \dual{U}}{V}$. Let us 
summarise  previous observations: 
\begin{align}
\label{cond-loopfree}
 &  \triv{B}=U\oplus V, \; B=\triv{B}\oplus
 \bbar{B}, \text{ and }
 \set{B,\dual{B}}= \set{U,\dual{U}} \oplus \set{V,\dual{V}} \oplus
\set{\bbar{B},\dual{\bbar{B}}} \text{ as dualizing pairs.}   \\
\label{eq.loopfreedecomp-mQ} &
\text{For bimodule morphism } \morph[/dist=2,/above={{\jtriv}}]{\triv{B}}{ \dual{B}} 
 \text{ we have: } \cderv\mfr{m}_{2}\circ\jtriv=\id_{\triv{B}},   \dual{U}\hspace{-0.4em}=\jtriv(V),\,
\dual{V}\hspace{-0.4em}=\jtriv(U).\hspace{-0.2em}
 \\[-0.5\baselineskip]
& \mfr{m}_{2} = \z[U\otimes V] =
\som{k=1}{p}{y_k\otimes x_k},\; \text{ where } \z[U\otimes
\dual{U}]=\som{k=1}{p}{y_k\otimes \dual{y_k}} \text{ and }
\z[\dual{V}\otimes V]=\som{k=1}{p}{\dual{x_k}\otimes x_k}.
\end{align}
Note in view of Lemma~\ref{lem.cyclequiv} that $\mfr{m}$ is
cyclically equivalent to a potential of the form
\mbox{$\z[U\otimes V] + S_1+S_2 + \bbar{\mfr{m}}_1$} with
\mbox{$S_1\in U \otimes \J[2]{\kQc}$}, \mbox{$S_2\in \J[2]{\what{\mrm{T}_{\K}(V\oplus \bbar{B})}}\otimes V\subset
\J[2]{\kQc}\otimes V$}, \mbox{$\bbar{\mfr{m}}_1 \in
\J[3]{\kQbc}$}, we can therefore write:
\\[-0.8\baselineskip]
\begin{equation}\label{eq.decomp-m}
\begin{aligned} & \mfr{m} \equiv_{\mrm{cyc}} \som{k=1}{p}{y_k\otimes x_k} +
\som{k=1}{p}{y_k\otimes v_k} + \som{k=1}{p}{u_k\otimes x_k} +
\bbar{\mfr{m}}_1, \\ &  \text {with } v_k=
\lderv(\dual{y_k}\otimes S_1),\, u_k=\rderv(S_2\otimes \dual{x_k})
\text{ for all } k \in \ninterv{p}.
\end{aligned}
\end{equation}
Letting  $\Seq{f:=(\cderv\mfr{m})\circ  \jtriv: \triv{B}\,\to \, \triv{B}\oplus \J[2]{\kQc}}$, we consider   the closed ideal
 $\Jtriv{\mfr{m}}:=\bbar{\ggen{f(\triv{B})}}$.

For a natural number $d\geq 1$,  a potential in the form
\eqref{eq.decomp-m} is   \emph{$d$-split} if
$u_k,v_k\in\J[d+1]{\kQc}$.

\begin{Lem} \label{lem.splitting-mQ} With the  assumption that the trivial part of $(\mQ,\mfr{m})$
splits, there exists a unitriangular automorphism
$\morph{\phi: \kQc}{\kQc}$ such that $\phi(\mfr{m})$ is cyclically
equivalent to a potential $\wtilde{\mfr{m}}$ in the form
\eqref{eq.decomp-m} with $u_k=0=v_k$ for all $k\in\ninterv{p}$ and such that: $\phi\restr{\kQbc}=\id_{\kQbc}$. Moreover, the map
$\pi_{\mfr{m}} =\rho \phi$ is a reduction on $(\mQ,\mfr{m})$  with $\Ker(\pi_{\mfr{m}})=\Jtriv{\mfr{m}}$  and
$\mfr{m}-\pi_{\mfr{m}}(\mfr{m}) \cyc{\equiv}
\phi^{-1}(\mfr{m}_{2}) \in (\Jtriv{\mfr{m}})^{2}$.
\end{Lem}

\begin{prv}
\Par{Claim} Suppose $S$ is a $d$-split potential
written in the form \eqref{eq.decomp-m}. Then there exists a
unitriangular automorphism $\morph{\vphi: \kQc}{\kQc}$ having
depth $d$, with $\vphi\restr{\kQbc}=\id_{\kQbc}$, such that
$\vphi(S)$ is cyclically equivalent to a $2d$-split potential
$S'$, with $S'-S\in \J[2d+2]{\kQc}$.

We write $S=\z[U\otimes V] + S_1+S_2 + \bbar{S}_1$ with $S_1\in U
\otimes \J[d+1]{\kQc}$, \mbox{$S_2\in \J[d+1]{\kQc}\otimes V$},
\mbox{$\bbar{S}_1 \in \J[3]{\kQbc}$}, keeping the
notations of \eqref{eq.decomp-m} for $S$. Then we  have a
unitriangular automorphism $\morph{\vphi: \kQc}{\kQc}$ having
depth $d$, defined by letting:   $\vphi\restr{\kQbc}=\id_{\kQbc}$,  $\morph{\vphi\restr{U} =\id_{U}-(\rderv S_2)\circ
\jtriv : U}{U \oplus \J[d+1]{\kQc}}$ and $\morph{\vphi\restr{V}= \id_{V} - (\lderv S_1)\circ \jtriv: V}{V \oplus \J[d+1]{\kQc} }$.
Thus, for all $k\in\ninterv{p}$ we have:
 $\vphi(y_k)=y_k - u_k$ and
$\vphi(x_k)=x_k - v_k$ with $v_k \in \J[d+1]{\kQc}$  and $u_k \in
\J[d+1]{\kQbc}$.   Since $\vphi$ has depth $d$, we have 
\mbox{$\vphi(u_k)=u_k+u_k'$} and \mbox{$\vphi(v_k)=v_k + v_k'$} for
some \mbox{$u_k',v_k'\in\J[2d+1]{\kQc}$}.
We get $\vphi(S)=\som{k=1}{p}{(y_k-u_k)(x_k-v_k) +
(y_k-u_k)(v_k+v_k') + (u_k+u_k')(x_k-v_k)} +  \bbar{S}_1
=\som{k=1}{p}{y_k\otimes x_k} + W + \bbar{S}_1$, where
\mbox{$W=\som{k=1}{p}{(y_k\otimes v_k'+u_k'\otimes x_k -u_k\otimes
v_k'-u_k'\otimes v_k -u_k\otimes v_k )}\in \J[2d+2]{\kQc}$} is a potential. Since
\mbox{$\kQc=L\oplus \kQbc$} with
\mbox{$L=\bbar{\kQc\mul\triv{B}\mul \kQc}$},  we can write
\mbox{$W=W' + \bbar{W}$}  for two potentials \mbox{$W'\in L\cap
\J[2d+2]{\kQc}$} and \mbox{$\bbar{W}\in\J[2d+2]{\kQbc}$}. But using again
 Lemma~\ref{lem.cyclequiv} and the fact that
$\triv{B}= U\oplus V$, we get that  $W'$ is cyclically equivalent to a sum
$W''=W_1''+W_2''$ of two potentials $W_1''=\som{k=1}{p}{y_k\otimes
v_k''} \in U \mul \J[2d+2]{\kQc}$ and
$W_2''=\som{k=1}{p}{u_k''\otimes x_k} \in  \J[2d+2]{\kQc}\mul V$
with \mbox{$u_k'' =\rderv[\dual{x_k}]W_2''\in \J[2d+1]{\kQbc}$}
and \mbox{$v_k''=\lderv[\dual{y_k}] W_1''\in \J[2d+1]{\kQc}$}. Hence,  \mbox{$W-(W''+ \bbar{W})$} lies in $
\J[2d+2]{\kQc}$ and  in the closed module
$\skcb{\Zck(\kQc_{(2)}),\Zck(\kQc_{(2)})}$ of skew commutators in
$\kQc$. Thus the desired potential $S'$ is  given by:
$S'=\som{k=1}{p}{y_k\otimes x_k} + S_1'+ S_2' + (\bbar{W} +
\bbar{S}_1)$  where $ S_1'=\som{k=1}{p}{y_k\otimes v_k''}$ and  $
S_2'= \som{k=1}{p}{u_k''\otimes x_k} $, with $u_k''
=\rderv[\dual{x_k}]S_2' \in \J[2d+1]{\kQbc}$ and
$v_k''=\lderv[\dual{y_k}] S_1'\in \J[2d+1]{\kQc}$.
This  completes the proof of our  claim.

Next,   starting  with a $1$-split potential $S_1$ in the form
\eqref{eq.decomp-m}  and using successively the above claim, one constructs a sequence of potentials $S_1,S_2,\dotsc,$
and a sequence of unitriangular automorphisms
$\phi_1,\phi_2,\dotsc,$ having the following properties:
\begin{enumerate}
    \item [$\msf{(p0)}$:] $\mfr{m} \equiv_{\mrm{cyc}} S_1$. \quad
    $\msf{(p1)}$:\, $S_d$ is $2^{d-1}$-split.  \quad
     $\msf{(p2)}$:\,  $\phi_d$ has depth $2^{d-1}$.
    \item [$\msf{(p3)}$:] $\phi_d(S_d) \equiv_{\mrm{cyc}} S_{d+1} $
    and  each element $C_d:=\phi_d(S_d)- S_{d+1}$ lies in $\J[2^{d}+2]{\kQc}\cap
    \skcb{\Zck(\kQc_{(2)}),\Zck(\kQc_{(2)})}$.
\end{enumerate}
Using $\msf{(p2)}$ we set  $\phi= \lim\limits_{l\rightarrow
\infty}\phi_l\phi_{l-1}\dotsm\phi_1$. By part $\msf{(a)}$
of Remark~\ref{rem.formelt-(A,JA)} and   
Proposition~\ref{prop.prolong}, $\phi$  is a well-defined
unitriangular  automorphism of $\kQc$ such that
$\phi\restr{\kQbc}=\id_{\kQbc}$. And by $\msf{(p3)}$,  each
element $C_d:=\phi_d(S_d)- S_{d+1}$ lies in $\J[2^{d}+2]{\kQc}\cap
\skcb{\Zck(\kQc_{(2)}),\Zck(\kQc_{(2)})}$. But by
Proposition~\ref{prop.morph-et-perm}, any unitriangular
automorphism sends  $\skcb{\Zck(\kQc_{(2)}),\Zck(\kQc_{(2)})}$ to  itself. Thus, 
using again $(\msf{p0})$, we get that
$\phi_l \phi_{l-1} \dotsm\phi_1(\mfr{m}) \equiv_{\mrm{cyc}} \phi_l
\phi_{l-1} \dotsm\phi_1(S_1) =S_{l+1} +\Som[l]{d=1}\phi_l
\phi_{l-1} \dotsm\phi_{d+1}(C_d)$ for all $l \geq 1$;
and passing to the limit as $l$ tends to $\infty$, we have
 $\phi(\mfr{m}) \equiv_{\mrm{cyc}} \phi(S_1) =
\lim\limits_{l\rightarrow \infty} S_l  + \phi\bigl(\Som{d\geq
1}(\phi_d \phi_{d-1} \dotsm\phi_{1})^{-1}(C_d)\bigr)$.
 \;  Letting $\wtilde{\mfr{m}}=\lim\limits_{l\rightarrow \infty}
S_l$, we get that $\phi(\mfr{m})$ is cyclically equivalent to
$\wtilde{\mfr{m}}$ and, in view of  $\msf{(p1)}$, $\wtilde{\mfr{m}}$ is in the form~\eqref{eq.decomp-m} with $u_k=0=v_k$. To complete the proof, we now let  $\pi_{\mfr{m}}:= \rho \phi$. Then   $\pi_{\mfr{m}}$ is clearly a split reduction with $\Ker(\pi_{\mfr{m}})= \vphi^{-1}(L)$. Let us   check that  $\vphi^{-1}(L)$ coincides with the  closed ideal  $\Jtriv{\mfr{m}}=\bbar{\ggen{f(\triv{B})}}$, where we   recall that
   $\Seq{f:=(\cderv \mfr{m})\circ  \jtriv: \triv{B}\,\to \, \triv{B}\oplus \J[2]{\kQc}}$ with $\cderv\mfr{m}_{2} \circ  \jtriv=\id_{\triv{B}}$. We have  $\wtilde{\mfr{m}}=\z[U\otimes V]
+\bbar{W}=\Som[p]{k=1}y_k\otimes x_k + \bbar{W}$ for some
potential $\bbar{W} \in \J[3]{\kQbc}$.
In view of \eqref{cond-loopfree} and \eqref{eq.loopfreedecomp-mQ} above, $B=\triv{B}\oplus \bbar{B}$ with $\triv{B}=U\oplus V$, $\dual{B}=\dual{\triv{B}} \oplus \dual{\bbar{B}}$ with $\dual{\triv{B}}=\dual{U}\oplus \dual{V}$, and  (using  again
Proposition~\ref{prop.dec.pot-mQ}-$\msf{(1)}$) the Casimir element in $B\otimes \dual{B}$ expresses as $\z[B\otimes \dual{B}]=\z[\triv{B}\otimes \dual{\triv{B}}]+\z[\bbar{B}\otimes \dual{\bbar{B}}]$. Write 
$\z[\triv{B}\otimes \dual{\triv{B}}]=\Som[r]{k=1}{z_k\otimes\dual{z_k}}$.   Since
$\phi \restr{\bbar{B}} =\id_{\bbar{B}}$, we deduce that
$\lgrad[\xi](\phi(\bbar{z}))=\lgrad[\xi](\bbar{z}) =0$ for all
\mbox{$\bbar{z} \in \bbar{B}$} and \mbox{$\xi\in \dual{\triv{B}}$}. Hence, the chain-rule~\eqref{eq-chainrule}
and   the fact that $\phi(\mfr{m}) \equiv_{\mrm{cyc}}
\wtilde{\mfr{m}}$ give the following conclusion: for all $\xi\in \dual{\triv{B}}$ we have
$\cderv[\xi]\wtilde{\mfr{m}}=\cderv[\xi]\phi(\mfr{m})=\Som[r]{k=1}\lgrad[\xi](\phi(z_k))\bilbox \phi(\cderv[\dual{z_k}] \mfr{m}) \in \phi(\Jtriv{\mfr{m}})$. Hence
$L\subseteq \phi(\Jtriv{\mfr{m}}) $, and
applying the inverse unitriangular automorphism $\phi^{-1}$ to
$\wtilde{\mfr{m}}$, we also have $\Jtriv{\mfr{m}} \subseteq\phi^{-1}(L)$, thus
$\phi(\Jtriv{\mfr{m}}) \subseteq \phi\phi^{-1}(L) =L$. Hence
$\phi(\Jtriv{\mfr{m}}) = L $, showing that $\Ker(\pi_{\mfr{m}})=\Jtriv{\mfr{m}}$.  Consequently, we get:
 $\phi(\mfr{m} - \pi_{\mfr{m}}(\mfr{m})) = \phi(\mfr{m} ) -\rho\phi(\mfr{m})
\equiv_{\mrm{cyc}} \wtilde{\mfr{m}} - \rho(\wtilde{\mfr{m}} ) =
\z[U\otimes V] \in L^2 $, so that $\mfr{m} -
\pi_{\mfr{m}}(\mfr{m}) \equiv_{\mrm{cyc}} \phi^{-1}(\z[U\otimes V]
) \in \bigl(\Jtriv{\mfr{m}}\bigr)^{(2)}$. This completes the proof
of Lemma~\ref{lem.splitting-mQ}  and  the first part in
Theorem~\ref{theo.red-mQp}.
\end{prv}

The rest of this section is consecrated to the proof of the second part of  Theorem~\ref{theo.red-mQp}.
Keeping the same assumptions on $(\mQ,\mfr{m})$, to each  direct sum decomposition  $B=\triv{B}\oplus \bbar{B}$ corresponds  a split reduction
$\morph[/dir=->>]{\pi_{\mfr{m}} : \kQc}{\kQbc}$ with kernel  denoted by $\Jtriv{\mfr{m}}$, such that $(\phi_{\mfr{m}})\restr{\kQbc}=\id_{\kQbc}$.
The second part of Theorem~\ref{theo.red-mQp} is given by the following lemma.
\begin{Lem}\label{lem.redwelldefined} Let $\morph[/isomark={\approx}]{\phi:
\kQc}{\kQpc}$ be a weak right-equivalence between $(\mQ,\mfr{m}) $
and a modulated quiver with potential $(\mQ',\mfr{m}')$ with
$\mQ'=(B',\K,\tr)$. Then   $\triv{(\mQ',\mfr{m}')}$
splits and,  for any split reduction $\morph[/dir=->>]{\vphi':=\pi_{\mfr{m}'}': \kQpc}{\kQbpc}$ corresponding to a direct sum decomposition $B'=\triv{B}'\oplus \bbar{B}'$  of $\K$-bimodules,  we have a weak right-equivalence
$\morph[/dist=3,/isomark={\approx}]{\psi:=(\vphi'\circ\phi)\restr{\kQbc}:
\kQbc}{\kQbpc}$ between  $\Red(\mQ,\mfr{m})$ and $\Red(\mQ',\mfr{m}')$.
\end{Lem}

\begin{prv} We  write  $\morph[/dir=->>]{\vphi=\pi_{\mfr{m}} :\kQc}{\kQbc}$
for the reduction defined by $\mfr{m}$ with respect to a direct
sum decomposition $B=\triv{B}\oplus \bbar{B}$.  Let us agree with
the following abbreviations: $J=\J{\mfr{m}},\, J':=\J{\mfr{m}'}$,
$\bbar{J}=\vphi(J),\, \bbar{J}'=\vphi'(J')$,
\mbox{$\what{J}=\J{\kQc} $} and \mbox{$ \what{J}'=\J{\kQpc}$}. 
Recall that the    trivial parts $\Jtriv{\mfr{m}}$ and
$\triv{B}$ are such that:
$\triv{B}=(J+ \what{J}^2 )\cap B=(\Jtriv{\mfr{m}} +
\what{J}^2)\cap B$. Similarly,  $\triv{B}'=(J'+  \what{J}'^2 )\cap B'=(\Ker(\vphi') + \what{J}'^2)\cap B'$.
 In view of
Proposition~\ref{prop.prolong}-$\msf{(a)}$, the degree-$1$
homogeneous  component of isomorphism
$\phi$    is a $\K$-bimodule isomorphism  $\morph[/isomark={\sim}]{\phi_1:
B=\triv{B}\oplus \bbar{B}}{B'}$, and
\mbox{$B'=\phi_1(\triv{B})\oplus \phi_1(\bbar{B})$}. But using
Theorem~\ref{theo.Jm-et-morph-Jtriv} and   the assumption on  $\phi$, we get
$\J{\phi(\mfr{m})}=\phi(\J{\mfr{m}})=\J{\mfr{m}'}=J'$.
 As a path algebra isomorphism, $\phi(\what{J}^l)=\what{J}'^l$ for all $l\in\N$. We then have:
 $\phi(\triv{B})=\phi((J + \what{J}^2)\cap B)=(\phi(J) + \phi(
\what{J}^2))\cap\phi(B)=(J' +  \what{J}'^2)\cap{\phi(B)}$, so that
$\phi_1(\triv{B})=(J' + \what{J}'^2)\cap{\phi_1(B)}=(J' +
\what{J}'^2)\cap B'=\triv{B}'$, implying  that   $\triv{B}'$   splits in $B'$.

 For the rest of the proof,  fix a direct sum decomposition $B'=\triv{B}'\oplus \bbar{B}'$ and let $\bbar{\mQ}'$ be the reduced modulated quiver associated with $\bbar{B}'$. Then let
 $\morph[/dir=->>]{\vphi':=\pi_{\mfr{m}'}':\kQpc}{\kQbpc}$ be the corresponding split reduction. Recall that $\vphi\restr{\kQbc}=\id_{\kQbc}$ and  $\vphi'\restr{\kQbpc}=\id_{\kQbpc}$, hence letting $\bbar{\mfr{m}}=\vphi(\mfr{m})$ and  $\bbar{\mfr{m}}'=\vphi'(\mfr{m}')$ for the corresponding reduced potentials, we derive the following conclusions:
 \\[-0.5\baselineskip]
 \begin{equation}\tag{$*$} \label{eq.prv-thred1}
 \begin{aligned}
 & \text{ for all } z\in\kQc \text{ and } z'\in \kQpc \text{ we have }   
 \vphi(z)-z\in \Jtriv{\mfr{m}}=\Ker(\vphi) \text{ and } \vphi'(z')-z' \in \Ker(\vphi'). \\
&  \bbar{J}=\J{\bbar{\mfr{m}}}= \vphi(\J{\mfr{m}}) \subseteq
 \J{\mfr{m}}= J \text{ and } \bbar{J}'=\J{\bbar{\mfr{m}}'}=
 \vphi'(\J{\mfr{m}'}) \subseteq  \J{\mfr{m}'} =J'.
 \end{aligned}
 \end{equation}
The algebra morphism $\morph{\psi:=(\vphi'\circ\phi)\restr{\kQbc}:
\kQbc}{\kQbpc}$ is already a path algebra morphism, that is, 
$\psi\restr{\K}=\id_{\K}$ and $\psi(\bbar{B}) \subset \J{\kQbpc}$.
In the previous paragraph, we proved that  $\phi_1(\triv{B})=\triv{B}'$.  Write the  bimodule isomorphism  $\morph{\phi_1 : B}{B'}$    in  a matrix form:
$\morph[/isomark={\sim}]{\phi_1=\Psmatr{\phi_{1,1}
& 0 \\ \phi_{2,1} & \phi_{2,2}}: \bbar{B}\oplus
\triv{B}}{\bbar{B}'\oplus \triv{B}'}$, with $\morph{\phi_{1,1} : \bbar{B} }{\bbar{B}'}$, $\morph{\phi_{2,1} : \bbar{B} }{\triv{B}'}$ and $\morph{\phi_{2,2} : \triv{B} }{\triv{B}'}$.  Hence,
  $\phi_{1,1}$ is also an isomorphism.   If we put
$\morph{\phi\restr{B}=\Psmatr{\phi_1 \\ \phi_{(2)}} : B}{B'\oplus
\what{J}'^{2}}=\prd{l\geq 1}{B'^{l}}$, then  for each $z\in
\bbar{B}$ we have: $\phi(z)=\phi_{1}(z) +
\phi_{(2)}(z)=\phi_{1,1}(z) + v$ with $v=\phi_{2,1}(z)+
\phi_{(2)}(z)$ belonging to $\triv{B}'\oplus \what{J}'^{2}$,  thus
$\psi(z)=\vphi'\phi(z)=\vphi'(\phi_{1,1}(z)) +
\vphi'(v)=\phi_{1,1}(z) + \vphi'(v)$ with $\vphi'(v)\in
\J[2]{\kQbpc}$. Then  the \mbox{degree-$1$} component of $\psi$ coincides with the isomorphism $ \phi_{1,1} : \bbar{B} \to \bbar{B}'$, implying  by
Proposition~\ref{prop.prolong}-$\msf{(a)}$  that
$\psi$ is a path algebra  isomorphism.
It now remains to check that $\psi(\bbar{J})=\bbar{J}'$. In view
of \eqref{eq.prv-thred1} above, we have:
$\psi(\bbar{J})=\vphi'\phi(\vphi(J))\subseteq\vphi'\phi(J)=\vphi'(J')=\bbar{J}'$,
 so that $\psi(\bbar{J})\subseteq \bbar{J}'$. Reciprocally, let
$z'\in\bbar{J}'$, then  $\psi$   being already an isomorphism we
have  $z'=\psi(z)$  for some $z\in\kQbc$; we  have to
check that  $z$  belongs to $\bbar{J}$. We have
$z'\in\bbar{J}'=\vphi'(J')=\vphi'(\phi(J))$,  so that
$z'=\vphi'\phi(x)$ for some $x\in J$ and
$\vphi'\phi(z)=\psi(z)=z'=\vphi'\phi(x)$. Thus
$\phi(z-x)=\phi(z)-\phi(x)\in\Ker(\vphi')\subset J'=\phi(J)$,
showing that $z-x\in J$. But then,   $x$ being  already  an
element in $J$ we  get $z\in J \cap \kQbc$, so that
$z=\vphi(z)\in \bbar{J}$. We conclude that
$\psi(\bbar{J})=\bbar{J}'$ and 
$ \psi $ is a weak right-equivalence  between  $\Red(\mQ',\mfr{m}')$ and $\Red(\mQ',\mfr{m}')$. 
\end{prv}

 When  the trivial part of $(\mQ,\mfr{m})$ does not split, reductions as described in Definition~\ref{defn.trivpart-reduction} may not exist. However,  examples from  section~\ref{secexample}   illustrate the fact that  reduction or  a notion a skew reduction still survive is some cases; but  weak right-equivalence   is still too restrictive to be a comparison tool between  non-split  reductions.


\section{Symmetric potentials} \label{sec:sympot}
 The main result of this section is that the study of modulated quivers with symmetric potentials mimics the simply-laced study of quivers with potentials; in particular the sophisticate  issue of skew permutations of general potentials is made easy  for symmetric potentials. As before,  $\mQ=(B,\K,\tr)$ is a fixed modulated quiver over $(\K,\tr)$;  the data
   $(\rmQ[1],\rdmQ[1])$ and $(\lmQ[1],\ldmQ[1])$ are respectively  
  right projective  and left projective bases associated with the dualizing pair  $\set{B,\dual{B};\bilf}$.

 Let
$\Som{s\in\Lambda}{e_s \otimes \dual{e_s}} \in\Zck(\K\tenss[\kk]
\K)$ be the  Casimir element of the  symmetric $\kk$-algebra $(\K,\tr)$, then the set $\set{e_s,\,\dual{e_s}\, : \,
s\in\Lambda}$ is a projective $\kk$-basis of $\K$ characterized by
identities~\eqref{casimir-K} which we  recall:
for all  $a\in \K$, $ \som{s\in\Lambda}{}{e_s
\tr(\dual{e_s}a)}= a = \som{s\in\Lambda}{}{\tr(ae_s) \dual{e_s}}$.
The enveloping algebra
$\Ke=\K\tenss[\kk] \op{\K}$ is endowed with the involution sending
each $ (a\otimes b)$ to $\op{(a\otimes b)}:=b\otimes a$. Each
$\K$-bimodule $M$ is naturally viewed as right and left $\Ke$-module; we have: $x\mul(a\otimes b)=b x a =
\op{(a\otimes b)}\mul x$ for all $x\in M$ and $a,b\in\K$. Consider
the $\Zc(\K)$-module $\kQc\tenss[\Ke] \K$ which is identified with
the $\Zc(\K)$-module $\Frac{\kQc}{[\K,\kQc]}$, where $[\K,\kQc]$
is the $\kk$-module generated by commutators $[a,v]:=a v- va$,
with $a\in\K$ and $v\in\kQc$. Indeed, writing $\bbar{v}$ for the
coset $v+[\K,\kQc]$ with $v\in\kQc$,  the map
$\morph{\kQc\tenss[\Ke] \K}{\Frac{\kQc}{[\K,\kQc]}: v\otimes a
\mapsto \bbar{v a}}$ yields a natural isomorphism with inverse
$\morph{\Frac{\kQc}{[\K,\kQc]}}{\kQc\tenss[\Ke] \K: \bbar{v}
\mapsto v \otimes 1}$.  In the sequel, each  $v\tenss[\Ke] 1
\in \kQc\tenss[\Ke] \K$  will be simply denoted by $v\tenss[e] 1$.
\begin{Lem} \label{lem.casimi-sympot}
\begin{itemize}
    \item [$\msf{(a)}$] We have a map 
$\morph{\zc: \kQc}{\Zck(\kQc):\, v \mapsto \Som{s\in\Lambda}{e_s v
\dual{e_s}} = \Som{s\in\Lambda}{\dual{e_s} v e_s } }$,  referred to as  \emph{Casimir
operator} for $\kQc$. It induces  a  $\Zc(\K)$-linear map $\morph[/dir=->>]{
\ztc: \kQc\tenss[\Ke]\K}{\zc(\kQc) :\, v\tenss[e] 1 \mapsto
\ztc(v\tenss[e] 1) := \zc(v)}$.
 \item [$\msf{(b)}$] The
$\Zc(\K)$-module $\kQc\tenss[\Ke]\K$ enjoys an ordinary cyclic
permutation operator
$\morph[/above={\cperm}]{\kQc\tenss[\Ke]\K}{\kQc\tenss[\Ke]\K}$
such that for every  $v=x_0\dotsm x_n\in B^{(n+1)}$ indexed over the cyclic group
$\mb{Z}_{n+1}:=\set{0,\dotsc,n}$, with $x_i\in B$ and with
corresponding ordinary cyclic permutation
$\cperm[\mrm{ord}]v=\Som[n]{i=0}{x_i\dotsm x_{i+n}}$,  we
have: $\cperm(v\tenss[e] 1) = (\cperm[\mrm{ord}] v)\tenss[e]1$.
\end{itemize}
\end{Lem}

\begin{prv} Part $\msf{(a)}$ readily follows  by the the
characterizing identities~\eqref{casimir-K} of the Casimir element
$\Som{s\in\Lambda}{e_s \otimes \dual{e_s}} \in\Zck(\Ke)$. And 
part $\msf{(b)}$ follows by the fact that $au\tenss[e] 1=ua\tenss[e]
1$ for all $u\in\kQc$ and $a\in \K$.
 \end{prv}

The module $\zc(\K)$ is called the \emph{Casimir
 ideal} of $\Zc(\K)$ (\cite[\S3.2]{Lorenz2010},  \cite[\S2]{Higman1955b}) and does not
depends on the choice of a projective $\kk$-basis
$\set{e_s,\dual{e_s}\; : \, s\in \Lambda}$ and, in view of part
$\msf{(a)}$ of  Lemma~\ref{lem.trace-dualizing-pair}, $\zc(\K)$ does
not depend on the  trace $\tr$
 chosen for $\K$.  We shall refer to
the elements of $\zc(\kQc_{(2)})$ as \emph{symmetric potentials},
and in view of point $\msf{(b)}$ of Lemma~\ref{lem.casimi-sympot}
we will also refer to the elements of $\kQc_{(2)}\tenss[\Ke] \K$ as
\emph{symmetric potentials}.

Next we consider the $\K$-bimodule   $\ddual{B}:=\HM{\Ke}{B,\Ke}$
where $B$ is regarded as right $\Ke$-module, thus the natural left
$\Ke$-module structure of $\ddual{B}$ is  such that  for all
$\xi\in \ddual{B}$ and $a,b\in\K$ we have $((a\otimes
b)\xi)(x)=(a\otimes b)\xi(x)$, with $x\in B$. We refer to
$\ddual{B}$ as the \emph{dual of $B$ as  bimodule}. 
For all  $a\in \Zc(\K)$,  $x\in B$ and $\alpha\in\ddual{B}$, we note that  $\alpha(ax)=\alpha(xa)$.  Hence, we naturally define a partial (left)  derivative operator
$\morph{\lderv: \ddual{B}\tenss[\Zc(\K)]\kQc\tenss[\Ke]\K}{\kQc}$
such that for all $\alpha \in \ddual{B}$,  $x\in B$ and
 $v\in\kQc$ we have:
 $\lderv[\alpha](xv\tenss[e]1)=\lderv(\alpha\otimes(xv\tenss[e]1))=\alpha(x)\mul
 v$.
Now, as in the simply-laced case, the cyclic derivative operator
$\morph{\cderv: \ddual{B}\tenss[\Zc(\K)]\kQc\tenss[\Ke]\K}{\kQc}$
acts on symmetric potentials  as follows:
let $v=x_0\dotsm x_n\in B^{(n+1)}$ be any  homogeneous tensor indexed
over  $\mb{Z}_{n+1}:=\set{0,\dotsc,n}$,
    with $x_i\in B$ and with corresponding ordinary cyclic permutation
    $\cperm[\mrm{ord}]=\Som[n]{i=0}{x_i\dotsm x_{i+n}}$, then
    \\[-1.0\baselineskip]
\begin{equation}\label{eq.cderv-ord}
\cderv[\alpha](v\tenss[e] 1):=\cderv(\alpha\otimes (v\tenss[e] 1))
= \lderv[\alpha](\cperm[\mrm{ord}] v\tenss[e] 1) =  
\Som[n]{i=0}{\alpha(x_i)\mul (x_{i+1}\dotsm x_{i+n})} 
\end{equation}
As  for general potentials, to a symmetric
potential $S\in\kQc\tenss[\Ke] \K$ is associated a \emph{Jacobian ideal}
 $\J{S}:=\bbar{\ggen{\im(\cderv S)}}$.

  The  next result  shows that the class of Jacobian ideals obtained from
symmetric potentials in $\kQc\tenss[\Ke] \K$ and the corresponding
ordinary cyclic derivative is exactly  the class of Jacobian
ideals obtained from symmetric potentials in $\zc(\kQc_{(2)})$ and
cyclic skew permutation and cyclic skew derivative.

\begin{Prop} \label{prop.dual-as_bim.sympot-cycderv}
 \begin{itemize}
 \item[$\msf{(1)}$] The trace    of the symmetric  algebra $(\K,\tr)$ yields
a bimodule isomorphism $\morph[/above={\sim}]{\htr: \ddual{B}}{\dual{B}}$ such that, for all $\alpha\in\ddual{B}$ we have:
 $\bilf(\htr(\alpha)\otimes\sdash)=(\id\otimes\tr)\circ
\alpha$ or equivalently, $\bilf(\sdash\otimes
\htr(\alpha))=(\tr\otimes \id)\circ \alpha$.  And for all $\xi\in \dual{B}$ and $x\in B$  we have:
$(\htr^{-1}(\xi))(x):=\Som{s\in\Lambda}\bilf(\xi\otimes e_s
x)\otimes \dual{e_s} =\Som{s\in\Lambda} e_s\otimes \bilf(
x\dual{e_s} \otimes x)$.
Thus  $\set{B,\ddual{B}}$ is a dualizing pair  naturally isomorphic to the pair
$\set{B,\dual{B};\bilf}$.
\item[$\msf{(2)}$] Let $S \in\kQc\tenss[\Ke] \K$ be a symmetric potential and
$\ztc(S)$ its  image in $\zc(\kQ_{(2)})$. Then
$\cperm(\ztc(S))=\ztc(\cperm S)$; and  for every
$\alpha\in\ddual{B}$ and $\xi=\htr(\alpha)\in\dual{B}$  we have $\cderv[\alpha]S=\cderv[\xi]\ztc(S)$.
\end{itemize}
\end{Prop}

\begin{prv} The proof of  part $\msf{(1)}$ follows by a direct application of the characterizing identities~\eqref{casimir-K} for the Casimir element
$\z[\Ke]=\Som{s\in\Lambda}e_s\otimes \dual{e_s}$.  Let us prove part $\msf{(2)}$.
Let $S=v\tenss[e] 1 \in\kQc\tenss[\Ke] \K$. To show that $\cperm(\ztc(S))=\ztc(\cperm S)$, we may assumed without lost of generality    that $v=x_0\dotsm x_n\in B^{(n+1)}$  is an homogeneous  element indexed over the cyclic group
$\mb{Z}_{n+1}:=\set{0,\dotsc,n}$, with $x_i\in B$.   We have
$\cperm(v\tenss[e] 1) = (\Som[n]{i=0}{x_i\dotsm x_{i+n}})\tenss[e] 1$. Writing
$\z[\dual{B}\otimes B]=\Som{x\in \lmQ[1]}\dual{x}\otimes x$, we
compute the skew left permutation of $\ztc(S)$ as follows:\\
$\begin{array}{l}
\hspace{-0.5em} \lperm(\ztc(S))\hspace{-0.2em}= \lperm(\Som{s\in\Lambda}e_s
x_0\mul x_1\dotsm x_n\dual{e_s}) = \Som{x\in
\lmQ[1]}\Som{s\in\Lambda} \bilf(\dual{x}\otimes e_s x_0)\mul
x_1\dotsm x_n\dual{e_s} \otimes x  \\
\hspace{1.5em} =  \Som{r\in\Lambda}
\Som{x\in \lmQ[1]}\Som{s\in\Lambda} e_r
\tr(\dual{e_r}\bilf(\dual{x}\otimes e_s x_0)) x_1\dotsm
x_n\dual{e_s} x   =  \Som{r\in\Lambda}e_r x_1\dotsm x_n \Bigl(
\Som{x\in \lmQ[1]}\bigl(\Som{s\in\Lambda} \dual{e_s}
\tr(\bilf(\dual{e_r}\dual{x}\otimes e_s x_0)) \bigr) x\Bigr) \\ 
\hspace{1.5em} =\hspace{-0.2em}\Som{r\in\Lambda}e_r x_1\dotsm x_n
 \Bigl(\Som{x\in \lmQ[1]}\bigl(\Som{s\in\Lambda} \dual{e_s} \tr(\bilf(e_s
x_0 \otimes \dual{e_r}\dual{x}))\bigr)  x \Bigr)\hspace{-0.2em}=\hspace{-0.2em} 
\Som{r\in\Lambda}e_r x_1\dotsm x_n
 \Bigl(\Som{x\in \lmQ[1]}\bigl(\Som{s\in\Lambda} \dual{e_s} \tr(e_s \bilf(
x_0\dual{e_r} \otimes \dual{x}))\bigr)  x\Bigr) \\
\hspace{1.5em} =  \Som{r\in\Lambda}e_r x_1\dotsm x_n
 \Bigl(\Som{x\in \lmQ[1]}\bilf(x_0\dual{e_r} \otimes \dual{x})
 x\Bigr)  = \Som{r\in\Lambda}e_r x_1\dotsm x_n
x_0\dual{e_r} = \ztc((x_1\dotsm x_n x_0)\tenss[e] 1).
 \end{array} $ \\
We deduce  that: $\cperm(\ztc(S)) =
\Som[n]{i=0}\lperm[i](\ztc(S)) =\ztc\bigl( (\Som[n]{i=0}{x_i\dotsm
x_{i+n} })\tenss[e] 1 \bigr) = \ztc(\cperm S)$.
Next, $S$ being assumed to be any general symmetric potential, write  $\cperm S
=(\Som[p]{k=1}x_k\mul v_k)\tenss[e] 1$ with $x_k\in B$ and $v_k\in
\kQc$ for each $k\in \ninterv{p}$.   For all  $\xi\in \dual{B}$
and $\alpha= \htr^{-1}(\xi)\in \ddual{B}$, we have: \\
 $\begin{array}{rl}\cderv[\alpha] S= &
\Som[p]{k=1}\alpha(x_k)\mul v_k = \Som[p]{k=1}\Som{s\in\Lambda}(
\bilf(\xi \otimes e_s x_k)\otimes \dual{e_s}) \mul v_k =
\Som[p]{k=1}\Som{s\in\Lambda}\bilf(\xi \otimes e_sx_k) \mul
v_k\mul \dual{e_s} =\lderv[\xi]\bigl(\Som{s\in\Lambda}
\Som[p]{k=1} e_s \mul x_k\mul v_k \dual{e_s} \bigr)   \\
=& \lderv[\xi] \bigl( \ztc(\cperm S) )= \lderv[\xi] \bigl( \cperm(\ztc(S)) \bigr) = \cderv[\xi] \bigl(
\ztc(S) \bigr),
\end{array}$ \\
completing the proof of part $\msf{(2)}$.
\end{prv}

\ParIt{A connection with the simply laced framework} To the modulated quiver $\mQ=(B,\K,\tr)$ is associated a \emph{$\kk$-quiver} $\mQk:=(B_{\kk},\kk^{\mfr{n}})$  described in the following way: chose a (finite) system $\set{\idl[i]\, :i\in\ninterv{\mfr{n}} }$ of central primitive orthogonal idempotents for $\K$, so that  $\K$ appears as direct product $\prd{1\leq i\leq \mfr{n}}{\kk_i}$ of indecomposable $\kk$-algebras with $\kk_i=\idl[i]\mul \K\mul\idl[i]$;  then   $\kk^{n}$ is the direct product $\prd{1\leq i\leq \mfr{n}}{\kk\mul\idl[i]}$;  the $\K$-bimodule $B$ is obviously a central $\kk^{\mfr{n}}$-bimodule, which we denote by $B_{\kk}$. Next,  the path algebra of the $\kk$-quiver  $\mQk$ is the tensor algebra $\kQk$ of the central $\kk^{\mfr{n}}$-bimodule $B_{\kk}$, we write $\kQkc$ for the complete path algebra of $\mQk$. The identity map $\id_{B}$  yields a natural surjective morphism of topological path algebras $\morph[dir=->>]{\pi: \kQkc}{\kQc}$. We may refer to  $\kQkc$ as  the \emph{simply laced counterpart} of $\kQc$. The simply-laced study of quivers with potentials applies to $\kQkc$, (the framework of \cite{DWZ} is obtained precisely when $\kk$ is a field).
Note that the  central $\Zc(\K)$-bimodule  $\kQc\tenss[\Ke] \K$  is  obviously  a  central $\kk^{\mfr{n}}$-bimodule. In view of  Lemma~\ref{lem.casimi-sympot} we get the following useful  connection with the simply-laced study.
\begin{Rem} \label{rem.connection-simplylaced} The natural surjective map of topological bimodules $\morph[/dir=->>]{\pi_{\mrm{c}}: \kQkc}{\kQc\tenss[\K^e]\K \,:x\mapsto \pi(x)\otimes 1}$  preserves permutations of tensors elements and cyclical equivalence, hence any property of a simply laced potential $w \in \kQkc$ with respect to cyclical equivalence is transferred  to the symmetric potential $w\tenss[e] 1  \in \kQc\tenss[\K^e] \K$.  In particular, when   $\kk$ is a field, the study of quivers with potentials with respect to cyclical equivalence applies to potentials in $\kQc\tenss[\Ke] \K$.
\end{Rem}

The above connection been made, we can derive the following useful result on symmetric potentials.
\begin{Lem}  \label{lem.cycequiv-IJ} Suppose  $\kk$ is a field.  Let $I$ be a closed ideal in $\J{\kQc}$ and $J$ the closure of an ideal generated by finitely many elements $m_1,\dotsc,m_p \in \J{\kQc}$.   Then any symmetric potential $S$ belonging to  $(I\mul J) \tenss[\Ke] 1$ is cyclically equivalent to a symmetric potential lying in $(I\mul m_1+\dotsm+ I\mul m_p) \tenss[\Ke] 1$, thus $\ztc(S)$ is cyclically equivalent to a symmetric potential $W$ lying in $(I\mul m_1+\dotsm+ I\mul m_p)$.
\end{Lem}

\begin{prv} Let $S=v\tenss[e]1$  with  $v\in IJ$. Then in view of Remark~\ref{rem.connection-simplylaced},  the fact that $S$ is cyclically equivalent to a symmetric potential lying in $(I\mul m_1+\dotsm+ I\mul m_p) \tenss[\Ke] 1$  is given by the corresponding simply-laced result  in $\kQkc$ (see  \cite[Lem~13.8]{DWZ}). By part $\msf{(2)}$ of Proposition~\ref{prop.dual-as_bim.sympot-cycderv}, $\cperm(\ztc(S))= \ztc(\cperm S)$, implying that the symmetric potential $\ztc(S)$ is cyclically equivalent to a symmetric potential $W$ lying in $(I\mul m_1+\dotsm+ I\mul m_p)$.
\end{prv}

Let us mention that  a direct proof of Lemma~\ref{lem.cycequiv-IJ} above, though a little bit  technical, is still possible. Indeed, the notions of ``$C$-space and $D$-space`` used in \cite[\S13]{DWZ} to prove the simply laced analogue of Lemma~\ref{lem.cycequiv-IJ}  are easily seen to be  special cases of   symmetrizable weakly dualizing pairs of bimodules  $\set{M,\dual{M}}$ where $M$ arises as a union $\Union[\infty]{n}{M_n}$ of  $\K$-bimodules $M_n$ which are finite-dimensional as free and semisimple $\kk$-modules, with $0=M_0\subset M_1\subset M_2\subset \dotsm$.

The assumption that  $\kk$ is a field is directly required only by  Lemma~\ref{lem.cycequiv-IJ} above, it enables us to quickly  establish the second main result of this work as follows.

\begin{Thm}\label{theo.sympot} Under the Casimir operator  $\morph[/dir=->>]{
\ztc: \kQc\tenss[\Ke]\K}{\zc(\kQc)}$ and  the natural isomorphism 
$\ddual{B}\cong \dual{B}$, ordinary permutations and cyclic derivatives of symmetric potentials from $\kQc\tenss[\Ke] \K$ agree with skew permutations and cyclic derivatives of their images in $\zc(\kQc)$, and when the Casimir ideal
$\zc(\K)$ coincides with the center of $\K$, all potentials on
$\mQ$ are symmetric. Moreover, over a  field $\kk$, the split reduction of modulated quivers with symmetric potentials $(\mQ,\mfr{m})$  such that the cyclic derivative
$\morph[/dir=->>,/above={\cderv \mfr{m}}]{\dual{B}}{\im(\cderv
\mfr{m})}$ also splits can be defined up to right-equivalences.
\end{Thm}

Before proving Theorem~\ref{theo.sympot}, the following question
retains our attention.
\begin{Quest} \label{Quest.sym-casimir-operSurj}
When does the Casimir ideal $\zc(\K)$ of  a symmetric algebras $\K$ coincide with the center  of $\K$?
\end{Quest}
Recall the following definition to be  compared  with
\cite[Defn~2.1]{Aguiar}, \cite[thm3.1]{KS} presenting  nine
equivalent characterizations of symmetrically separable algebras.
\begin{Defn}  A $\kk$-algebra $A$ is  \emph{symmetrically separable} (or \emph{strongly
separable}) if there exists a $\kk$-linear trace $\tau$ on $A$ such
that $(A,\tau)$ is  symmetric  and the associated Casimir
element
$\Som[\mfr{r}]{s=1}\varepsilon_s\otimes\dual{\varepsilon_s}$ is a
(symmetric) separability idempotent for $A$: in particular
$\Som[\mfr{r}]{s=1}\varepsilon_s\dual{\varepsilon_s} = 1=
\Som[\mfr{r}]{s=1}\dual{\varepsilon_s} \varepsilon_s$.
\end{Defn}

\begin{enumerate}
\item[$\trr$] If $ \zc(\K)=\Zc(\K)$,     then for all $\K$-bimodule $M$ we also have $\zc(M)=\Zck(M)$.
 \item[$\trr$] Part $\msf{(a)}$ of
 Lemma~\ref{lem.trace-dualizing-pair} shows that  if $\K$ is symmetrically separable algebra then $ \zc(\K)=\Zc(\K)$.
\end{enumerate}
 Thanks to a well-known Higman's Theorem,
  Question~\ref{Quest.sym-casimir-operSurj} is completely solved when the ground ring is a field.
\begin{Rem}\label{rem.product-of-separableAlg-field}
\begin{enumerate}
 \item[$\trr$] By Higman's Theorem \cite[thm~10]{Higman1955a}
  (or \cite[thm~1]{Higman1955b}), separable algebras over a field are exactly those
  symmetric algebras $\K$ such that $\zc(\K)=\Zc(\K)$.
 \item[$\trr$] By a well-known result (see  P.M. Cohn {\cite[Cor~11.6.8]{Cohn2003}}),
 the tensor product over a field   of  two separable algebras is again 
separable  and hence   semisimple.
\end{enumerate}
\end{Rem}

\begin{Cor} \label{cor.sympot-perfectfield} If $\K$ is separable over a  field $\kk$,
then so is the enveloping algebra $\Ke$, $\zc(\K)=\Zc(\K)$, all
potentials on $\mQ$ are symmetric  and can treated
symmetrically using ordinary cyclic permutation and ordinary
cyclic derivative from $\kQc\tenss[\Ke] \K$, and the reduction of
every $\kk$-modulated quiver with potential  is well-defined  up to right-equivalences. This is in particular the case  when   $\kk$ is a perfect  field.
\end{Cor}

\subsubsection*{Proof of Theorem~\ref{theo.sympot}.}
The first part  of Theorem~\ref{theo.sympot}~is a direct consequence of  Proposition~\ref{prop.dual-as_bim.sympot-cycderv}. We dedicate the rest of this section to establish the last part of Theorem~\ref{theo.sympot}.
Thus $(\K,\tr)$ is   symmetric  over a field $\kk$,  and   $(\mQ,\mfr{m})$ is a
modulated quiver  with symmetric potential having a
split trivial part such that 
$\morph[/dir=->>]{\cderv \mfr{m}: \dual{B}}{\im(\cderv
\mfr{m})}$ also splits. We need  to construct  and manipulate
 unitriangular automorphisms.  In the following   the assumption that $\kk$ is a field is not used.

\begin{Lem} \label{lem.cycequiv-unitriang}  Let $\morph[/dir=->>]{f: \dual{B}}{M}$ be a split bimodule epimorphism   with right inverse $f'$,   $W\in
M'\mul M$ a potential with $M'\subset \J[2]{\kQc}$ and $M\subset \J{\kQc}$, and put  $W'=(\id_{M'} \otimes
f')(W)$. Then the following assertions hold.
\begin{itemize}
\item [$\msf{(1)}$]  There is  a bimodule morphism $\morph{\alpha: B}{M'}$ such that  $W=\Som{y\in\rmQ[1]}{b_y\otimes f(\dual{y})}$ with $b_y:=\alpha(y)$ for all $y\in \rmQ$. There is a  unitriangular automorphism $\morph{\phi: \kQc}{\kQc}$ with
\mbox{$\phi\restr{B} = \id_{B} + \alpha$}.
 \item [$\msf{(2)}$] Let $S$ be a reduced potential on $\mQ$ and $\phi$   the unitriangular automorphism above. Then
 $\phi(S)-S-\Som{y\in\rmQ[1]}{b_y\cderv[\dual{y}]S}$ is cyclically
 equivalent to a potential $S'$ in $\J{\kQc}\mul I^2$ where $I$ is the closed
 ideal   given by: $I=\bbar{\ggen{\set{b_y\, :\,
 y\in\rmQ[1]} }} =\bbar{\ggen{\set{b_x\, :\,
 x\in\lmQ[1]} }}$ with $b_x=\rderv(W'\otimes x)$ for each
 $x\in\lmQ[1]$. Moreover, if $S$ is symmetric, then so is $S'$.
\end{itemize}
\end{Lem}

\begin{prv} Let us prove $\msf{(1)}$. The element $W'=(\id_{M'} \otimes
f')(W) \in \J[2]{\kQc}\otimes \dual{B}$ is   $\K$-central
 since $W$ is,  hence $W'$ is a potential on
$\mQ \oplus \dual{\mQ}$ and we get a bimodule
 morphism  $\morph{\alpha:=\rderv W':
B}{M'}$. In view of
\eqref{eq.dervdualbasesKQc} we have
$W'=\Som{y\in\rmQ[1]}{\rderv(W'\otimes y)\otimes\dual{y}}$, and since $f\circ f'=\id_{M}$  we get
 $W=(\id_{M'} \otimes  f)(W')=\Som{y\in\rmQ[1]}{b_y\otimes f(\dual{y})}$
 with $b_y=\rderv(W'\otimes y)=\alpha(y)$ for each $y\in\rmQ[1]$.
By Proposition~\ref{prop.prolong}, the bimodule morphism 
$\morph{\id_{B}+ \alpha: B}{B \oplus M'\subset B\oplus  \J[2]{\kQc}}$ induces  a unitriangular automorphism $\phi$ of $\kQc$ with $\phi\restr{B}=\id_B+\alpha$.

We now turn to the proof of part $\msf{(2)}$. We have  $\set{b_y\, :\,
 y\in\rmQ[1]}\mul \K = \im(\rderv W') =\K\mul \set{b_x\, :\,
 x\in\lmQ[1]}$ with $b_y=\rderv(W'\otimes y)$ and $b_x=\rderv(W'\otimes x)$
 for each $(x,y)\in\lmQ[1]\times \rmQ[1]$. We start with the case of an homogeneous potential
$S=\Som[n]{i=1}{u_{i,0}\mul u_{i,1}\dotsm \mul u_{i,d} }$ with
$d\geq 2$ and $u_{i,r} \in B$ for all $(i,r)\in \ninterv{n} \times \ninterv[0]{d}$. As in the statement of the cyclic
Leibniz rule~\eqref{eq.cyclicLeibnizrule}, for each $(i,r)\in
\ninterv{n}\times \ninterv[0]{d}$ we  write: \mbox{$u_{i,<r}=\prod\limits_{k=0}^{r-1}u_{i,k}$},
\mbox{$u_{i,>r}=\prod\limits_{k=r+1}^{d}u_{i,k}$} and
\mbox{$u_{i,\geq r}=\prod\limits_{k=r}^{d}u_{i,k}$}, where the
empty products coincide with $\idl\in\K$. Then,
expanding $\phi(S)$ we write:
$\phi(S)=\Som[n]{i=1}{(u_{i,0}+\alpha(u_{i,0}))\mul
(u_{i,1}+\alpha(u_{i,1}))\dotsm \mul (u_{i,d}+\alpha(u_{i,d}))} =
S + S_1 + S_{(2)}$, where
$S_1=\Som[d]{r=0}\Som[n]{i=1}{u_{i,<r}\mul \alpha(u_{i,r}) \mul
u_{i,>r}}$ is a potential while the potential $ S_{(2)}$ is  the
sums of all the rest of the terms in the expansion of $\phi(S)$
containing at least two occurrences of the form $\alpha(u_{i,k})$. We show
that $S_1 \cyc{\equiv} \Som{y\in\rmQ[1]}{b_y\cderv[\dual{y}]} S$.
For each $r\in \ninterv[0]{d}$,  the term
$S_{1,r}:=\Som[n]{i=1}{u_{i,<r}\mul \alpha(u_{i,r}) \mul
u_{i,>r}}$ is an homogeneous potential, which is then cyclically
equivalent to the left permutation $\lperm[r] S_{1,r}$ and we
have:\\
 $ \begin{array}{rl}
\lperm[r] S_{1,r} = &  \Som{x\in
\lmQ[r]}\Som[n]{i=1}{\bilfh(\dual{x}\otimes u_{i,<r})\mul
\alpha(u_{i,r}) \mul u_{i,>r}} = \Som{x\in
\lmQ[r]}\Som[n]{i=1}{\alpha(\bilfh(\dual{x}\otimes u_{i,<r})\mul
u_{i,r}) \mul u_{i,>r}} \\
  = &  \Som{x\in \lmQ[r]}\Som[n]{i=1}{\alpha\bigl(
  \Som{y\in \rmQ[1]} y\mul \bilf(\dual{y}\otimes \bilfh(\dual{x}\otimes u_{i,<r})\mul
u_{i,r})\bigr) \mul u_{i,>r}}  \\
= &   \Som{y\in \rmQ[1]} \alpha(y) \mul \Som{x\in
\lmQ[r]}\Som[n]{i=1}{
  \bilf(\dual{y} \otimes \bilfh(\dual{x}\otimes u_{i,<r})\mul
u_{i,r}) \mul u_{i,>r}} \\
  = &  \Som{y\in \rmQ[1]} \alpha(y)
\mul \lderv[\dual{y}]\bigl( \Som{x\in \lmQ[r]}\Som[n]{i=1}{
  \bilfh(\dual{x}\otimes u_{i,<r})\mul
u_{i,\geq r}}  \bigr)  = \Som{y\in \rmQ[1]} \alpha(y)
\mul \lderv[\dual{y}]\bigl( \lperm[r] S \bigr) \\
 = &   \Som{y\in \rmQ[1]} b_y \mul
\lderv[\dual{y}]\bigl( \lperm[r] S \bigr).
\end{array} $ \\
Thus, $S_1 \cyc{\equiv} \Som[d]{r=0} \lperm[r] S_{1,r} =
 \Som[d]{r=0}\Som{y\in \rmQ[1]} b_y \mul \lderv[\dual{y}]\bigl(
\lperm[r] S \bigr) = \Som{y\in \rmQ[1]} b_y \mul
\lderv[\dual{y}]\bigl( \Som[d]{r=0} \lperm[r] S \bigr) = \Som{y\in
\rmQ[1]} b_y \mul \lderv[\dual{y}]\bigl(\cperm S \bigr) =
 \Som{y\in \rmQ[1]} b_y \mul \cderv[\dual{y}] S$. Next,   write $S_{(2)}$ as  sum of  potentials of the form
 $S_{(2),r}=\Som[n]{i=1}v_r\alpha(u_{i,r})u_{i,>r}$, with
 $r\in\ninterv{d}$ and $v_r\in \J[r+1]{\kQc}\cap I$. As above, we get that the right
 permutation $\rperm[d-r] S_{(2),r}$ belongs to
$\Som{x\in \lmQ[1]} \J{\kQc}\mul (\J[d]{\kQc} \cap I)b_x$,
implying   that $S_{(2)}$ is cyclically equivalent to an element
of $\Som{x\in \lmQ[1]} \J{\kQc}\mul (\J[d]{\kQc} \cap I)b_x$.

Now, for a general potential $S\in \J[3]{\kQc}$ written as  sum
of homogeneous potentials, the previous discussion shows that
$\phi(S)-S-\Som{y\in\rmQ[1]}{b_y\cderv[\dual{y}]S}$ is cyclically
 equivalent an element of the form $\Som{x\in \lmQ[1]}{c_x b_x}$,
 where $c_x=\Som{y\in \rmQ[1]} y \Som{l\geq 3}c_{x,y,l}$ with $c_{x,y,l} \in
 \J[l-1]{\kQc} \cap I$. Since $I$ is closed, each $c_x$ is a
 well-defined element of $\J{\kQc}I$. Thus   $\phi(S)-S-\Som{y\in\rmQ[1]}{b_y\cderv[\dual{y}]S}$
 is cyclically  equivalent to a potential $S'$ in $\J{\kQc}I^2$.
 Finally, if  $S$ is  symmetric, then part $\msf{(2)}$ of Proposition~\ref{prop.dual-as_bim.sympot-cycderv} implies
  that  $S'$ is also symmetric.
\end{prv}

As in the simply-laced case, the fact  that the reduction of
$(\mQ,\mfr{m})$ is  defined up to right-equivalences  will be derived as 
consequence of the following result whose proof relies on
Lemma~\ref{lem.cycequiv-unitriang}  and Lemma~\ref{lem.cycequiv-IJ}.
\begin{Prop} \label{prop.rightequiv} Let $S$ and $S'$ be reduced
potentials on $\mQ$ such that $S'-S$ is cyclically
equivalent to a potential   $S''$ in $(\J{S})^{2}$. Then the following statements hold.
\begin{enumerate}
    \item [$\msf{(1)}$] $\J{S}= \J{S'}$.
    \item [$\msf{(2)}$] Suppose $S$ and $S''$ are symmetric
 and   the cyclic derivative $\morph[/dir=->>,/above={\cderv
S}]{\dual{B}}{\im(\cderv S)}$ splits. Then $(\mQ,S)$ is
right-equivalent to  $(\mQ, S')$; more precisely  there exists a
unitriangular automorphism  $\phi$ of $\kQc$ such that $\phi(S)$
is cyclically equivalent to $S'$ and $\phi(u)-u\in \J{S}$ for all
$u\in \kQc$.
\end{enumerate}
\end{Prop}

\begin{prv} Put $\what{J}=\J{\kQc}$.  We have $\J{S} \subset \what{J}^2$ because 
$(\mQ,S)$ is reduced. Part $\msf{(1)}$ is an easy consequence of the
cyclic-Leibniz rule~\eqref{eq.cyclicLeibnizrule}: indeed, for all $\xi\in
\dual{B}$ we have
$\cderv[\xi]S' - \cderv[\xi]S =\cderv[\xi](S'-S)=\cderv[\xi] S''
\in \what{J} \J{S} + \J{S} \what{J}$,
 implying in view of part $\msf{(b)}$ of Lemma~\ref{lem.trivIdeal-closedIdeal} that  $\J{S'}= \J{S}$.

For part  $\msf{(2)}$, by assumption we have a split epimorphism   $\morph{f:
\dual{B}}{M }$ with  $M=\cderv(\dual{B} \otimes S)$ and $f(\xi)=\cderv[\xi] S$ for all $\xi\in \dual{B}$.   Let $f'$ be a right inverse for $f$.  Using induction on $n$, we
will construct a sequence of unitriangular automorphisms
$\morph{\phi_n: \kQc}{\kQc}$, with $n\geq 0$ and $\phi_0=\id_{\kQc}$, taking each
generator $y\in \rmQ[1]$ to $y+b_{y,n}$ and having the following
properties:
\begin{enumerate}
    \item [$\msf{(i)}$] $b_{y,n} \in \what{J}^{n+ 1} \cap\J{S}$ for all $y\in
    \rmQ$  and   all $n\in \N$.
    \item [$\msf{(ii)}$] The sum $ \Som{y\in\rmQ} b_{y,n} \cderv[\dual{y}]S$ is a symmetric potential
    and $S'$ is cyclically equivalent to the symmetric
    potential $\phi_0\phi_1\dotsm\phi_{n-1}\bigl(S + \Som{y\in\rmQ} b_{y,n} \cderv[\dual{y}]S     \bigr)$ for all $n\geq 1$.
    \end{enumerate}
    The existence of $\phi_1$ with the desired properties follows by
 part $\msf{(1)}$ of Lemma~\ref{lem.cycequiv-unitriang} and by
    Lemma~\ref{lem.cycequiv-IJ} in which we take $I=\J{S}=J$ and
$m_y=f(\dual{y})=\cderv[\dual{y}]S$.

 Now assume that, for some $n\geq 1$, we have already defined
 $\phi_1,\dotsc, \phi_n$ having the
desired properties. We then want to construct a unitriangular
automorphism $\phi_{n+1}$ such that properties $\msf{(i)}\sdash
\msf{(ii)}$ are satisfied with $n$ replaced by $n+1$.   By part
$\msf{(2)}$ of Lemma~\ref{lem.cycequiv-unitriang},
\mbox{$(S+\Som{y\in\rmQ[1]}{b_{y,n}\cderv[\dual{y}]S})-\phi_n(S)$}
is cyclically  equivalent to a symmetric potential $W_1$ belonging
to $\what{J}(\what{J}^{n+1} \cap \J{S})^2$. In particular,
observing that $\what{J}(\what{J}^{n+1} \cap \J{S})^2 \subseteq
(\what{J}^{n+2} \cap \J{S})\J{S}$, we deduce that $S-\phi_n(S)$ is
cyclically equivalent to a symmetric potential in $(\J{S})^2$.
Thus combining part $\msf{(1)}$ of Theorem~\ref{theo.Jm-et-morph-Jtriv}
together with the already proved part $\msf{(1)}$  of
Proposition~\ref{prop.rightequiv}, we conclude that
$\phi_n(\J{S})=\J{\phi_n(S)}=\J{S}$. It follows that the symmetric
potential
$(S+\Som{y\in\rmQ[1]}{b_{y,n}\cderv[\dual{y}]S})-\phi_n(S)$ is in
fact cyclically  equivalent to a symmetric potential $\phi_n(W_2)$
belonging to $\phi_n\bigl( (\what{J}^{n+2} \cap \J{S})\J{S}
\bigr)$, where $W_2=\phi_n^{-1}(W_1)$ is a symmetric potential in
$(\what{J}^{n+2} \cap \J{S})\J{S}$. But then applying
Lemma~\ref{lem.cycequiv-IJ} to $I=\what{J}^{n+2}\cap \J{S}$ and
$J=\J{S}$, we get that $W_2$ is cyclically equivalent to a
symmetric potential $W$ lying in $(\what{J}^{n+2} \cap \J{S})\mul M$, and part $\msf{(1)}$ of Lemma~\ref{lem.cycequiv-unitriang} yields a unitriangular automorphism $\morph{\phi_{n+1}: \kQc}{\kQc}$  taking each $y\in\rmQ[1]$ to an element $y+ b_{y,n+1}$ with $b_{y,n+1}\in \what{J}^{n+2} \cap \J{S}$ and such that  $W=\Som{y\in\rmQ[1]}b_{y,n+1}
\cderv[\dual{y}] S$. Now, the fact that $(S+\Som{y\in\rmQ[1]}{b_{y,n}\cderv[\dual{y}]S})-\phi_n(S)$ is cyclically equivalent to $\phi_n(W)=\phi_n\bigl( \Som{y\in\rmQ[1]}{b_{y,n+1}\cderv[\dual{y}]S} \bigr)$ shows that \mbox{$(S+\Som{y\in\rmQ[1]}{b_{y,n}\cderv[\dual{y}]S})$} is cyclically equivalent to
\mbox{$\phi_n\bigl(S+ \Som{y\in\rmQ[1]}{b_{y,n+1}\cderv[\dual{y}]S} \bigr)$}, thus the assumption  that $S'$ is   cyclically equivalent to \mbox{$\phi_0\phi_1\dotsm\phi_{n-1}\bigl(S + \Som{y\in\rmQ} b_{y,n} \cderv[\dual{y}]S  \bigr)$} shows that $S'$ is cyclically equivalent to \mbox{$\phi_0\phi_1\dotsm\phi_{n-1}\phi_n\bigl(S + \Som{y\in\rmQ} b_{y,n+1} \cderv[\dual{y}]S \bigr)$}. We have therefore constructed a unitriangular automorphism
$\phi_{n+1}$ such that properties $\msf{(i)}\sdash \msf{(ii)}$ are
satisfied with $n$ replaced by $n+1$, completing the induction
step.

Now, in view of property~$\msf{(i)}$, letting
$\phi=\lim\limits_{n\rightarrow \infty} \phi_1\dotsm\phi_n$, we
get a well-defined unitriangular automorphism of $\kQc$ such that
$\phi(u)-u \in \J{S}$ for all $u\in \kQc$. And letting $n$ tends
to $\infty$ in property~$\msf{(ii)}$, we conclude that $S'$ is
cyclically equivalent to $\phi(S)$, completing the proof of part
$\msf{(2)}$ of Proposition~\ref{prop.rightequiv}.
\end{prv}

Now,  using Proposition~\ref{prop.rightequiv}, we will show that the reduction of
$(\mQ,\mfr{m})$ can be defined up to right-equivalences. Thus, let
$\morph{\phi: \kQc}{\kQpc}$ be a right-equivalence between
$(\mQ,\mfr{m})$ and a modulated quiver with potential $
(\mQ',\mfr{m}')$ where $\mQ'=(B',\K,\tr)$. Since $\phi$ is
obviously a weak right-equivalence and   the reduction is
 defined up to weak right-equivalences by the reduction
Theorem~\ref{theo.red-mQp}, we derive the following conclusions:
$(\mQ',\mfr{m}')$  has a split trivial part, and  keeping the notations of \eqref{cond-loopfree} and
\eqref{eq.loopfreedecomp-mQ} we write: $\mQ=\triv{\mQ} \oplus
\mQb$, $\mQ'=\triv{\mQ}' \oplus \mQb'$ as direct sums of naturally
induced  modulated quivers where $\triv{\mQ}=(\triv{B},\K,\tr)$
with $\triv{B}=\cderv(\dual{B}\otimes \mfr{m}_{2})$ and
$\mfr{m}_{2}$ denotes the degree-$2$ component of $\mfr{m}$,
$B=\triv{B}\oplus \bbar{B}$, $\mQb=(\bbar{B}, \K, \tr )$; and
similarly $\triv{\mQ}'=(\triv{B}',\K,\tr)$  with
$\triv{B}'=\cderv(\dual{B'}\otimes {\mfr{m}}_{2}')$,  the
degree-one component $\morph[/above={\sim}]{\phi_{1}:B}{B'} $ of
$\phi$ is an isomorphism with $\phi_{1}(\triv{B})=\triv{B}'$ and
$B'=\triv{B}' \oplus \bbar{B}'$, $\mQb'=(\bbar{B}', \K, \tr )$.
Still by   Theorem~\ref{theo.red-mQp}, consider the
reduction $\morph[/dir=->>]{\pi_{\mfr{m}}: \kQc}{\kQbc}$ from
$(\mQ,\mfr{m})$ into  
$\Red(\mQ,\mfr{m})=(\mQb,\bbar{\mfr{m}})$ and
reduction $\morph[/dir=->>]{\pi_{\mfr{m}'}: \kQpc}{\kQbpc}$ from
$(\mQ',\mfr{m}')$ into $\Red(\mQ,\mfr{m})=(\mQb',\bbar{\mfr{m}}')$. Recall that
\begin{center}
$\bbar{\mfr{m}}=\pi_{\mfr{m}}(\mfr{m})$,
$\bbar{\mfr{m}}'=\pi_{\mfr{m}'}(\mfr{m}')$,
$\pi_{\mfr{m}}\restr{\kQbc}=\mrm{id}$ and
$\pi_{\mfr{m}'}\restr{\kQbpc}=\mrm{id}$.
\end{center}
 Next,  $\phi$ induces a weak right-equivalence
$\morph{\psi:=\pi_{\mfr{m}'}\phi\restr{\kQbc}:
\kQbc}{\kQbpc}$ 
 between reduced modulated quivers with potentials $(\mQb,\bbar{\mfr{m}})$ and
$(\mQb',\bbar{\mfr{m}}')$, and we have:
\begin{equation} \label{eq.prvthmsym} \tag{$\ast$}
\phi(\J{\mfr{m}})=\J{\mfr{m}'}; \; \,
\pi_{\mfr{m}}(\J{\mfr{m}})=\J{\bbar{\mfr{m}}} \text{ and }
\J{\psi(\bbar{\mfr{m}})}=\psi(\J{\bbar{\mfr{m}}})=\J{\bbar{\mfr{m}}'}
=\pi_{\mfr{m}'}(\J{\mfr{m}'}).
\end{equation}
To show that $(\mQb,\bbar{\mfr{m}})$ and
$(\mQb',\bbar{\mfr{m}}')$ are   right-equivalent, by
Proposition~\ref{prop.rightequiv}  it suffices to prove that:
\begin{enumerate}
    \item [$(\alpha)$] The potential $S:=\psi(\bbar{\mfr{m}})$ is symmetric  and
    $\bbar{\mfr{m}}'-S$ is cyclically equivalent to a
    symmetric potential in $(\J{S})^2$.
    \item [$(\beta)$] The cyclic derivative $\morph[/dir=->>]{\cderv
    S: \dual{\bbar{B}'}}{\im(\cderv S)}$ splits.
\end{enumerate}

For  $(\alpha)$, we have  $S=\psi(\bbar{\mfr{m}})=
\pi_{\mfr{m}'}\phi(\pi_{\mfr{m}}(\mfr{m}))$, showing that $S$ is  symmetric
 since $\mfr{m}$ is. The last part of Lemma~\ref{lem.splitting-mQ}  shows that
$\mfr{m}-\pi_{\mfr{m}}(\mfr{m}) $ is cyclically  equivalent to a
symmetric potential $W $ lying in $ (\Ker(\pi_{\mfr{m}}))^{2}$. Since
reductions and weak right-equivalences send cyclically equivalent
potentials to cyclically equivalent ones, we
deduce that: \\
$ \begin{array}{r l l} \bbar{\mfr{m}}'-S & =&
\pi_{\mfr{m}'}(\mfr{m}')
-\pi_{\mfr{m}'}\phi(\pi_{\mfr{m}}(\mfr{m})) =
\pi_{\mfr{m}'}(\mfr{m}' -
\phi(\pi_{\mfr{m}}(\mfr{m}))) \\
 & \cyc{\equiv} &  \pi_{\mfr{m}'}\bigl( \phi(\mfr{m}) -
\phi\pi_{\mfr{m}}(\mfr{m})\bigr), \qquad\,  (\text{since by the
definition of a right-equivalence, } \mfr{m}'
 \cyc{\equiv} \phi(\mfr{m}) \,) \\
 & \cyc{\equiv} &
  \pi_{\mfr{m}'}\phi\bigl(\mfr{m}-\pi_{\mfr{m}}(\mfr{m})\bigr) \,
  \cyc{\equiv}\,   \pi_{\mfr{m}'}\phi\bigl(W \bigr) \in
    \pi_{\mfr{m}'}\phi\bigl( (\Ker(\pi_{\mfr{m}}))^{2}\bigr) \subset
  \pi_{\mfr{m}'}\phi\bigl( (\J{\mfr{m}})^{2}\bigr).
\end{array}$ \\
 But, using  \eqref{eq.prvthmsym} above we get
$\pi_{\mfr{m}'}\phi(\J{\mfr{m}})=\pi_{\mfr{m}'}(\J{\mfr{m}'})=\J{\bbar{\mfr{m}}'}
=\psi(\J{\bbar{\mfr{m}}}) =\J{S}$, so that
$\pi_{\mfr{m}'}\phi\bigl( (\J{\mfr{m}})^{2} \bigr) = (\J{S})^2$.
Hence,  $\bbar{\mfr{m}}'-S$ is cyclically equivalent to a
symmetric potential in $(\J{S})^2$, completing the proof of
$(\alpha)$.

We now turn to the proof of $(\beta)$. Observe that the map
$\morph[/dir=->>,/above={\cderv S}]{\dual{\bbar{B}'}}{\im(\cderv
S)}$ splits whenever  its kernel is a direct summand in
$\dual{\bbar{B}'}$, but since $S = \psi(\bbar{\mfr{m}})$ and $\psi$ is a path algebra isomorphism, applying part $\msf{(1)}$ of
Theorem~\ref{theo.Jm-et-morph-Jtriv} we get that, $\Ker(\cderv
\psi(\bbar{\mfr{m}}))$ is a direct summand in $\dual{\bbar{B}'}$
if and only if $\Ker(\cderv\bbar{\mfr{m}})$ is a direct summand in
$\dual{\bbar{B}}$. Hence, we have to show that the cyclic
derivative $\morph[/dir=->>,/above={\cderv
\bbar{\mfr{m}}}]{\dual{\bbar{B}}}{\im(\cderv \bbar{\mfr{m}})}$
splits. By the reduction Theorem~\ref{theo.red-mQp} (or 
Lemma~\ref{lem.splitting-mQ}), there is a unitriangular
automorphism $\morph{\vphi: \kQc}{\kQc}$ such that
$\vphi(\mfr{m})$ is cyclically equivalent to $\z[U\otimes V] +
\bbar{\mfr{m}}$ with $\triv{B}=U\oplus V$ for an 
induced dualizing pair  $\set{U,V}$. And
By by assumption the cyclic derivative
$\morph[/dir=->>,/above={\cderv \mfr{m}}]{\dual{B}}{\im(\cderv
\mfr{m})}$ splits, implying  that the cyclic derivative
$\morph[/dir=->>,/above={\cderv
\vphi(\mfr{m})}]{\dual{B}}{\im(\cderv \vphi(\mfr{m}))}$ splits, so
that:
 \begin{equation} \label{eq.prvthmsym+} \tag{$\ast\ast$}
\text{the cyclic derivative } \,
\morph[/dir=->>]{\cderv(\z[U\otimes V] + \bbar{\mfr{m}}) :
\dual{B}}{\im(\cderv(\z[U\otimes V] + \bbar{\mfr{m}}))}\,
\text{splits}.
\end{equation}
Recall   that $B=\triv{B} \oplus \bbar{B}$, $\triv{B}=U\oplus V =\im(\cderv(\z[U\otimes V]))$, $\bbar{\mfr{m}} \in \J[3]{\kQbc}$,
$\dual{B}=\dual{\triv{B}} \oplus \dual{\bbar{B}}$ where
$\dual{\bbar{B}}$   is the kernel of the cyclic
derivative $\morph[/dir=->>]{\cderv(\z[U\otimes V]) :
\dual{B}}{\im(\cderv(\z[U\otimes V]))}$,\,
 Therefore,
$\Ker\bigl(\cderv(\z[U\otimes V] + \bbar{\mfr{m}})
\bigr)=\Ker(\z[U\otimes V])\cap \Ker(\Seq{\cderv\bbar{\mfr{m}} :\dual{B} \to[/dist=1]  \J[2]{\kQbc}}) \subset \dual{\bbar{B}}$ and $\dual{\triv{B}}\subset \Seq{\cderv\bbar{\mfr{m}} :\dual{B} \to[/dist=1]  \J[2]{\kQbc}} $, implying by  \eqref{eq.prvthmsym+} above   that the cyclic
derivative $\morph[/dir=->>,/above={\cderv
\bbar{\mfr{m}}}]{\dual{\bbar{B}}}{\im(\cderv \bbar{\mfr{m}})}$
also splits and completing the proof of $(\beta)$. Hence,  the proof
of the last part of Theorem~\ref{theo.sympot} is finished.

\section{Some examples in the inseparable context}\label{secexample}
Here we illustrate the fact that the reduction of a modulated quiver with potential $(\mQ,\mfr{m})$ may still be carried even if the trivial part of $(\mQ,\mfr{m})$ does not split.

Let $\kk=\mb{F}_2(u)=\mb{F}_2[u,u^{-1}]$ be the non perfect
  function field of one variable over   the prime
field $\mb{F}_2$ of characteristic $2$; $\msf{E}=\mb{F}_2(u^{\frac{1}{2}}) =\kk\mul\{1, u^{\frac{1}{2}}\}$
and $\msf{F}=\mb{F}_2(u^{\frac{1}{4}})
=\kk\mul\{1,u^{\frac{1}{4}},u^{\frac{1}{2}}, u^{\frac{3}{4}}\}
=\msf{E}\mul\{1,u^{\frac{1}{4}}\}$. Then $\msf{E}$ and $\msf{F}$ are  finite-dimensional inseparable extensions of  $\kk$. Let $\K=\kk_1\times \kk_2\times \kk_3$ with
$\kk_1=\msf{F}=\kk_2$ and $\kk_3=\kk$;  each
$\kk_i$ is viewed as  subfield in $\K$ with unity $\idl[i]$, thus the unity of $\K$ is $\idl[1]+\idl[2] + \idl[3]$.  For  $\lambda\in\set{\frac{1}{4},\frac{1}{2},\frac{3}{4}}$ we write $u_i^{\lambda}=\idl[i]\mul u^{\lambda} \in \kk_i$. We have a symmetric $\kk$-algebra $(\K,\tr)$  where  $\tr$ is the natural $\kk$-linear trace induced by its restriction on $\msf{F}$  as follows:  $\tr(1)=1$ and
 $\tr(u^{\frac{1}{4}})=\tr(u^{\frac{1}{2}})=\tr(u^{\frac{3}{4}})=0$.
 We have Casimir elements 
 $\z[{\msf{F}\tenss[\kk] \msf{F}}]
 =1\otimes 1+ u^{\frac{1}{4}} \otimes u^{-\frac{1}{4}} +  u^{\frac{1}{2}} \otimes u^{-\frac{1}{2}}
 +  u^{\frac{3}{4}} \otimes u^{-\frac{3}{4}}$ and
$ \z[\Ke] = \Som[9]{s=1}{e_s\otimes \dual{e_s}}$ with $(e_s)_{1\leq s\leq 9}=(\idl[1],\idl[2],\idl[3],u_1^{\frac{1}{4}},u_1^{\frac{1}{2}},
 u_1^{\frac{3}{4}},u_2^{\frac{1}{4}},u_2^{\frac{1}{2}},
 u_2^{\frac{3}{4}})$ and $(\dual{e_s})_{1\leq s\leq 9}=(\idl[1],\idl[2],\idl[3],u_1^{-\frac{1}{4}},u_1^{-\frac{1}{2}},
 u_1^{-\frac{3}{4}},u_2^{-\frac{1}{4}},u_2^{-\frac{1}{2}},
 u_2^{-\frac{3}{4}})$. The canonical element $\Som[8]{s=1}{e_s\dual{e_s}}$ is
 equal to   $\idl[3]$ and   the Casimir ideal $\zc(\K)$ coincides with $\kk_3$,  showing that $\zc(\K)\subsetneqq \Zc(\K)=\K$. Let us  pose some useful notations: for $M\in\set{\msf{F},\msf{F}\tenss[\msf{E}]\msf{F}}$ and  $i,j\in\set{1,2,3}$, let  $\indds{i}{M}{j}$  have the natural $\kk_i\sdash\kk_j$-bimodule structure on $M$; when viewed as element in  $\indds{i}{M}{j}$ each $x\in M$ may still be written as $x$ or be subscripted as $\elts{x}{i}{j}$ if more precision is needed. In particular put
$\elts{z}{i}{j}:=\z[{\indds{i}{\msf{F}\tenss[\msf{E}]\msf{F}}{j}}]=
 1\otimes 1 +   u^{\frac{1}{4}}\otimes u^{{-\frac{1}{4}}} \in {\indds{i}{\msf{F}\tenss[\msf{E}]\msf{F}}{j}}$ with $(i,j)=(1,2),(2,1)$. 
\\[-1.0\baselineskip]
Consider the modulated quiver  
\mbox{$\mQ=(B,\K,\tr) :$} $ \diagram[/bline=-1.7em]{ {} \&  \kk_3 \&  {} \\[-0.5em]
 \kk_1 \&  {} \&  \kk_2 }{
 \path[<-] ($(m-2-1.east) + (0.em,+0.25em)$) edge node[above] {$\indds{2}{\msf{F}}{1}$}
($(m-2-3.west)+(0.em,+0.25em)$);
 \path[->] ($(m-2-1.east) +
(0.em,-0.25em)$) edge node[below] {$\indds{1}{\msf{F}\tenss[\msf{E}] \msf{F}}{2}$}($(m-2-3.west)+(0.em,-0.25em)$);
 \path[->] (m-2-3) edge node[above,sloped]
 {$\indd{2}{\msf{F}}{3}$} (m-1-2)
 (m-1-2) edge node[above,sloped]
 {$\indd{3}{\msf{F}}{1}$} (m-2-1);
 }$ with   \mbox{$ B=\indds{1}{\msf{F}\tenss[\msf{E}] \msf{F}}{2}\oplus \indds{2}{\msf{F}}{1} \oplus \indds{2}{\msf{F}}{3} \oplus
\indds{3}{\msf{F}}{1}$}  and 
\mbox{$\dual{B}=\indds{2}{\msf{F}\tenss[\msf{E}] \msf{F}}{1}\oplus \indds{1}{\msf{F}}{2} \oplus \indds{3}{\msf{F}}{2} \oplus \indds{1}{\msf{F}}{3}$}.
The bilinear form in the symmetrizable dualizing pair
 $\set{B,\dual{B};\bilf}$  is   induced by   $\tr$,
 the multiplication   of  $\msf{F}$ and   the   projection $\morph{\msf{p_E}: \msf{F}}{\msf{E}}$ taking each $ a + b u^{\frac{1}{4}}\in \msf{E}$ to $a \in \msf{F}$. Thus for   \mbox{$x,x',y,y' \in \msf{F}$} we have $\bilf\bigl((x\otimes y)\otimes (x'\otimes y')\bigr)=  \bilf\bigl(x\mul\bilf(y\otimes x') \otimes
 y'\bigr)=x\mul\msf{p_E}(yx')y'$;   for $i\!=\!1,2$,
 $x\in \indd{i}{\msf{F}}{3}$ and \mbox{$x'\in
 \indd{3}{\msf{F}}{i}$} we have $\bilf(x\otimes x')=xx'$
 and $\bilf(x'\otimes x)=\tr(x'x)$.  Note that  $\bilf(\elts{z}{1}{2}\otimes \elts{z}{2}{1})=\idl[1]+\idl[1]=0$. The Casimir  operator  $\morph{\zc: \kQc}{\Zck(\kQc),x\mapsto \Som[9]{s=1}{e_s\mul x\mul \dual{e_s}}}$  vanishes on $(\indds{1}{\msf{F}\tenss[\msf{E}] \msf{F}}{2})\otimes\indds{2}{\msf{F}}{1} \oplus \indds{2}{\msf{F}}{1}\otimes( \indds{1}{\msf{F}\tenss[\msf{E}] \msf{F}}{2})$.  For  Casimir element $\z[\dual{B}\otimes B]=\Som{x\in\lmQ}{\dual{x}\otimes x}$ and  $\z[B\otimes \dual{B}]=\Som{y\in\rmQ}{y\otimes \dual{y}}$, take
  \begin{equation} \label{eq.explprjb1-2}
  \begin{aligned} 
  &   \lmQ=\set{\elts{z}{1}{2}, (\elts{1}{1}{}\otimes 1_{2}),\elts{1}{2}{1}, \elts{1}{2}{3},\elts{1}{3}{1},  \inddsu{3}{u}{1}{\frac{1}{4}}, \inddsu{3}{u}{1}{\frac{1}{2}}, \inddsu{3}{u}{1}{\frac{3}{4}}},   \ldmQ=\set{\elts{z}{2}{1}+(\elts{1\otimes 1}{2}{1}),\elts{z}{2}{1},\elts{1}{1}{2},
   \elts{1}{3}{2},\elts{1}{1}{3},\inddsu{1}{u}{3}{-\frac{1}{4}},
   \inddsu{1}{u}{3}{-\frac{1}{2}}, \inddsu{1}{u}{3}{-\frac{3}{4}}} \\
 & \rmQ=\set{\elts{z}{1}{2}, (\elts{1}{1}{}\otimes 1_{2}),\elts{1}{2}{1},\elts{1}{2}{3},  \inddsu{2}{u}{3}{\frac{1}{4}}, \inddsu{2}{u}{3}{\frac{1}{2}}, \inddsu{2}{u}{3}{\frac{3}{4}}, \elts{1}{3}{1} },
  \rdmQ=\set{\elts{z}{2}{1}+(\elts{1\otimes 1}{2}{1}),\elts{z}{2}{1},\elts{1}{1}{2},
    \elts{1}{3}{2},\inddsu{3}{u}{2}{-\frac{1}{4}},
    \inddsu{3}{u}{2}{-\frac{1}{2}}, \inddsu{3}{u}{2}{-\frac{3}{4}},\elts{1}{1}{3}}. 
 \end{aligned}
 \end{equation}

\ParIt{Example \ref{secexample}.1: A nonsymmetric potential of degree $2$}  $W:=\elts{z}{1}{2}\otimes \elts{1}{2}{1}$ is  a  nonsymmetric potential on  $\mQ$,   the   subbimodule $U=\msf{F}_1\mul\set{\elts{z}{1}{2}}\mul\msf{F}_2$  is one-dimensional on both sides and   is not a direct summand in  $\indds{1}{\msf{F}\tenss[\msf{E}] \msf{F}}{2}$.  We compute: $\cderv[\dual{\elts{z}{1}{2}}]W=\elts{1}{2}{1}$ and $\cderv[{\dual{(\elts{1}{2}{1})}}]W=\elts{z}{1}{2}$. Thus for the modulated quiver with potential  $(\mQ,W)$ we have  $ \triv{B}=U \oplus \indd{2}{\msf{F}}{1}$ and     $\bbar{B}:=\red{B}=\bbar{U} \oplus \indd{2}{\msf{F}}{3} \oplus \indd{3}{\msf{F}}{1}$ where
 $\bbar{U}=\Frac{(\indds{1}{\msf{F}\tenss[\msf{E}] \msf{F}}{2})}{U}$. We let $\bbar{\mQ}=(\bbar{B},\K,\tr)$ be the corresponding reduced modulated quiver.  We get that $\J{W}$ coincides with the closed ideal $L=\bbar{\ggen{\triv{B}}}$ and the natural projection $\morph{\rho: \kQc}{\kQbc}$ is a non-split reduction from $(\mQ,W)$ to $(\bbar{\mQ},0)$.

 \ParIt{Example~\ref{secexample}.2:  A family ($\mfr{m}_n)_{n\in\N}$ of  nonsymmetric potentials} Keep the  notations from Example~\ref{secexample}.1. For each $n\in \N$,  we  consider the  nonsymmetric potential below, with the convention that $W^0$  is the unity of $\K$, \begin{center}
$\mfr{m}:=W  + S\mul W^n$
where $ S:=(\elts{1}{1}{}\otimes 1_{2}) \otimes (\alpha u_1^{\frac{1}{4}}) +
   (\elts{u^{\frac{1}{4}}\otimes 1}{1}{2}) \otimes \alpha \in B^{(3)} $ and
   $\alpha:= (\elts{1}{2}{3}\otimes \elts{1}{3}{1} + \inddsu{2}{u}{3}{{\frac{1}{2}}}\otimes \inddsu{3}{u}{1}{{-\frac{1}{2}}}) \in \Zc[\msf{\E}](\indd{2}{\msf{F}}{3}\otimes \indd{3}{\msf{F}}{1})$.
 \end{center}
Thus $S$  is a symmetric potential, we get
$S=\elts{z}{1}{2}\otimes u_2^{\frac{1}{4}}\mul\alpha +
(\elts{1\otimes 1}{1}{2})\otimes(\alpha\mul u_1^{\frac{1}{4}} +u_2^{\frac{1}{4}}\mul \alpha)$  since $(\elts{u^{\frac{1}{4}}\otimes 1}{1}{2})=\elts{z}{1}{2}\mul u_2^{\frac{1}{4}} + (\elts{1\otimes 1}{1}{2})\mul u_2^{\frac{1}{4}} $.  The degree-$2$ component of $\mfr{m}$ is   $\mfr{m}_2=W=\elts{z}{1}{2}\otimes \elts{1}{2}{1} \in \indds{1}{\msf{F}\tenss[\msf{E}] \msf{F}}{2}\otimes \indd{2}{\msf{F}}{1}$ and, as in Example~\ref{secexample}.1 above,  $(\mQ,\mfr{m}_n)$ is a modulated quiver with potential  such that  $\triv{B}=U \oplus \indd{2}{\msf{F}}{1}$ and  $\bbar{B}=\bbar{U} \oplus \indd{2}{\msf{F}}{3} \oplus \indd{3}{\msf{F}}{1}$.

We  will  compute the  cyclic derivatives $\cderv[\dual{\elts{z}{1}{2}}]\mfr{m}_n$ and
  $\cderv[\dual{\elts{1}{2}{1}}]\mfr{m}_n$.  To compute  $\cderv[\dual{\elts{z}{1}{2}}](S\mul W^n)=\lderv[\dual{\elts{z}{1}{2}}](\cperm (S\mul W^n))$ we  need  the following permutations: $S\mul W^n,\rperm[2](S\mul W^n),\dots,\rperm[2n](S\mul W^n)$. We have $\cderv[\dual{\elts{z}{1}{2}}](S)=u_2^{\frac{1}{4}}\mul\alpha$ and by definition,\\
  $\myeqar{\rperm[2](S\mul W)= & \elts{z}{1}{2}\otimes \elts{1}{2}{1} \otimes S\bilf^2(W\otimes(\dual{\elts{1}{2}{1}}\otimes \dual{\elts{z}{1}{2}})) +
  (\elts{1\otimes 1}{1}{2})\otimes \elts{1}{2}{1} \otimes S\bilf^2(W\otimes(\dual{\elts{1}{2}{1}}\otimes \dual{(\elts{1\otimes 1}{1}{2})})) \\ =& \elts{z}{1}{2}\otimes \elts{1}{2}{1} \otimes S +0=W\mul S.
  }$\\
 Thus  $\cderv[\dual{\elts{z}{1}{2}}](S\mul W^n)=u_2^{\frac{1}{4}}\mul\alpha\mul W^n + \elts{1}{2}{1}\mul (S\mul W^{n-1} +W \mul S\mul W^{n-2}+\dots+W^{n-1}S )$ for all $n\geq 1$.  Apply a similar argument to compute $\cderv[{\dual{\elts{1}{2}{1}}}](S\mul W^n)$. We deduce that
 \\[-1.5\baselineskip]
  \begin{align}\label{eq.explcderv1}
   & \cderv[\dual{\elts{z}{1}{2}}]\mfr{m}_n=\elts{1}{2}{1} +
  u_2^{\frac{1}{4}}\mul\alpha\mul W^n + \elts{1}{2}{1}\mul \Som[n-1]{r=0}{W^r}\mul S\mul W^{n-1-r}   \\[-0.8\baselineskip]  \label{eq.explcderv2}
  &   \cderv[\dual{\elts{1}{2}{1}}]\mfr{m}_n=\elts{z}{1}{2}+\bigl(\Som[n-1]{r=0}{W^r\mul S\mul W^{n-1-r}}\bigr)\mul \elts{z}{1}{2}=\bigl(1+\Som[n-1]{r=0}{W^r\mul S\mul W^{n-1-r}}\bigr)\mul \elts{z}{1}{2}.
  \end{align}
  Let $J_0$ be  the closed ideal in $\kQc$ generated by  $\cderv[\dual{\elts{z}{1}{2}}]\mfr{m}_n$ and $\cderv[\dual{\elts{1}{2}{1}}]\mfr{m}_n$. Let $\Seq{\rho :\kQc \, \to \,  \kQbc}$ be the natural projection. The right $\K$-linear map $\morph{\bbar{U}}{\indds{1}{\msf{F}\tenss[\msf{E}] \msf{F}}{2}}$, taking the coset $\bbar{1\otimes 1 }$ of $1\otimes 1$ to $1\otimes 1$, yields a right $\K$-linear map   $\morph{\rho': \bbar{B}}{B}$ such that $\rho\restr{B}\circ\rho'=\id_{\bbar{B}}$ , inducing a bimodule morphism $\morph{\bbar{\rho}'=\pi\circ \rho': \bbar{B}}{\Frac{(B+J_0)}{J_0}}$ where $\morph{\pi: \kQc}{\kQc/J_0}$ is the natural projection. Thus  $J_0=\Ker(\pi)$ satisfies condition $\msf{(1.ii)}$ of  Definition~\ref{defn.trivpart-reduction} for trivial parts  of Jacobian ideals (that is, kernels of reductions). Now we have the two following cases.
 \begin{itemize}
 \item[$\msf{(a)}$] The case $n\geq 1$. By \eqref{eq.explcderv1} above we get  that $\elts{z}{1}{2}\in J_0$ because the element $(1+\Som[n-1]{r=0}{W^r\mul S\mul W^{n-1-r}}\bigr)$ is invertible  in $\kQc$, and next, $\elts{1}{2}{1}\in J_0$, implying that $J_0$  coincides with the closed ideal $L=\bbar{\ggen{\triv{B}}}$. Therefore the natural projection $\morph{\rho: \kQc}{\kQbc}$ is a non-split reduction from $(\mQ,\mfr{m}_n)$  to   $(\bbar{\mQ},0)$.
 \item[$\msf{(b)}$] The case of potential $\mfr{m}_0=W+S$. Here,   $\cderv[\dual{\elts{1}{2}{1}}]\mfr{m}_0=\elts{z}{1}{2}\in J_0$ and  $\cderv[\dual{\elts{z}{1}{2}}]\mfr{m}_0=\elts{1}{2}{1} +
  u_2^{\frac{1}{4}}\mul\alpha\in J_0$. However,  $\elts{1}{2}{1}\in \indds{2}{\msf{F}}{1}$ is $\msf{F}$-central  while the $\msf{E}$-central element $  u_2^{\frac{1}{4}}\mul\alpha$ is  not $\msf{F}$-central, and we note that the $\K$-bimodule generated by the set $\set{\cderv[\dual{\elts{1}{2}{1}}]\mfr{m}_0,\,\cderv[\dual{\elts{z}{1}{2}}]\mfr{m}_0}$ contains an element in $\J[2]{\kQc}$, namely the element $\gamma=(\alpha u_1^{\frac{1}{4}} + u_2^{\frac{1}{4}}  \alpha)=(\elts{1}{2}{1} +
  u_2^{\frac{1}{4}}\mul\alpha) + u_2^{-\frac{1}{4}}(\elts{1}{2}{1} +
  u_2^{\frac{1}{4}}\mul\alpha)u_1^{\frac{1}{4}}$. Thus $J_0$ does not satisfy condition $\msf{(1.i)}$ in Definition~\ref{defn.trivpart-reduction}. Indeed, one can check that there is no $\K$-bimodule morphism  $\morph{f=\Psmatr{ \id_{\triv{B}} \\ f'}: \triv{B}}{ \triv{B} \oplus \J[2]{\kQc}}$  such that $f(\triv{B})\subset (\cderv \mfr{m}_0)(\dual{B})$, showing that there is no reduction on $(\mQ,\mfr{m}_0)$   as described by  Definition~\ref{defn.trivpart-reduction}. However, (as in the proof of part $\msf{(2)}$ of Theorem~\ref{theo.Jm-et-morph-Jtriv}) the right $\K$-linear map   $\morph{\rho': \bbar{B}}{B}$ extends to a right $\K$-linear morphism $\morph{\rho': \kQbc}{\kQc}$ such that   the map  $\morph{\bbar{\rho}'=\pi\circ \rho': \kQbc}{\Frac{\kQc}{J_0}}$  is  a  surjective morphism of topological algebras; one checks that $\Ker(\bbar{\rho}')$ is the closed ideal  $I_0$  generated by the element $\gamma$ above. We have an epimorphism of  path  algebras $\Seq{\phi: \kQc\, \to \, \Frac{\kQbc}{I_0}}$ defined on the bimodule $B=\indds{1}{\msf{F}\tenss[\msf{E}] \msf{F}}{2}\oplus \indds{2}{\msf{F}}{1} \oplus \indds{2}{\msf{F}}{3} \oplus \indds{3}{\msf{F}}{1}$ as follows:
  \begin{itemize}
   \item[$(\star)$] $\phi$ is induced by  $\rho$ over $ \indds{1}{\msf{F}\tenss[\msf{E}] \msf{F}}{2} \oplus \indds{2}{\msf{F}}{3} \oplus \indds{3}{\msf{F}}{1} $, and over $\indds{2}{\msf{F}}{1}$  we have: $\phi(\elts{1}{2}{1})=-u_2^{\frac{1}{4}}\mul\alpha + I_0=u_2^{\frac{1}{4}}\mul\alpha + I_0$.
   \end{itemize}
   Thus $\Ker(\phi) =J_0$  (and we note   that $\phi(\mfr{m}_0)=0$). Let  $\bbar{\mfr{m}_0}=\rho(\mfr{m}_0) \in \kQbc$,  then $\bbar{\mfr{m}_0}=\bbar{\indds{1}{1\otimes 1}{2}}\otimes (\alpha u_1^{\frac{1}{4}} + u_2^{\frac{1}{4}}  \alpha) \in \bbar{U} \otimes \indd{2}{\msf{F}}{3} \otimes \indd{3}{\msf{F}}{1}$.
  For  the reduced quiver with potential $(\bbar{\mQ},\bbar{\mfr{m}_0})$, the morphism of topological  algebras $\phi$ above yields an isomorphism of Jacobian algebras $\Jc(\mQ,\mfr{m}_0)\cong \Jc(\bbar{\mQ},\bbar{\mfr{m}_0})$.
    \end{itemize}

\ParIt{Skew reductions} In point $\msf{(b)}$ above, $\phi$ is an instance of what we may refer to as \emph{skew reduction}.

 With previous observations , it is not difficult to derive  the following  consequence of Theorem~\ref{theo.Jm-et-morph-Jtriv}.
 \begin{Cor} \label{Cor-red-nonsplit} Let $(\mQ,\mfr{m})$ be a modulated quiver with potential,  with $\triv{B}=U\oplus V$ and $\mfr{m}_{2}=\z[U\otimes V]$ and $\bbar{\mQ}$   the  corresponding reduced modulated quiver. Suppose  that $B=V\oplus B_1$  for some subbimodule $B_1$   containing   $U$   such that $(\mfr{m}-\mfr{m}_{2}) \in \what{\kQ^{(1)}}$  where $\mQ^{(1)}=(B_1,\K,\tr)$. Then there is  a reduction or a skew reduction from $(\mQ, \mfr{m})$ to a reduced modulated quiver with potential $(\bbar{\mQ},\rho(\mfr{m}))$.
  \end{Cor}


\section{Mutations of modulated quivers with potentials}\label{sec:mutationsmqp}
Let us mention that a motivation for lifting mutations of quivers with potentials to mutations of modulated quivers with potentials  comes from  a successful  non simply-laced generalization of  cluster structures  on $2$-Calabi-Yau
categories over arbitrary fields.
As before, we fix a modulated quiver with potential
$(\mQ,\mfr{m})$ over a symmetric algebra $(\K,\tr)$, with
$\mQ=(B,\K,\tr)$; where $B$ is part of a symmetrizable dualizing
pair  $\set{B,\dual{B}; \bilf}$. Also, 
 $(\rmQ[1],\rdmQ[1])$ and $(\lmQ[1],\ldmQ[1])$ stand respectively for a chosen right    projective basis and   left projective basis for $B$ over $\K$.

Note that we can write $\K=\prd{i\in\I}{\kk_i}$ as direct product of indecomposable $\kk$-algebras $\kk_i$, each $\kk_i$
 viewed as  subalgebra in $\K$ with $\idl[i]$.  The unity of $\K$
is $\idl=\som{i\in\I}{}{\idl[i]}$, the set $\set{\idl[i]\ :\;
i\in\I}$ is a system of central primitive orthogonal idempotents for $\K$, referred to as \emph{set of points of $\mQ$}. We  fix a point 
$e$   of $\mQ$ and  write $\bbar{e}=\idl -e$, such that:
\begin{equation}
\label{eq.idempoint}
\begin{aligned} & \text{The idempotent } e \text{ is
loop-free and 2-loop free,  that is, }
  e\mul B \mul e = 0 \text{ and }  (B \mul e ) \cap (  e\mul B
)=0. \\
 & \text{Replacing if necessary $\mfr{m}$ by
a cyclically equivalent potential, we have: }  e \mfr{m}=0=\mfr{m}e.
\end{aligned}
\end{equation}
In view of the first part of \eqref{eq.idempoint}, we derive  the following relations:
\begin{equation}\label{eq.idempoint+decomp}
\begin{aligned} 
& B= B \mul e \oplus e \mul B  \oplus \bbar{e} \mul B \mul  \bbar{e}
 \dual{B}= e\mul \dual{B} \oplus
\dual{B}\mul e \oplus \bbar{e} \mul \dual{B} \mul  \bbar{e}, \text{ and } \\
& \set{B,\dual{B}}=\set{B \mul e ,
e\mul \dual{B}}\oplus \set{e \mul B, \dual{B}\mul e} \oplus \set{\bbar{e} \mul B \mul  \bbar{e},\bbar{e} \mul \dual{B}\mul
\bbar{e}} \text{ as naturally induced dualizing pairs}. 
\end{aligned} 
\end{equation}

\begin{Defn}\label{defn.mutmqpot1} Whenever
\eqref{eq.idempoint} is satisfied, the \emph{semi-mutation} of
$(\mQ,\mfr{m})$ at  point $e$  is the modulated quiver with potential
$\wtilde{\mu}_{e}(\mQ,\mfr{m})=(\wtilde{\mu}_{e}(\mQ),\wtilde{\mfr{m}})$  described as follow:
\begin{itemize}
\item [$\msf{(i)}$] $\wtilde{\mu}_{e}(\mQ)=(\widetilde{B},\K,\tr)$
with $\widetilde{B}=[BeB] \oplus e\mul \dual{B} \oplus
\dual{B}\mul e \oplus \bbar{e} \mul B \mul  \bbar{e}$, where  $[BeB]$ is still the bimodule $BeB$  regarded as being part of the arrow bimodule of $\wtilde{\mu}_{e}(\mQ)$, and letting  $\dual{\widetilde{B}}=[e\dual{B}\mul \dual{B}e] \oplus B \mul e
\oplus e\mul B  \oplus \bbar{e} \mul \dual{B} \mul  \bbar{e}$, the data  $\set{\widetilde{B},\dual{\widetilde{B}}}$ is a symmetrizable dualizing pair canonically induced by the pair $\set{B,\dual{B};\bilf}$.  Each tensor
element in $[BeB]$ may be written as $[xy]$ or $[x\otimes y]$ with
$x\in Be$ and $y\in e B$.
 \item [$\msf{(ii)}$]  $\wtilde{\mfr{m}}=[\mfr{m}]+\z[e]$  with $\z[e]=\z[{[BeB]\otimes\dual{(BeB)}}]
= \Som{y\in\rmQ[1]} \Som{z\in\rmQ[1]} [yez]\otimes
\dual{z}e\otimes e\dual{y}$, where $[\mfr{m}]$ coincides with $\mfr{m}$ but is
regarded as an element in the complete path algebra of
$\wtilde{\mu}_{e} (\mQ)$.
\end{itemize}
\end{Defn}

Observe (using \eqref{eq.idempoint}) that $\wtilde{\mfr{m}}$ is necessarily  $2$-loop free, so that $(\wtilde{\mu}_{e}(\mQ),\wtilde{\mfr{m}})$ is indeed a modulated quiver with potential. Using  part $\msf{(2)}$ of
 Lemma~\ref{lem.cyc-stable} we  obtain:
 \begin{Rem}\label{rem.ze-cycequiv} $
\rperm(\z[e])\hspace{-0.2em}=\hspace{-0.5em}\Som{x\in\lmQ[1]}\,\Som{z\in\rmQ[1]}{
\hspace{-0.2em} e\dual{x}\otimes [x \mul e\mul z]\otimes \dual{z}e}= [\z[e\dual{B}\otimes Be] \otimes \z[e B\otimes
\dual{B}e] ]$, thus  $\widetilde{\mfr{m}} \cyc{\equiv}  [ \mfr{m} + \z[e\dual{B}\otimes Be] \mul \z[e B\otimes
\dual{B}e] ]$.
\end{Rem}
\begin{Rem} Let  $(\mQ', \mfr{m}')$ be another   modulated quiver with
potential  with $(\mQ'=(B',\K,\tr))$ such that $e\mul B'=0=B'\mul
e$. Then $\wtilde{\mu}_{e}(\mQ\oplus
\mQ',\mfr{m}+\mfr{m}')=\wtilde{\mu}_{e}(\mQ,\mfr{m})\oplus (\mQ',
\mfr{m}')$.
\end{Rem}

\begin{Thm}\label{theo.mtmq-ceqf} For each  modulated quiver with potential $(\mQ,\mfr{m})$ satisfying condition
\eqref{eq.idempoint}, the  right-equivalence
class  of the semi-mutation $\wtilde{\mu}_{e}(\mQ,\mfr{m})=(\wtilde{\mu}_{e}(\mQ),\wtilde{\mfr{m}})$
is determined by  that of $(\mQ,\mfr{m})$.
\end{Thm}

\begin{prv}  Let $\shdual{\mQ}=(\shdual{B},\K,\tr)$
with $\shdual{B}=B \oplus e\mul \dual{B} \oplus \dual{B} \mul e$.
Clearly,  $\shdual{B}$ is  part of a naturally induced
symmetrizable dualizing pair $\set{\shdual{B},\shdual{{\dual{B}}}}$  with $\shdual{{\dual{B}}}=\dual{B}
\oplus B\mul e \oplus e\mul B$. Then, the natural embedding
$\morph[/dir=right hook->]{B}{\shdual{B}}$ identifies $\kQc$ with
a closed subalgebra in $\kQshc$. We also have  a natural embedding
$\morph[/dir=right hook->]{\widetilde{B}}{\kQshc}$ sending each
degree-$1$ element $[xy]$ of $\widetilde{B}$ to the tensor element
$xy$ in $\kQshc$, allowing us to  identify $\kQtc$ with a closed
subalgebra in $\kQshc$. Under this identification, with
$\widetilde{\mfr{m}}$ now viewed as as element in $\kQshc$,
Remark~\ref{rem.ze-cycequiv}  states that $\widetilde{\mfr{m}}$ is
cyclically equivalent to the potential $\mfr{m} + \bigl(\z[e\dual{B}\otimes Be] \mul \z[e B\otimes \dual{B}e] \bigr) $. Taking
the latter into account, Theorem~\ref{theo.mtmq-ceqf}
becomes a direct consequence of the following lemma.
\end{prv}

\begin{Lem} Every automorphism $\phi$ of $\kQc$ extends to an
automorphism $\shdual{\phi}$ of $\kQshc$ such that: for all
$\xi\in  \dual{B}$ we have
$\shdual{\phi}(e\xi)=\Som{x\in \lmQ[1]} e\dual{
x}\phi\bigl(\rderv[\xi]\phi^{-1}(x\mul e) \bigr)$ and
$\shdual{\phi}(\xi e)=\Som{y\in
\rmQ[1]}\phi\bigl(\lderv[\xi]\phi^{-1}(e y) \bigr) \dual{y}e$.
Consequently,
\begin{equation}
\label{eq.semi-ext1}  \shdual{\phi}(\z[(e\dual{B})\otimes (B\mul
e)])= \z[(e\dual{B})\otimes (B\mul e)], \; \quad
\shdual{\phi}(\z[(e\mul B)\otimes (\dual{B}\mul e)])= \z[(e\mul
B)\otimes (\dual{B}\mul e)], \text{ and }
 \shdual{\phi}(\kQtc) =\kQtc.
\end{equation}
\end{Lem}

\begin{prv}  We first check that $\shdual{\phi}$ is indeed a morphism of
$\K$-bimodules on $e \dual{B}$ and $\dual{B}\mul e$. Thus let
$\xi\in \dual{B}$ and $a,b\in \K$; we have:
$\shdual{\phi}(a e\xi b)=\Som{x\in \lmQ[1]} e\dual{
x}\phi\bigl(\rderv(\phi^{-1}(x\mul e) \otimes a \xi b
\bigr)=\Som{x\in \lmQ[1]} e\dual{ x}\phi\bigl(\rderv(\phi^{-1}(x a
\mul e) \otimes \xi \bigr)b$.  Using
identities~\eqref{eq.dualbases} in the sequel we write
$xa=\Som{z\in \lmQ[1]}\bilf(xa\otimes
\dual{z})\mul z$ , thus \\
 $\begin{array}{rl}\shdual{\phi}(a e \xi b)= & \Som{z\in
\lmQ[1]}e\bigl(\Som{x\in \lmQ[1]} \dual{x } \bilf(x\otimes
a\dual{z})\bigr) \phi\bigl(\rderv(\phi^{-1}( z \mul e) \otimes \xi
) \bigr)b =\Som{z\in \lmQ[1]}e a\dual{z}
\phi\bigl(\rderv(\phi^{-1}( z \mul e) \otimes \xi ) \bigr)b \\
 = & a\Som{z\in \lmQ[1]}e \dual{z} \phi\bigl(\rderv(\phi^{-1}( z \mul
e) \otimes \xi ) \bigr)b = a\shdual{\phi}(e\xi)b.
 \end{array}$ \\
Similarly, it is easily checked that $\shdual{\phi}(a\xi e
b)= a\shdual{\phi}(\xi e)b$.  Let us show that the degree-$1$ component
$\morph{\shdual{\phi_1}: B \oplus e\mul\dual{B} \oplus
\dual{B}\mul e }{B \oplus e\mul\dual{B} \oplus \dual{B}\mul e}$ of
$\shdual{\phi}$ is an automorphism. 
$\morph{\shdual{\phi_1}\restr{B}=\phi_1 : B }{B} $   and  by 
part $\msf{(a)}$ of Proposition~\ref{prop.prolong}, $\phi_1$ is
an automorphism of $B$. For all $\xi\in \dual{B}$, we have: \\
 $\begin{array}{rl}\shdual{\phi_1}(e\xi)= & \Som{x\in \lmQ[1]} e\dual{
x}\phi\bigl(\rderv[\xi]\phi_1^{-1}(x\mul e) \bigr) = \Som{x\in
\lmQ[1]} e\dual{ x}\phi\bigl(\bilf(\phi_1^{-1}(x\mul e) \otimes
\xi) \bigr)= \Som{x\in \lmQ[1]} e\dual{ x}\phi\bigl(\bilf(x\mul e
\otimes {\dual{\phi_1}}^{-1}(\xi)) \bigr) \\
=& \Som{x\in \lmQ[1]} e\dual{ x} \bilf(x \otimes
{\dual{\phi_1}}^{-1}(e\xi))={\dual{\phi_1}}^{-1}(e\xi).
\end{array}$ \\
 Thus,   $\morph[/above={\sim}]{ \shdual{\phi_1} \restr{e\mul\dual{B}} = {\dual{\phi_1}}^{-1} \restr{e\mul\dual{B}} :
e\mul\dual{B}}{e\mul\dual{B}}$.  Similarly, $\morph[/above={\sim}]{ \shdual{\phi_1} \restr{\dual{B}\mul e} = {\dual{\phi_1}}^{-1} \restr{ \dual{B}\mul e} :
 \dual{B}\mul e }{ \dual{B} \mul e}$. Therefore,
 $\shdual{\phi_1}$ is an automorphism and  part $\msf{(a)}$ of Proposition~\ref{prop.prolong} yields that
 $\shdual{\phi}$ is an automorphism of $\kQshc$
extending $\phi$.

We now prove the first identity in \eqref{eq.semi-ext1}.  In view of identities~\eqref{eq.dualbasesKQc}, for each
$u\in B$ we know that $\phi^{-1}(u) = \Som{x\in \lmQ[1]}
\bigl(\rderv[\dual{x}]\phi^{-1}(u) \bigr) \otimes x $, so that $u \phi(\phi^{-1}(u))= \Som{x\in\lmQ[1]} \phi \bigl(\rderv[\dual{x}]\phi^{-1}(u) \bigr) \otimes
\phi(x) $. Then,  we compute: \\
$\begin{array}{rl} \z[(e\dual{B})\otimes (B\mul e)] = & \Som{z
\in\lmQ[1]}e\dual{z}\otimes  ze  = \Som{z\in \lmQ[1]} \Som{x\in
\lmQ[1]}  e\dual{z} \phi\bigl(\rderv( \phi^{-1}(ze)e\otimes
\dual{x})\bigr) \otimes \phi( x e)     = \Som{x\in \lmQ[1]}
\shdual{\phi}(e\dual{x}) \phi(x e)   =  \shdual{\phi}(\z[(e\dual{B})\otimes (B\mul
 e)]).
\end{array}$\\
The second identity of \eqref{eq.semi-ext1} is established in the
same way and  the last identity follows by the definition of $\shdual{\phi}$.
\end{prv}

The reduction Theorem~\ref{theo.red-mQp} together with Theorem~\ref{theo.mtmq-ceqf} above yield the following result. 
\begin{Cor}\label{cor.mutclequivf}
 Suppose \eqref{eq.idempoint} holds and the trivial part $\triv{(\wtilde{\mu}_{e}(\mQ), \wtilde{\mfr{m}} )}$ splits. Then the weak right-equivalence class of
$\Red(\wtilde{\mu}_{e}(\mQ,\mfr{m}))$ is determined by  that of
$(\mQ,\mfr{m})$. \finprv
\end{Cor}

\begin{Defn} With the assumptions of Corollary~\ref{cor.mutclequivf}, the
\emph{mutation} of $(\mQ,\mfr{m})$ at   point $e$ is
the reduced modulated quiver with potential
$\Red(\wtilde{\mu}_{e}(\mQ),\wtilde{\mfr{m}})$, unique up to weak
right-equivalence:  we write
$\mu_e(\mQ,\mfr{m})=\Red(\wtilde{\mu}_{e}(\mQ),\wtilde{\mfr{m}})$.
\end{Defn}

Another important  result of this work establishes that, every mutation is an
involution.
\begin{Thm}\label{theo.mut-involutif} The mutation $\mu_e$   is an  involution over the set of weak right-equivalence classes of the modulated quivers with
potentials $(\mQ,\mfr{m})$ satisfying  \eqref{eq.idempoint} and whose semi-mutations have a split trivial part. If moreover $\mfr{m}$ is
a symmetric  and  $\morph{\cderv\mfr{m}: \dual{B}}{\cderv(\dual{B}\otimes
\mfr{m})}$ also splits, then $\mu_e$ is an involution up to
right-equivalences.
\end{Thm}

\begin{prv}  Suppose that $(\mQ,\mfr{m})$ is a reduced modulated quiver with potential  satisfying the assumptions of the  theorem.
 Then write: $\widetilde{\mu}_{e}^{2}(\mQ,\mfr{m})= \widetilde{\mu}_{e}(\widetilde{\mQ},\widetilde{\mfr{m}})
 =(\widetilde{\widetilde{\mQ}},\widetilde{\widetilde{\mfr{m}}})$. In
 view of  reduction Theorem~\ref{theo.red-mQp}  and Theorem~\ref{theo.sympot} (for the symmetric case), it suffices to show that  $(\widetilde{\widetilde{\mQ}},\widetilde{\widetilde{\mfr{m}}})$ is
 right-equivalent to $(\mQ,\mfr{m}) \oplus (\triv{\widetilde{\widetilde{\mQ}}},W)$ where
 $W$ is cyclically equivalent to the degree-$2$ component  $\widetilde{\widetilde{\mfr{m}}}_{2}$ of $\widetilde{\widetilde{\mfr{m}}}$.
 By definition,  $\widetilde{B}=[BeB] \oplus e\mul \dual{B}
\oplus \dual{B}\mul e \oplus \bbar{e} \mul B \bbar{e}$.  We have \eqref{eq.mut2B} and \eqref{mu2-pot1}  below where \eqref{eq.mut2B}  uses the assumption that  \eqref{eq.idempoint} holds. 
\begin{equation}\label{eq.mut2B}
\widetilde{\widetilde{B}}=[e\dual{B} \dual{B}e] \oplus Be \oplus e
B \oplus [BeB] \oplus \bbar{e}B\bbar{e}= B\oplus ([e\dual{B}
\dual{B}e] \oplus [BeB]).
\end{equation}
\\[-1.0\baselineskip]
\begin{equation}\label{mu2-pot1}
\begin{aligned}
\widetilde{\widetilde{\mfr{m}}} &= [[\mfr{m}]]
+[\z[{[BeB]\otimes\dual{(BeB)}}]] + \z[{[e\mul\dual{B}\dual{B}\mul
e]\otimes (e\mul B B \mul e) }]  \\
&= [\mfr{m}] + \hspace{-0.4em}\Som{y\in\rmQ[1]}\Som{z\in\rmQ[1]}
\hspace{-0.4em}[yez]\mul [\dual{z}e\otimes e\dual{y}] + \Som{x\in
\lmQ[1]} \Som{x'\in \lmQ[1]} [e \dual{x'} \mul \dual{x}e]\mul e x
\mul x' e.
\end{aligned}
\end{equation}
But   using  part $\msf{(2)}$ of
 Lemma~\ref{lem.cyc-stable} we know  that: 
\mbox{$\lperm(\z[{[e\mul\dual{B}\dual{B}\mul e]\otimes }]) =
\z[{(e\mul B B \mul e) \otimes [e\mul\dual{B}\dual{B}\mul e] }] =
\hspace{-0.4em}\Som{y\in\rmQ[1]}\Som{z\in\rmQ[1]} \hspace{-0.4em}
ye\mul  ez [\dual{z}e \otimes e\dual{y}]$.} Hence,
\begin{equation}\label{mu2-pot2}
\widetilde{\widetilde{\mfr{m}}} \cyc{\equiv}\; \mfr{m}' :=
[\mfr{m}] + \hspace{-0.4em}\Som{y\in\rmQ[1]}\Som{z\in\rmQ[1]}
\hspace{-0.4em}([ye\otimes ez] + ye \otimes e z )\mul
[\dual{z}e\otimes e\dual{y}].
\end{equation}

We then consider the trivial modulated quiver with potential
$(\triv{\widetilde{\widetilde{\mQ}}},W)$ with
$\triv{\widetilde{\widetilde{\mQ}}}=[BeB] \oplus \dual{[BeB]}=[Be
\otimes eB] \oplus [\dual{B}e\otimes e\dual{B}]$ and
 $W=\z[{[BeB] \otimes \dual{[BeB]}}]=\hspace{-0.4em}\Som{y\in\rmQ[1]}\Som{z\in\rmQ[1]}
\hspace{-0.4em}[ye\otimes ez]\mul [\dual{z}e\otimes e\dual{y}]$;
(note that $W$ is of course cyclically equivalent to
$\widetilde{\widetilde{\mfr{m}}}_{2})$.  Now, to prove
Theorem~\ref{theo.mut-involutif} it suffices to show that the
modulated quiver with potential $(\widetilde{\widetilde{\mQ}},
\mfr{m}' )$ is right-equivalent to $(\mQ,\mfr{m}) \oplus
(\triv{\widetilde{\widetilde{\mQ}}},W)$, here  $\mfr{m}'$ is
given by \eqref{mu2-pot2} above.
The term $S:=\z[{(Be B)\otimes
[\dual{B}e\dual{B}]}]=\hspace{-0.4em}\Som{y\in\rmQ[1]}\Som{z\in\rmQ[1]}
\hspace{-0.4em} ye \otimes e z \mul [\dual{z}e\otimes e\dual{y}]$
of $\mfr{m}'$ is a potential and the right derivative morphism
$\morph[/above={\rderv S}]{[BeB]}{Be \otimes eB}$ is a bimodule 
isomorphism taking each $[ye\mul ez]$ to $ye
\otimes ez$ for all $y,z\in\rmQ[1]$. Whence the following
unitriangular automorphism $\morph{\vphi:
\widehat{{\text{$\kk$}}\widetilde{\widetilde{\mQ}}} }{
\widehat{{\text{$\kk$}}\widetilde{\widetilde{\mQ}}} }$ whose
restriction on the bimodule $\widetilde{\widetilde{B}}= B\oplus
([e\dual{B} \dual{B}e] \oplus [BeB])$ is given by:
\begin{equation*}
\begin{aligned}& \vphi\restr{B\oplus  [\dual{B}e\mul e\dual{B}]}=\id_{B\oplus
[\dual{B}e\mul e\dual{B}]} \text{ and } \morph{\vphi\restr{[BeB]}
= \id_{[BeB]} - \rderv S : [BeB] }{[BeB] \oplus Be\otimes eB}, \\[0.2em]
& \text{ thus } \vphi([ye\mul ez])=[ye\mul ez] -ye \otimes ez
\text{ for all } y,z\in\rmQ[1].
\end{aligned}
\end{equation*}
Since  $[\mfr{m}]$ is obtained $\mfr{m}$ by substituting $[xex']$ for each tensor element
\mbox{$xe\otimes e x' \in Be\otimes e B$} occurring  the expansion of $\mfr{m}$, we deduce the following.
\begin{center}
$\vphi(\mfr{m}')=\mfr{m} +
\hspace{-0.4em}\Som{y\in\rmQ[1]}\Som{z\in\rmQ[1]}
\hspace{-0.4em}[ye\otimes ez]\mul [\dual{z}e\otimes e\dual{y}] +
S'$,
\end{center} where $S'$ is a potential lying  in the closed ideal in
$ \widehat{{\text{$\kk$}}\widetilde{\widetilde{\mQ}}}$ generated
by  $[BeB]$ which is a direct summand  in  $\widetilde{\widetilde{B}}$. Then by
Lemma~\ref{lem.cyclequiv}, $S'$ is cyclically equivalent to a
potential $S''\in [BeB]\mul \widehat{\widetilde{\widetilde{J}}}^2$
where  we write $\widehat{\widetilde{\widetilde{J}}}$ for the
closed arrow ideal in $
\widehat{{\text{$\kk$}}\widetilde{\widetilde{\mQ}}}$.  Now let $\morph{f=\lderv S'':
 [\dual{B}e\mul e\dual{B}]}{\widehat{\widetilde{\widetilde{J}}}^2}$, applying
\eqref{eq.dervdualbasesKQc} we know that
$S''=\hspace{-0.4em}\Som{y\in\rmQ[1]}\Som{z\in\rmQ[1]}
\hspace{-0.4em}[ye\otimes ez]\mul f([\dual{z}e\otimes
e\dual{y}])$, hence
\begin{center}
$\vphi(\mfr{m}') \cyc{\equiv} \; \mfr{m}'' := \mfr{m} +
\hspace{-0.4em}\Som{y\in\rmQ[1]}\Som{z\in\rmQ[1]}
\hspace{-0.4em}[ye\otimes ez]\mul \bigl([\dual{z}e\otimes
e\dual{y}]+ f([\dual{z}e\otimes e\dual{y}])\bigr)$.
\end{center}

Next, we  deduce  the unitriangular automorphism $\vphi'$ of
$\widehat{{\text{$\kk$}}\widetilde{\widetilde{\mQ}}} $ with
\begin{equation*}
\begin{aligned}& \vphi'\restr{B\oplus [BeB]}=\id_{B\oplus [BeB]}
\text{ and } \morph{\vphi'\restr{ [\dual{B}e\mul e\dual{B}]} =
\id_{ [\dual{B}e\mul e\dual{B}]} - f \,:\, [\dual{B}e\mul
e\dual{B}] }{ [\dual{B}e\mul e\dual{B}] \oplus
\widehat{\widetilde{\widetilde{J}}}^2}, \\[0.2em]
& \text{thus } \vphi'([\dual{z}e\otimes
e\dual{y}])=[\dual{z}e\otimes e\dual{y}] -f([\dual{z}e\otimes
e\dual{y}]) \text{ for all } y,z\in\rmQ[1].
\end{aligned}
\end{equation*}
We get that $\vphi'(\mfr{m}')= \mfr{m}+W$. Hence, letting
$\phi=\vphi'\circ\vphi$, we obtain a right-equivalence from
$(\widetilde{\widetilde{\mQ}},\widetilde{\widetilde{\mfr{m}}})$ to
 $(\mQ,\mfr{m}) \oplus (\triv{\widetilde{\widetilde{\mQ}}},W)$, completing the proof of
 Theorem~\ref{theo.mut-involutif}.
\end{prv}


\section{Examples of mutations in the mutation class of Dynkin type $\msf{F}_4$} \label{sec:explmutationreduction}
 Consider the $\mb{R}$-algebra  $\K=\kk_1\times \kk_2\times \kk_3\times \kk_4$, with $\kk_1=\mb{R}_1=\mb{R},\kk_2=\mb{R}_2=\mb{R},\kk_3=\mb{C}_3=\mb{C}$ and  $  \kk_4=\mb{C}_4=\mb{C}$, regarded  as $\mb{R}$-subalgebras of $\K$ with respective unities   $1_1,1_2,1_3,1_4$. An  $\mb{R}$-basis of $\K$  set $\mrm{S}:=\set{1_1,1_2,1_3,\ii_3,1_4,\ii_4}$ where for $s=3,4$, the element  $\ii_s\in \mb{C}_s$  corresponds to the complex number $\ii \in \mbm{C}$.  We have the canonical trace $\morph{\tr: \K}{\mb{R}}$ with $\tr(1_s)=1$ for each $s\in \ninterv{4}$ and $\tr(\ii_3)=\tr(\ii_4)=0$. Then $(\K,\tr)$ is a symmetric and separable $\mbm{R}$-algebra with  Casimir element 
  $\z[\Ke]=\som{s=1}{4}{1_s\otimes \dual{1_s}}+\ii_3\otimes \dual{\ii_3} + \ii_4\otimes \dual{\ii_4}
= \som{s=1}{4}{1_s\otimes 1_s}-\ii_3\otimes \ii_3 - \ii_4\otimes \ii_4          $.
Let $M$ be a $\K$-bimodule. Recall that  $\Zc(M)$ is the central $\Zc(\K)$-bimodule consisting  of $\K$-central elements in $M$, the associated Casimir operator is
$\zc: M \to \Zc(M); x\mapsto \zc(x)=\som{s=1}{4}{1_s \mul x \mul \dual{1_s}+ \ii_3 \mul x \mul \dual{\ii_3}+\ii_4 \mul x \mul \dual{\ii_4} }=
\som{s=1}{4}{1_s \mul x \mul 1_s - \ii_3 \mul x \mul \ii_3-\ii_4 \mul x \mul \ii_4 }  $.
Also recall by Corollary~\ref{cor.sympot-perfectfield} that, since $\mbm{R}$ is a perfect field,  any potential on an $\mbm{R}$-modulated quiver is  symmetric, cyclic (left or right) permutation  mimics the simply laced case: it is  obtained as the image under the Casimir operator of the corresponding ordinary cyclic (left or right)  permutation.

As in  the illustrative section~\ref{secexample}, we fix some notations for some useful $\K$-bimodules here. For all $(M,i,j)\in \cal{X}(\K):=\set{(\mb{R},1,2),(\mb{R},2,1),(\mb{C},1,3),(\mb{C},3,1),(\mb{C},2,3),(\mb{C},3,2),(\mb{C},1,4),(\mb{C},4,1),(\mb{C},3,4),(\mb{C},4,3)}$, we write $\indd{i}{M}{j}$ for the  natural $\kk_i\sdash\kk_j$-bimodule structure on $M$; and  when viewed as element in $\indd{i}{M}{j}$  each $x\in M$ is still written as $x$ (if this can be easily inferred from the context), otherwise $x$ is subscripted as $\elts{x}{i}{j}$.  For each $(\mb{C},s,t)\in \cal{X}(\K)$,  
the \emph{conjugate natural bimodule}
$\indds{s}{\mb{C}}{\bbar{t}}$  is   obtained by conjugating the right module structure of the natural bimodule $\indds{s}{\mb{C}}{t}$: thus we have  $z\mul x\mul z':=zx\bbar{z'}$ for all $z\in \mb{\kk}_s,z'\in \mb{C}_t$ and $x\in \indds{s}{\mb{C}}{\bbar{t}}$. Also put $\elts{\!1}{s}{\bbar{t}}=1$ and $\elts{\ii}{s}{\bbar{t}}=\ii$ as elements of $\indds{s}{\mb{C}}{\bbar{t}}$, and observe for example that $\ii_3\mul\elts{\!1}{3}{\bbar{4}}=-\elts{\!1}{3}{\bbar{4}}\mul \ii_4=\elts{\ii}{3}{\bbar{4}}$.

Below, each symmetrizable dualizing pair over $(\K,\tr)$ is (naturally isomorphic) to one of the following.
\begin{itemize}
  \item [$\bullet$] The self-dual pairs $\set{\indds{i}{\mb{R }}{j},\indds{j}{\mb{ R }}{i}}$ or $\set{\indds{s}{\mb{C }}{t},\indds{t}{\mb{ C }}{s}}$, with $(i,j)=(1,2),(2,1)$ and $(s,t)=(2,3),(3,2),(3,4),(4,3)$; here associated bilinear forms are given by the ordinary multiplication.
  \item [$\bullet$] The pairs $\set{\indds{1}{\mb{C }}{s},\indds{s}{\mb{C }}{1}}$ with $s=3,4$ and with associated bilinear forms given by the ordinary multiplication $\morph{\indds{s}{\mb{C }}{1} \otimes \indds{1}{\mb{C }}{s}}{\mb{C}_s}$ and the map
      $\morph{\indds{1}{\mb{C }}{s} \otimes \indds{s}{\mb{C }}{1}}{\mb{R}_1:\, (z\otimes z')\mapsto \scprd{z\otimes z'}=\tr(zz')}$.
  \item [$\bullet$] The conjugate pairs $\set{\indds{s}{\mb{C }}{\bbar{t}},\indds{t}{\mb{C }}{\bbar{s}}}$ with $(s,t)=(3,4),(4,3)$ and with associated bilinear forms induced  by conjugating the second argument of the ordinary multiplication:
      $\morph{\indds{s}{\mb{C }}{\bbar{t}} \otimes \indds{t}{\mb{C }}{\bbar{s}}}{\mb{C}_s:\, (z\otimes z')\mapsto \scprd{z\otimes z'}=z\bbar{z'}}$.
\end{itemize}

Now start with the modulated quiver with zero potential $\Seq{\mQ^{(0)}: \mb{R}_1 \to[/above={\indd{\mb{R}}{\mb{R}}{\mb{R}}}] \mb{R}_2
\to[/above={\indd{\mb{R}}{\mb{C}}{\mb{C}}}] \mb{C}_3 \to[/above={\indd{\mb{C}}{\mb{C}}{\mb{C}}}] \mb{C}_4}$. First observe the following picture of successive mutations of the underlying valued quivers $\lin{\msf{F}_4}$ of $\mQ^{(0)}$:

$\Seq{\lin{\msf{F}_4}: 1\to 2\to[/above={1,2}] 3\to 4\; \;  \to[/dir=<->,/dist=1.2,/above={\mu_2}] \! \graphI[/bline=1em,/dir1=<-,/scale=0.8]{\NodeI{1}{2}{3}{4}}{}{1,2}{1,2}{}\!
 \to[/dir=<->,/dist=1.2,/above={\mu_3}] \; \cyclFd[/bline=1em,]{1}{2}{3}{4}{1,2}{1,2}{1,2} \;
\graph[/clsep=2]{ {} \& {} }{\path[<->] (m-1-1) edge node[sloped,above] {$\mu_4$} ($ (0,0.25)+(m-1-2) $) edge node [sloped,below] {$\mu_1$} ($ (0,-0.25)+(m-1-2) $) } \!\!\!\!
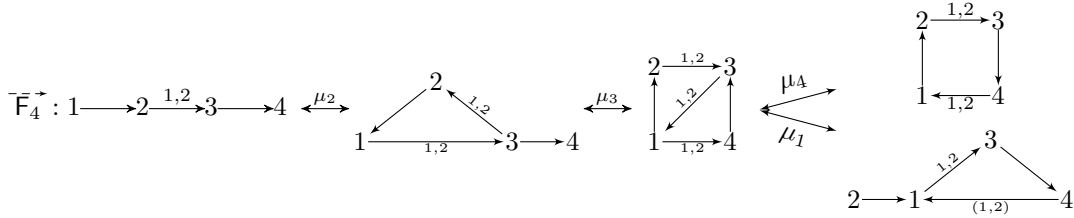
\begin{tabular}{c} \cyclF{1}{2}{3}{4}{1,2}{1,2}  \\
\graphIi[/scale=0.8]{\NodeIi{2}{1}{3}{4}}{}{1,2}{}{(1,2)}  \end{tabular}} $

The first mutation below is clear from the definition of mutation, where as usual a tensor element $x\otimes y$ is also written as $x\mul y$ or   $xy$ and where the bimodule $[\indd{1}{\mb{R}}{2}\otimes \indd{2}{\mb{C}}{3}]$ is naturally identified with  $\indd{1}{\mb{C}}{3}$: \\[-1.5\baselineskip]
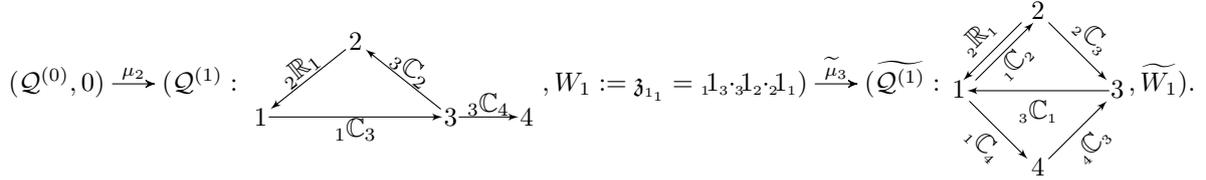
\begin{center}
$\Seq{(\mQ^{(0)},0)\; \to[/dir=->,/dist=1,/above={\mu_2}] \;
( \mQ^{(1)}: \graphI[/bline=1em,/dir1=<-]{
\NodeI{1}{2}{3}{4}}{\textstyle{\indd{2}{\mb{R}}{1}}}{
\textstyle{\indd{3}{\mb{C}}{2}}}{
\textstyle{\indd{1}{\mb{C}}{3}}}{\textstyle{\indd{3}{\mb{C}}{4}}}, W_1:=\z[1_1]=\elts{\!1}{1}{3}\mul\elts{\!1}{3}{2}\mul\elts{\!1}{2}{1})\;
 \to[/dir=->,/dist=1,/above={\wtilde{\mu}_3}]\;
( \wtilde{\mQ^{(1)}}:
 \graph[/bline=0.08,/dist=1.8,/clsep=1.8]{{} \& 2 \& {} \\ 1\& \&  3 \\ {} \& 4 \& {} }{
 \path[->] ($ (0.08,0)+(m-2-1.north east) $) edge node[sloped,below=-0.08] {$\indds{1}{\mb{C}}{2}$} ($ (0.08,0)+(m-1-2.south west) $)  ($ (-0.08,0)+ (m-1-2.south west)  $) edge  node[above,sloped] {$\indds{2}{\mb{R}}{1}$} ($ (-0.08,0)+ (m-2-1.north east)  $)
 (m-1-2) edge node[sloped,above] {$\indds{2}{\mb{C}}{3}$}  (m-2-3)
  (m-2-3) edge node[below] {$\indds{3}{\mb{C}}{1}$} (m-2-1)
 (m-2-1) edge node[sloped,below] {$\indds{1}{\mb{C}}{4}$} (m-3-2)
 (m-3-2) edge node[sloped,below] {$\indds{4}{\mb{C}}{3}$} (m-2-3);
 }, \wtilde{W_1}). }$
\end{center}
For the semi-mutation $\wtilde{\mu}_3$ above,  we naturally identify   $[\indds{1}{\mb{C}}{3}\otimes_{\mb{C}_3} \indds{3}{\mb{C}}{2}]$ and
$[\indds{1}{\mb{C}}{3}\otimes_{\mb{C}_3} \indds{3}{\mb{C}}{4}]$  with   $\indds{1}{\mb{C}}{2}$  and $\indds{1}{\mb{C}}{4}$ respectively. We  have
$\wtilde{W_1}=[W_1]+\z[1_3]$. Here the Casimir element $\z[1_3]$ is the sum of two Casimir elements:
$\z[1_3]=\z[{\indds{1}{\mb{C}}{2}\otimes (\indds{2}{\mb{C}}{3}\mul \indds{3}{\mb{C}}{1})}]+\z[{\indds{1}{\mb{C}}{4}\otimes (\indds{4}{\mb{C}}{3}\mul \indds{3}{\mb{C}}{1})}]$, with
 $\z[\indds{1}{\mb{C}}{2}\otimes (\indds{2}{\mb{C}}{3}\mul \indds{3}{\mb{C}}{1})]=\elts{\!1}{1}{2}\mul\elts{\!1}{2}{3}\mul\elts{\!1}{3}{1}+
\elts{\ii}{1}{2}\mul (-\elts{\ii}{2}{3}\mul\elts{\!1}{3}{1})$
 and $\z[{\indds{1}{\mb{C}}{4}\otimes (\indds{4}{\mb{C}}{3}\mul \indds{3}{\mb{C}}{1})}]=\elts{\!1}{1}{4}\mul\elts{\!1}{4}{3}\mul\elts{\!1}{3}{1}$.
 Thus
 $\wtilde{W_1}= \elts{\!1}{1}{2}\mul\elts{\!1}{2}{1} + \elts{\!1}{1}{2}\mul\elts{\!1}{2}{3}\mul\elts{\!1}{3}{1}-
\elts{\ii}{1}{2}\mul\elts{\ii}{2}{3}\mul\elts{\!1}{3}{1} + \elts{\!1}{1}{4}\mul\elts{\!1}{4}{3}\mul\elts{\!1}{3}{1} $ with $2$-cyclic component $\wtilde{W_1}_{,2}=\elts{\!1}{1}{2}\mul\elts{\!1}{2}{1} \in \indds{1}{\mb{C}}{2}\mul \indds{2}{\mb{R}}{1}$.

 At point $1=\dual{(\elts{\!1}{2}{1})}\in \dual{(\indds{2}{\mb{R}}{1})}$, we  compute:
$\cderv[\dual{(\elts{\!1}{2}{1})}](\wtilde{W_1})=\elts{\!1}{1}{2}\in \indds{1}{\mb{C}}{2}$. In view of the Casimir element $\z[{\indd{1}{\mb{C}}{2}\otimes \indd{2}{\mb{C}}{1}}]=\elts{\!1}{1}{2}\mul \elts{\!1}{2}{1}+\elts{\ii}{1}{2}\mul(-\elts{\ii}{2}{1})$,
  we have $\cderv[{\dual{(\elts{\!1}{1}{2})}}](\wtilde{W_1})=\elts{\!1}{2}{1} + \elts{\!1}{2}{3}\mul\elts{\!1}{3}{1} \in \indds{2}{\mb{R}}{1} \oplus \indds{2}{\mb{C}}{3}\mul \indds{3}{\mb{C}}{1}$. Thus  the trivial bimodule for $(\wtilde{\mQ^{(1)}},\wtilde{W_1})$  is $\indds{2}{\mb{R}}{1}\oplus \mb{R}\mul \elts{\!1}{1}{2}$ (and is of course   a direct summand of the arrow bimodule of $\wtilde{\mQ^{(1)}}$),  the corresponding reduced bimodule is $\mb{R}\mul\elts{\ii}{1}{2}\oplus \indds{2}{\mb{C}}{3}\oplus \indds{3}{\mb{C}}{1}\oplus\indds{1}{\mb{C}}{4} \oplus \indds{4}{\mb{C}}{3}$; the closed ideal $\Jtriv{\wtilde{W_1}}=\bbar{\K\mul \set{\elts{\!1}{1}{2},2\mul\elts{\!1}{2}{1}+\elts{\!1}{2}{3}\mul\elts{\!1}{3}{1}}}$  (of the complete path algebra of $\mQ^{(1)}$) is the kernel of a   reduction $\pi:\what{\mb{R}\wtilde{\mQ^{(1)}}}\to \what{\mb{R}\Red{(\wtilde{\mQ^{(1)}})}}$ which fixes the reduced arrow bimodule and  such that:   $\pi(\elts{\!1}{1}{2})=0$, $\pi( \elts{\!1}{2}{1} + \elts{\!1}{2}{3}\mul\elts{\!1}{3}{1})=0$ so that $\pi(\elts{\!1}{2}{1}) =-\elts{\!1}{2}{3}\mul\elts{\!1}{3}{1}$. Thus the  reduced potential   is given by
$\pi(\wtilde{W_1})= -\elts{\ii}{1}{2}\mul \elts{\ii}{2}{3}\mul\elts{\!1}{3}{1} +\mul\elts{\!1}{1}{4}
\mul\elts{\!1}{4}{3}\mul\elts{\!1}{3}{1} $. Naturally identifying $\mb{R}\mul\elts{\ii}{1}{2}$ with  $\indds{1}{\mb{R}}{2}$,  the previous details are summarized in the following picture:
\\[-1.5\baselineskip]
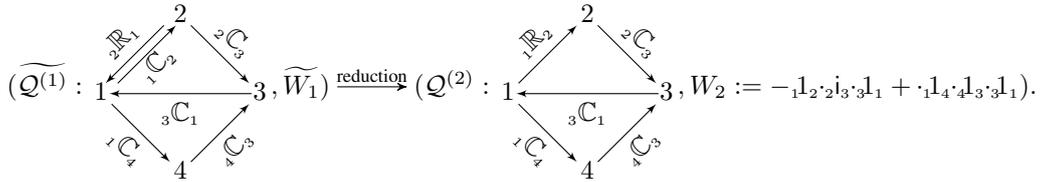
\begin{center}
$ \Seq{ ( \wtilde{\mQ^{(1)}}:
 \graph[/bline=0.08,/dist=1.8,/clsep=1.8]{{} \& 2 \& {} \\ 1\& \&  3 \\ {} \& 4 \& {} }{
 \path[->] ($ (0.08,0)+(m-2-1.north east) $) edge node[sloped,below=-0.08] {$\indds{1}{\mb{C}}{2}$} ($ (0.08,0)+(m-1-2.south west) $)  ($ (-0.08,0)+ (m-1-2.south west)  $) edge  node[above,sloped] {$\indds{2}{\mb{R}}{1}$} ($ (-0.08,0)+ (m-2-1.north east)  $)
 (m-1-2) edge node[sloped,above] {$\indds{2}{\mb{C}}{3}$}  (m-2-3)
  (m-2-3) edge node[below] {$\indds{3}{\mb{C}}{1}$} (m-2-1)
 (m-2-1) edge node[sloped,below] {$\indds{1}{\mb{C}}{4}$} (m-3-2)
 (m-3-2) edge node[sloped,below] {$\indds{4}{\mb{C}}{3}$} (m-2-3);
 }, \wtilde{W_1} )  \;  \to[/dir=->,/dist=2,/above={\text{reduction}}]\;
   ( \mQ^{(2)}:
 \graph[/bline=0.08,/dist=1.8,/clsep=1.8]{{} \& 2 \& {} \\ 1\& \&  3 \\ {} \& 4 \& {} }{
 \path[->] (m-2-1) edge node[sloped,above] {$\indds{1}{\mb{R}}{2}$} (m-1-2)   (m-1-2) edge node[sloped,above] {$\indds{2}{\mb{C}}{3}$}  (m-2-3)
  (m-2-3) edge node[below] {$\indds{3}{\mb{C}}{1}$} (m-2-1)
 (m-2-1) edge node[sloped,below] {$\indds{1}{\mb{C}}{4}$} (m-3-2)
 (m-3-2) edge node[sloped,below] {$\indds{4}{\mb{C}}{3}$} (m-2-3);
 }, W_2:= -\elts{\!1}{1}{2}\mul \elts{\ii}{2}{3}\mul\elts{\!1}{3}{1} +\mul\elts{\!1}{1}{4} \mul\elts{\!1}{4}{3}\mul\elts{\!1}{3}{1}).} $
\end{center}

We can perform more mutations as shown is the following picture, where one should notice the presence of  the conjugate natural bimodule $\indds{3}{\mb{C}}{\bbar{4}}$ in the last modulated quiver and   $ W_4:=\frac{1}{2}(\elts{\!1}{4}{1}\mul \elts{\!1}{1}{3} + \elts{\ii}{4}{1}\mul \elts{\ii}{1}{3})\elts{\!1}{3}{\bbar{4}}$.
\\[-1.0\baselineskip]
\begin{center}
$ \Seq{ (\mQ^{(2)},W_2)\;
\graph[/clsep=3]{ {} \& {} }{\path[->] (m-1-1) edge node[sloped,above] {$\wtilde{\mu}_4$} ($ (0,0.5)+(m-1-2) $) edge node [sloped,below] {$\wtilde{\mu}_1$} ($ (0,-0.5)+(m-1-2) $) } 
\hspace{-3.5em}
\begin{array}{l} ( \wtilde{\mu}_4(\mQ^{(2)}): \graph[/bline=0.08,/dist=1.8,/clsep=1.8]{{} \& 2 \& {} \\ 1\& \&  3 \\ {} \& 4 \& {} }{
 \path[->] (m-2-1) edge node[sloped,above] {$\indds{1}{\mb{R}}{2}$} (m-1-2)   (m-1-2) edge node[sloped,above] {$\indds{2}{\mb{C}}{3}$}  (m-2-3)
  ($ (0,0.08)+(m-2-3.west) $) edge node[above] {$\indds{3}{\mb{C}}{1}$}
  ($ (0,0.08)+(m-2-1.east) $)
  ($ (0,-0.08)+(m-2-1.east) $) edge node[below] {$\indds{1}{\mb{C}}{3}$}
  ($ (0,-0.08)+(m-2-3.west) $)
 (m-2-1) edge[<-] node[sloped,below] {$\indds{4}{\mb{C}}{1}$} (m-3-2)
 (m-3-2) edge[<-] node[sloped,below] {$\indds{3}{\mb{C}}{4}$} (m-2-3);
 }, \wtilde{\mu}_4(W_2)) \; \to[/dir=->,/dist=2,/above={\text{reduction}}] \;  \Bigl( \mQ^{(3)}: \graph[/bline=0.08,/dist=1.8,/clsep=1.8]{{} \& 2 \& {} \\ 1\& \&  3 \\ {} \& 4 \& {} }{
 \path[->] (m-2-1) edge node[sloped,above] {$\indds{1}{\mb{R}}{2}$} (m-1-2)   (m-1-2) edge node[sloped,above] {$\indds{2}{\mb{C}}{3}$}  (m-2-3)
  (m-2-3) edge node[sloped,below] {$\indds{3}{\mb{C}}{4}$} (m-3-2)
 (m-3-2) edge node[sloped,below] {$\indds{4}{\mb{C}}{1}$} (m-2-1);
 },   W_3:=\elts{\!1}{1}{2} \mul\elts{\ii}{2}{3} \mul\elts{\!1}{3}{4} \mul\elts{\!1}{4}{1} \Bigr);   \\[-1.8\baselineskip]
 (\wtilde{\mu_1}(\mQ^{(2)}): \graph[/bline=0.08,/dist=2.5,/clsep=2.8]{{} \& 2 \& {} \\ 1\& \&  3 \\ {} \& 4 \& {} }{
 \path[->] (m-2-1) edge[<-] node[sloped,above] {$\indds{2}{\mb{R}}{1}$} (m-1-2)   ($ (0.08,0)+(m-1-2.south east) $) edge node[sloped,above] {$\indds{2}{\mb{C}}{3}$}  ($ (0.08,0)+(m-2-3.north west) $)
 ($ (-0.08,0)+(m-1-2.south east) $) edge[<-] node[sloped,below] {$\indds{3}{\mb{C}}{2}$}  ($ (-0.08,0)+(m-2-3.north west) $)
   (m-2-3) edge[<-] node[above,pos=0.6] {$\indds{1}{\mb{C}}{3}$} (m-2-1)
 (m-2-1) edge[<-] node[sloped,below] {$\indds{4}{\mb{C}}{1}$} (m-3-2)
($ (-0.04,0.04)+(m-3-2.north east) $) edge node[sloped,above] {$\indds{4}{\mb{C}}{3}$}  ($ (-0.04,0.04)+(m-2-3.south west) $)
 ($ (0.04,-0.04)+(m-3-2.north east) $) edge[<-] node[sloped,below] {$[\indds{3}{\mb{C}}{1}\mul \indds{1}{\mb{C}}{4}]$}  ($ (0.04,-0.04)+(m-2-3.south west) $);
 },\wtilde{\mu_1}(\lperm(W_2))) \; \to[/dir=->,/dist=2,/above={\text{reduction}}] \;  \Bigl( \mQ^{(4)}: \graph[/bline=0.08,/dist=1.8,/clsep=1.8]{{} \& 2 \& {} \\ 1\& \&  3 \\ {} \& 4 \& {} }{
  \path[<-] (m-2-1) edge node[sloped,above] {$\indds{2}{\mb{R}}{1}$} (m-1-2)
  (m-2-3) edge node[above] {$\indds{1}{\mb{C}}{3}$} (m-2-1)
 (m-2-1) edge node[sloped,below] {$\indds{4}{\mb{C}}{1}$} (m-3-2)
 (m-3-2) edge node[sloped,below] {$\indds{3}{\mb{C}}{\bbar{4}}$} (m-2-3);
 },   W_4 \Bigr).
\end{array}
}$
\end{center}
The details for the semi-mutation $\wtilde{\mu}_4$ and the reduction in the first row of the above diagram are similar (in their form) to the semi-mutation and corresponding reduction from the previous paragraph. Let us shed more light on how the second row of the above diagram is obtained. Write $B$ for the arrow bimodule of $\mQ^{(2)}$. We have naturally identified the   bimodule $[\indds{3}{\mb{C}}{1}\mul \indds{1}{\mb{R}}{2}]$ in $\wtilde{\mu}_1(\mQ^{(2)})$ with  $\indds{3}{\mb{C}}{2}$. The arrow bimodule  $\wtilde{\mu}_1(\mQ^{(2)})$ is $\wtilde{B}=(\indds{3}{\mb{C}}{2}\oplus \indds{2}{\mb{C}}{3}) \oplus ([\indds{3}{\mb{C}}{1}\mul \indds{1}{\mb{C}}{4}]\oplus \indds{4}{\mb{C}}{3}) \oplus B_1'$ with  $B_1':= \indds{2}{\mb{R}}{1}\oplus  \indds{1}{\mb{C}}{3} \oplus \indds{4}{\mb{C}}{1}$. Note that each component of  $W_2=-\elts{\!1}{1}{2}\mul \elts{\ii}{2}{3}\mul\elts{\!1}{3}{1} +\elts{\!1}{1}{4} \mul\elts{\!1}{4}{3}\mul\elts{\!1}{3}{1}$  at point $\idl[1]$, in order to perform the semi-mutation $\wtilde{\mu}_1$ we  must replace $W_2$ by a cyclically equivalent potential $W_2'$ not starting at point $\idl[1]$. We can take   $W_2':=\lperm(W_2)$. Using the Casimir operator  $\morph{\zc:\mb{R}\mQ^{(2)}}{\Zc(\mb{R}\mQ^{(2)})}$,
we know that $W_2':=\lperm(-\elts{\!1}{1}{2}\mul \elts{\ii}{2}{3}\mul\elts{\!1}{3}{1} +\elts{\!1}{1}{4} \mul\elts{\!1}{4}{3}\mul\elts{\!1}{3}{1})=
\zc(-\elts{\ii}{2}{3}\mul\elts{\!1}{3}{1}\mul \elts{\!1}{1}{2} +
\elts{\!1}{4}{3}\mul\elts{\!1}{3}{1}\mul\elts{\!1}{1}{4})=
-\elts{\ii}{2}{3}\mul\elts{\!1}{3}{1}\mul \elts{\!1}{1}{2} +
(\elts{\!1}{4}{3}\mul\elts{\!1}{3}{1}\mul\elts{\!1}{1}{4}
-\ii_4\mul \elts{\!1}{4}{3}\mul\elts{\!1}{3}{1}\mul\elts{\!1}{1}{4} \ii_4)=-\elts{\ii}{2}{3}\mul\elts{\!1}{3}{1}\mul \elts{\!1}{1}{2} +
\elts{\!1}{4}{3}(\elts{\!1}{3}{1}\mul\elts{\!1}{1}{4}
-\elts{\ii}{3}{1}\mul\elts{\ii}{1}{4})$.
Now $\wtilde{\mu}_1(W_2')=[W_2']+\z[1_1]$. Let $S$ be the potential obtained from the latter by replacing the Casimir element $\z[1_1]=\z[{[B1_1B]\otimes {\dual{(B1_1B)}}}]$ with the Casimir element
$\z[1_1]'=\z[{\dual{(B1_1B)}\otimes [B1_1B]}]$ (since $\z[1_1]'$ is the common value of the left and right permutation  $\z[1_1]$). We have
\[ \myeqar{
\z[1_1]'= & \z[{(\indd{2}{\mb{R}}{1}\mul \indd{1}{\mb{C}}{3})\otimes \indd{3}{\mb{C}}{2}}] + \z[{(\indd{4}{\mb{C}}{1}\mul \indd{1}{\mb{C}}{3})\otimes [\indds{3}{\mb{C}}{1}\mul \indds{1}{\mb{C}}{4}]}] =\elts{\!1}{2}{1}\mul \elts{\!1}{1}{3}\mul\elts{\!1}{3}{2} +(\elts{\!1}{4}{1}\mul \elts{\!1}{1}{3}\mul[\elts{\!1}{3}{1}\mul \elts{\!1}{1}{4}] -\elts{\ii}{4}{1}\mul \elts{\!1}{1}{3}\mul [\elts{\!1}{3}{1}\mul \elts{\ii}{1}{4}])  \text{ and } \\
S= & -\elts{\ii}{2}{3}\mul \elts{\!1}{3}{2} + \elts{\!1}{4}{3}([\elts{\!1}{3}{1}\mul \elts{\!1}{1}{4}]-[\elts{\ii}{3}{1}\mul \elts{\ii}{1}{4}]) + (\elts{\!1}{2}{1}\mul \elts{\!1}{1}{3}\mul\elts{\!1}{3}{2}+
\elts{\!1}{4}{1}\mul \elts{\!1}{1}{3}\mul[\elts{\!1}{3}{1}\mul \elts{\!1}{1}{4}] -\elts{\ii}{4}{1}\mul \elts{\!1}{1}{3}\mul [\elts{\!1}{3}{1}\mul \elts{\ii}{1}{4}]). } \]
Let $u:=[\elts{\!1}{3}{1}\mul \elts{\!1}{1}{4}]-[\elts{\ii}{3}{1}\mul \elts{\ii}{1}{4}]=[1\otimes 1-\ii \otimes \ii ], v:=[\elts{\!1}{3}{1}\mul \elts{\!1}{1}{4}]+[\elts{\ii}{3}{1}\mul \elts{\ii}{1}{4}]=[1\otimes 1+\ii \otimes \ii ]\in [\indds{3}{\mb{C}}{1}\mul \indds{1}{\mb{C}}{4}]$, note that $\ii_3\mul u=u\mul \ii_4$  while $\ii_3 \mul v=-v\mul \ii_4$. Moreover,   $[\indds{3}{\mb{C}}{1}\mul \indds{1}{\mb{C}}{4}]=U\oplus V$ where $U$ is the subbimodule generated by $u$ and $V$ is the subbimodule generated by $v$, and we have natural isomorphisms $U\cong \indds{3}{\mb{C}}{4} $ and $V \cong \indds{3}{\mb{C}}{\bbar{4}}$. Since $1\otimes 1=\frac{1}{2}(u+v)$ and $ \ii\otimes \ii=\frac{1}{2}(-u+v)$, we have:
\begin{center}
$\myeqar{S= & -\elts{\ii}{2}{3}\mul \elts{\!1}{3}{2} + \elts{\!1}{4}{3}\mul u + \elts{\!1}{2}{1}\mul \elts{\!1}{1}{3}\mul\elts{\!1}{3}{2}+ \frac{1}{2}\mul
\elts{\!1}{4}{1}\mul \elts{\!1}{1}{3}\mul (u+v) + \frac{1}{2}\mul\elts{\ii}{4}{1}\mul \elts{\ii}{1}{3}\mul (-u+v) \\
=& -\elts{\ii}{2}{3}\mul \elts{\!1}{3}{2} + \elts{\!1}{4}{3}\mul u + \elts{\!1}{2}{1}\mul \elts{\!1}{1}{3}\mul\elts{\!1}{3}{2}+ \frac{1}{2}
(\elts{\!1}{4}{1}\mul \elts{\!1}{1}{3} -\elts{\ii}{4}{1}\mul \elts{\ii}{1}{3})u + \frac{1}{2}(\elts{\!1}{4}{1}\mul \elts{\!1}{1}{3}+\elts{\ii}{4}{1}\mul \elts{\ii}{1}{3})v . }$
\end{center}
The $2$-cyclic component of $S$ is $\triv{S}:=-\elts{\ii}{2}{3}\mul \elts{\!1}{3}{2} + \elts{\!1}{4}{3}\mul u$. At points $\dual{(\elts{\!1}{3}{2})}\in \dual{(\indds{3}{\mb{C}}{2})}$ and $\dual{(\elts{\ii}{2}{3})}\in \dual{(\indds{2}{\mb{C}}{3})}$ we  get: $\cderv[{\dual{(\elts{\!1}{3}{2})}}]S=-\elts{\ii}{2}{3}+\elts{\!1}{2}{1}\mul
\elts{\!1}{1}{3}$,  $\cderv[{\dual{(\elts{\ii}{2}{3})}}]S=-\elts{\!1}{3}{2}$.
For $\dual{(\elts{\!1}{4}{3})}\in \dual{(\indds{4}{\mb{C}}{3})}$ we have:
$\cderv[{\dual{(\elts{\!1}{4}{3})}}]S=u$. The Casimir element in the product $\dual{([\indds{3}{\mb{C}}{1}\mul \indds{1}{\mb{C}}{4}])}\otimes [\indds{3}{\mb{C}}{1}\mul \indds{1}{\mb{C}}{4}]$ is given by the element
$ \dual{u}\otimes u+ \dual{v} \otimes v$ with $\dual{([\indds{3}{\mb{C}}{1}\mul \indds{1}{\mb{C}}{4}])}=[\indds{4}{\mb{C}}{1}\mul \indds{1}{\mb{C}}{3}]$, (here we note that $\dual{u}=[\elts{\!1}{4}{1}\mul \elts{\!1}{1}{3}]-[\elts{\ii}{4}{1}\mul \elts{\ii}{1}{3}], \dual{v}=[\elts{\!1}{4}{1}\mul \elts{\!1}{1}{3}]+[\elts{\ii}{4}{1}\mul \elts{\ii}{1}{3}]$). We compute the cyclic derivative: $\cderv[\dual{u}]S=\elts{\!1}{4}{3}+\frac{1}{2}
(\elts{\!1}{4}{1}\mul \elts{\!1}{1}{3} -\elts{\ii}{4}{1}\mul \elts{\ii}{1}{3})$. We  now deduce that  $\triv{\wtilde{B}}=(\indds{2}{\mb{C}}{3}\oplus \indds{3}{\mb{C}}{2})\oplus (\indds{4}{\mb{C}}{3}\oplus U)$ and the corresponding reduced bimodule is $\red{\wtilde{B}}=V\oplus B_1'\cong \indds{3}{\mb{C}}{\bbar{4}}\oplus B_1'$. Associated with the decomposition $\wtilde{B}=\triv{\wtilde{B}} \oplus \red{\wtilde{B}}$, we have a reduction $\pi$ fixing $\red{\wtilde{B}}$ and such that $\pi(\elts{\!1}{3}{2})=0=\pi(u)$, thus the reduced potential is
$\pi(S)=\frac{1}{2}(\elts{\!1}{4}{1}\mul \elts{\!1}{1}{3}+\elts{\ii}{4}{1}\mul \elts{\ii}{1}{3})v$, naturally identified with $\frac{1}{2}(\elts{\!1}{4}{1}\mul \elts{\!1}{1}{3}+\elts{\ii}{4}{1}\mul \elts{\ii}{1}{3})\elts{\!1}{3}{\bbar{4}}$ under the identification $V\cong \indds{3}{\mb{C}}{\bbar{4}}$.
Hence  the details for  the second row of the above   are   complete.

\begin{Rem} If  instead of a perfect field we  consider a non perfect field, then all the sequences of mutations and reductions above can till be performed, provided, in view of Corollary~\ref{Cor-red-nonsplit}, skew reductions are also allowed.
\end{Rem}


\section{Graded modulated quiver with potentials}
\label{sec:gradedcontext}  
This section is motivated by \cite[\S~6.2]{AO} about graded quivers with potentials.  We will quickly explain why the results of preceding sections holds in the graded context.    
We fix an abelian group $G$ which should be $\mb{Z}$, $\mb{Z}/p\mb{Z}$ or $\mb{Z}^n$ for some $n,p\in \N$. Let $\C$ be an additive category. A \emph{$G$-graded object} in $\C$   is just a family $X=(X^p)_{p\in \mb{Z}}$ of objects of $\C$; the \emph{degree-$p$  component}  of $X$ is $X^p$. A \emph{graded morphism $\morph{f: X}{Y}$ of degree $n\in G$}  between two graded objects consists of a family of morphisms $\morph {f^p: X^p}{Y^p{p+n}},\, p\in \mb{Z}$.  Graded morphisms of degree $0$ are simply referred to as \emph{graded morphisms}.    A \emph{complex} (or a \emph{\dg~ (differential graded) object}  in $\C$ consists of a graded object  $X=(X_p)_{p\in \mb{Z}}$  together with a \emph{differential} $d=d_X$, the latter is  a graded morphism of degree $1$  such that '$d\circ  d=0$' (that is, $d^p\circ d^{p-1}=0$ for all $p \in G$).   When the category $\C$ has all direct sums, we identify each graded object $X$ with the direct sum $\somd{p\in G}{X^p}$.  
Giving a $\kk$-algebra $\Lambda$, denote by  $\Grd(\Lambda)$   the category of $G$-graded (right) $\Lambda$-modules and graded morphisms (of degree $0$). Let $M=\somd{p\in G}{M^p}$ in  $\Grd(\Lambda)$. The   \emph{$G$-graded left   $\Lambda$-module  $M'=\HM{A}{M,A}$}  has components 
$M'^p=\HM{A}{M^{-p},A}$, $p\in G$; with this $G$-grading,  $\HM{A}{M,A}$ is called the \emph{dual} of the $G$-graded  $\Lambda$-module $M$. For $n\in G$, the \emph{$n$-shift} of $M$ is the graded module $M[n]$ with components $(M[n])^p=M^{p+n}$ for all $p\in G$.  The \emph{tensor product} of a $G$-graded left $\Lambda$-module  $L$ by $M$ is the   $G$-graded $\kk$-module $M\otimes_{\Lambda}M$  with components: $(M\otimes_{\Lambda}M)^n=\somd{p+q=n}L^p\otimes_{\Lambda}M^q, n\in G$.

\ParIt{Graded modulated quivers and their complete path algebras}  A modulated quiver  $\mQ=(B,\K,\tr)$ is \emph{$G$-graded} if the (finitely generated) arrow $\K$-bimodule $B$ is $G$-graded (thus, $B=\somd{p\in G}{B^p}$ and only finitely many components $B^p$ are nonzero). Assume that $\mQ$  is $G$-graded. 
The path  algebra $\kQ$ is  a topological $G$-graded algebra with respect to $\J{\kQ}$-adic topology on $\kQ$, with  grading induced by that of $\mQ$ and with $\K$ lying in degree $0$. For  $p\in G$, the degree-$p$ component $B^p$ of $B$ should not be confused with the notation $B^{(l)}$ for $l\in \N$, the latter being the $l$-fold tensor product of  $B$ over $\K$; in particular $(B^{(l)})^p$ is the degree-$p$  component of $B^{(l)}$.  The \emph{complete path algebra  $\kQc$} of   $\mQ$  is the completion of $\kQ$ with respect to $\J{\kQ}$-adic topology in the category $\Grd(\kk)$  of $G$-graded $\kk$-modules;  thus  $\kQc$ coincides with the projective limit (in $\Grd(\kk)$) of the natural inverse system 
\begin{center}
$\Seq{\K=\kQ/\J{\kQ}\, \to[/dir=<-] \K\oplus B\cong \kQ/\J[2]{\kQ} \dotsm \to[/dir=<-] \somd{0\leq d < l}{B^{(d)}} \cong \kQ/\J[l]{\kQ} \,  \to[/dir=<-] \somd{0\leq d < l+1}{B^{(d)}}\cong  \kQ/\J[l+1]{\kQ} \dotsm }$
\end{center}
As $G$-graded $\kk$-module, we get: $\kQc=\somd{p\in G}{\kQc^p}$ with degree-$p$ component $\kQc^p=\prd{l\geq 0}{(B^{(l)})^p}$.

Recall that $B$ is part of a symmetrizable dualizing pair $\set{B,\dual{B};\bilf}$.  Since the dual $\dual{B}$ is canonically isomorphic to the $\kk$-dual $\HM{\kk}{B,\kk}$, we get the following observation. 
\begin{Rem} Endowing $\dual{B}$ with the dual $G$-grading induced by that of $B$,  the dualizing pair $\set{B,\dual{B};\bilf}$ arises as direct sum of induced dualizing pairs $\set{B^p,(\dual{B})^{-p}}=\set{B^p,\dual{(B^p)}}, p\in G$,  the bilinear form $\bilf$ is $G$-graded of degree $0$ and   vanishes on $B^p\otimes (\dual{B})^q$ and $(\dual{B})^q \otimes B^p$ for all $p,q\in G$ with $q\neq -p$.
\end{Rem}


In the sequel, let $n\in G$ and $(\mQ,\mfr{m})$ be a $G$-graded  modulated quiver with potential homogeneous of degree $n$, ($\mfr{m}$ needs not be homogeneous with respect to path-grading). Each component $\mfr{m}_d\in B^{(d)}$ of $\mfr{m}$ (with respect to path grading) is therefore homogeneous of degree $n\in G$  and the following lemma is an easy observation.
\begin{Lem} \label{lem.graded-derv} For each  $0<d \in \N$, the   derivative morphisms $\morph{\lderv\mfr{m},\rderv\mfr{m}:(\dual{B})^{(d)}}{\kQc}$ and  $\morph{\cderv\mfr{m}:\dual{B}}{\kQc}$ are $G$-graded  morphisms of degree $n\in G$. In particular the trivial part $\triv{\mQ}=(\triv{B},\K,\tr)$  and  the reduced part $\mQb=\red{\mQ}=(\bbar{B},\K,\tr)$  are naturally  $G$-graded with  $G$-gradings induced by that of $B$,   the natural projection $\morph{\rho:\kQc}{\kQbc}$ is  a $G$-graded morphism,  left and right $\K$-linear right inverses to the projection $\morph{\cderv\mfr{m}_2: \dual{B}}{\triv{B}}$ can be chosen as $G$-graded morphisms of degree $-n$. 
\end{Lem}

Note that two cyclically equivalent potentials (of homogeneous degree with respect to $G$-grading)  have the same degree.  We  adapt the notion of (weak) right-equivalence and reductions to the graded context. Below,    $\mQ'=(B',\K,\tr)$ is another graded modulated quiver.
\begin{Defn} \label{defn.reduction-graded} Let $\morph{\phi:\kQc}{\kQpc}$ be a graded algebra morphism. Then $\phi$ is a \emph{path algebra morphism} if $\phi\restr{\K}=\id_{\K}$ and $\phi(B)\subset \J{\kQpc}$, in this case $\phi\restr{B}=(\phi_l)_{l\geq1}: B\to \J{\kQpc}$ with $\morph{\phi_l:B}{B'^{(l)}},\,l\geq 1$. We call  $\phi$   a \emph{reduction on $(\mQ,\mfr{m})$} if  $\phi$ is a path algebra epimorphism satisfying the graded version of  properties $\msf{(1.i)}$ and $\msf{(1.i)}$ from Definition~\ref{defn.trivpart-reduction}, namely:
 \begin{itemize}
 \item [$\msf{(1.i)'}$] $\Ker(\phi)$  is the closed graded ideal in  $\kQc$ generated by the image of a graded $\K$-bimodule morphism  $\morph{f=\Psmatr{\id \\ f'}: \triv{B}}{\triv{B}\oplus \J[2]{\kQc}}$ with $\im(f)\subset(\cderv \mfr{m})(\dual{B})$.
 \item [$\msf{(1.ii)'}$] Let $\morph{\pi: \kQc}{\kQc/\Ker(\phi)}$ be the natural projection, $\morph{\rho:B}{\bbar{B}=B/\triv{B}}$ the natural projection, and $\morph{\bbar{\rho}: \Frac{(B+\Ker(\phi))}{\Ker(\phi)}}{\bbar{B}}$ the $\K$-bimodule epimorphism  with $\rho=\bbar{\rho}\circ \pi\restr{B}$. Then $\bbar{\rho}$  has a right inverse   $\morph{ \bbar{\rho}':  \bbar{B}}{(B+\Ker(\phi))/\Ker(\phi)}$ which lifts to a left (respectively, right) graded $\K$-linear  map $\morph{ \rho':  \bbar{B}}{B}$  such that $ \phi\circ \rho' =\phi_1\circ \rho'$. 
\end{itemize}
 Similarly, \emph{(weak) right equivalences} between G-graded modulated quivers with potentials are defined.
\end{Defn}

As   direct consequence of Lemma~\ref{lem.graded-derv} above, we have the following. 
\begin{Cor}  Unitriangular automorphisms appearing in Lemma~\ref{lem.splitting-mQ}  and the proof of 
Proposition~\ref{lem.cycequiv-unitriang} can be constructed as  $G$-graded   algebra morphisms.
\end{Cor}

Applying the previous discussion and Definition~\ref{defn.reduction-graded} above, we get that the results from the first section to the fifth (as well as  the setting of \cite[\S2-6]{DWZ}) generalize to the graded context. We therefore  state the following.

\begin{Thm} \label{thm.grreduction} The reduction Theorem~\ref{theo.red-mQp} and its symmetric version Theorem~\ref{theo.sympot}   holds for $G$-graded modulated quivers with potentials of degree 
$n\in G$.  In particular,   the reduced potential is also   of degree $n$.
\end{Thm}

\ParIt{Mutation of graded modulated quiver with potentials}

Let $e\in \K$ be a point in $\mQ$ ($e$   belongs to a system of central primitive orthogonal idempotents for $\K$) satisfying \eqref{eq.idempoint}.  We then adapt  Definition~\ref{defn.mutmqpot1} to the graded context as follows (compare with \cite[\S~6.2]{AO}).

\begin{Defn}\label{defn.mutgrmqpot1} The \emph{left semi-mutation} of
$(\mQ,\mfr{m})$ at  point $e$  is the graded modulated quiver with potential
$\wtilde{\mu}_{e}^{\lft}(\mQ,\mfr{m})$  whose
  underlying non $G$-graded modulated quiver with potential   is the semi-mutation
 $(\wtilde{\mu}_{e}(\mQ),\wtilde{\mfr{m}})$ with arrow bimodule $\widetilde{B}=[BeB] \oplus e\mul \dual{B} \oplus \dual{B}\mul e \oplus \bbar{e} \mul B \mul  \bbar{e}$, and the $G$-graded arrow bimodule of $\wtilde{\mu}_{e}^{\lft}(\mQ,\mfr{m})$  is
 $\wtilde{\mu}_{e}^{\lft}(B)=[BeB] \oplus (e\mul \dual{B})[n] \oplus \dual{B}\mul e \oplus \bbar{e} \mul B \mul  \bbar{e}$ where the $G$-grading of $[BeB]$ is induced by that of the tensor product 
$B\otimes B$.   Similarly, the \emph{right semi-mutation} $\wtilde{\mu}_{e}^{\rgt}(\mQ,\mfr{m})$ is defined by letting $\wtilde{\mu}_{e}^{\rgt}(B)=[BeB] \oplus e\mul \dual{B} \oplus 
(\dual{B}\mul e)[n] \oplus \bbar{e} \mul B \mul  \bbar{e}$.
\end{Defn}

In the above definition,  the potential $\wtilde{\mfr{m}}$ is  homogeneous of the same  degree $n\in G$ as $\mfr{m}$.  We obtain the following graded version of Theorem~\ref{theo.mtmq-ceqf} and Corollary~\ref{cor.mutclequivf}, where Theorem~\ref{thm.grreduction} above is also used.

\begin{Thm}\label{theo.mtgrmq-ceqf} The  right-equivalence classes  of the left semi-mutation $\wtilde{\mu}_{e}^{\lft}(\mQ,\mfr{m})$ and the right semi-mutation $\wtilde{\mu}_{e}^{\rgt}(\mQ,\mfr{m})$ are  determined by  that of $(\mQ,\mfr{m})$. Thus, if the trivial part $\triv{(\wtilde{\mu}_{e}(\mQ), \wtilde{\mfr{m}} )}$ splits, then the weak right-equivalence classes of $\Red(\wtilde{\mu}_{e}^{\lft}(\mQ,\mfr{m}))$  and $\Red(\wtilde{\mu}_{e}^{\rgt}(\mQ,\mfr{m}))$ are determined by  that of $(\mQ,\mfr{m})$.
\end{Thm}

Let $\mrm{s}\in \set{\lft,\rgt}$, in the situation of previous theorem, the reduced $G$-graded modulated quiver with potential $\red{\wtilde{\mu}_{e}^{\mrm{s}}(\mQ,\mfr{m})}$ is the \emph{left (or right) mutation} of $(\mQ,\mfr{m})$ at point $e$,   it is unique up to weak right-equivalence (or right equivalence if $\kk$ is a perfect field).  We also deduce the following result.

\begin{Thm} In the graded context, left and right mutation are again involutive on the set of weak right-equivalence classes of $G$-graded modulated quivers with potentials of homogeneous degree in $G$.
\end{Thm}


\section{The cluster category of a graded modulated quiver with
potential}\label{sec:clct-mqp}

Here the abelian group $G$ is $\mb{Z}$ and we keep the notions of graded objects as defined in the previous section.  We let $n\in G$ and  $(\mQ,\mfr{m})$  be a  graded modulated quiver with potential homogeneous of degree $n-3 \in G$, where as before  $\mQ=(B,\K,\tr)$  and $(\K,\tr)$ is a symmetric  $\kk$-algebra. $(\mQ,\mfr{m})$  is   \emph{Jacobian-finite} whenever
$\mfr{m}\in \kQ$ and 
the Jacobian algebra $\Jc(\mQ,\mfr{m})$ is finitely generated as  $\kk$-module.

\subsection{Complete Ginzburg \dg-algebra and the generalized cluster category}
Refer to \cite{Keller06}   for  concepts about  \emph{differential graded
categories} and \emph{differential graded algebras} (in short,
\dg-categories, \dg-algebras). Simply-laced   Ginzburg  \dg-algebra  appears in
\cite[sec~4.2]{Ginzburg} (for $\mQ$ concentrated in degree $0$ and $n=3$), see also \cite[\S~2.5]{Pla2010b} and \cite[\S~6.2]{Keller08c}.  To  $\mQ$ we associate a 
\emph{graded modulated quiver} \mbox{$\mQh=\mQ\oplus\dual{\mQ}[n-2]\oplus \indd{\K}{\K}{\K}[n-1]$} where the $G$-graded modulated quiver  $\dual{\mQ}=(\dual{B},\K,\tr)$ is the dual of $\mQ$; thus the  $G$-graded arrow bimodule of $\mQh$ is  $\what{B}=B\oplus \dual{B} \oplus \indd{\K}{\K}{\K}$ with   $\indd{\K}{\K}{\K}$ concentrated in degree $0$. 

 \begin{Defn}\label{defn.dg-aldeginsburg} The \emph{complete
 Ginzburg \dg-algebra} $\what{\Gamma}_n=\what{\Gamma}_n(\mQ,\mfr{m})$ is the complete path algebra $\kQhc$ of the graded modulated quiver $\mQh$. 
$\what{\Gamma}_n$ is endowed with the unique continuous differential  $\morph[/above={\df}]{\kQh}{\kQh}$,   satisfying  the  Leibniz rule: $\df(uv)=\df(u)v + (-1)^p u\mul \df(v)$   for all  $u\in \what{\Gamma}_n^p$, given on   $\what{B}$ as follows:
\begin{enumerate}
\item[$\trr$]   $\df$ vanishes on $B$,  and  $\morph{\df\restr{\dual{B}}=\cderv\mfr{m}: \dual{B}}{\J{\kQc},\ \xi\mapsto \df(\xi)=\cderv[\xi]\mfr{m}}$.
\item[$\trr$] The restriction of $\df$ on the self-dual natural
bimodule $\indd{\K}{\K}{\K}$  is  the Casimir morphism
\begin{center}
$\morph{\z[B\otimes\dual{B}]-\z[\dual{B}\otimes B]:
\indd{\K}{\K}{\K}}{ (B\otimes\dual{B} ) \oplus (\dual{B}\otimes
B), }$ 
\end{center}
thus for all $a \in \K$ we have:
$\df(a)=a \bigl(\Som{y\in\rmQ[1]} y\mul
\dual{y} -\Som{x\in\lmQ[1]} \dual{x}\mul
x\bigl)= \bigl(\Som{y\in\rmQ[1]} y\mul
\dual{y}-\Som{x\in\lmQ[1]} \dual{x}\mul
x\bigl)a$,
 where $(\rmQ[1],\rdmQ[1])$ and $(\lmQ[1],\ldmQ[1])$ are
 respectively right and left projective bases for $B$   defined by  the
 Casimir morphisms $\z[B\otimes\dual{B}]$ and $\z[\dual{B}\otimes
 B]$.
\end{enumerate}
In case $\mfr{m}$ lies in $\kQ$,  the \emph{non-complete Ginzburg \dg-algebra $\Gamma_n$} is the path algebra  $\kQh$ endowed  with the differential defined above.
\end{Defn}

\begin{Rem}[{\cite[Lem~2.8]{KY2011}} for the simply laced case] If $\mQ$ is concentrated in degree $0$ and $n=3$, then $\Jc(\mQ,\mfr{m})$
coincides with the $0$-homology $\mrm{H}^0\what{\Gamma}_3$  of the
differential graded algebra $\what{\Gamma}_3$.
\end{Rem}

 Let $\D\what{\Gamma}_n$ be the derived category  of $\what{\Gamma}_n$ and  view $\what{\Gamma}_n$  as  object of $\D\what{\Gamma}_n$.  The\emph{perfect derived category} of $\what{\Gamma}_n$ is
the smallest full triangulated subcategory $\msf{per}\what{\Gamma}_n$ of  $\D\what{\Gamma}_n$
generated by $\what{\Gamma}_n$ and closed under taking direct summands.
Denote by $\D[fd]\what{\Gamma}_n$  the subcategory of $\D\what{\Gamma}_n$
consisting of \dg~modules $M$ with finite-length total homology, that is, the homology 
$\mrm{H}(M)=\somd{p\in\mb{Z}}{\mrm{H}^p(M)}$ has  finite length over $\kk$.  For $n=3$, it is shown in the simply-laced framework  that $\D[fd]\what{\Gamma}_n$  is a triangulated subcategory of
$\msf{per}\what{\Gamma}_n$ \cite[\S2.15,2.18,4]{KY2011}, and $\D[fd]\Gamma_n$ enjoys a relative
$n$-Calabi-Yau property in $\D\Gamma_n$  \cite[Lem~4.1]{Keller2008}
and \cite[Thm~6.3]{Keller08c}.  
 \begin{Defn} When $\mQ$ is concentrated in degree $0$ and $n=3$, the \emph{cluster category} $\C_{(\mQ,\mfr{m})}$ associated with
$(\mQ,\mfr{m})$ is   the  idempotent completion of the
triangulated quotient \mbox{$\Frac{\msf{per}\what{\Gamma}_3}{\Db\what{\Gamma}_3}$}. 
 \end{Defn}

The following questions  arise naturally since  Calabi-Yau property is fundamental in cluster theory.
\begin{enumerate}
\item[$\msf{(a)}$] Does the relative Calabi-Yau property of
$\D[fd]\what{\Gamma}_n$ in $\D\what{\Gamma}_n$ survive in the non simply-laced
framework?
 \item[$\msf{(b)}$] Suppose  $\kk$ is a field and $(\mQ,\mfr{m})$ Jacobian-finite. Is  $\C_{(\mQ,\mfr{m})}$ $\Hm$-finite, $(n-1)$-Calabi-Yau? Is $\what{\Gamma}_n$ a
cluster-tilting object in $\C_{(\mQ,\mfr{m})}$?
\end{enumerate}

A \dg-algebra $A$ is \emph{homologically smooth}  if $A\in\msf{per}(A^{e})$
where $A^e=\op{A}\otimes A$ is the enveloping  \dg~$\kk$-algebra
of $A$. And    $A$ is \emph{$n$-Calabi-Yau  as  bimodule}  if in  $\D(A^e)$
there is a bimodule isomorphism     $\msf{R}\HM{A^e}{A,A^e} \to[/above={\simeq}] A[-n]$. 
A notion of  \emph{topological and  homological  smoothness} is defined for bilaterally pseudocompact \dg-algebras \cite[\S7.11]{KY2011}.

We expect the following result due to  Bernhard Keller to hold  in the general framework.
  
\begin{Thm}[{\cite[Thm~6.3]{Keller08c}, \cite[Thm~7.17]{KY2011}}]\label{theo.keller-homlisse}
The   non-complete (resp. complete)  Ginzburg \dg~algebra (or \dg~ category)  of a    quiver with potential     is (topologically) homologically smooth and $3$-Calabi-Yau  as bimodule.
\end{Thm}

\begin{Conj}\label{conj.homglisse} Generalized Ginzburg \dg-algebras (\dg-categories) For $n=3$,  $\Gamma_n$ and  $\what{\Gamma}_n$ are (topologically) homologically smooth and $n$-Calabi-Yau as bimodules, at least when the symmetric
algebra $\K$ is separable over a field.
\end{Conj}
 
In the sequel, suppose $\mQ$ is concentrated in degree $0$.
With exactly the same argument as in \cite[Thm~3.6]{Amiot2009} and \cite[\S~7.20]{KY2011}, we derive the following.
\begin{Thm}[{\cite[Thm~3.6]{Amiot2009}, \cite[\S~7.20]{KY2011}  for simply-laced case}]\label{Cqmpot}  Suppose  Conjecture~\ref{conj.homglisse} holds and
$\kk$ is a field. Then  the  generalized cluster category $\C_{(\mQ,\mfr{m})}$
of a Jacobian-finite modulated quiver with symmetric
potential is  $\Hm$-finite $2$-Calabi-Yau and the image $T$
of the free module $\Gamma$ into $\C_{(\mQ,\mfr{m})}$ is a cluster
tilting object such that $\End_{\C_{(\mQ,\mfr{m})}}(T)$ coincides
with the Jacobian algebra $\Jc(\mQ,\mfr{m})$.  
\end{Thm}


Recall the following  interesting characterization of cluster categories  inside the context of $2$-Calabi-Yau categories.
   \begin{Thm}[Keller-Reiten\cite{KR}]\label{theo.KR-cartCatAm} Assume $\kk$ is
 a perfect field. Let $\C$ be the stable category of  a Frobenius category  
such that $\C$  is $2$-Calabi-Yau; let $\T\subset\C$ be a
cluster tilting subcategory. Then, if the  category $\md \T$ of
finite presented modules over $\T$ is hereditary then $\C$ is
exactly equivalent to the cluster category $\C_{\T}=\Frac{\Db(\md
\T)}{(\tau^{-1}[1])^{\mb{Z}}}$.
\end{Thm}

\begin{Cor}\label{cor.KR-cartC_mQ} If Conjecture~\ref{conj.homglisse} holds and  $\kk$ is a perfect field,  then  for an acyclic $\mQ$ the category  $\C_{(\mQ,0)}$
  is exactly equivalent to the cluster category  $\C_{\mQ}$ of \cite{BMRRT}.
\end{Cor}
\begin{prv} The argument of the proof is the same  as in the simply-laced case. When $\kk$
is a perfect  field and Conjecture~\ref{conj.homglisse} holds, it follows by Theorem~\ref{Cqmpot}
that $\C_{(\mQ,0)}$ is  $2$-Calabi-Yau, admitting a  cluster
tilting object $T$ such that $\Endc{T}=\kk\mQ$,  so that we have
the expected result in view of   Keller-Reiten
Theorem~\ref{theo.KR-cartCatAm}.
\end{prv}

\section*{Acknowledgements}
\addcontentsline{toc}{section}{Acknowledgements}

We are grateful to A. Zelevinsky for accepting to read two earlier versions of this work and for useful remarks and critics that he made to improve the quality of the paper.  An earlier attempt to part of the present  work started around middle
2008, after a talk Zelevinsky gave in the University
of Sherbrooke some months ago about quivers with potentials. We
are also acknowledging useful discussions  we had with Derksen
(about quivers with potentials) and Iyama (mostly about cluster
structures for non simply-laced $2$-Calabi-Yau categories) during
2008 international conference on cluster algebras and related
topic.   Also, one of the motivations to study modulated quivers with potentials comes from the need to characterize non simply-laced cluster tilted algebras of finite
representation type in terms of modulated quivers with relations.

\bibliographystyle{elsarticle-num}

{\small{
\begin{tabular}{l}
\texttt{Bertrand Nguefack}\\
Department of Mathematics, Univ. of Yaounde I, \\
P.O. Box: 812 Yaounde, Cameroon\\
e-mail: \texttt{b.nguefack@uy1.uninet.cm}
\end{tabular}
}}

\end{document}